\documentclass[11pt,oneside,a4paper,titlepage]{book}
\usepackage[all]{xy}
\usepackage[pdftex]{graphicx}  
\usepackage{amsmath,amssymb,amsthm}
\usepackage{fancyhdr}
\fancypagestyle{plain}{
 \fancyhead{} 
 }

\makeatletter
\def\cleardoublepage{\clearpage\if@twoside \ifodd\c@page\else%
         \hbox{}%
     \thispagestyle{empty}
     \newpage%
     \if@twocolumn\hbox{}\newpage\fi\fi\fi}
\makeatother

\newcommand{\A}{{\mathbb A}}
\newcommand{\B}{{\mathbb B}}
\newcommand{\C}{{\mathbb C}}
\newcommand{\D}{{\mathbb D}}
\newcommand{\E}{{\mathbb E}}
\newcommand{\F}{{\mathbb F}}
\newcommand{\N}{{\mathbb N}}
\newcommand{\X}{{\mathbb X}}
\newcommand{\Z}{{\mathbb Z}}
\newcommand{\PP}{{\mathbb P}}
\newcommand{\ES}{{\mathbb S}}
\newcommand{\Q}{{\mathbb Q}}

\newcommand{\K}{{\mathbb K}}
\newcommand{\I}{{\mathbb{I}}}
\newcommand{\Id}[1]{\raisebox{-.9ex}{\xy
(0,2.4)*{\mathbb{I}};
(0,0)*{\scriptscriptstyle ({#1})}
\endxy}}
\newcommand{\x}[1]{\raisebox{-.5ex}{\xy
(0,2)*{\times};
(0,0)*{\scriptscriptstyle ({#1})}
\endxy}}
\newcommand{\eq}[1]{\ \raisebox{-.5ex}{\xy
(0,2)*{\scriptscriptstyle ({#1})};
(0,0)*{=}
\endxy}\ }

\newtheorem{Theorem}{Theorem}[chapter]
\newtheorem{Lemma}[Theorem]{Lemma}
\newtheorem{Proposition}[Theorem]{Proposition}
\newtheorem{Definition}[Theorem]{Definition}

\newtheorem{Corollary}[Theorem]{Corollary}

\theoremstyle{remark}\newtheorem{Remark}[Theorem]{Remark}
\theoremstyle{plain}\newtheorem{Claim}[Theorem]{Claim}

\newtheorem{Defprop}[Theorem]{Definition/Proposition}
\newtheorem{UP}[Theorem]{Universal Property}
\newtheorem{Notation}[Theorem]{Notation}
\usepackage[final,expansion=true,protrusion=true]{microtype}

\begin{document}
\frontmatter

\def \tm{\!\times\!}

\newenvironment{changemargin}[2]{\begin{list}{}{
\setlength{\topsep}{0pt}
\setlength{\leftmargin}{0pt}
\setlength{\rightmargin}{0pt}
\setlength{\listparindent}{\parindent}
\setlength{\itemindent}{\parindent}
\setlength{\parsep}{0pt plus 1pt}
\addtolength{\leftmargin}{#1}\addtolength{\rightmargin}{#2}
}\item}{\end{list}}

\begin{titlepage}

\begin{center}
\Large
{UNIVERSIT\`A DEGLI STUDI DI MILANO \\
\vspace{-,4cm}\rule{9truecm}{.4pt}\\
Facolt\`a di Scienze Matematiche Fisiche e Naturali \\
Dipartimento di Matematica ``F. Enriques''}\\

\vspace{1truecm}

\begin{figure}[h]
\centering
\end{figure}

\vspace{0.5truecm}

Corso di Dottorato di Ricerca in Matematica, XX ciclo\\

\vspace{0.7truecm}

Tesi di Dottorato di Ricerca\\

\vspace{0.5truecm}

\LARGE{\textsc{The Ziqqurath \\of Exact Sequences of $n$-Groupoids}}

\vspace{0.5truecm}

\Large{MAT$\backslash$02}
\end{center}

\vspace{0.7truecm}

\large{
\begin{center}
\begin{tabular}{lcl}
Relatori & \hspace{1cm} & Candidato\\
Prof. Enrico M. Vitale & \hspace{1,5cm} & Dott. Giuseppe Metere\\
Prof. Stefano Kasangian& & \\
& & \\
& & \\
Coordinatore di Dottorato & & \\
Prof. Antonio Lanteri & &
\end{tabular}
\end{center}}

\vspace{1,1cm}

\vfill
\begin{flushbottom}
\begin{center}
\rule{9truecm}{.4pt}\\
\small{Anno Accademico 2006--2007}
\end{center}
\end{flushbottom}
\end{titlepage}

\newpage
\thispagestyle{empty}
\verb"   "
\newpage
\thispagestyle{empty}

Ai miei genitori Lalla e Saverio,

che hanno sempre creduto in me.

A  Valentina,

senza di lei questa Tesi non sarebbe mai stata scritta.

\newpage
\thispagestyle{empty}
\verb"   "
\newpage
\thispagestyle{empty}
\begin{center}
    \textbf{Acknowledgments}
\end{center}

I would like to thank my supervisors Stefano Kasangian and Enrico Vitale,
for their encouragement and the pleasant time spent together working at this project.

I would also like to thank George Janelidze and Marco Grandis for their helpful suggestions and
remarks.

I acknowledge the members of the Category Seminar of the University of Milano for many fruitful
discussions, and I am especially grateful to Sandra Mantovani for her tireless support.

I am in debt with my wife Valentina Guarino, who  fought in vain against my
Ital-English while revising this thesis.

\newpage
\thispagestyle{empty}
\verb"   "
\tableofcontents

\mainmatter

\chapter{Introduction}

\emph{Higher Dimensional Categories}  are showing relevant implications in
Algebraic Topology (classification of homotopy $n$-types, operads, cobordism), and in Algebraic Geometry (Grothedieck's
$n$-stacks, non-abelian cohomology), not to mention recent
applications in Mathematical Physics (TQFT, higher order gauge
theory) and Computer Science.

Nevertheless basic algebraic tools, in order to further develop the
theory, are far from being established. A step forward towards  this
direction would be having an essential understanding of the notion
of exactness for Higher Dimensional Categories, and of the limits involved
in defining this notion.

\section{Summary}

In 1970 Ronald Brown published a paper \cite{brown70ja} on an
uprising area of mathematical research: the theory of groupoids.

According to the category-theorist, a \emph{groupoid} is a category
with all of its morphisms invertible w.r.t. composition
\cite{MacLane,Higgins71}. In fact the notion of groupoid was introduced earlier
as a generalization of the notion of \emph{group}, where
the binary operation is only partially defined.

In studying connections with algebraic topology and non-abelian
cohomology, Brown showed that, given a fibration $F$ of pointed
groupoids and its (strict) kernel $\K_s$:
\begin{equation}\label{sixterm}
\xymatrix{ \K_s\ar[r]^K&\B\ar[r]^F&\C}
$$
it was possible to obtain a 6-term exact sequence
$$
\xymatrix{ \pi_1\K_s \ar[r]^{\pi_1 K} & \pi_1\B \ar[r]^{\pi_1 F} &
\pi_1\C \ar@{.>}[r]^{\delta} & \pi_0\K_s \ar[r]^{\pi_0 K} & \pi_0\B
\ar[r]^{\pi_0 F} & \pi_0\C}
\end{equation}
of groups and pointed sets, where $\pi_0$ is the functor giving set
of isomorphism classes of object, and $\pi_1$
 the group of endomorphisms of the point.\\

Still, the condition of $F$ being a fibration was motivated by
``the analogy [\dots] with topological situations'' \cite{brown70ja}:
 this influenced the kind of limit considered, {\em i.e.} a (strict) kernel.

However it suffices to consider the homotopy kernel instead of the
strict kernel to remove the need of restricting to fibrations, so that construction above
still holds for a generic groupoid morphism.

These ideas developed further in \cite{HKK2002JPAA,Vitale2004fic},
where the authors generalized  Brown's result to 2-groupoids.
The 2-groupoids considered are weakly invertible 2-categories.
Duskin, Kieboom and Vitale showed in \cite{HKK2002JPAA} that, given
a morphism of 2-groupoids $G:\C\to\D$, it is still possible to get a
6-term sequence as in (\ref{sixterm}) of (strict) categorical
groups and pointed groupoids, where $\pi_0$ is the
classifying-functor, and $\pi_1$ gives the cat-group of
endomorphisms of the point. This sequence is exact in a suitable
sense (see \cite{Vitale2002jpaa}). Further, since
$\pi_0\circ\pi_1=\pi_1\circ\pi_0$, applying $\pi_1$ and $\pi_0$
again we get two 6-term exact sequences that can be pasted together
in a 9-term exact sequence, where the left-most three terms are
abelian groups, the central three  are groups and the
right-most three are just pointed sets.\\

The purpose of this thesis is to extend these results to a
n-dimensional context.

The setting is the sesqui-category \cite{street1996cs} of
$n$-groupoids, strict $n$-functors and lax $n$-transformations. We
consider a notion of $n$-groupoid equivalent to that of
\cite{MR1130401}, {\em i.e.} a weakly invertible $n$-category,
 but our approach is genuinely recursive.

Being more precise,  $n$-categories and $n$-functors are given by
means of the standard enrichment in the category of $(n-1)$-categories and
$(n-1)$-functors, w.r.t its cartesian bi-closed structure.

Differently, $n$-transformations considered come in a lax version,
being a direct generalization of those of \cite{B94HANDBOOK} for
2-categories.
In fact a notion of strict natural $n$-transformation is also sketched, but that
has shown inadequate in developing the theory for all morphisms, and not just for
fibrations.

The lax $n$-transformations introduced are equivalent to that considered
by Crans in \cite{Crans95}, yet inductive definition allows to deal with
such morphisms directly, without the complications of the theory of
pasting schemes \cite{Johnson89} (as in \cite{Crans95}).

The $n$-groupoids we consider are  $n$-categories which are locally
$(n-1)$-groupoids and whose 1-cells are equivalences. In this
setting we introduce a straightforward generalization of $\pi_0$
(lower ``$*$'' stays for pointed):
$$
\pi_0^{(n)}:n\mathbf{Gpd}_*\to(n-1)\mathbf{Gpd}_*
$$

Further, standard $h$-pullbacks \cite{Mather76} are introduced,
together with their 1-dimensional, i.e. ordinary, universal property. This allows to
deal with $h$-kernels, and to define a \emph{loop}-functor
$$
\Omega:n\mathbf{Cat}\to n\mathbf{Cat},\qquad
\C\mapsto\raisebox{5ex}{
\xymatrix{\Omega(\C)\ar@{.>}[r]\ar@{.>}[d]&\I\ar[d]^{[*]}\ar@{}[dl]|(.35){}="1"|(.65){}="2"
\\
\I\ar[r]_{[*]} &\C\ar@{=>}"2";"1"}}
$$
which  gives the functor
$$
\pi_1^{(n)}:n\mathbf{Gpd}_*\to(n-1)\mathbf{Gpd}_*,\qquad
\mathrm{with\ }\pi_1^{(n)}(\C)=\pi_0^{(n)}(\Omega(\C))
$$
To extend the last to a (contra-variant) sesqui-functor, it is
necessary to consider $h$-pullbacks with their 2-dimensional
universal property, this involving 3-morphisms between $n$-natural
transformations,  dimension-raising horizontal composition of
2-morphisms and the necessary algebra for those. To this end, an appropriate notion
of sesqui${}^2$category has been introduced.\\

Finally, we set up a notion of exactness for  a triple
$(K,\varphi,F)$ in $n\mathbf{Gpd}_*$
$$
\xymatrix@!C=8ex{ \K\ar[r]_{K}\ar@/^7ex/[rr]^{0}="a1"
&\B\ar@{}|{}="a2"\ar[r]_{F}&\C
\ar@{}"a1";"a2"|(.3){}="A1"|(.75){}="A2"\ar@{=>}"A1";"A2"^{\varphi}}
$$

and we prove the\\

\textsc{Main Result} \vskip3pt For any natural number $n$,
\begin{enumerate}
  \item Sesqui-functors $\pi_0^{(n)}$ and $\pi_1^{(n)}$ preserve exactness, the last up to
  reversing the directions of 2-morphisms.
  \item Sesqui-functors $\pi_0^{(n)}$ and $\pi_1^{(n)}$ commute, {i.e.}
  $$\pi_1^{(n)}\circ\pi_0^{(n-1)}=\pi_0^{(n)}\circ\pi_1^{(n-1)}.$$
  \item Given a morphism of pointed $n$-groupoids $F:\B\to\C$, with $h$-kernel $\K$,
  there exists a canonical   morphism  $\Delta$ of pointed $(n-1)$-groupoids, and 2-morphism $\delta$,
  such that the sequence below is exact:
$$
\xymatrix@!C=5ex{
\pi_1\K\ar[r]^{\pi_1K}\ar@/_5ex/[rr]_{0}="a1"
&\pi_1\B\ar@{}|{}="a2"\ar[r]^{\pi_1F}\ar@/^5ex/[rr]^{0}="d1"
&\ar@{}|{}="d2"
\pi_1\C\ar[r]_{\Delta}
\ar@/_5ex/[rr]_{0}="c1"
\ar@{}"a2";"a1"|(.3){}="A1"|(.75){}="A2"\ar@{=>}"A1";"A2"^{\pi_1\varphi}
\ar@{}"d1";"d2"|(.3){}="D1"|(.75){}="D2"\ar@{=>}"D1";"D2"^{\delta}
&
\pi_0\K\ar[r]_{\pi_0K}\ar@/^5ex/[rr]^{0}="xa1"\ar@{}|{}="c2"
&\pi_0\B\ar@{}|{}="xa2"\ar[r]_{\pi_0F}&\pi_0\C
\ar@{}"c2";"c1"|(.3){}="C1"|(.65){}="C2"\ar@{=}"C1";"C2"
\ar@{}"xa1";"xa2"|(.3){}="xA1"|(.75){}="xA2"\ar@{=>}"xA1";"xA2"^{\pi_0\varphi}}
$$
\end{enumerate}
Applying $\pi_0^{(n-1)}$ and $\pi_1^{(n-1)}$,  we get two six-term
sequences, exact by (1) above. Those can be pasted by (2) above, in a
nine-term exact sequence of $(n-2)$-groupoids (cells to be pasted
are dotted in the diagram):
$$
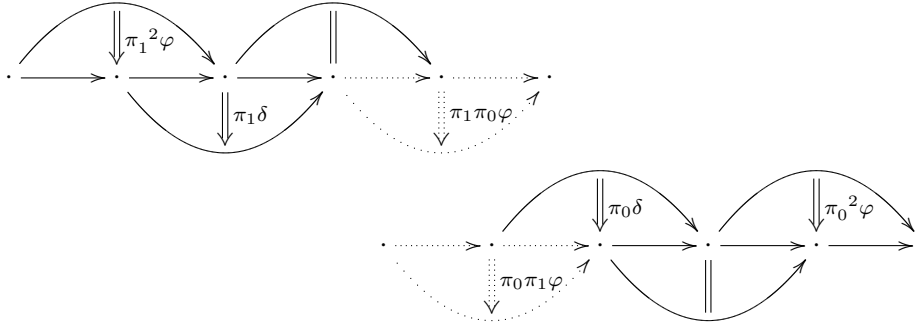
\raisebox{14ex}{\xymatrix@!C=3.5ex{
\cdot\ar[r]\ar@/^6ex/[rr]^{}="A2"
&\cdot\ar[r]\ar@{}|{}="A1"\ar@{=>}"A2";"A1"^{\pi_1^{\ 2}\varphi}
\ar@/_6ex/[rr]^{}="K2"
&\cdot\ar@{}|{}="K1"\ar@/^6ex/[rr]^{}="H2"
\ar[r]
&
\cdot\ar@{.>}[r]\ar@{..>}@/_6ex/[rr]^{}="B2"\ar@{}|{}="H1"
&\cdot\ar@{.>}[r]\ar@{}|{}="B1"\ar@{::>}"B1";"B2"^{\pi_1\pi_0\varphi}
&\cdot
\ar@{=>}"K1";"K2"^{\pi_1\delta}
\ar@{=}"H2";"H1"}}
\!\!\!\!\!\!\!\!\!\!\!\!\!\!\!\!\!\!\!\!\!\!\!\!\!\!\!\!\!\!\!\!\!\!\!\!\!\!\!
\xymatrix@!C=3.5ex{
\cdot\ar@{.>}[r]\ar@{..>}@/_6ex/[rr]_{}="A2"
&\cdot\ar@{.>}[r]\ar@{}|{}="A1"\ar@{::>}"A1";"A2"^{\pi_0\pi_1\varphi}
\ar@/^6ex/[rr]^{}="K2"
&\cdot\ar@{}|{}="K1"\ar@/_6ex/[rr]^{}="H2"
\ar[r]
&
\cdot\ar[r]\ar@/^6ex/[rr]^{}="B2"\ar@{}|{}="H1"
&\cdot\ar[r]\ar@{}|{}="B1"\ar@{=>}"B2";"B1"^{\pi_0^{\ 2}\varphi}
&\cdot
\ar@{=>}"K2";"K1"^{\pi_0\delta}
\ar@{=}"H1";"H2"}
$$

Iterating the process we finally obtain a $9\cdot n$ exact sequence
of pointed sets.

Furthermore, since the sesqui-functors $\pi_0$ force $h$-groups
structures, the first $9\cdot (n-1)$ terms are structured as an
exact sequence of groups, abelian up to the $9\cdot
(n-2)^{\mathrm{th}}$.

Similarly, the last-but-one step produces a $9\cdot (n-1)$ term
sequence of pointed groupoids, $9\cdot (n-2)$ of which are indeed
categorical groups \cite{Vitale2002jpaa}, finally $9\cdot (n-3)$ are
commutative.\\

Arranging these sequences from top (dimension $n$) to bottom (dimension $0$) we unveil
the shape of  a \emph{Ziqqurath}\footnote{Ziqquraths (or Ziggurats) were a type of step pyramid temples
 common to the inhabitants of  ancient Mesopotamia \cite{Ziqqu1,Ziqqu2}.}, in which each level is an exact sequence of $k$-groupoids ($k=0,\dots,n$).
In the relations between contiguous levels are nested classification properties of $n$-groupoids and their morphisms,
many of them are still to be investigated.

\section{Further developments}
The sesqui-categorical setting  presented here yields a fruitful
 perspective in the study of $n$-dimensional categorical structures.
In fact this is a general fact, and it permeates categorical investigations from its very beginning:

\begin{quote}
``{categories are \emph{two} steps away from naturality, the concept they were designed to formalize.
[\dots]
From the very study of the established practice of routinely specifying morphisms along with each mathematical
structure, we were presented, in the 1940's, with an extra dimension: morphisms between morphisms. We were
naturally led by naturality to objects, arrows \emph{and} 2-cells.} \cite{street1996cs}''
\end{quote}

Moreover the inductive approach permits to deal with a sesqui-categorical environment in any dimension,
carrying along the constructions while ascending  the dimensional ladder.\\

Two main points are currently being  investigated by the author.\\

First, the native  sesqui-categorical setting offers the chance to study weak structures inductively,
reducing most of coherence issues to (inductively nested) planar diagrams.
This would make  possible to describe some lax $n$-dimensional structures more easily
in a pure diagrammatic way, in terms of cells and of compositions given explicitly. Namely,
it would be interesting to consider an inductively defined sesqui-categories of $n$-groupoids with
weak units, in order to compare it with the special connected 3-dimensional case of \cite{JK2007}.
Similarly for other semi-strict  versions, as in \cite{Paoli2007}.  \\

Second, an application to abelian chain complexes is announced.
In fact, most authors prefer globally
defined versions of $n$-categorical structures, as it makes the internalization process easier.
In this way is usually proved the equivalence between $\omega$-categories in \textbf{Ab} (the category of
abelian groups) and non-negatively graded chain complexes of abelian groups.

 Yet it is possible to deal with this equivalence in our context
too. In fact, the category of length-$n$ abelian chain complexes is equivalent to the category of abelian group objects
in $n$\textbf{Cat} (\cite{LeinsterBook}, Example 1.4.11). More interestingly
this equivalence extends to homotopies and natural $n$-transformation.

Now, abelian group objects in $n$\textbf{Cat} behave well with respect to the $h$-structure of $n$\textbf{Cat}, similarly
to what is shown for monoid objects in a similar situation  by Grandis in \cite{MR1478526}.
 That is the reason why
it is worth studying  in this new setting how the theory extends to.\\

\section{Conventions}
The purpose of this section is to make life easier to the reader, stating the notational conventions (and other)
adopted. Nevertheless exceptions to these conventions are not rare, although always pointed out.\\

Doubled capitals as $\mathbb{A},\B,\C$ are used for $n$-categories and $n$-groupoids, capitals as $F,G,H$
for their morphisms, lower case Greek letters $\alpha, \beta, \gamma$ for their 2-morphism. Finally capital
Greek letters as $\Sigma, \Lambda$ are reserved to 3-morphisms.\\

Objects of a $n$-category ($n$-groupoid) are denoted by lower case Latin letters, subscripted by the number
``$0$''. Objects have often the same letter as the big doubled capital denoting the whole structure, using
primes (or other modifiers) for different objects, {\em e.g.} $a_0,a_0',a_0''$ are objects of $\A$.\\

Cells follow a similar convention, where the number represent the dimension of the cell, as the
 1-cell $b_1:\xymatrix@C=3ex{b_0\ar[r]& b_0'}$ of $\B$. 2-Cells are often
represented by double arrows ($b_2:\xymatrix@C=3ex{b_1\ar@2[r]& b_1'}$). As in higher dimension it
 would not be quite  practical,  for representing a $k$-cell we label the arrow itself with the number $k$:
$b_k:\xymatrix{b_{k-1}\ar[r]|{k}& b_{k-1}'}$.\\

 In order to avoid confusion with the name \emph{morphisms} (reserved for the sorts of the environment sesqui-category)
  the sorts of our $n$-categories ($n$-groupoids) are always named \emph{cells}.\\

A major exception to these rules is the Chapter on \emph{Sesqui-Categories}. Indeed it uses its own
notational conventions, explained therein.\\

Compositions are dealt with different symbols in order to distinguish the (internal) compositions
of cells from the (external) compositions of morphisms.\\

We adopt for cells-composition the empty circle superscripted by a number that represents the dimension of the
intersection cell. Example: $c_h\circ^m c_k$ means that the $h$-cell $c_h$ and the $k$-cell $c_k$
are composed along their common boundary, that is a $m$-cell $c_m$. In this case we will use the terms
$m$-domain and $m$-codomain. Let us point out that  superscripted $m$ is often omitted, specially when
it is $0$.\\

For the compositions of morphisms we use the filled circle superscripted by a number that represents the
dimension of the intersection morphisms. In the present work we will use only $0$-compositions and
$1$-compositions of morphism, hence the symbols $\bullet^0$ and $\bullet^1$. They come with a lower-scripted
$L$ or $R$ if they are left or right whiskering, respectively. Moreover dimension raising $0$-composition
of 2-morphisms is denoted $*$.\\

All composition-symbols are omitted when clear from the context.\\

The compositions of cells  morphisms will be written in algebraic order, {\em e.g.} for
$c_1:c_0\to c_0'$ and $c_1':c_0'\to c_0''$ we will write $c_1\circ^0 c_1':c_0\to c_0''$.
The other order is considered as \emph{evaluation}, so parentheses will be used, {\em e.g.} for
$F:\A\to\B$ and $G:\B\to\C$ we will write $G(F(-)):\A\to\C$.\\

\newpage
\section{Synopsis}
The thesis is organized as follows:\\

the rest of the introduction is dedicated to analyze low-dimensional cases, that inspired this
generalization;\\

the second chapter gives the sesqui-categorical-theoretical framework: different characterizations
are compared, finite products and $h$-pullbacks are introduced with their universal properties;\\

strict $n$-categories are defined in the third chapter, together with their morphisms ($n$-functor)
and their 2-morphisms (lax $n$-transformations); moreover finite products and standard $h$-pullbacks
are constructed  explicitly;\\

the fourth chapter introduces $n$-groupoids, $h$-surjective morphisms and
equivalences: these are necessary to formulate a notion of exactness for the sesqui-categories
of pointed $n$-groupoids; moreover we extend to them some classical result as the adjunction
discrete/iso-classes functor, and one point suspension;\\

in order to deal with 3-morphisms of $n$-categories, a new  framework is
defined in the fifth chapter, namely \emph{sesqui}${}^2$-\emph{categories};
lax $n$-modifications are introduced thereafter, together with other whiskering
and compositions; in fact a dimension raising 0-composition of 2-morphisms is given
and many useful algebraic properties are proved;\\

with the machinery developed in the previous ones, the sixth chapter
presents the main result: the construction of (a \emph{Ziqqurath} of) exact sequences in
any lower dimension from a given morphism of $n$-groupoids; this is achieved in few
steps, the starting point being a 2-dimensional property of pullbacks that $h$-pullbacks
in $n$\textbf{Cat} are proved to satisfy;\\

finally the appendix contains a comparison with the globular approach, the groupoid condition and
the choice of inverses.

\newpage

\section{Case study: dimension one }
In this section we recall, for reader's convenience, the construction set up in \cite{brown70ja}. \\
\subsection{Browns'result}

Let us suppose we are given a functor
$$
F:\B\to\C
$$
between two groupoids.\\

For every fixed object $b_0$ in $\B$, $F$ induces a map
$$
\mathrm{St}_F(b_0):\mathrm{St}_{\B}(b_0)\to\mathrm{St}_{\C}(Fb_0)
$$
where $\mathrm{St}_{\B}(b_0)=\bigcup_{b_0'\in \B_0} \B_1(b_0,b_0')$.

Suppose now that for every $b_0$, the map $\mathrm{St}_F(b_0)$ is surjective. Such $F$ is called
\emph{star-surjective}, or \emph{fibration}.

We can consider its (strict) kernel w.r.t. an object $b$ of $\B$:

\begin{equation}\label{dgm:brown1}
\xymatrix{\K_s\ar[r]^{K}&\B\ar[r]^F&\C}
\end{equation}

Here the groupoid $\K_s$ is just the strict fiber over the object $Fb$ of $\C$, \emph{i.e.} the groupoid
with objects $b_0$ of $\B$ such that $Fb_0=Fb$, and arrows $b_1$ of $\B$ such that $Fb_1=1_{b}$. Finally
$K$ is the natural inclusion.

Diagram (\ref{dgm:brown1}) above  can be \emph{restricted} to automorphisms groups over fixed objects,
 thus giving the exact sequence
$$
\xymatrix{1\ar[r]&\K_s(b,b)\ar[r]^{K^{b,b}}&\B(b,b)\ar[r]^{F^{b,b}}&\C(Fb,Fb)}
$$
\begin{proof}
Exactness in $\K_s(b,b)$ for $K$ injective, in $\B(b,b)$ for $\K_s= F^{-1}(Fb)$.
\end{proof}
Furthermore (\ref{dgm:brown1}) gives also an exact sequence of pointed sets, when we apply
the \emph{isomorphism-classes-functor}
$$
\pi_0:\mathbf{Gpd}\to\mathbf{Set}
$$
that sends a groupoid in the set of  classes of isomorphic objects.

We obtain the diagram
$$
\xymatrix{\pi_0\K_s\ar[r]^{\pi_0K}&\pi_0\B\ar[r]^{\pi_0F}&\pi_0\C}
$$
exact in $\pi_0\B$.
\begin{proof}
Clearly $\mathrm{Im}(\pi_0K)\subseteq \mathrm{Ker}(\pi_0F)$. Suppose then $\pi_0F(\{b_0\})=\{Fb_0\}=\{Fb\}$, where brackets
means iso-class. Then the hom-set $\C(Fb,Fb_0)$ is nonempty, containing an element $c_1$, say.
Star-surjectivity
in $b_0$  implies that there is a an arrow $b_1:b_0\to b'$ such that $F(b_1)={c_1}^{-1}$, but this means
$\{b_0\}=\{b'\}$. Since $Fb'=Fb$ the proof is complete.
\end{proof}

Finally we define a morphism of pointed sets
$$
\delta:\xymatrix{\C(Fb,Fb)\ar[r]&\pi_0\K_s}
$$
in the following way: given the arrow $c_1:Fb\to Fb$ star-surjectivity yields  a $b_1:b\to b'$ such that
$Fb_1=c_1$. Then we let $\delta(c_1)=\{b'\}$. Clearly this map is well defined, since for a different lifting
of $c_1$, its codomain is isomorphic to $b'$. Moreover it is obviously pointed by the identity.\\

Now the new sequence connected by $\delta$ is everywhere exact

\begin{equation}\label{eqn:Brown}
\xymatrix{1\ar[r]&\K_s(b,b)\ar[r]^{K^{b,b}}&\B(b,b)\ar[r]^{F^{b,b}}&\C(Fb,Fb)\ar[r]^{\delta}&\pi_0\K_s\ar[r]^{\pi_0K}&\pi_0\B\ar[r]^{\pi_0F}&\pi_0\C}
\end{equation}

where the last three terms are pointed sets, the other are groups.
\begin{proof}
It remains to prove the exactness in $\C(Fb,Fb)$  and in   $\pi_0\K_s$.

As for the first, let $b_1:b\to b$ be given.  Then among the liftings of $Fb_1$ there is $b_1$ itself, hence
$\delta(Fb_1)=\{b\}$. Conversely let $c_1:Fb\to Fb$ in the kernel of $\delta$. This means that for a lifting
$b_1:b\to b'$ of $c_1$, $\{b'\}=\{b\}$ in the fiber, {\em i.e.} there is a $b_1':b'\to b$ such that
$F(b_1')=1_{b}$.
Then $F(b_1\circ b_1')=F(b_1)\circ F( b_1')=F(b_1)\circ 1_{b}=F(b_1)=c_1$, with $b_1\circ b_1'\in \B(b,b)$.

For the second, let a $c_1:Fb\to Fb$ be given. Then $\delta(c_1)=\{b'\}$ for a lifting
$b_1:b\to b'$ of $c_1$. Hence $\{b'\}=\{b\}$ in $\pi_0 \B$. Conversely let $\{b_0\}$ in the kernel of
$\pi_0K$. This means $\{b_0\}=\{b\}$ in $\pi_0\B$, {\em i.e.} there exists a $b_1:b\to b_0$
in $\B$, which implies $\delta(Fb_1)=\{b_0\}$.
\end{proof}

A first attempt in extending this to higher dimensional structures has been done by
Hardie, Kamps and Kieboom in \cite{HKK2002JPAA}, where they obtain a similar result for
fibrations of bigroupoids.

Nevertheless the necessity to consider only fibrations is not just a limitation
in the choice of  morphisms, but it introduces also serious difficulties in trying to
further extend the result to $n$-groupoids.

On this lines Duskin, Kieboom, and Vitale proposed in \cite{Vitale2004fic} the different setting given by considering
homotopy kernels instead of kernels. Consequently a new notion of exactness was introduced.

\subsection{Brown's result revisited}\label{sec:Brown_rev}
In order to fully understand the generalization of \cite{Vitale2004fic}, we start by considering the
one dimensional construction w.r.t. homotopy kernels instead of strict kernels. Notice that here we
provide only constructions  since proofs are consequence of our general result.\\

Let us consider a functor
$$
F:\B\to\C
$$
between two groupoids. We define the $h$-fiber over a chosen object $Fb$ of $\C$
$$
\xymatrix@!C=8ex{ \K\ar[r]_{K}\ar@/^7ex/[rr]^{[Fb]}="a1"
&\B\ar@{}|{}="a2"\ar[r]_{F}&\C
\ar@{}"a1";"a2"|(.3){}="A1"|(.75){}="A2"\ar@{=>}"A1";"A2"^{\varphi}}
$$
where
\begin{itemize}
\item $[Fb]$ is the constant functor;
\item $\K$ is the comma-groupoid with objects the pairs $(b_0,c_1:\xymatrix{Fb\ar[r]&Fb_0})$.
An arrow $(b_0,c_1)\to(b_0',c_1')$ is a pair $(b_1,=)$ where $b_1:b_0,b_0'$, and the ``$=$''
stays for the equality $c_1=Fb_1\circ c_1'$;
\item $K:\K\to\B$ is the faithful functor defined by
$$
K\big((b_0,c_1)\big)=b_0,\qquad K\big((b_1,=)\big)=b_1;
$$
\item $\varphi:\xymatrix@C=3ex{[Fb]\ar@2[r]&KF}$ is the natural isomorphism with components
$$
\varphi_{(b_0,c_1)}=c_1:Fb\to Fb_0.
$$
\end{itemize}
Again we can get an exact sequence of groups and pointed sets

\begin{equation}\label{eqn:Brown_rev}
\xymatrix{1\ar[r]&\K(b,b)\ar[r]^{K^{b,b}}&\B(b,b)\ar[r]^{F^{b,b}}&\C(Fb,Fb)\ar[r]^{\delta}&\pi_0\K\ar[r]^{\pi_0K}&\pi_0\B\ar[r]^{\pi_0F}&\pi_0\C}
\end{equation}

where the connecting map $\delta$ is defined in a natural way
$$
\delta(c_1:\xymatrix@C=3ex{Fb\ar[r]&Fb})=(b,c_1)
$$
Now Brown's result can be seen as a corollary, and this gives a conceptual insight about
 the relation between these two different settings.

In fact, there exists a fully faithful functor $I_b:\K_s\to\K$; this is given explicitly
by letting $I_b(b_0)=(b_0,1_{Fb})$, indeed it is provided by the universal property defining the homotopy
kernel as a comparison functor.

Clearly $F$ is a fibration of groupoids if, and only if, for each object $b$ of $\B$ the functor
$I_b$ is essentially surjective on objects. Then when $F$ is a fibration, one can replace $\K(b,b)$ and
$\pi_0\K$ by $\K_s(b,b)$ and $\pi_0\K_s$, and obtain Brown's exact sequence (\ref{eqn:Brown}) from (\ref{eqn:Brown_rev}).

\section{Case study: dimension two}
In \cite{Vitale2004fic} the authors prove a similar result for morphisms of $2$-groupoids, {\em i.e.} weakly invertible strict
$2$-categories, and they claim that it easily extends to bi-groupoids. In the present work we will keep close to the
first setting in order to generalize it to $n$-groupoids (weakly invertible strict $n$-categories).\\

\subsection{Homotopy fibers}

Let us suppose we are given a morphism ($2$-functor) of $2$-groupoids
$$
F:\B\to\C
$$
that is a $2$-functor between two 2-categories in which every arrow is an equivalence and every 2-cell is an
isomorphism. If we fix an object $b$ of $\B$, the homotopy fiber $\F=\F_{F,Fb}$ of $F$ over $Fb$ is
the following $2$-groupoid:
\begin{itemize}
  \item objects are pairs $(b_0,\xymatrix@C=3ex{Fb_0\ar[r]^{c_1}&Fb})$;
  \item an arrow $(b_0,c_1)\to(b_0',c_1')$ is a pair $(b_1,c_2)$ as in the diagram below
$$
\xymatrix{
b_0\ar[rr]^{b_1}
&&b_0'
\\
Fb_0\ar[rr]^{Fb_1}
\ar[ddr]_{c_1}^{}="1"
&&Fb_0'\ar@{}|{}="2"
\ar[ddl]^{c_1'}
\\
\\
&Fb
\ar@{}"1";"2"|(.25){}="a1"|(.6){}="a2"
\ar@2"a1";"a2"^{c_2}}
$$
  \item a 2-cell $(b_1,c_2)\Rightarrow(b_1',c_2')$ is a pair $(b_2,\equiv)$ as in the diagram below
$$
\xymatrix{
Fb_1\ar@2[rr]^{Fb_2}
&&Fb_1'
\\
Fb_1\circ c_1'\ar@2[rr]^{Fb_2\circ id_{c_1'}}\ar@{}|{}="2"
&&Fb_1'\circ c_1'
\\
\\
&c_1
\ar@2[uur]_{c_2'}^{}="1"
\ar@2[uul]^{c_2}
\ar@{}"1";"2"|(.25){}="a1"|(.55){}="a2"
\ar@{3-3-}"a1";"a2"}
$$
\end{itemize}
The homotopy fiber comes with a 2-functor that embeds it in $\B$
$$
K:\F_{F,Fb}\to \B
$$
which sends the 2-cell
$$
(b_2,\equiv):(b_1,c_2)\Rightarrow(b_1',c_2'):(b_0,c_1)\to(b_0',c_1')
$$
to $b_2:b_1\Rightarrow b_1'$.

\subsection{$2$-Exact sequences}
Before going further with our description, we urge to introduce a notion of exactness suitable for a
2-dimensional context.

A notion of \emph{2-exactness} has been introduced  by Vitale in \cite{Vitale2002jpaa}, in order to study some classical exact sequences of abelian groups associated with a
morphism of commutative unital rings from sequences of pointed groupoids and categorical groups.

We report the definition in the context of pointed groupoids, as for categorical groups it applies plainly with
no changes\footnote{As a guiding analogy, do consider that exactness in the category of groups may be defined on the
underlying pointed sets.}.\\

Let us consider the 2-category of pointed groupoids $\mathbf{Gpd}_*$, where the  morphisms are
functors that preserves the base points, 2-morphisms are natural isomorphisms whose component at the base point
is the identity on the point.

For a given morphism $F\B\to\C$, we define its $h$-kernel as the triple $(\K,K:\K\to\B,\varphi:[*]\Rightarrow KF)$
($[*]$ denoting the constant $0$-functor) satisfying the following universal property
\begin{UP}[$h$-kernels]
  For any other triple $(\K',K',\varphi')$ there exists a unique $L:\K'\to\K$ such that $K'=LK$.
$$
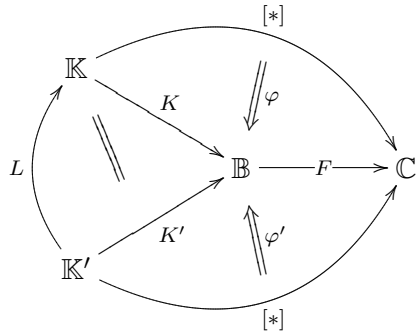
\xymatrix@!C=8ex{
\K\ar@{}|{}="k2"
\ar[dr]^{K}\ar@/^7ex/[drr]^{[*]}="a1"
\\
&\B\ar@{}|{}="a2"\ar[r]|{F}&\C
\\
\K'\ar@/^3.5ex/[uu]^{L}
\ar[ur]_{K'}^(.4){}="k1"\ar@/_7ex/[urr]_{[*]}="b1"
\ar@{}"a1";"a2"|(.3){}="A1"|(.75){}="A2"
\ar@2"A1";"A2"^{\varphi}
\ar@{}"b1";"a2"|(.3){}="B1"|(.75){}="B2"
\ar@2"B1";"B2"_{\varphi'}
\ar@{}"k1";"k2"|(.3){}="K1"|(.7){}="K2"
\ar@{=}"K1";"K2"
}
$$
\end{UP}
This universal property defines the ($h$-)kernel up to isomorphism.

\begin{Remark}
 Let us notice that last universal property uses only whiskerings of morphisms with a 2-morphism, and does not use
the full horizontal composition of 2-morphisms available in a 2-category.
Hence a step forward towards a full generalization of Brown's result to weak $n$-structures can be accomplished
by developing a theory that deals with these $1.5$-universal properties (a.k.a. \textit{sesqui}-universal, a.k.a. $h$-universal).
\end{Remark}

\begin{Remark}
In dimension $1$, our $h$-kernel satisfies also a universal property  of a bi-limit, and this is indeed
a point of view well considered in  \cite{Vitale2004fic}. Nevertheless, last \emph{Remark} motivate
the choice to restrict our attention to the $h$-limits considered.
\end{Remark}

Finally we are able to give the following
\begin{Definition}
A triple $(E,\varphi,F)$
$$
\xymatrix@!C=8ex{ \A\ar[r]_{E}\ar@/^7ex/[rr]^{[*]}="a1"
&\B\ar@{}|{}="a2"\ar[r]_{F}&\C
\ar@{}"a1";"a2"|(.3){}="A1"|(.75){}="A2"\ar@{=>}"A1";"A2"^{\varphi}}
$$
in $\mathbf{Gpd}_*$ is called \emph{exact} if the comparison with the 2-kernel is full and essentially
 surjective on objects.
\end{Definition}
Let us notice that in the above definition the 2-kernel can be replaced by the $h$-kernel, since fullness and
essential surjectivity are preserved by equivalences.

We leave to the conscious reader the deepening of the theory of such exactness for the $2$-dimensional
context of pointed groupoids and categorical groups, w.r.t. cohomology, extensions, chain-complexes,
(\cite{Vitale2002jpaa, Vitale2004fic, MR1983014, MR2072344, KV2004TAC, MR2259264, MR2205216, MR2110526, MR2262378, MR2141594, MR2116324}).

\subsection{Lowering the dimension: first step}
Back to \cite{Vitale2004fic}, let us consider a morphism of 2-groupoids $F:\B\to\C$ and its homotopy
kernel $\K\to\B$. They define indeed an exact sequence in the sesqui-category $2\mathbf{Gpd}_*$, but this
point of view is not analyzed explicitly in \cite{Vitale2004fic}.
Instead the authors consider the diagram
$$
\xymatrix{\K\ar[r]^{K}&\B\ar[r]^{F}&\C}
$$
and they apply to that hom-of-the-point functor $[-]_1(*,*):2\mathbf{Gpd}_*\to\mathbf{Gpd}_*$
the classifying functor $\mathcal{C}\ell:2\mathbf{Gpd}_*\to\mathbf{Gpd}_*$.\\

The first yields a 2-exact sequence
$$
\xymatrix@!C=8ex{ \K_1(*,*)\ar[r]^{K_1^{*,*}}
&\B_1(*,*)\ar[r]^{F_1^{*,*}}&\C_1(*,*)
}
$$

The second assigns to a 2-groupoid, the groupoid with the same set of objects, and whose
arrows are 2-isomorphism
classes of arrows \cite{MR0220789}, thus providing the 2-exact sequence
$$
\xymatrix@!C=8ex{ \mathcal{C}\ell\K\ar[r]^{\mathcal{C}\ell K}
&\mathcal{C}\ell\B\ar[r]^{\mathcal{C}\ell F}&\mathcal{C}\ell \C
}$$
Moreover it is possible to define a connecting functor
$$
\delta : \C_1(*,*)\to\mathcal{C}\ell\K
$$
such that the 6-term  sequence obtained this way is everywhere 2-exact.
The functor $\delta$ is defined as follows:
\begin{itemize}
  \item (\emph{on objects}) given an object $c_1:*\to *$ in the domain,
  $$
\delta(c_1)=(*,c_1);
  $$
  \item (\emph{on arrows}) given an arrow $c_2:c_1\Rightarrow c_1'$ in the domain,
  $$
\delta(c_2)=\{(1_*,c_2)\},
  $$
  where brackets denote $\mathcal{C}\ell$asses (it is well defined).
\end{itemize}

\subsection{Lowering the dimension: second step}
The 6-term exact sequence
$$
\xymatrix@!C=8ex{ \K_1(*,*)\ar[r]^{K_1^{*,*}}
&\B_1(*,*)\ar[r]^{F_1^{*,*}}&\C_1(*,*)
\ar[r]^{\delta}
&
\mathcal{C}\ell\K\ar[r]^{\mathcal{C}\ell K}
&\mathcal{C}\ell\B\ar[r]^{\mathcal{C}\ell F}&\mathcal{C}\ell \C}
$$
with obvious transformations is such that the left-most three terms underly a strict monoidal structure
given by the (former) 0-composition. Moreover, since we started with (weakly) invertible strict 2-categories,
they are indeed categorical groups.\\

Now, if we denote by $\pi_1$ the functor $\mathbf{Gpd}_*\to\mathbf{Set}_*$ that assigns to a pointed groupoid
the pointed set of the isomorphism classes of its objects, it is possible to show that it preserves exactness,
{\em i.e.} it sends 2-exact sequences of pointed groupoids (categorical groups), to exact-sequences of
pointed sets (groups). The same can be said of the functor $\pi_0$.

Moreover $\pi_1(\mathcal{C}\ell(-))=\pi_0([-]_1(*,*))$, hence we get a 9-term exact sequence
\begin{changemargin}{-10ex}{-10ex}
$$
\xymatrix@C=4ex{
\pi_1(\K_1(*,*))\ar[r]
&\pi_1(\B_1(*,*))\ar[r]
&\pi_1(\C_1(*,*))\ar[dll]
\\
\pi_1(\mathcal{C}\ell\K)=\pi_0(\K_1(*,*))\ar[r]
&\pi_1(\mathcal{C}\ell\B)=\pi_0(\B_1(*,*))\ar[r]
&\pi_1(\mathcal{C}\ell \C)=\pi_0(\C_1(*,*))\ar[dll]
\\
\pi_0(\mathcal{C}\ell\K)\ar[r]
&\pi_0(\mathcal{C}\ell\B)\ar[r]
&\pi_0(\mathcal{C}\ell \C)
}
$$
\end{changemargin}
where the three left-most terms are abelian groups, the three central terms are groups, the three right-most
terms are pointed sets. The reason why the three leftmost terms are abelian follows from a general
fact of strict $n$-categories for homs over an identity cell (see \cite{simpson-1998}), that is another
\textit{variazione} on the classical Eckmann-Hilton argument.

\chapter{Basics on sesqui-categories}\label{cha:sesqui}
\section{Sesqui-categories}

The notion of sesqui-category is due to Ross Street \cite{street1996cs}. The term \emph{sesqui} comes from the latin
\textit{semis-que}, that means (one and) a half. Hence a sesqui-category is something in-between a category and a 2-category.
More precisely
\begin{Definition}
A sesqui-category $\mathcal{C}$ is a category $\lfloor\mathcal{C}\rfloor$ with a lifting of the \mbox{hom-functor} to \textbf{Cat}, such that
the following diagram of categories and functors commutes, $\mathbf{obj}$ being the functor that forgets the morphisms:
\begin{equation}\label{def:sesquicat}
\xymatrix@C=15ex@R=10ex{
&\mathbf{Cat}\ar[d]^{\mathbf{obj}}\\
\lfloor\mathcal{C}\rfloor^{\mathrm{op}}\times\lfloor\mathcal{C}\rfloor
\ar[ur]^{\mathcal{C}(-,-)}
\ar[r]_{\lfloor\mathcal{C}\rfloor(-,-)}
&\mathbf{Set}
}
\end{equation}
Objects  and morphisms of $\lfloor\mathcal{C}\rfloor$ are also objects and 1-cells of $\mathcal{C}$, while  morphisms of
$\mathcal{C}(A,B)$'s
(with $A$ and $B$ running in $\mathbf{obj}(\lfloor\mathcal{C}\rfloor)$) are the 2-cells of $\mathcal{C}$.
\end{Definition}

We first observe that the definition above induces a 2-graph structure on $\mathcal{C}$, whose underlying graph underlies
the category $\mathcal{C}$. Besides, the functor $\mathcal{C}(-,-)$ provides hom-sets of the category $\lfloor\mathcal{C}\rfloor$
with a category structure, whose composition is termed \emph{vertical composition} (or 1-composition) of 2-cells.
Finally, condition expressed by diagram (\ref{def:sesquicat}) on the lifting $\mathcal{C}(-,-)$ gives a reduced horizontal
composition, or \emph{whiskering} (or 0-composition), compatible with 1-cell composition and with the 2-graph structure of
$\mathcal{C}$.

In fact, for $\xymatrix{A'\ar[r]^a&A}$ and $\xymatrix{B\ar[r]^b&B'}$ in $\lfloor\mathcal{C}\rfloor$, the functor
$$
\xymatrix{\mathcal{C}(a,b):\mathcal{C}(A,B)\ar[r] &\mathcal{C}(A',B')}
$$
gives explicitly such a composition: for a 2-cell $\xymatrix@C=3ex{\alpha:f\ar@{=>}[r]&g:A\ar[r]&B}$,
it whiskers the diagram
$$
\xymatrix{
A'\ar[r]^a
&A\ar@/^3ex/[r]^f_{}="1"\ar@/_3ex/[r]_g^{}="2"
&B\ar[r]^b
&B'
\ar@{=>}"1";"2"^{\alpha}
}
$$
to get the 2-cell
$$
\xymatrix@C=16ex{
A'\ar@/^3ex/[r]^{a\bullet f\bullet b}_{}="1"\ar@/_3ex/[r]_{a\bullet g\bullet b}^{}="2"
&B'
\ar@{=>}"1";"2"^{a\bullet \alpha\bullet  b}
}
$$
where $a\bullet \alpha\bullet  b$ is just a concise form for $\mathcal{C}(a,b)(\alpha)$.\\

By functoriality of whiskering, the operation may also be given in a \emph{left-and-right} fashion. In fact
it suffices to identify
$$
a\bullet_L\alpha = a\bullet \alpha\bullet 1_B, \quad \alpha\bullet_R b= 1_A \bullet \alpha \bullet b
$$
This fact can be made precise, and gives a more tractable definition,  by the following characterization
(see, for example \cite{grandis1994hah,stell1994mtr}):
\begin{Proposition}\label{prop:sesqui-categories}
Let $\mathcal{C}$ be a reflexive 2-graph
$$
\xymatrix{
\mathcal{C}_2\ar@<.7ex>^{s}[r]\ar@<-.7ex>_{t}[r]
&\ar[l]|e\mathcal{C}_1\ar@<.7ex>^{s}[r]\ar@<-.7ex>_{t}[r]
&\ar[l]|e\mathcal{C}_0
}
$$
whose underlying graph $\lfloor\mathcal{C}\rfloor=\xymatrix{
\mathcal{C}_1\ar@<.7ex>^{s}[r]\ar@<-.7ex>_{t}[r]
&\ar[l]|e\mathcal{C}_0
}$ has a category structure.
Then $\mathcal{C}$ is a sesqui-category precisely when the following conditions hold:
\begin{enumerate}
  \item for every pair of objects $A, B$ of $\mathcal{C}_0$, the  graph $\mathcal{C}(A,B)$ has a category
  structure, called the hom-category of $A,B$.
  \item (partial) reduced horizontal compositions are defined, {\em i.e.} for every $A',A, B$ and $B'$ objects of
  $\mathcal{C}_0$, composition in $\lfloor\mathcal{C}\rfloor$ extends to binary operations
\begin{eqnarray}
\bullet_L &:& \xymatrix{: \lfloor\mathcal{C}\rfloor(A',A)\times\mathcal{C}(A,B)\ar[r]& \mathcal{C}(A',B)} \\
\bullet_R &:& \xymatrix{\mathcal{C}(A,B)\times\lfloor\mathcal{C}\rfloor(B,B')\ar[r]& \mathcal{C}(A,B')},
\end{eqnarray}
that satisfy equations below, whenever the composites are defined:
\end{enumerate}
$$
\xymatrix{A''\ar[r]^{a'}&A'\ar[r]^a & A\ar@/^3ex/[r]^f_{}="1"\ar@/_3ex/[r]_g^{}="2"& B\ar[r]^b&B'\ar[r]^{b'}&A''\ar@{=>}"1";"2"^{\alpha}}
$$
\begin{changemargin}{-10ex}{-10ex}
$$
  \begin{array}{crclcrcl}
    (L1) & 1_A\bullet_L \alpha &=& \alpha & (R1) &\alpha\bullet_R 1_B&=&\alpha \\
    (L2)&a'a\bullet_L\alpha &=& a'\bullet_L(a\bullet_L\alpha)& (R2)& \alpha\bullet_R bb'&=&(\alpha\bullet_R b)\bullet_R b' \\
    (L3)& a\bullet_L 1_f&=&1_{af}& (R3) &1_f\bullet_R b &=& 1_{fb}\\
    (L4)& a\bullet_L(\alpha\cdot\beta)&=&(a\bullet_L \alpha)\cdot(a\bullet_L \beta)& (R4) &(\alpha\cdot\beta)\bullet_R b &=& (\alpha\bullet_R b)\cdot(\beta\bullet_R b)\\
    (LR5)& (a\bullet_L \alpha)\bullet_R b&=&a\bullet_L (\alpha\bullet_R b)&\\
  \end{array}
$$
\end{changemargin}
In these equations, $1_A$ and $1_B$ are identity 1-cells, while $1_f$, $1_{af}$ and $1_{fb}$ are identity 2-cells, and $\cdot$
is the (vertical) composition inside the hom-categories. Axiom (LR5) will be also called \emph{whiskering axiom}.
\end{Proposition}
\begin{proof}
Let $\mathcal{C}$ be a sesqui-category.
The fact that $\mathcal{C}(A,B)$ are categories is clear from the definition, hence $1$ is satisfied.
Now, define for chosen $A,A',B$ and $B'$
$$
-_1\bullet_L-_2 = \mathcal{C}(-_1,1_B)(-_2)\quad -_1\bullet_R-_2 = \mathcal{C}(1_A,-_2)(-_1)
$$
Then for reduced left composition axioms, we have:\\

(L1)
$$
1_A\bullet_L\alpha=\mathcal{C}(1_A,1_B)(\alpha)=\alpha
$$
by functoriality w.r.t. units of $\mathcal{C}(-,-)$.\\

(L2)
$$
(a'a)\bullet_L\alpha=\mathcal{C}(a'a,1_B)(\alpha)=\mathcal{C}((a',1_B)(a,1_B))(\alpha)=
\mathcal{C}(a',1_B)(\mathcal{C}(a,1_B)(\alpha))
$$
by  functoriality w.r.t. composition of $\mathcal{C}(-,-)$. Notice the contravariance on the first component.\\

(L3)
$$
a\bullet_L 1_f = \mathcal{C}(a,1_B)(1_f)=1_{\mathcal{C}(a,1_B)(f)}=1_{af}
$$
by  functoriality w.r.t. units of $\mathcal{C}(a,1_B)$.\\

(L4)
$$
a\bullet_L (\alpha\cdot\beta)= \mathcal{C}(a,1_B)(\alpha\cdot\beta) =\mathcal{C}(a,1_B)(\alpha)\cdot\mathcal{C}(a,1_B)(\beta)=
(a\bullet_L\alpha)\cdot (a\bullet_L\beta)
$$
by  functoriality w.r.t. composition of $\mathcal{C}(a,1_B)$.\\

Analogous proofs hold for (R1) to (R4). Finally\\

(LR5)
$$
(a\bullet_L \alpha)\bullet_R b= \mathcal{C}(1_{A'},b)(\mathcal{C}(a,1_B)(\alpha))=\mathcal{C}((1_{A'},b)(a,1_B))(\alpha)=
\mathcal{C}(a,b)(\alpha)=
$$
$$
=\mathcal{C}((a,1_{B'})(1_A,b))(\alpha)=\mathcal{C}(a,1_{B'})(\mathcal{C}(1_A,b)(\alpha))
=a\bullet_L (\alpha\bullet_R b)
$$
by  functoriality w.r.t. composition of $\mathcal{C}(-,-)$.\\

Conversely, supposing  we are given a reflexive 2-graph $\mathcal{C}$ underlying a category $\lfloor\mathcal{C}\rfloor$,
hom-categories (1) and left/right-compositions satisfying conditions above (2).  We show $\mathcal{C}$ is a sesqui-category.

To this end we define a functor
$$\xymatrix@C=15ex{
\lfloor\mathcal{C}\rfloor^{\mathrm{op}}\times\lfloor\mathcal{C}\rfloor
\ar[r]^{\mathcal{C}(-,-)}
&\mathbf{Cat}
}$$
On objects, this is given by condition (1); on arrows, for
$$
\xymatrix{(a,b):(A,B)\ar[r]&(A',B')}
$$
({\em i.e.}\/ $\xymatrix{a:A'\ar[r]&A}$ and $\xymatrix{b:B\ar[r]&B'}$), equation (2.5) helps us to define
$$
\xymatrix{\mathcal{C}(a,b):\mathcal{C}(A,B)\ar[r]&\mathcal{C}(A',B')}
$$
by
$$
\mathcal{C}(a,b)(-)= a\bullet_L - \bullet_R b
$$
These assignments give indeed functors and are functorial. \\

\textit{(i)} $\mathcal{C}(a,b)$ is functor w.r.t. units
$$
\mathcal{C}(a,b)(1_f) = a\bullet_L 1_f \bullet_R b = 1_{afb}= 1_{\mathcal{C}(a,b)(f)}
$$
by (L3) and (R3)\\

\textit{(ii)} $\mathcal{C}(a,b)$ is functor w.r.t. composition
$$
\mathcal{C}(a,b)(\alpha\cdot\beta) = a\bullet_L (\alpha\cdot\beta) \bullet_R b =
((a\bullet_L \alpha) \cdot (a\bullet_L \beta))\bullet_R b =
$$
$$
=(a\bullet_L \alpha \bullet_R b)\cdot(a\bullet_L \beta \bullet_R b)
=\mathcal{C}(a,b)(\alpha)\cdot\mathcal{C}(a,b)(\beta)
$$
by (L4) and (R4).\\

\textit{(iii)} $\mathcal{C}(-,-)$ is functor w.r.t. units
$$
\mathcal{C}(1_A,1_B)(\alpha)= 1_A \bullet_L \alpha \bullet_R 1_B= \alpha
$$
by (L1) and (R1).\\

\textit{(vi)} $\mathcal{C}(-,-)$ is functor w.r.t. composition
$$
\mathcal{C}(a'a,bb')(\alpha)= (a'a)\bullet_L \alpha \bullet_R (bb') = ((a'\bullet_L (a\bullet_L\alpha))\bullet_R b)\bullet_R b' =
$$
$$
= (a'\bullet_L \mathcal{C}(a,b)(\alpha))\bullet_R b' =   \mathcal{C}(a',b') (\mathcal{C}(a,b)(\alpha))
$$
by (L2) and (R2).\\

That $\mathcal{C(-,-)}$ makes (\ref{def:sesquicat}) commute is immediate from its definition.
\end{proof}

Notice that reduced horizontal left/right composition will be often denoted simply by $\bullet$, when this does not cause
ambiguity.
\section{Morphisms of sesqui-categories}
Morphisms between sesqui-categories are termed sesqui-functors. More precisely a sesqui-functor
$\xymatrix{\mathcal{F}:\mathcal{C}\ar[r]&\mathcal{D}}$ is a 2-graph morphism such that
\begin{itemize}
  \item $\xymatrix{\lfloor\mathcal{F}\rfloor:\lfloor\mathcal{C}\rfloor\ar[r]&\lfloor\mathcal{D}\rfloor}$
is a functor,
  \item for every $A,B$ in $\mathcal{C}_0$,
$$\xymatrix{\mathcal{F}^{A,B}:\mathcal{C}(A,B)\ar[r]&\mathcal{D}(\mathcal{F}(A),\mathcal{F}(B))}$$
are  functors component of a  natural transformation $\mathfrak{F}$
\begin{equation}
\xymatrix{
\lfloor\mathcal{C}\rfloor^{op}\times\lfloor\mathcal{C}\rfloor\ar[dr]^{\mathcal{C}(-,-)}_{}="1"
\ar[dd]_{\lfloor\mathcal{F}\rfloor^{op}\times\lfloor\mathcal{F}\rfloor}\\
&\mathbf{Cat}\\
\lfloor\mathcal{D}\rfloor^{op}\times\lfloor\mathcal{D}\rfloor\ar[ur]_{\mathcal{D}(-,-)}
\ar@{}"1";|(.2){}="3"|(.6){}="4"
\ar@{=>}"3";"4"_{\mathfrak{F}}}
\end{equation}
that lifts
$\xymatrix@C=3ex{\lfloor\mathfrak{F}\rfloor: \lfloor\mathcal{C}\rfloor(-,-)\ar@{=>}[r]
&(\lfloor\mathcal{F}\rfloor^{op}\times\lfloor\mathcal{F}\rfloor)\cdot}\lfloor\mathcal{D}\rfloor(-,-)$.
\end{itemize}
\begin{Remark}
Notice that every functor between categories gives rise to such a natural transformation as $\lfloor\mathfrak{F}\rfloor$
for $\lfloor\mathcal{F}\rfloor$. From this point of view, the last condition may be re-formulated saying that
\emph{a sesqui-functor is the lifting of a functor between the underlying categories}.
\end{Remark}

We can translate the definition of sesqui-functor in terms of left/right compositions:

\begin{Proposition}
Let $\mathcal{C}$ and $ \mathcal{D}$ be sesqui-categories, and let
$$
\xymatrix{
\mathcal{F}:\mathcal{C}\ar[r]&\mathcal{D}}
$$
be a 2-graphs homomorphism, whose underlying graph homomorphism
$$
\xymatrix{\lfloor\mathcal{F}\rfloor:\lfloor\mathcal{C}\rfloor\ar[r]&\lfloor\mathcal{D}\rfloor}
$$
is a functor.

Then $\mathcal{F}$ is a sesqui-functor precisely when the following conditions hold:
\begin{enumerate}
  \item for every pair of objects $A, B$ of $\mathcal{C}_0$, the  graph homomorphism
  $$
  \xymatrix{\mathcal{F}^{A,B}:\mathcal{C}(A,B)\ar[r]&\mathcal{D}(\mathcal{F}(A),\mathcal{F}(B))}
  $$ is a functor,
  called the hom-functor at $A,B$.
  \item (partial) horizontal reduced compositions are preserved, {\em i.e.} for every diagram
  \end{enumerate}
$$
\xymatrix{A'\ar[r]^a & A\ar@/^3ex/[r]^f_{}="1"\ar@/_3ex/[r]_g^{}="2"& B\ar[r]^b&B'\ar@{=>}"1";"2"^{\alpha}}
$$
in $\mathcal{C}_0$, equations below hold:
$$
  \begin{array}{crclcrcl}
  (L6)&\mathcal{F}(a\bullet_L \alpha)&=&\mathcal{F}(a)\bullet_L\mathcal{C}(\alpha)&\quad(R6)&\mathcal{F}(\alpha\bullet_R b)&=&\mathcal{F}(\alpha)\bullet_R \mathcal{F}(b)
  \end{array}
$$
\end{Proposition}
\begin{proof}
Let $\xymatrix{\mathcal{F}:\mathcal{C}\ar[r]&\mathcal{D}}$ be a sesqui-functor. Then $\mathcal{F}$ is \textit{a fortiori}\/ a
homomorphism of 2-graphs, underlying a functor $\lfloor\mathcal{F}\rfloor$. Furthermore, for every choice of $A$ and $B$ in
$\mathcal{C}_0$, the $\mathcal{F}^{A,B}$ are functors too. What remains to prove is that $\mathcal{F}$ preserves
left/right-compositions in the sense of (L6) and (R6), and this follows easily from naturality of $\mathfrak{F}$.
In fact, for (L6)
\begin{eqnarray*}
\mathcal{F}(a\bullet_L\alpha) &=& \mathcal{F}(\mathcal{C}(a,1_B)(\alpha)) \mathrm{\ by\ definition}\\
&=&\mathcal{D}(\mathcal{F}(a),\mathcal{F}(1_B))(\mathcal{F}(\alpha)) \mathrm{\ by\ naturality}\\
&=&\mathcal{D}(\mathcal{F}(a),1_{\mathcal{F}(B)})(\mathcal{F}(\alpha)) \mathrm{\ by\ functoriality}\\
&=&\mathcal{F}(a)\bullet_L\mathcal{F}(\alpha)\mathrm{\ by\ definition}
\end{eqnarray*}
(R6) from a similar calculation.\\

Conversely, suppose we are given two sesqui-categories $\mathcal{C}$ and $\mathcal{D}$,
together with   a 2-graph homomorphism $\mathcal{F}$ satisfying
conditions (1) and (2) above.

We will prove naturality of $\mathfrak{F}$, \emph{i.e.} for every
$\xymatrix{A'\ar[r]^a&A}$ and $\xymatrix{B\ar[r]^b&B'}$
in $\mathcal{C}$, the following is a commutative diagram in \textbf{Cat}:
$$
\xymatrix@C=16ex{
\mathcal{C}(A,B)\ar[r]^{\mathcal{C}(a,b)}\ar[d]_{\mathcal{F}^{A,B}}
&\mathcal{C}(A',B') \ar[d]^{\mathcal{F}^{A',B'}}\\
\mathcal{D}(\mathcal{F}(A),\mathcal{F}(B))\ar[r]_{\mathcal{D}(\mathcal{F}(a),\mathcal{F}(b))}
&\mathcal{D}(\mathcal{F}(A'),\mathcal{F}(B'))}
$$
That this diagram commutes on objects (\emph{i.e.} on 1-cells of $\mathcal{C}$ and $\mathcal{D}$) is clear from the fact that
$\lfloor\mathcal{F}\rfloor$ is a functor and that left/right-compositions extend 1-cell-compositions.
Finally, for a 2-cell $\alpha$ as above,
\begin{eqnarray*}
\mathcal{F}(\mathcal{C}(a,b)(\alpha))&=&\mathcal{F}(a\bullet_L\alpha\bullet_R b)\\
&=&\mathcal{F}(a\bullet_L\alpha)\bullet_R \mathcal{F}(b)\\
&=&\mathcal{F}(a)\bullet_L\mathcal{F}(\alpha)\bullet_R \mathcal{F}(b)\\
&=&\mathcal{D}(\mathcal{F}(a),\mathcal{F}(b))(\mathcal{F}(\alpha))
\end{eqnarray*}
follows from (L6) and (R6).\\
\end{proof}

\begin{Remark}
In the following we will need the notion of $2$\emph{-contravariant sesqui-functor}. This is simply a sesqui-functor as
above, such that the functors component $\mathcal{F}^{A,B}$ are usual contravariant functors.

Of course, characterization above still holds, \emph{mutatis mutanda}: {\em e.g.} if
$$
\alpha:f\Rightarrow g:A\to B,
$$
then
$$
\mathcal{F}(\alpha):\mathcal{F}(g)\Rightarrow \mathcal{F}(f):\mathcal{F}(A)\to \mathcal{F}(B).
$$
\end{Remark}

\section{2-Natural transformation of sesqui-functors}

\begin{Definition} [strict sesqui-transformations]
Let two parallel sesqui-functors
$$
\mathcal{F,G}:\mathcal{C}\to\mathcal{D}
$$
be given, and let be given a 2-graph transformation $\Delta:\mathcal{F}\Rightarrow\mathcal{G}$ whose
underlying 1-transformation
$$
\lfloor\Delta\rfloor:\lfloor\mathcal{F}\rfloor\Rightarrow\lfloor\mathcal{G}\rfloor
$$
is a natural transformation of functors. Then $\Delta$ is a (strict) natural transformation
of sesqui-functors when, for every $\alpha:f\Rightarrow g:A\to B$ in $\mathcal{C}$,
$$
\mathcal{F}(\alpha) \bullet_R \Delta_B = \Delta_A\bullet_L \mathcal{G}(\alpha)
$$
$$
\xymatrix@C=15ex@R=12ex{
\mathcal{F}(A)\ar[r]^{\Delta_A}
\ar@/_4ex/[d]_{\mathcal{F}(f)}="a1"
\ar@/^4ex/[d]^{\mathcal{F}(g)}="a2"
&
\mathcal{G}(A)
\ar@/_4ex/[d]_{\mathcal{G}(f)}="b1"
\ar@/^4ex/[d]^{\mathcal{G}(g)}="b2"
\\
\mathcal{F}(B)\ar[r]_{\Delta_B}
&
\mathcal{G}(B)
\ar@{}"a1";"a2"|(.4){}="ax1"|(.6){}="ax2"
\ar@{}"b1";"b2"|(.4){}="bx1"|(.6){}="bx2"
\ar@{=>}"ax1";"ax2"^{\mathcal{F}(\alpha)}
\ar@{=>}"bx1";"bx2"^{\mathcal{G}(\alpha)}
}
$$
\end{Definition}
Notice that while vertical composition of (strict) natural transformation of sesqui-functors
can be easily defined, the same is not true for horizontal composition. Therefore the category $\mathrm{Sesqui}\mathbf{CAT}$
of sesqui-categories, regardless of size issues, is indeed a sesqui-category itself.\\

The notion of (strict) natural transformation of sesqui-functors is essentially of a categorical nature. Namely
the ``functor''
$$
\lfloor-\rfloor: \mathrm{Sesqui}\mathbf{CAT}\to\mathbf{CAT}
$$
is also a ``sesqui-functor'', when we consider the 2-category $\mathbf{CAT}$ as a sesqui-category.

Therefore those are just usual natural transformations that behave nice with respect to reduced left and right compositions.
For the same reason the notions of adjunction and  equivalence of sesqui-categories (w.r.t. strict transformations)
are straightforward generalization of their categorical analogues.\\

Generalizing sesqui-categories (Chapter \ref{cha:three}) we will need a further  notion of sesqui-transformation,
whose definition follows

\begin{Definition}[lax sesqui-transformations]\label{def:lax_sesqui_transformation}
Let two parallel sesqui-functors
$$
\mathcal{F,G}:\mathcal{C}\to\mathcal{D}
$$
be given, and let be given a 2-graph transformation $\Gamma:\mathcal{F}\Rightarrow\mathcal{G}$.

 Then a  lax natural transformation $\Gamma:\ \mathcal{G}\Rightarrow\mathcal{G}$ is given by the following data:\\

$\bullet$ For every object $A$ of $\mathcal{C}$, an arrow
$$
\Gamma_{A}:\ \mathcal{F}(A)\rightarrow\mathcal{G}(A)
$$
\vskip4ex

$\bullet$ {\em ( naturality w.r.t. 1-cells)} For every arrow $f:\ A\to B$ of $\mathcal{C}$, a 2-cell

$$
\Gamma_f:\ \Gamma_A\bullet \mathcal{G}(f)\Rightarrow \mathcal{F}(f)\bullet \Gamma_B
$$
$$
\xymatrix@C=15ex@R=12ex{
\mathcal{F}(A)\ar[r]^{\Gamma_A}
\ar[d]_{\mathcal{F}(f)}
&
\mathcal{G}(A)
\ar[d]^{\mathcal{G}(f)}
\ar@{}[dl]|(.3){}="1"|(.7){}="2"
\ar@2"1";"2"_{\Gamma_f}
\\
\mathcal{F}(B)\ar[r]_{\Gamma_B}
&
\mathcal{G}(B)
}
$$
\vskip4ex

$\bullet$ {\em (naturality w.r.t. 2-cells)} For every 2-cell $\alpha:f\Rightarrow g:A\to B$ in $\mathcal{C}$, an equation
$$
\xymatrix@C=12ex@R=12ex{
\Gamma_{A}\bullet \mathcal{G}(f)
\ar@2[r]^{\Gamma_{A}\bullet_L \mathcal{G}(\alpha)}
\ar@2[d]_{\Gamma_{f}}
&
\Gamma_{A}\bullet \mathcal{G}(g)
\ar@2[d]^{\Gamma_g}
\\
\mathcal{F}(f)\bullet \Gamma_B
\ar@2[r]_{\mathcal{F}(\alpha)\bullet_R}
&
\mathcal{F}(g)\bullet \Gamma_B
}
$$
\vskip4ex

Those data have to satisfy the following {\em functoriality} axioms:\\

$\bullet$ For every object $A$ of $\mathcal{C}$
$$
\raisebox{8ex}{\xymatrix@C=13ex@R=12ex{
\mathcal{F}(A)\ar[r]^{\Gamma_A}
\ar[d]_{\mathcal{F}(1_A)}
&
\mathcal{G}(A)
\ar[d]^{\mathcal{G}(1_A)}
\ar@{}[dl]|(.3){}="1"|(.7){}="2"
\ar@2"1";"2"_{\Gamma_{1_A}}
\\
\mathcal{F}(A)\ar[r]_{\Gamma_A}
&
\mathcal{G}(A)
}}
\ =\
\raisebox{8ex}{
\xymatrix@C=13ex@R=12ex{
\mathcal{F}(A)\ar[r]^{\Gamma_A}
\ar@{-}[d]_{1_{\mathcal{F}(A)}}
&
\mathcal{G}(A)
\ar@{-}[d]^{1_{\mathcal{G}(A)}}
\ar@{}[dl]|(.3){}="1"|(.7){}="2"
\ar@{=}"1";"2"_{1_{\Gamma_{A}}}
\\
\mathcal{F}(A)\ar[r]_{\Gamma_A}
&
\mathcal{G}(A)
}}
$$

$\bullet$ For every composable pair $\xymatrix{A\ar[r]^f&B\ar[r]^h&C}$ in $\mathcal{C}$
$$
\raisebox{10ex}{\xymatrix@C=13ex@R=6ex{
\mathcal{F}(A)\ar[r]^{\Gamma_A}
\ar[d]_{\mathcal{F}(f)}
&
\mathcal{G}(A)
\ar[d]^{\mathcal{G}(f)}
\ar@{}[dl]|(.3){}="a1"|(.7){}="a2"
\ar@2"a1";"a2"_{\Gamma_f}
\\
\mathcal{F}(B)\ar[r]|{\Gamma_B}
\ar[d]_{\mathcal{F}(h)}
&
\mathcal{G}(B)
\ar[d]^{\mathcal{G}(h)}
\ar@{}[dl]|(.3){}="b1"|(.7){}="b2"
\ar@2"b1";"b2"_{\Gamma_h}
\\
\mathcal{F}(C)\ar[r]_{\Gamma_C}
&
\mathcal{G}(C)
}}
\ =\
\raisebox{10ex}{
\xymatrix@C=13ex@R=16ex{
\mathcal{F}(A)\ar[r]^{\Gamma_A}
\ar[d]_{\mathcal{F}(fh)}
&
\mathcal{G}(A)
\ar[d]^{\mathcal{G}(fh)}
\ar@{}[dl]|(.3){}="1"|(.7){}="2"
\ar@2"1";"2"_{\Gamma_{fh}}
\\
\mathcal{F}(C)\ar[r]_{\Gamma_C}
&
\mathcal{G}(C)
}
}
$$
\end{Definition}
\begin{Remark}
In general, a lax sesqui-transformation is \emph{not} a natural transformation of the functors underlying domain and co-domain
 sesqui-functors.
\end{Remark}

\section{Sesqui-categories and 2-categories}

That a sesqui-category induces a category structure on the underlying graph is clear from the very definition of
sesqui-categories.

Hence, the question that naturally arises concerns \emph{when} a sesqui-category is also a 2-category.
In fact, given a sesqui-category $\mathcal{C}$, this underlies a 2-category precisely
when, for every diagram of the kind
$$
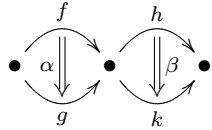
\xymatrix{\bullet\ar@/^3ex/[r]^f_{}="1" \ar@/_3ex/[r]_g_{}="2" & \bullet\ar@/^3ex/[r]^h_{}="3"\ar@/_3ex/[r]_k^{}="4"&\bullet
\ar@{=>}"1";"2"_{\alpha}\ar@{=>}"3";"4"^{\beta}}
$$
the following equation is satisfied:
\begin{equation}\label{eqn:1-comp}
(f\bullet_L\beta)\cdot(\alpha\bullet_R k)=(\alpha\bullet_R h)\cdot(g\bullet_L \beta)
\end{equation}
In this situation, the two composites are denoted $\alpha\bullet\beta$, and termed \emph{horizontal composition} of
 $\alpha$ and $\beta$.

In other terms, it is possible to show that such a composition defines a family of functors
$$
\xymatrix@C=16ex{\mathcal{C}(A,B)\times\mathcal{C}(B,C)\ar[r]^(.6){\bullet^{A,B,C}}&\mathcal{C}(A,C)}
$$
indexed by triples $(A,B,C)$ of objects of  $\mathcal{C}$, satisfying 2-categorical axioms.

There follows an interchange law for horizontal and vertical composition holds: for every four 2-cells
$\alpha, \beta, \gamma, \delta$
$$
\xymatrix{
\bullet \ar@/^4ex/[r]_{}="a1"\ar[r]^{}="a2"^{}="a3"\ar@/_4ex/[r]^{}="a4"&\bullet
\ar@/^4ex/[r]_{}="b1"\ar[r]^{}="b2"^{}="b3"\ar@/_4ex/[r]^{}="b4"&\bullet
\ar@{=>}"a1";"a2"_{\alpha}
\ar@{=>}"a3";"a4"_{\gamma}
\ar@{=>}"b1";"b2"^{\beta}
\ar@{=>}"b3";"b4"^{\delta}
}
\quad
(\alpha\cdot\gamma)\bullet(\beta\cdot\delta)=(\alpha\bullet\beta)\cdot(\gamma\bullet\delta)
$$

Even when equation (\ref{eqn:1-comp}) does not hold, some pasting operations of 2-cells are
 still available. We show this  with an example.

Consider the diagram:
$$
\xymatrix{
\bullet\ar[r]^a \ar[d]_d
&\bullet\ar[r]^b \ar[d]_g\ar@{}[dl]|(.3){}="a1"|(.7){}="a2"
&\bullet\ar[d]^c \ar@{}[dl]|(.3){}="b1"|(.7){}="b2"
\\
\bullet\ar[r]_e
&\bullet\ar[r]_f
&\bullet
\ar@{=>}"a1";"a2"_{\alpha}
\ar@{=>}"b1";"b2"_{\beta}
}
$$
Since intersection between $\alpha$ and $\beta$ is one dimensional, it is unambiguous to define the \emph{pasting}
$$
(\alpha|\beta)= (a\bullet \beta)\cdot(\alpha\bullet f)
$$
In the present work, we will not go any further into this subject.

\section{Finite products in a sesqui-category}

In the sesqui-categorical context we will refer to binary products according to the following 2-dimensional
universal property
\begin{Definition}\label{def:sesqui-product}\label{UP:products}
Let $\mathcal{C}$ be a sesqui-category, $A$ and $B$ two objects of $\mathcal{C}$. A product of $A$ and $B$
is a triple $(A\times B,\pi_A:A\times B\rightarrow A,\pi_B:A\times B\rightarrow B)$
satisfying  the following\\

\textbf{Universal property}\\

For every object $Q$ of $\mathcal{C}$ and  2-cells
$$
\xymatrix@C=3ex{\alpha:a\ar@{=>}[r]&a':Q\ar[r] &A},\quad
\xymatrix@C=3ex{\beta :b\ar@{=>}[r]&b':Q\ar[r] &B}
$$
there exists a unique 2-cell
$$
\xymatrix@C=3ex{\gamma:q\ar@{=>}[r]&q':Q\ar[r] & A\times B}
$$
with $\gamma \bullet \pi_A = \alpha$ and $\gamma \bullet \pi_B = \beta$.\\

We will write $\gamma=\langle\alpha,\beta\rangle$
\end{Definition}
The situation may be visualized on the diagram below
$$
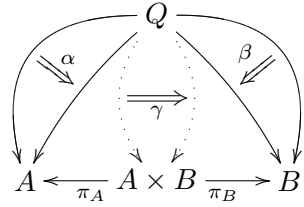
\xymatrix@R=10ex{
&Q\ar@/_6ex/[dl]_{}="a1"\ar@/_1ex/[dl]^{}="a2"
\ar@/_3ex/@{.>}[d]^{}="q1"\ar@/^3ex/@{.>}[d]_{}="q2"
\ar@/^6ex/[dr]_{}="b1"\ar@/^1ex/[dr]^{}="b2"
\\
A&\ar[l]^{\pi_A} A\times B\ar[r]_{\pi_B}&B
\ar@{=>}"a1";"a2"^{\alpha}
\ar@{=>}"q1";"q2"_{\gamma}
\ar@{=>}"b1";"b2"_{\beta}
}
$$
Such a product satisfies also the universal property defining categorical products. It suffices to choose
$\alpha=1_a$ and $\beta=1_b$: the unique $\gamma$ of the property satisfies $\gamma\bullet\pi_A=1_a$
and $\gamma\bullet\pi_B=1_b$. Hence, taking domains and codomains, we get
$$
q\pi_A=a, \quad q\pi_B=b, \quad q'\pi_A=a, \quad q'\pi_B=b
$$
This gives
\begin{eqnarray*}
  1_q\bullet\pi_A &=& 1_{q\pi_A }\mathrm{\ by\ axiom\ (R3)}\\
    &=& 1_a\\
  1_q\bullet\pi_B &=& 1_{q\pi_B }\mathrm{\ by\ axiom\ (L3)}\\
    &=& 1_b
\end{eqnarray*}
Finally, uniqueness forces $\gamma=1_q$, and in turn, there exists a unique $q\ (=q')$ such that $q\pi_A=a$ and $q\pi_B=b$.

\begin{Definition}\label{def:sesqui-terminal}
Let $\mathcal{C}$ be a sesqui-category. A terminal object is an object $I$ of $\mathcal{C}$ satisfying the
following universal property\\

\emph{(UP)} for every other object $X$ of $\mathcal{C}$, there exists a unique 2-cell
$$
\xi:x\Rightarrow x':X\rightarrow I
$$
\end{Definition}
With a calculation similar to that of products, this universal property is equivalent to the existence of a unique $!_X:X\to I$,
henceforth $\xi$ is indeed the identity 2-cell on $!_{X}$.\\

Products and terminals defined  this way are determined up to isomorphism.
Furthermore finite products  and canonical isomorphisms are defined as in the categorical case.

\section{Product interchange rules}\label{Product interchange rules}
In the previous section we were concerned with properties of products in a sesqui-category that specialize in classical
(viz. categorical)  ones.\\

Now we focus our attention on products of 2-cells. What we recapture is the idea of independence of the components of a product, and a
sort of commutativity that arises.\\

Consider the 2-cells $\alpha:f\Rightarrow g:A\rightarrow B$ and $\beta:h\Rightarrow k:C\rightarrow D$ in a sesqui-category $\mathcal{C}$.
A 2-cell
$$
\alpha\times\beta:f\times h\Rightarrow g\times k:A\times C\rightarrow B\times D
$$
is uniquely determined by the universal property: $\alpha\times\beta=\langle\pi_A\bullet\alpha,\pi_C\bullet\beta\rangle$.

Notice that this induces a kind of commutative horizontal composition of 2-cells, provided they are on different product-components.

In fact, we need the following
\begin{Lemma}\label{lemma:prod_interchange}
For $\alpha$ and $\beta$ as above,
$$
((1_A\tm \beta)\bullet (f\tm 1_D))\cdot((1_A\tm k)\bullet(\alpha\tm 1_D))=
((1_A\tm h)\bullet(\alpha\tm 1_D))\cdot((1_A\tm\beta)\bullet( g\tm 1_D))
$$
\end{Lemma}
\begin{proof} By universality of products, they are both equal to
$\alpha\times\beta$, because they have the same composite
with projections. In fact we prove just the left side, the right
side being analogous.
\begin{eqnarray*}
  & &\!\!\!\!\!\!\!\!\!\!\!\!\!((1_A\tm \beta)\bullet (f\tm 1_D))\cdot((1_A\tm k)\bullet(\alpha\tm 1_D))\bullet \pi_D =\\
  &=& ((1_A\tm \beta)\bullet (f\tm 1_D)\bullet\pi_D)\cdot(((1_A\tm k)\bullet(\alpha\tm 1_D))\bullet\pi_D) \\
  &=& ((1_A\tm \beta) \bullet\pi_D)\cdot((1_A\tm k)\bullet(\alpha\tm 1_D)\bullet\pi_D) \\
  &=& ((1_A\tm \beta) \bullet\pi_D)\cdot((1_A\tm k)\bullet\pi_D) \\
  &=&  (\pi_C\bullet\beta) \cdot (\pi_C\, k) =\pi_C\bullet\beta
\end{eqnarray*}
\begin{eqnarray*}
  & &\!\!\!\!\!\!\!\!\!\!\!\!\!((1_A\tm \beta)\bullet (f\tm 1_D))\cdot((1_A\tm k)\bullet(\alpha\tm 1_D))\bullet \pi_B =\\
  &=& ((1_A\tm \beta)\bullet (f\tm 1_D)\bullet\pi_B)\cdot(((1_A\tm k)\bullet(\alpha\tm 1_D))\bullet\pi_B) \\
  &=& ((1_A\tm \beta) \bullet\pi_A\,f)\cdot((1_A\tm k)\bullet(\alpha\tm 1_D)\bullet\pi_B) \\
  &=& (\pi_A\,f)\cdot((1_A\tm k)\pi_A\bullet\alpha) \\
  &=& (\pi_A\,f) \cdot (\pi_A\bullet \alpha)=\pi_A\bullet \alpha
\end{eqnarray*}
\end{proof}
\newpage
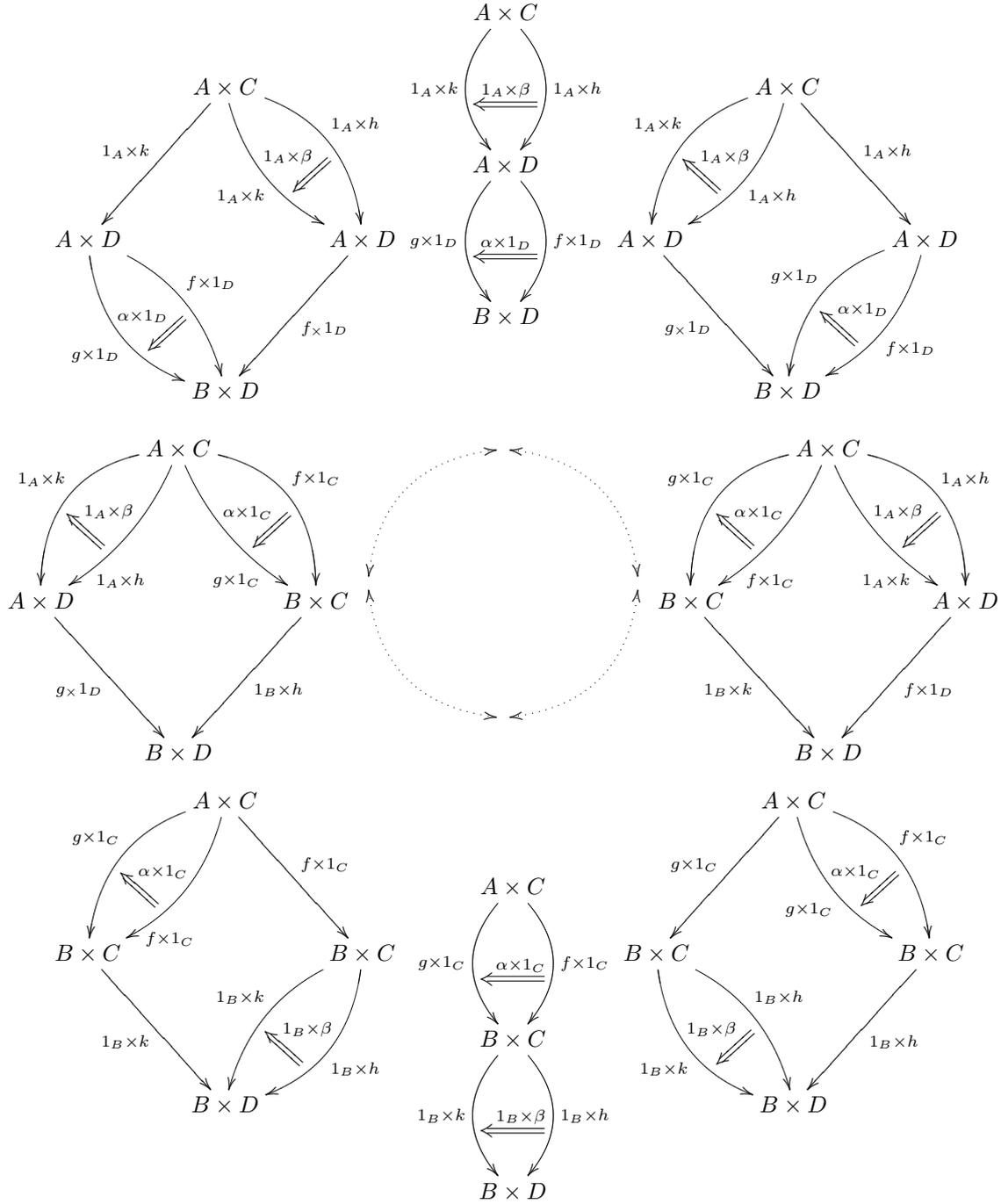
\begin{changemargin}{-10ex}{8ex}
{\small
$$
\qquad\xymatrix@R=12ex{
&A\times C
\ar@/_3ex/[dr]_{1_A\times k}^(.6){}="b2"
\ar@/^4ex/[dr]^(.5){1_A\times h}_(.6){}="b1"
\ar[dl]_{1_A\times k}
\\
A\times D
\ar@/_4ex/[dr]_{g\times 1_D}^(.6){}="a2"
\ar@/^3ex/[dr]^(.5){f\times 1_D}_(.6){}="a1"
&&
A\times D
\ar[dl]^{f_\times 1_D}
\\
& B\times D
\ar@{=>}"a1";"a2"_(.4){\alpha\times 1_D}
\ar@{=>}"b1";"b2"_(.4){1_A\times \beta}
}
\raisebox{8ex}{\xymatrix@R=12ex{
A\times C
\ar@/_4ex/[d]_{1_A\times k}^(.6){}="b2"
\ar@/^4ex/[d]^(.5){1_A\times h}_(.6){}="b1"
\\A\times D
\ar@/_4ex/[d]_{g\times 1_D}^(.6){}="a2"
\ar@/^4ex/[d]^(.5){f\times 1_D}_(.6){}="a1"
\\B\times D
\ar@{=>}"a1";"a2"_(.5){\alpha\times 1_D}
\ar@{=>}"b1";"b2"_(.5){1_A\times \beta}
}}
\xymatrix@R=12ex{
&A\times C
\ar@/_4ex/[dl]_{1_A\times k}^(.6){}="b2"
\ar@/^3ex/[dl]^(.5){1_A\times h}_(.6){}="b1"
\ar[dr]^{1_A\times h}
\\
A\times D
\ar[dr]_{g_\times 1_D}
&&
A\times D
\ar@/_4ex/[dl]_{g\times 1_D}^(.6){}="a2"
\ar@/^3ex/[dl]^(.5){f\times 1_D}_(.6){}="a1"
\\
& B\times D
\ar@{=>}"a1";"a2"_(.6){\alpha\times 1_D}
\ar@{=>}"b1";"b2"_(.6){1_A\times \beta}
}
$$
$$
\xymatrix@R=12ex{
&A\times C
\ar@/_5ex/[dl]_{1_A\times k}^(.6){}="b2"
\ar@/^2ex/[dl]^(.7){1_A\times h}_(.6){}="b1"
\ar@/^5ex/[dr]^{f\times 1_C}^(.6){}="a1"
\ar@/_2ex/[dr]_(.7){g\times 1_C}_(.6){}="a2"\\
A\times D
\ar[dr]_{g_\times 1_D}
&&
B\times C
\ar[dl]^{1_B\times h}
\\
& B\times D
\ar@{=>}"a1";"a2"_(.4){\alpha\times 1_C}
\ar@{=>}"b1";"b2"_(.6){1_A\times \beta}}
\ \xymatrix@R=12ex@C=12ex{
&{}\ar@/_4ex/@{<.>}[dl]\ar@/^4ex/@{<.>}[dr]\\
{}\ar@/_4ex/@{<.>}[dr]&&{}\ar@/^4ex/@{<.>}[dl]\\
&{}}
\ \xymatrix@R=12ex{
&A\times C
\ar@/_2ex/[dr]_(.7){1_A\times k}^(.6){}="b2"
\ar@/^5ex/[dr]^{1_A\times h}_(.6){}="b1"
\ar@/^2ex/[dl]^(.7){f\times 1_C}^(.6){}="a1"
\ar@/_5ex/[dl]_{g\times 1_C}_(.6){}="a2"\\
B\times C
\ar[dr]_{1_B\times k}
&&
A\times D
\ar[dl]^{f\times 1_D}
\\
& B\times D
\ar@{=>}"a1";"a2"_(.6){\alpha\times 1_C}
\ar@{=>}"b1";"b2"_(.4){1_A\times \beta}}
$$
$$
\qquad\xymatrix@R=12ex{
&A\times C
\ar@/^3ex/[dl]^(.7){f\times 1_C}^(.6){}="a1"
\ar@/_4ex/[dl]_{g\times 1_C}_(.6){}="a2"
\ar[dr]^{f\times 1_C}
\\
B\times C
\ar[dr]_{1_B \times k}
&&
B\times C
\ar@/_3ex/[dl]_{1_B\times k}^(.6){}="b2"
\ar@/^4ex/[dl]^(.5){1_B\times h}_(.6){}="b1"
\\
& B\times D
\ar@{=>}"a1";"a2"_(.6){\alpha\times 1_C}
\ar@{=>}"b1";"b2"_(.6){1_B\times \beta}
}
\raisebox{-8ex}{
\xymatrix@R=12ex{
A\times C
\ar@/_4ex/[d]_{g\times 1_C}^(.6){}="b2"
\ar@/^4ex/[d]^(.5){f\times 1_C}_(.6){}="b1"
\\B\times C
\ar@/_4ex/[d]_{1_B\times k}^(.6){}="a2"
\ar@/^4ex/[d]^(.5){1_B\times h}_(.6){}="a1"
\\B\times D
\ar@{=>}"b1";"b2"_(.4){\alpha\times 1_C}
\ar@{=>}"a1";"a2"_(.4){1_B\times \beta}
}}
\xymatrix@R=12ex{
&A\times C
\ar@/_3ex/[dr]_{g\times 1_C}^(.6){}="b2"
\ar@/^4ex/[dr]^(.5){f\times 1_C}_(.6){}="b1"
\ar[dl]_{g\times 1_C}
\\
B\times C
\ar@/_4ex/[dr]_{1_B\times k}^(.6){}="a2"
\ar@/^3ex/[dr]^(.5){1_B\times h}_(.6){}="a1"
&&
B\times C
\ar[dl]^{1_B\times h}
\\
& B\times D
\ar@{=>}"b1";"b2"_(.4){\alpha\times 1_C}
\ar@{=>}"a1";"a2"_(.4){1_B\times \beta}
}
$$
}
\end{changemargin}
\emph{Lemma \ref{lemma:prod_interchange}} allows us to define a
horizontal composition of this kind of 2-cells
$$
(1_A\times \beta) \bullet (\alpha\times1_D) =\alpha\times\beta
=(\alpha\times 1_C)\bullet(1_B\times\beta)
$$
and to prove diagram equalities, such as the one above. These kind of
diagrammatic equations will be called \emph{product interchange rules}.

\section{$h$-Pullbacks}

We introduce here a notion of standard \textit{h}-pullback suitable for our purposes. This notion has been formalized by
 Michael Mather in \cite{MR0402694}, for generic categories of spaces, with (eventually pointed) topological spaces in mind.
It has been further generalized to \textit{h}-categories\footnote{A \textit{h}-category is a weaker notion than that of
a sesqui-category, see \cite{grandis1994hah}.} by Marco Grandis in \cite{grandis1994hah}. We (ab)use the term $h$-pullback,
instead of that of \emph{comma-square} because we will work mainly in a n-groupoidal context, with 2-morphisms being weakly
invertible.\\

\begin{Definition}\label{def:h-pullbacks}\label{UP:h-pullbacks}
Consider the following diagram in a sesqui-category $\mathcal{C}$
$$
\xymatrix{&C\ar[d]^g\\A\ar[r]_f&B}
$$
An $h$-pullback of $f$ and $g$ is a four-tuple $(P(f,g), p,q,\varepsilon)$
$$
\xymatrix{
P\ar[r]^q\ar[d]_p
&C\ar[d]^g
\\
A\ar[r]_f\ar@{}[ur]|(.3){}="1"|(.7){}="2"
&B
\ar@{=>}"1";"2"^{\varepsilon}}
$$
where $P=P(f,g)$, that satisfies the following\\

\textbf{Universal Property}

For any other four-tuple $(X, m,n,\lambda)$ as in
$$
\xymatrix{
X\ar[r]^n\ar[d]_m
&C\ar[d]^g
\\
A\ar[r]_f\ar@{}[ur]|(.3){}="1"|(.7){}="2"
&B
\ar@{=>}"1";"2"^{\lambda}}
$$

there exists a unique $\ell:X\rightarrow P$ such that
\begin{enumerate}
  \item $\ell p=m$
  \item $\ell q=n$
  \item $\ell\bullet_L \varepsilon=\lambda$
\end{enumerate}

\end{Definition}

\begin{Lemma}
{\em Universal Property \ref{UP:h-pullbacks}} defines $h$-pullbacks up to isomorphisms.
\end{Lemma}
\begin{proof}
Let $(P,p,q,\varepsilon)$ be a $h$-pullback, according to definition above, and let $(P',p',q',\varepsilon')$ be another
four-tuple satisfying the universal property. Then, since the first is a $h$-pullback, there exists $\ell:P'\to P$
\begin{equation*}
    \ell p=p',\quad \ell q=q',\quad \ell\varepsilon = \varepsilon'
\end{equation*}
and, since the second is a $h$-pullback, there exists $\ell':P\to P'$
\begin{equation*}
    \ell'p'=p,\quad \ell'q'=q,\quad \ell'\varepsilon' = \varepsilon
\end{equation*}
Now, applying the universal property of the first one to itself, $\ell'\ell$ and $id_{P}$ satisfy the same equations, and by uniqueness are equal. Similarly
applying the universal property of the second one to itself, $\ell \ell'$ and $id_{P'}$ satisfy the
 same equations, and by uniqueness they are equal too.
Hence $\ell$ and $\ell'$ are isomorphisms.
\end{proof}

\begin{Lemma}[Pullback of $h$-projections.]\label{lemma:h-pb_h-proj}
 In the sesqui-category $\mathcal{C}$, let be given the diagram below, where the left-hand
 square is commutative and the right-hand square $\varepsilon$ is a $h$-pullback
$$
\xymatrix{
R\ar[r]^s\ar[d]_r
&P\ar[r]^q\ar[d]_p
&D\ar[d]^g
\\
A\ar[r]_e
&B\ar[r]_f\ar@{}[ur]|(.3){}="1"|(.7){}="2"\ar@2"1";"2"_{\varepsilon}
&C}
$$
then the composition $s\bullet_L \varepsilon$ is a $h$-pullback if, and only if, the left hand square is a pullback.
\end{Lemma}
\begin{proof}
We will show that the four-tuple $(R,r,sq,s\bullet_L \varepsilon)$ satisfies the universal property \ref{UP:2h-pullbacks}.

Let the four-tuple $(\X,y,z,\xi)$ be given as in the diagram below
$$
\xymatrix{
X\ar[rr]^z\ar[d]_y
&&D\ar[d]^g
\\
A\ar@{}[urr]|(.3){}="1"|(.7){}="2"\ar@2"1";"2"_{\xi}
\ar[r]_e
&B\ar[r]_f
&C}
$$
Since $P$ is an $h$-pullback, there exists a unique $\ell:X\to P$ such that
$$
(i)\ \ell p=ye,\quad (ii)\ \ell q=z, \quad (iii)\ \ell\bullet_L \varepsilon = \xi.
$$
Yet since $R$ is a pullback, condition $(i)$ is equivalent to:\\

there exists a unique
$x:X\to R$ such that
$$
(iv)\ xr=y,\quad (v)\ xs=\ell.
$$
Substituting, there exists a unique $x:X\to R$ such that
\begin{eqnarray*}
  &(i)'& xr\eq{iv}y     \\
  &(ii)'& xsq\eq{v}\ell q\eq{ii} z \\
  &(iii)'&x\bullet_L(s\bullet_L \varepsilon)\eq{L2} xs\bullet_L\varepsilon \eq{v}\ell\bullet_L \varepsilon\eq{iii}\xi
\end{eqnarray*}
\end{proof}
\begin{Remark}
This Lemma still holds in a mere $h$-category (\cite{grandis1994hah} {\em Lemma 2.2.}).
\end{Remark}

\chapter{Strict $n$-categories}\label{cha:strictncat}

We give an inductive definition of the sesqui-category $n$\textbf{Cat},
whose objects are (strict and small) $n$-categories, morphisms are
(strict) $n$-functors and $2$-morphisms are (lax)
$n$-transformations. Furthermore, $n$\textbf{Cat} has sesqui-categorical finite products.\\

In the next three sections, we recall a standard inductive construction of $n$\textbf{Cat},
well known in literature, recalled for instance in \cite{MR920944}.
This is in fact a notion based upon a more general and influential theory of enrichment,
developed by Gregory Maxwell Kelly  (see \cite{MR2177301}).

This new perspective arises in the sesqui-categorical structure
developed thereafter.\\

More precisely our notion of lax natural $n$-transformation generalizes the  ordinary 2-dimensional
version, as  recalled for example in \cite{B94HANDBOOK}.
This coincides with the inductive definition given in the internal abelian case in \cite{MR1075339}.  Moreover it is equivalent to the global definition  given in \cite{Crans95} (see
{\em Definition 9.1}, {\em Lemma 9.2} for a comparison), closely related to that of {\em$m$-fold homotopies}
 of $\omega$-groupoids in \cite{MR906402}.

For $n=0$, $n$\textbf{Cat} is \textbf{Set}, the sesqui-category of
sets and maps with trivial (\textit{i.e.} identity)
transformations. Cartesian product provides the required sesqui-categorical product.\\

For $n=1$, the $2$-category \textbf{Cat} of categories, functors and
natural transformations has a underlying sesqui-categorical
structure, when we consider only reduced horizontal composition of
natural transformations with
functors. Categorical product gives again the required sesqui-categorical product.\\

\section{$n$\textbf{Cat}: the data}
\subsection*{$n$-categories}

For given integer $n>1$, a (strict) n-category $\C$ consists of the following data:\\

$\bullet$ a set of objects $\C_0$;\\

$\bullet$ for every pair of objects $c_0, c_0'$ of $\C_0$, a
$(n-1)$category
$$
\C_1(c_0,c_0')
$$
called \textit{hom$(n-1)$category} over $c_0$ and $c_0'$, and sometimes written $[c_0,c_0']$ in order
to simplify notation;\\

$\bullet$ for every object $c_0$ of $\C_0$, a morphism of
$(n-1)$categories
$$
\xymatrix@C=10ex{\Id{n-1}\ar[r]^(.35){u^0(c_0)}&\C_1(c_0,c_0)}
$$
called the \textit{$0$-identity} of $c_0$;

$\bullet$ for every triple of objects $c_0, c_0', c_0''$ of $\C_0$,
a morphism of $(n-1)$categories
$$
\xymatrix@C=12ex{\C_1(c_0,c_0')\x{n-1}
\C_1(c_0',c_0'')\ar[r]^(.6){\circ^0_{c_0, c_0',
c_0''}}&\C_1(c_0,c_0'')}
$$
called \textit{$0$-composition}, following the \textsl{dimension-intersection} convention.\\

Here, $\x{n-1}$ stays for the binary product in
$(n-1)$\textbf{Cat}, while $\Id{n-1}$ is the \mbox{$0$-ary} product
in
$(n-1)$\textbf{Cat}. Subscripts will be usually omitted, unless this causes confusion.

All these data must satisfy the following axioms, expressed by commutative diagrams in $(n-1)$\textbf{Cat}:\\

$\bullet$ (\emph{associativity axiom}) for every four-tuple of objects
$c_0, c_0', c_0'', c_0'''$ of $\C_0$ {\small
\begin{equation}\label{nCat:associativity}
\xymatrix@C=-2ex{ (\C_1(c_0,c_0')\times \C_1(c_0',c_0''))\times
\C_1(c_0'',c_0''') \ar[rr]^{\alpha}_{\sim}
\ar[d]_{\circ^0_{c_0,c_0',c_0''}\times id}
&& \C_1(c_0,c_0')\times (\C_1(c_0',c_0'')\times \C_1(c_0'',c_0'''))\ar[d]^{id \times \circ^0_{c_0',c_0''.c_0'''}}\\
\C_1(c_0,c_0'')\times
\C_1(c_0'',c_0''')\ar[dr]_{\circ^0_{c_0,c_0''.c_0'''}}
&& \C_1(c_0,c_0')\times \C_1(c_0',c_0''')\ar[dl]^{\circ^0_{c_0,c_0',c_0'''}}\\
&\C_1(c_0,c_0''')}
\end{equation}
} where
$$
\alpha=\alpha_{\C_1(c_0,c_0'), \C_1(c_0',c_0''), \C_1(c_0'',c_0''')}
$$
is the usual associator given by universal property of product;\\

$\bullet$ (\emph{left and right unit axioms}) for every pair of
objects $c_0, c_0'$ of $\C_0$ {\small
\begin{equation}\label{nCat:units}
\xymatrix@C=10ex{ \I\times \C_1(c_0,c_0')\ar[d]_{u^0(c_0)\times id}
&
\C_1(c_0,c_0')\ar[l]_{\lambda}^{\sim}\ar@{=}[d]|{id}\ar[r]^{\rho}_{\sim}
& \C_1(c_0,c_0')\times \I\ar[d]^{id \times {u^0(c_0)}}
\\
\C_1(c_0,c_0)\times
\C_1(c_0,c_0')\ar[r]_(.6){\circ^0_{c_0,c_0,c_0'}} &\C_1(c_0,c_0')
&\C_1(c_0,c_0')\times\C_1(c_0',c_0')\ar[l]^(.6){\circ^0_{c_0,c_0',c_0'}}}
\end{equation}
} where
$$
\lambda=\lambda_{\C_1(c_0,c_0')}\quad \mathrm{and}\quad
\rho=\rho_{\C_1(c_0,c_0')}
$$
are the usual left and right unit isomorphisms given by the universal property of the product.\\

\subsection*{Morphisms of $n$-categories}

For a given integer $n>1$, and given $n$-categories $\C$ and $\D$, a
(strict) n-functor
$$
\xymatrix{F:\C\ar[r]&\D}
$$
is a pair $(F_0,F_1)$ where:\\

$\bullet$ $\xymatrix{F_0:\C_0\ar[r]&\D_0}$ is a map;\\

$\bullet$ for every pair of objects $c_0, c_0'$ of $\C_0$
$$
\xymatrix{F_1^{c_0,c_0'}:\C_1(c_0,c_0')\ar[r]&\D_1(F_0c_0,F_0c_0')}
$$
is a morphism of $(n-1)$categories.\\

These data must satisfy the following axioms, expressed by commutative diagrams in $(n-1)$\textbf{Cat}:\\

$\bullet$ (\emph{functoriality w.r.t. composition}) for every triple
of objects $c_0, c_0', c_0''$ of $\C_0$
{\small\begin{equation}\label{nfunct:composition} \xymatrix@C=10ex{
\C_1(c_0,c_0')\times \C_1(c_0',c_0'')\ar[r]^{^{\C}\circ^0}
\ar[d]_{F_1^{c_0,c_0'}\times F_1^{c_0',c_0''}}
& \C_1(c_0,c_0'')\ar[d]^{F_1^{c_0,c_0''}}\\
\D_1(F_0c_0,F_0c_0')\times
\D_1(F_0c_0',F_0c_0'')\ar[r]_(.6){^{\D}\circ^0} & \D_1(F_0 c_0,F_0
c_0'') }
\end{equation}}
$\bullet$ (\emph{functoriality w.r.t. units}) for every triple of
object $c_0$ of $\C_0$ {\small\begin{equation}\label{nfunct:units}
    \xymatrix@C=14ex{
    \I \ar[r]^{^\C u^0(c_0)}\ar[dr]_{^\D u^0(F_0c_0)}
    & \C_1(c_0,c_0)\ar[d]^{F_1^{c_0,c_0}}\\
    & \D_1(F_0c_0,F_0c_0)}
\end{equation}}

\subsection*{2-Morphisms of $n$-categories}

For a given integer $n>1$, and given $n$-functors
$\xymatrix{F,G:\C\ar[r]&\D}$, a lax natural $n$-transformation
$$
{
\def \xa{\C}
\def \xb{F}
\def \xc{G}
\def \xd{\D}
\def \xe{\alpha}
\xymatrix@C=15ex{ \xa \ar@/^3ex/[r]^{\xb}_{}="x1"
\ar@/_3ex/[r]_{\xc}^{}="x2" & \xd
\ar@{}"x1";"x2"|(.2){}="xx1"|(.8){}="xx2" \ar@{=>}"xx1";"xx2"^{\xe}
} }
$$
is a pair $\alpha=(\alpha_0,\alpha_1)$, where\\

$\bullet$ $\xymatrix{\alpha_0:\C_0\ar[r]& {}\substack{\coprod
\\c_0\in \C_0} [\D_1(F_0c_0,G_0c_0)]_0}$ is a
 map such that, for every $c_0$ in $\C_0$,
$\xymatrix{\alpha_0(c_0):F_0(c_0)\ar[r]&G_0(c_0)}$;\\

$\bullet$ (\emph{n-naturality}) for every pair $c_0$, $c_0'$ of $\C$,
$\alpha_1^{c_0,c_0'}$ is a 2-morphism of $(n-1)\mathbf{Cat}$, as in
the diagram below
$$
{
\def \xxa{\C_1(c_0,c_0')}
\def \xxb{\D_1(F_0c_0,F_0c_0')}
\def \xxc{\D_1(G_0c_0,G_0c_0')}
\def \xxd{\D_1(F_0c_0,G_0c_0')}
\def \xxe{F_1^{c_0,c_0'}}
\def \xxf{G_1^{c_0,c_0'}}
\def \xxg{- \circ^0 \alpha_0c_0'}
\def \xxh{\alpha_0c_0 \circ^0 -}
\def \xxi{\alpha_1^{c_0,c_0'}}
\xymatrix@C=-2ex{
& \xxa \ar[dl]_{\xxe}\ar[dr]^{\xxf}\\
\xxb\ar@{}[rr]|(.4){}="xxb"|(.6){}="xxc"\ar[dr]_{\xxg}
&&\xxc\ar[dl]^{\xxh}\\
&\xxd \ar@{=>}"xxc";"xxb"_{\xxi}}}
$$

In order to keep notation lighter we will often write $\alpha_{c_0}$
instead of $\alpha_0(c_0)$.

These data must satisfy functoriality axioms expressed by the following equations of diagrams in $(n-1)$\textbf{Cat}:\\

$\bullet$ (\emph{functoriality w.r.t. composition}) for every triple of
objects $c_0, c_0', c_0''$ of $\C_0$,
\begin{changemargin}{-10ex}{-10ex}\small
\begin{equation}\label{ntrans:composition}
\def \xya{\C_1(c_0,c_0')\times\C_1(c_0',c_0'')}
\def \xyb{\C_1(c_0,c_0')\times\D_1(F_0c_0',F_0c_0'')}
\def \xyc{\D_1(G_0c_0,G_0c_0')\times\C_1(c_0',c_0'')}
\def \xyd{\C_1(c_0,c_0')\times\D_1(G_0c_0',G_0c_0'')}
\def \xye{\D_1(F_0c_0,F_0c_0')\times\C_1(c_0',c_0'')}
\def \xyf{\C_1(c_0,c_0')\times\D_1(F_0c_0',G_0c_0'')}
\def \xyg{\D_1(F_0c_0,G_0c_0')\times\C_1(c_0',c_0'')}
\def \xyh{\D_1(F_0c_0,F_0c_0')\times\D_1(F_0c_0',G_0c_0'')}
\def \xyi{\D_1(F_0c_0,G_0c_0')\times\D_1(G_0c_0',G_0c_0'')}
\def \xyj{id\times \alpha_1^{c_0',c_0''}}
\def \xyk{\alpha_1^{c_0,c_0'}\times id}
\def \xyl{\D_1(F_0c_0,G_0c_0'')}
\def \xym{id\times F_1^{c_0',c_0''}}
\def \xyn{G_1^{c_0,c_0'}\times id}
\def \xyo{id\times(-\circ\alpha_0c_0'')}
\def \xyp{id\times G_1^{c_0',c_0''}}
\def \xyq{F_1^{c_0,c_0'}\times id}
\def \xyr{(\alpha_0c_0\circ -)\times id}
\def \xys{id\times(\alpha_0c_0'\circ -)}
\def \xyt{(-\circ\alpha_0c_0')\times id}
\def \xyu{F_1^{c_0,c_0'}\times id}
\def \xyv{id\times G_1^{c_0',c_0''}}
\def \xyw{\circ^0}
\def \xyz{\circ^0}
\xymatrix@C=-12ex{
&&\xya\ar[dll]_{\xym}\ar[ddl]|{\xyp}\ar[ddr]|{\xyq}\ar[drr]^{\xyn}\ar@{}[ddddd]|{\equiv}\\
\xyb\ar[dd]_{\xyo}\ar@{}[dr]|(.35){}="2"|(.65){}="1"&&&&\xyc\ar[dd]^{\xyr}\ar@{}[dl]|(.35){}="3"|(.65){}="4"\\
&\xyd\ar[dl]^{\xys}&&\xye\ar[dr]_{\xyt}\\
\xyf\ar[d]_{\xyu}&&&&\xyg\ar[d]^{\xyv}\\
\xyh\ar[drr]_{\xyw}&&&&\xyi\ar[dll]^{\xyz}\\
&&\xyl \ar@{=>} "1";"2"_{\xyj}\ar@{=>} "3";"4"_{\xyk} }
\end{equation}
\end{changemargin}
{ \small
$$
\def \xxa{\C_1(c_0,c_0'')}
\def \xxb{\D_1(F_0c_0,F_0c_0'')}
\def \xxc{\D_1(G_0c_0,G_0c_0'')}
\def \xxd{\D_1(F_0c_0,G_0c_0'')}
\def \xxe{F_1^{c_0,c_0''}}
\def \xxf{G_1^{c_0,c_0''}}
\def \xxg{- \circ \alpha_0c_0''}
\def \xxh{\alpha_0c_0 \circ -}
\def \xxi{\alpha_1^{c_0,c_0''}}
=\raisebox{15ex}{ \xymatrix@C=-2ex{
&\C_1(c_0,c_0')\times\C_1(c_0',c_0'')\ar[d]^{\circ^0}\\
& \xxa \ar[dl]_{\xxe}\ar[dr]^{\xxf}\\
\xxb\ar@{}[rr]|(.4){}="xxb"|(.6){}="xxc"\ar[dr]_{\xxg}
&&\xxc\ar[dl]^{\xxh}\\
&\xxd \ar@{=>}"xxc";"xxb"_{\xxi}}}
$$
}

$\bullet$ (\emph{functoriality w.r.t. units}) for every object $c_0$ of
$\C_0$, {\small
\begin{equation}\label{ntrans:units}
\def \xxa{\C_1(c_0,c_0)}
\def \xxb{\D_1(F_0c_0,F_0c_0)}
\def \xxc{\D_1(G_0c_0,G_0c_0)}
\def \xxd{\D_1(F_0c_0,G_0c_0)}
\def \xxe{F_1^{c_0,c_0}}
\def \xxf{G_1^{c_0,c_0}}
\def \xxg{- \circ \alpha_0c_0}
\def \xxh{\alpha_0c_0 \circ -}
\def \xxi{\alpha_1^{c_0,c_0}}
\raisebox{13ex}{\xymatrix@C=-2ex{
&\I\ar[d]^{u^0(c_0)}\\
& \xxa \ar[dl]_{\xxe}\ar[dr]^{\xxf}\\
\xxb\ar@{}[rr]|(.4){}="xxb"|(.6){}="xxc"\ar[dr]_{\xxg}
&&\xxc\ar[dl]^{\xxh}\\
&\xxd \ar@{=>}"xxc";"xxb"_{\xxi}}} =
\raisebox{13ex}{\xymatrix@R=26ex{ \I
\ar@/^3ex/[d]^{[\alpha_0c_0]}="1" \ar@/_3ex/[d]_{[\alpha_0c_0]}="2"
\\
\D_1(F_0c_0,G_0c_0) \ar@{}"1";"2"|(.4){}="3"|(.6){}="4"
\ar@{=>}"3";"4"_{id}}}
\end{equation}}
\begin{Remark}
In defining transformations, we used expressions such as $-\circ \alpha_{c_0'}$ or $\alpha_{c_0}\circ -$ to denote
composition (n-1)functors. In fact, given a n-category $\C$ and a 1-cell $c_1:c_0\to c_0'$, it is always possible
to define a pair of (n-1)functors for every other chosen objects $\bar{c}_0$ of $\C$:
$$
[-\circ^0 c_1]_{\bar{c}_0} : \C_1(\bar{c}_0,c_0)\to\C_1(\bar{c}_0,c_0')
$$
$$
[c_1 \circ^0 - ]_{\bar{c}_0} : \C_1(c_0',\bar{c}_0\to\C_1(c_0,\bar{c}_0)
$$
As a matter of fact,  they are just restrictions of 0-composition functors in $\C$. The assignment is (mutually) natural
in $\bar{c}_0$.

In fact, whiskering makes the following diagram commute, for $\bar{c}_1:\bar{\bar{c}}_0\to\bar{c}_0$:
$$
\xymatrix{
\C_1(\bar{c}_0,c_0) \ar[r]^{-\circ c_1}\ar[d]_{\bar{c}_1\circ-}
&\C_1(\bar{c}_0,c_0') \ar[d]^{\bar{c}_1\circ-}
\\
\C_1(\bar{\bar{c}}_0,c_0)\ar[r]_{-\circ c_1}
&\C_1(\bar{\bar{c}}_0,c_0')}
$$
\end{Remark}
A natural $n$-transformation of $n$-functors $\alpha:F\Rightarrow G:\C\to\D$ is called strict when
 for every pair of objects $c_0,c_0'$ of $\C$, $\alpha_1^{c_0,c_0'}$ is an identity.

\section{$n$\textbf{Cat}: the underlying category}
$n$\textbf{Cat} has the underlying category denoted $\lfloor n\mathbf{Cat}\rfloor$,
whose description follows.
Notice that in this section we will assume $n$ be an integer greater than $1$.\\

Given $n$-functors
$$
\xymatrix{\C\ar[r]^F&\D\ar[r]^G&\E}
$$
their composition $F\bullet^0 G$ (or simply $FG$) is the morphism with
$$
[F\bullet^0  G]_0=F_0 G_0
$$
in \textbf{Set}, and, for every pair of objects $c_0,c_0'$ of
$\C_0$,
$$
[F\bullet^0  G]_1^{c_0,c_0'}= F_1^{c_0,c_0'}\bullet^0  G_1^{F_0c_0, F_0 c_0'}
$$
in $(n-1)$\textbf{Cat}.

The pair $([F\bullet^0  G]_0,[F\bullet^0  G]_1)$ defines indeed a morphism of
\mbox{$(n-1)$\textbf{Cat}}. This is clear by pasting  the commutative
diagrams below, for every triple of objects $c_0, c_0', c_0''$ of
$\C_0$. They ensure functoriality w.r.t. composition and units of
(\ref{nfunct:composition}) and (\ref{ntrans:units}).
\begin{changemargin}{-10ex}{10ex}
{\small
$$
\xymatrix@C=10ex{ [c_0,c_0'] \times
[c_0',c_0'']\ar[r]^{^{\C}\circ^0} \ar[d]_{F_1^{c_0,c_0'}\times
F_1^{c_0',c_0''}}
& [c_0,c_0'']\ar[d]^{F_1^{c_0,c_0''}}\\
[Fc_0,Fc_0']\times [Fc_0',F_0c_0'']\ar[r]_(.6){^{\D}\circ^0}
\ar[d]_{G_1^{Fc_0,Fc_0'}\times G_1^{Fc_0',Fc_0''}}
& [F c_0,F_0 c_0'']\ar[d]^{G_1^{Fc_0,Fc_0''}}\\
[G(Fc_0),G(Fc_0')]\times
[G(Fc_0'),G(F_0c_0'')]\ar[r]_(.6){^{\E}\circ^0} & [G(F c_0),G(F_0
c_0'')]} \quad\xymatrix@C=14ex{ \I \ar[r]^{^\C u^0(c_0)}\ar[dr]|{^\D
u^0(Fc_0)}\ar[ddr]_{^\E u^0(G(Fc_0))}
& [c_0,c_0]\ar[d]^{F_1^{c_0,c_0}}\\
& [Fc_0,Fc_0]\ar[d]^{G_1^{Fc_0,Fc_0}}\\
& [G(Fc_0),G(Fc_0)]}
$$
}
\end{changemargin}

Furthermore, for every $n$-category $\C$, an identity functor
$$
\xymatrix{\C\ar[r]^{id_{\C}}&\C}
$$
is defined by the pair $([id_{\C}]_0],[id_{\C}]_1])$, where
$$
[id_{\C}]_0=id_{\C_0}
$$
in \textbf{Set}, and, for every pair of objects $c_0,c_0'$ of $\C$
$$
[id_{\C}]_1^{c_0,c_0'} = id_{\C_1(c_0,c_0')}
$$
in $(n-1)$\textbf{Cat}.

Notice that $([id_{\C}]_0],[id_{\C}]_1])$ satisfies trivially
functoriality diagrams (\ref{nfunct:composition}) and
(\ref{ntrans:units}).

\begin{Proposition}
(small and strict) $n$-categories and (strict) $n$-functors define a
category: $\lfloor n\mathbf{Cat}\rfloor$
\end{Proposition}
\begin{proof}
Firstly, for every pair $\C$, $\D$ of (small) $n$-categories,
$$
Hom(\C,\D)=\{n\mathrm{-functors}\  F:\C\to\D\}
$$
is a set, since $\D$ is small.

Hence we will show that composition is associative and that identities are neutral.\\

Concerning associativity, let a composable triple of morphisms  be given:
$$
\xymatrix{\C\ar[r]^F&\D\ar[r]^G&\E\ar[r]^H&\F}
$$
We want to prove $(FG)H = F(GH)$.

On objects, let us consider the
following  equalities in the category \textbf{Set}:
\begin{eqnarray*}
  [(FG)H]_0 &=& [FG]_0H_0 =\\
   &=& (F_0G_0)H_0 =\\
   &=& F_0(G_0 H_0) =\\
   &=& F_0 [GH]_0 = [F(GH)]_0
\end{eqnarray*}
Besides, for every pair of objects $c_0$, $c_0'$ of $\C$,
$(n-1)$-associativity implies:
\begin{eqnarray*}
  [(FG)H]_1^{c_0,c_0'} &=& [FG]_1^{c_0,c_0'}H_1^{G(Fc_0),G(Fc_0')} =\\
   &=& \left( F_1^{c_0,c_0'}G_1^{Fc_0,Fc_0'}\right) H_1^{G(Fc_0),G(Fc_0')} =\\
   &=& F_1^{c_0,c_0'}\left( G_1^{Fc_0,Fc_0'} H_1^{G(Fc_0),G(Fc_0')}\right)  =\\
   &=& F_1^{c_0,c_0'} [GH]_1^{Fc_0,Fc_0'} = [F(GH)]_1^{c_0,c_0'}
\end{eqnarray*}

Turning now to identities, let us consider the situation:
$$
\xymatrix{\C\ar[r]^{id_{\C}}&\C\ar[r]^F&\D\ar[r]^{id_{\D}}&\D}
$$
We want to prove $id_{\C}F=F=F id_{\D}$. On objects,
\begin{eqnarray*}
  [id_{\C}F]_0 &=& [id_{\C}]_0 F_0 =  \\
   &=& id_{\C_0} F_0 = \\
   = F_0 &=& F_0 id_{\D_0} =\\
   &=& F_0 [id_{\D}]_0 = [F_0 id_{\D}]_0
\end{eqnarray*}
and, for every pair of objects $c_0$, $c_0'$ of $\C$, neutral
$(n-1)$-identities imply:
\begin{eqnarray*}
  [id_{\C}F]_1^{c_0,c_0'} &=& [id_{\C}]_1^{c_0,c_0'} F_1^{c_0,c_0'} =  \\
   &=& id_{\C_1(c_0,c_0')} F_1^{c_0,c_0'} = \\
= F_1^{c_0,c_0'} &=& F_1^{c_0,c_0'} id_{\D_1(Fc_0,Fc_0')} =\\
   &=& F_1^{c_0,c_0'} [id_{\D}]_1^{Fc_0,Fc_0'} = [F_0 id_{\D}]_1^{c_0,c_0'}
\end{eqnarray*}
\end{proof}

\begin{Proposition}\label{prop:nCat_has_products}
The category $\lfloor n\mathbf{Cat}\rfloor$ has finite products.
\end{Proposition}
\begin{proof}
We will show that $\lfloor\mathbf{Cat}\rfloor$ has a terminal object and binary products.

Given $n$-categories $\C$ and $\D$, their (standard) product is
defined as follows:
$$
[\C\times\D]_0= \C_0\times\D_0
$$
and, for every pair $(c_0,d_0)$,  $(c_0',d_0')$ in $[\C\times\D]_0$,
$$
[\C\times\D]_1\left((c_0,d_0),(c_0',d_0')\right)=
\C_1(c_0,c_0')\times\D_1(d_0,d_0')
$$
Composition is  defined by means of universality of products in
$(n-1)$\textbf{Cat}: for every triple $(c_0,d_0)$,
$(c_0',d_0')$ and $(c_0'',d_0'')$ in  $[\C\times\D]_0$, the dotted
arrow in the diagram below gives composition: {\small$$
\xymatrix@C=7ex{
[\C\times\D]_1((c_0,d_0),(c_0',d_0'))\times[\C\times\D]_1((c_0',d_0'),(c_0'',d_0''))
\ar@{=}[d]_{id}\ar@{.>}[r]^(.65){{}^{\C\times \D}\circ}
&[\C\times\D]_1((c_0,d_0),(c_0'',d_0''))\ar@{=}[dd]^{id}\\
\Big(\C_1(c_0,c_0')\times
\D_1(d_0,d_0')\Big)\times\Big(\C_1(c_0',c_0'')\times
\D_1(d_0',d_0'')\Big)
\ar[d]_{\tau}\\
\Big(\C_1(c_0,c_0')\times
\C_1(c_0',c_0'')\Big)\times\Big(\D_1(d_0,d_0')\times
\D_1(d_0',d_0'')\Big) \ar[r]_(.65){{}^{\C}\circ \times{\
}^{\D}\circ} &\C_1(c_0,c_0'')\times \D_1(d_0,d_0'')}
$$}
where the twist isomorphism $\tau=\tau_{\C_1(c_0,c_0'),
\D_1(d_0,d_0'),\C_1(c_0',c_0''),\D_1(d_0',d_0'')}$
is given by products properties in $(n-1)$\textbf{Cat}.\\

Identities are defined in the same way. For every object $(c_0,d_0)$
in $\C\times\D$, by the dotted arrow in the diagram below:
{\small$$ \xymatrix@C=20ex{
\I\ar@{.>}[r]^(.4){{}^{\C\times\D}u^0((C_0,d_0))}\ar[d]_{iso}
&[\C\times\D]_1((c_0,d_0),(c_0,d_0))\ar@{=}[d]^{id}\\
\I\times\I\ar[r]_(.4){{}^{\C}u^0(c_0)\times {}^{\D}u^0(d_0)}
&\C_1(c_0,c_0)\times \D_1(d_0,d_0) }
$$}

Product projections
$$
\xymatrix{ \C &\ar[l]_{\Pi_\C} \C\times \D \ar[r]^{\Pi_\D} &\D}
$$
are given respectively  by projections
$$
[\Pi_{\C}]_0 = \pi^{\C_0\times \D_0}_{\C_0},\quad [\Pi_{\D}]_0 =
\pi^{\C_0\times \D_0}_{\D_0}
$$
and by the following compositions in $(n-1)$\textbf{Cat}, for every
pair $(c_0,d_0)$ and  $(c_0',d_0')$:
$$
\xymatrix@C=-0ex{ [\C\times \D]_1((c_0,d_0),(c_0',d_0'))
\ar@{.>}[rr]^(.6){[\Pi_\C]_1^{(c_0,d_0),(c_0',d_0')}}
\ar[dr]_{id\times !}
&&\C_1(c_0,c_0')\\
&\C_1(c_0,c_0')\times \I \ar[ur]_{iso} },
$$
$$
\xymatrix@C=-0ex{ [\C\times \D]_1((c_0,d_0),(c_0',d_0'))
\ar@{.>}[rr]^(.6){[\Pi_\D]_1^{(c_0,d_0),(c_0',d_0')}}
\ar[dr]_{id\times !}
&&\D_1(d_0,d_0')\\
&\D_1(d_0,d_0')\times \I \ar[ur]_{iso} }
$$

It is a matter of repeated use of the universal  property of products in
$(n-1)$\textbf{Cat} to prove that all these data define
a $n$-category and two $n$-categories morphisms, and that $(\C\times\D,\Pi_\C,\Pi_\D)$ is a product in
$n$\textbf{Cat}.\\

Concerning the terminal $n$-category, a standard construction
follows. A terminal $\Id{n}$ is given by the pair
$$
[\Id{n}]_0= \{\ast\}, \qquad [\Id{n}]^{\ast,\ast}_0=\Id{n-1};
$$
with composition
$\xymatrix{\Id{n-1}\times\Id{n-1}\ar@<+1ex>[r]^{\sim}&\Id{n-1}}$.

Classical constructions of categorical limits help in defining
$n$-ary products and canonical isomorphisms
\begin{eqnarray}
  \alpha_{\A,\B,\C} &:&\xymatrix{(\A\times\B)\times\C\ar[r]^{\sim}&\A\times(\B\times\C)}\\
  \rho_{\A} &:& \xymatrix{\A\ar[r]^(.4){\sim}&\A\times\Id{n}}\\
  \lambda_{\A} &:& \xymatrix{\A\ar[r]^(.4){\sim}&\Id{n}\times\A }\\
  \tau_{\A,\B,\C,\D} &:&\xymatrix{(\A\times\B)\times(\C\times\D)\ar[r]^{\sim}&(\A\times\C)\times(\B\times\D)}
\end{eqnarray}
\end{proof}

\section{$n$\textbf{Cat}: the hom-categories}\label{section:homcats}
In this section we describe, hom-categories $n$\textbf{Cat}$(\C,\D)$, once $n$-categories $\C$ and $\D$ are fixed.\\

\subsection{Vertical composition}
Given the diagram:
$${
\def \xaa{\C}
\def \xab{\D}
\def \xac{E}
\def \xad{F}
\def \xae{G}
\def \xaf{\omega}
\def \xag{\alpha}
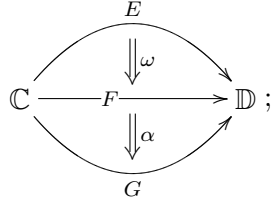
\xymatrix@C=15ex{ \xaa \ar@/_6ex/[r]_{\xae}^(.5){}="1a"
\ar[r]|(.4){\xad}^(.5){}="2a"_(.5){}="1b"
\ar@/^6ex/[r]^{\xac}_(.5){}="2b" & \xab
\ar@{}"1a";"2a"|(.2){}="1aa"|(.8){}="2aa"
\ar@{}"1b";"2b"|(.2){}="1bb"|(.8){}="2bb" \ar@{=>}"2aa";"1aa"^{\xag}
\ar@{=>}"2bb";"1bb"^{\xaf}}; }$$
one defines a (vertical, or $1$-)composition $\xymatrix@C=3ex{\omega\bullet^1\alpha:E\ar@{=>}[r]& G}$ in the following way:\\

$\bullet$ for every object $c_0$ in $\C$,
$$
[\omega \alpha]_0 (c_0)= \omega_0 c_0 \circ^0 \alpha_0 c_0 :
\xymatrix{E c_0\ar[r]& G c_0}
$$
$\bullet$ for every pair of objects $c_0, c_0'$ in $\C$, the diagram
below describes
$$
[\omega\bullet^1\alpha]_1^{c_0,c_0'}=\left(\omega c_0\circ\alpha_1^{c_0,c_0'}\right)\bullet^1\left(\omega_1^{c_0,c_0'}\circ\alpha c_0'\right)
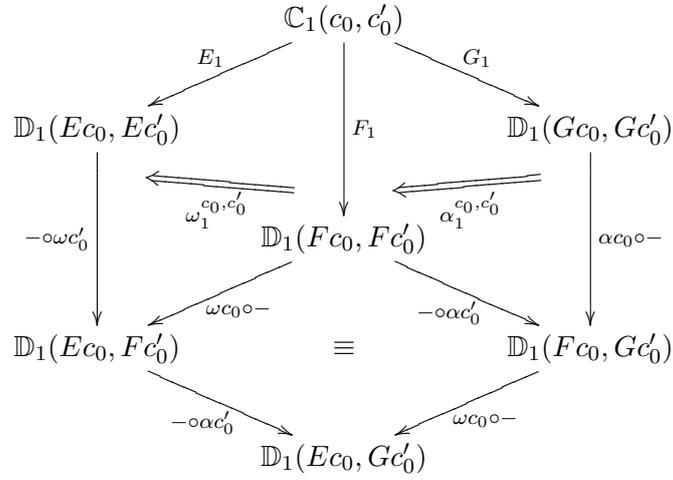
$$
$$
{\def \xwa{\C_1(c_0,c_0')}
\def \xwb{\D_1(Ec_0,Ec_0')}
\def \xwc{\D_1(Gc_0,Gc_0')}
\def \xwd{\D_1(Fc_0,Fc_0')}
\def \xwe{\D_1(Ec_0,Fc_0')}
\def \xwf{\D_1(Fc_0,Gc_0')}
\def \xwg{\D_1(Ec_0,Gc_0')}
\def \xwh{E_1}
\def \xwi{F_1}
\def \xwl{G_1}
\def \xwm{-\circ \omega c_0'}
\def \xwn{\omega c_0 \circ -}
\def \xwo{- \circ \alpha c_0'}
\def \xwp{\alpha c_0 \circ -}
\def \xwq{- \circ \alpha c_0'}
\def \xwr{\omega c_0 \circ -}
\def \xws{\omega_1^{c_0,c_0'}}
\def \xwt{\alpha_1^{c_0,c_0'}}
\xymatrix{ &\xwa
\ar[dl]_{\xwh} \ar[dd]^{\xwi}|(.8){}="2" \ar[dr]^{\xwl}\\
\xwb \ar[dd]_{\xwm}|(.2){}="3"
&&\xwc\ar[dd]^{\xwp}|(.2){}="1"\\
&\xwd \ar[dl]^{\xwn} \ar[dr]_{\xwo}\\
\xwe\ar[dr]_{\xwq} &\equiv&\xwf\ar[dl]^{\xwr}\\
&\xwg \ar@{} "1";"2"|(.2){}="a1"|(.8){}="a2" \ar@{}
"2";"3"|(.2){}="b1"|(.8){}="b2" \ar@{=>} "a1";"a2"^{\xwt}  \ar@{=>}
"b1";"b2"^{\xws} }}
$$
To prove that these data give indeed a $2$-morphism, unit functoriality
(\ref{ntrans:units}) and composition functoriality
(\ref{ntrans:composition}) equations must hold.\\

To this end, chose an object $c_0$ of $\C$ and consider the
following chain of diagrams equalities: {\small
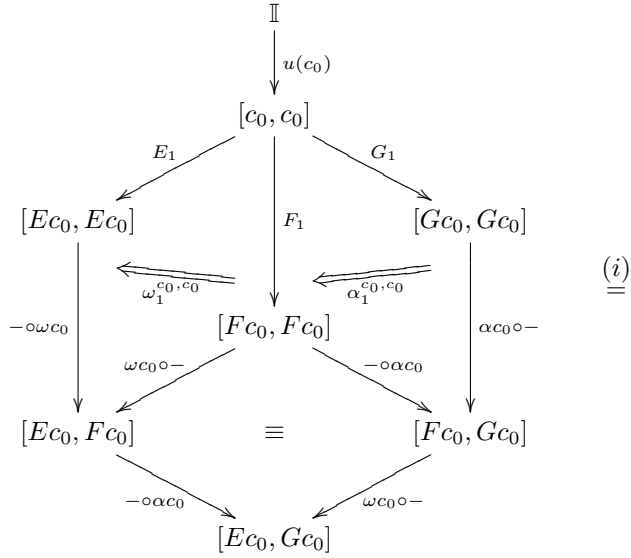
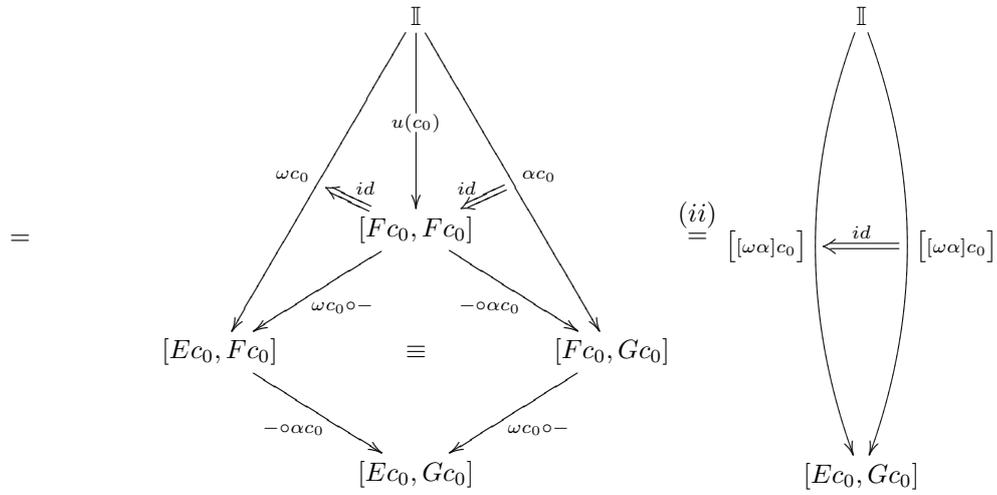
$$
\raisebox{25ex}{\def \xwa{[c_0,c_0]}
\def \xwb{[Ec_0,Ec_0]}
\def \xwc{[Gc_0,Gc_0]}
\def \xwd{[Fc_0,Fc_0]}
\def \xwe{[Ec_0,Fc_0]}
\def \xwf{[Fc_0,Gc_0]}
\def \xwg{[Ec_0,Gc_0]}
\def \xwh{E_1}
\def \xwi{F_1}
\def \xwl{G_1}
\def \xwm{-\circ \omega c_0 }
\def \xwn{\omega c_0\circ -}
\def \xwo{-\circ \alpha c_0}
\def \xwp{\alpha c_0\circ -}
\def \xwq{-\circ \alpha c_0}
\def \xwr{\omega c_0\circ -}
\def \xws{\omega_1^{c_0,c_0}}
\def \xwt{\alpha_1^{c_0,c_0}}
\xymatrix{
&\I\ar[d]^{u(c_0)}\\
&\xwa
\ar[dl]_{\xwh} \ar[dd]^{\xwi}|(.8){}="2" \ar[dr]^{\xwl}\\
\xwb \ar[dd]_{\xwm}|(.2){}="3"
&&\xwc\ar[dd]^{\xwp}|(.2){}="1"\\
&\xwd \ar[dl]_{\xwn} \ar[dr]^{\xwo}\\
\xwe\ar[dr]_{\xwq} &\equiv&\xwf\ar[dl]^{\xwr}\\
&\xwg \ar@{} "1";"2"|(.2){}="a1"|(.8){}="a2" \ar@{}
"2";"3"|(.2){}="b1"|(.8){}="b2" \ar@{=>} "a1";"a2"^{\xwt}  \ar@{=>}
"b1";"b2"^{\xws} }} \qquad\xy(0,3)*{(i)};(0,0)*{=}\endxy
$$
$$
=\raisebox{20ex}{\def \xwa{[c_0,c_0]}
\def \xwb{[Ec_0,Ec_0]}
\def \xwc{[Gc_0,Gc_0]}
\def \xwd{[Fc_0,Fc_0]}
\def \xwe{[Ec_0,Fc_0]}
\def \xwf{[Fc_0,Gc_0]}
\def \xwg{[Ec_0,Gc_0]}
\def \xwn{\omega c_0\circ -}
\def \xwo{-\circ \alpha c_0}
\def \xwq{-\circ \alpha c_0}
\def \xwr{\omega c_0\circ -}
\def \xws{\omega_1^{c_0,c_0}}
\def \xwt{\alpha_1^{c_0,c_0}}
\xymatrix@R=7ex{ &\I\ar[dd]|{u(c_0)}
\ar[dddl]_{\omega c_0}^{}="1" \ar[dddr]^{\alpha c_0}_{}="2"\\
\\
&\xwd \ar[dl]^{\xwn} \ar[dr]_{\xwo}\ar@{=>};"1"_(.6){id}\ar@{=>}"2";_(.4){id}\\
\xwe\ar[dr]_{\xwq} &\equiv&\xwf\ar[dl]^{\xwr}\\
&\xwg }} \xy(0,3)*{(ii)};(0,0)*{=}\endxy
\raisebox{20ex}{\xymatrix@R=37ex{\I
\ar@/_4ex/[d]_{\big[[\omega\alpha]c_0\big]}^{}="1"
\ar@/^4ex/[d]^{\big[[\omega\alpha]c_0\big]}_{}="2"
\\
[Ec_0,Gc_0] \ar@{=>}"2";"1"_{id}}}
$$
} where $(i)$ follows for units functoriality of $\omega$ and $\alpha$,
while $(ii)$ from functoriality of constant functors. This proves unit functoriality.\\

Concerning composition functoriality, take three objects $c_0$, $c_0'$
and $c_0''$ in $\C$, and consider the following diagram:

{
\def \objectstyle{\scriptstyle}
\def \labelstyle{\scriptstyle}
\def \xa{[c_0,c_0']\tm [c_0',c_0'']}
\def \xb{[c_0,c_0']\tm[Ec_0',Ec_0'']}
\def \xc{[c_0,c_0']\tm[Fc_0',Fc_0'']}
\def \xd{[Fc_0,Fc_0']\tm[c_0',c_0'']}
\def \xe{[Gc_0,Gc_0']\tm[c_0',c_0'']}
\def \xf{[c_0,c_0']\tm[Ec_0',Fc_0'']}
\def \xg{[c_0,c_0']\tm[Fc_0',Gc_0'']}
\def \xh{[Ec_0,Fc_0']\tm[Fc_0',Gc_0'']}
\def \xk{[Ec_0,Fc_0']\tm[c_0',c_0'']}
\def \xq{[Fc_0,Gc_0']\tm[c_0',c_0'']}
\def \xl{[c_0,c_0']\tm[Ec_0',Gc_0'']}
\def \xm{[Ec_0,Ec_0']\tm[Ec_0',Gc_0'']}
\def \xn{[Ec_0,Gc_0'']}
\def \xo{[Ec_0,Gc_0']\tm[Gc_0',Gc_0'']}
\def \xp{[Ec_0,Gc_0']\tm[c_0',c_0'']}
\def \xr{[c_0,c_0']\tm[Gc_0',Gc_0'']}
\def \xs{[Ec_0,Ec_0']\tm[c_0',c_0'']}
\def \xab{id\tm E_1}
\def \xac{id\tm F_1}
\def \xar{id\tm G_1}
\def \xas{E_1\tm id}
\def \xad{F_1\tm id}
\def \xae{G_1\tm id}
\def \xbf{id\tm (-\circ \omega c_0'')}
\def \xcf{id\tm ( \omega c_0'\circ -)}
\def \xcg{id\tm (-\circ \alpha c_0'')}
\def \xrg{id\tm ( \alpha c_0'\circ -)}
\def \xsk{(-\circ \omega c_0')\tm id}
\def \xdk{(\omega c_0 \circ -)\tm id}
\def \xdq{(-\circ\alpha c_0')\tm id}
\def \xeq{(\alpha c_0\circ-)\tm id}
\def \xfl{id\tm (-\circ \alpha_{c_0''})}
\def \xgl{id\tm ( \omega_{c_0'}\circ -)}
\def \xgh{(E_1(-)\circ \omega c_0')\tm id}
\def \xhn{\circ}
\def \xkh{id \tm (\alpha c_0' \circ G_1(-))}
\def \xkp{(-\circ\alpha c_0')\tm id}
\def \xqp{(\omega c_0 \circ -)\tm id}
\def \xlm{E_1\tm id}
\def \xmn{\circ}
\def \xon{\circ}
\def \xpo{id \tm G_1}
\def \xcb{id\times \omega_1^{c_0',c_0''}}
\def \xrc{id\times \alpha_1^{c_0',c_0''}}
\def \xds{\omega_1^{c_0,c_0'}\times id}
\def \xed{\alpha_1^{c_0,c_0'}\times id}
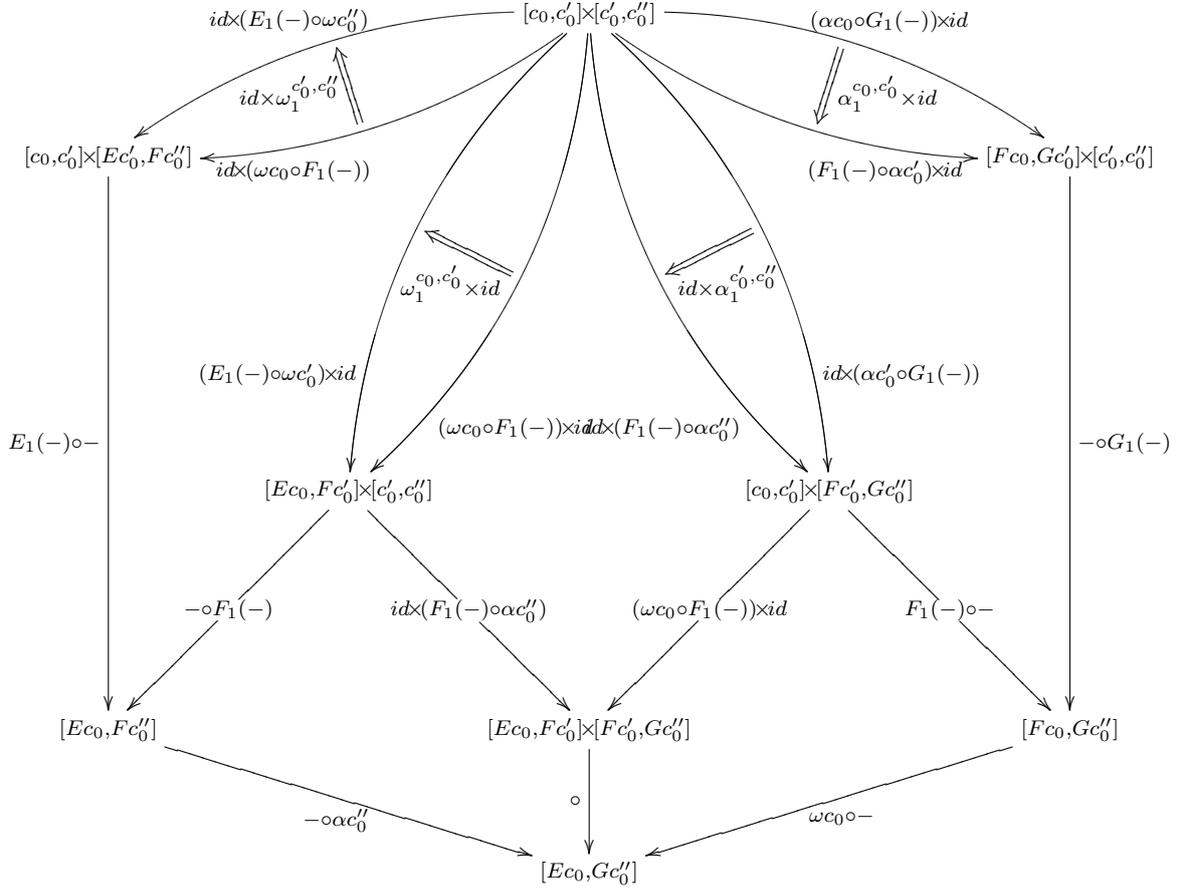
\begin{changemargin}{-15ex}{0ex}
$$
\xy 0;/r.27pc/: (0,60)*+{\xa}="a"; (-45,40)*+{\xb}="b";
(-30,25)*+{\xc}="c"; (30,25)*+{\xd}="d"; (45,40)*+{\xe}="e";
(-60,0)*+{\xf}="f"; (-40,-20)*+{\xg}="g"; (0,-35)*+{\xh}="h";
(40,-20)*+{\xk}="k"; (60,0)*+{\xq}="q"; (-45,-40)*+{\xl}="l";
(-25,-55)*+{\xm}="m"; (0,-60)*+{\xn}="n"; (25,-55)*+{\xo}="o";
(45,-40)*+{\xp}="p"; (-13,10)*+{\xr}="r"; (13,10)*+{\xs}="s";
{\ar_{\xab}"a";"b"}; {\ar|{\xac}"a";"c"}; {\ar|{\xar}"a";"r"};
{\ar|{\xas}"a";"s"}; {\ar|{\xad}"a";"d"}; {\ar^{\xae}"a";"e"};
{\ar_{\xbf}"b";"f"}; {\ar|{\xcf}"c";"f"}; {\ar|{\xcg}"c";"g"};
{\ar^{\xrg}"r";"g"}; {\ar|{\xsk}"s";"k"}; {\ar|{\xdk}"d";"k"};
{\ar|{\xdq}"d";"q"}; {\ar^{\xeq}"e";"q"}; {\ar_{\xfl}"f";"l"};
{\ar^{\xgl}"g";"l"}; {\ar@{..>}^{\xgh}"g";"h"};
{\ar@{..>}_{\xkh}"k";"h"}; {\ar_{\xkp}"k";"p"}; {\ar^{\xqp}"q";"p"};
{\ar_{\xlm}"l";"m"}; {\ar_{\xmn}"m";"n"}; {\ar@{..>}_{\xhn}"h";"n"};
{\ar^{\xpo}"p";"o"}; {\ar^{\xon}"o";"n"};
{\ar@{}"c";"b"|(.25){}="cb1"|(.75){}="cb2"};{\ar@{=>}_{\xcb}"cb1";"cb2"};
{\ar@{}"r";"c"|(.25){}="rc1"|(.75){}="rc2"};{\ar@{=>}_{\xrc}"rc1";"rc2"};
{\ar@{}"d";"s"|(.25){}="ds1"|(.75){}="ds2"};{\ar@{=>}_{\xds}"ds1";"ds2"};
{\ar@{}"e";"d"|(.25){}="ed1"|(.75){}="ed2"};{\ar@{=>}_{\xcb}"ed1";"ed2"};
\endxy
$$
\end{changemargin}}

After applying the product interchange (see section \ref{Product interchange rules})
to 2-morphisms $\omega_1^{c_0,c_0'}$
and $\alpha_1^{c_0',c_0''}$, the  last diagram becomes

\begin{changemargin}{-10ex}{-10ex}
{
\def \objectstyle{\scriptstyle}
\def \labelstyle{\scriptstyle}
\def \xa{[c_0,c_0']\tm[c_0',c_0'']}
\def \xb{[c_0,c_0']\tm[Ec_0',Fc_0'']}
\def \xc{[Fc_0,Gc_0']\tm[c_0',c_0'']}
\def \xd{[Ec_0,Fc_0']\tm[c_0',c_0'']}
\def \xe{[c_0,c_0']\tm[Fc_0',Gc_0'']}
\def \xf{[Ec_0,Fc_0']\tm[Fc_0',Gc_0'']}
\def \xg{[Ec_0,Fc_0'']}
\def \xh{[Fc_0,Gc_0'']}
\def \xl{[Ec_0,Gc_0'']}
\def \xabt{id\tm (E_1(-)\circ \omega c_0'')\quad}
\def \xabs{\quad id\tm (\omega c_0\circ F_1(-) )}
\def \xact{(F_1(-)\circ \alpha c_0')\tm id\quad}
\def \xacs{\quad(\alpha c_0 \circ G_1(-))\tm id}
\def \xadt{(E_1(-)\circ \omega c_0')\tm id}
\def \xads{(\omega c_0 \circ F_1(-))\tm id}
\def \xaet{id\tm(F_1(-)\circ \alpha c_0'')}
\def \xaes{id\tm(\alpha c_0' \circ G_1(-))}
\def \xbg{E_1(-)\circ -}
\def \xch{- \circ G_1(-)}
\def \xdg{-\circ F_1(-)}
\def \xdf{id\tm (F_1(-)\circ \alpha c_0'')}
\def \xef{(\omega c_0 \circ F_1(-))\tm id}
\def \xeh{F_1(-)\circ-}
\def \xfl{\circ}
\def \xgl{-\circ \alpha c_0''}
\def \xhl{\omega c_0 \circ -}
\def \abii{id \times \omega_1^{c_0',c_0''}}
\def \acii{\alpha_1^{c_0,c_0'}\times id}
\def \adii{\omega_1^{c_0,c_0'}\times id}
\def \aeii{id \times \alpha_1^{c_0',c_0''}}
$$
\xy 0;/r.3pc/: (0,60)*+{\xa}="a"; (-50,45)*+{\xb}="b";
(50,45)*+{\xc}="c"; (-25,10)*+{\xd}="d"; (25,10)*+{\xe}="e";
(0,-15)*+{\xf}="f"; (-50,-15)*+{\xg}="g"; (50,-15)*+{\xh}="h";
(0,-30)*+{\xl}="l"; {\ar@/_1.5pc/_(.5){\xabt}"a";"b"^{}="abt"};
{\ar@/^1.5pc/^(.7){\xabs}"a";"b"_{}="abs"};
{\ar@/_1.5pc/_(.7){\xact}"a";"c"^{}="act"};
{\ar@/^1.5pc/^(.5){\xacs}"a";"c"_{}="acs"};
{\ar@/_1.8pc/_(.8){\xadt}"a";"d"^{}="adt"};
{\ar@/^1.8pc/^(.8){\xads}"a";"d"_{}="ads"};
{\ar@/_1.8pc/_(.8){\xaet}"a";"e"^{}="aet"};
{\ar@/^1.8pc/^(.8){\xaes}"a";"e"_{}="aes"}; {\ar_{\xbg}"b";"g"};
{\ar^{\xch}"c";"h"}; {\ar|{\xdg}"d";"g"}; {\ar|{\xdf}"d";"f"};
{\ar|{\xef}"e";"f"}; {\ar|{\xeh}"e";"h"}; {\ar_{\xfl}"f";"l"};
{\ar_{\xgl}"g";"l"}; {\ar^{\xhl}"h";"l"};
{\ar@{=>}^{\abii}"abs";"abt"}; {\ar@{=>}^{\acii}"acs";"act"};
{\ar@{=>}^(.4){\adii}"ads";"adt"};
{\ar@{=>}^(.6){\aeii}"aes";"aet"};
\endxy
$$}
\end{changemargin}

By functoriality of 2-morphisms in
$(n-1)$\textbf{Cat}  the two sides of the diagram get

{
\def \objectstyle{\scriptstyle}
\def \labelstyle{\scriptstyle}
\def \xa{[c_0,c_0']\tm[c_0',c_0'']}
\def \xb{[c_0,c_0'']}
\def \xc{[c_0,c_0'']}
\def \xd{[Ec_0,Fc_0'']}
\def \xe{[Fc_0,Gc_0'']}
\def \xf{[Ec_0,Gc_0'']}
\def \xbdt{E_1(-)\circ \omega c_0''}
\def \xbds{\omega c_0\circ F_1(-)}
\def \xcet{F_1(-)\circ \alpha c_0''}
\def \xces{\alpha c_0\circ G_1(-)}
\def \xab{\circ}
\def \xac{\circ}
\def \xdf{-\circ \alpha c_0''}
\def \xef{\omega c_0 \circ -}
\def \bdii{\omega_1^{c_0,c_0''}}
\def \ceii{\alpha_1^{c_0,c_0''}}
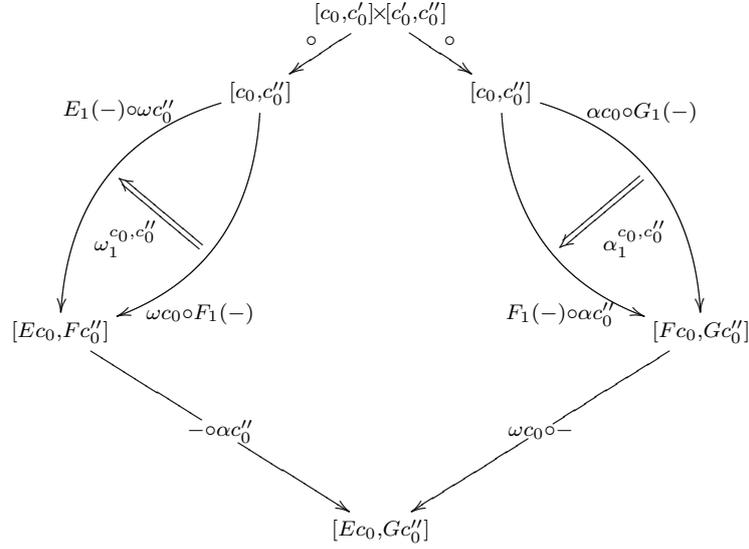
$$
\xy 0;/r.25pc/: (0,30)*+{\xa}="a"; (-15,20)*+{\xb}="b";
(15,20)*+{\xc}="c"; (-40,-10)*+{\xd}="d"; (40,-10)*+{\xe}="e";
(0,-35)*+{\xf}="f"; {\ar@/_2pc/_(.3){\xbdt}"b";"d"^{}="bdt"};
{\ar@/^2pc/^(.7){\xbds}"b";"d"_{}="bds"};
{\ar@/_2pc/_(.7){\xcet}"c";"e"^{}="cet"};
{\ar@/^2pc/^(.3){\xces}"c";"e"_{}="ces"}; {\ar_{\xab}"a";"b"};
{\ar^{\xac}"a";"c"}; {\ar|{\xdf}"d";"f"}; {\ar|{\xef}"e";"f"};
{\ar@{=>}^{\bdii}"bds";"bdt"}; {\ar@{=>}^{\ceii}"ces";"cet"};
\endxy
$$}

that is exactly $[\omega\bullet^1\alpha]_1^{c_0,c_0''}$, and this concludes
the proof.

\subsection{Units}

Given a morphism of $n$-categories $\xymatrix{F:\C \ar[r]^{F}&\D}$,
it is possible to define the \emph{unit 2-cell of $F$}, This is
denoted $id_F$, with
$$
  [id_F]_0 (c_0) = {\xymatrix{id_{Fc_0}:Fc_0 \ar[r]&Fc_0}}
$$
and
$$
[id_F]_1^{c_0,c_0'}  = id_{F_1^{c_0,c_0'}}
$$
since in the diagram
$$
{
\def \xxa{\C_1(c_0,c_0')}
\def \xxb{\D_1(Fc_0,Fc_0')}
\def \xxc{\D_1(Fc_0,Fc_0')}
\def \xxd{\D_1(Fc_0,Fc_0')}
\def \xxe{F_1^{c_0,c_0'}}
\def \xxf{F_1^{c_0,c_0'}}
\def \xxg{id_{Fc_0}\circ -}
\def \xxh{-\circ id_{Fc_0'}}
\def \xxi{id_{F_1^{c_0,c_0'}}}
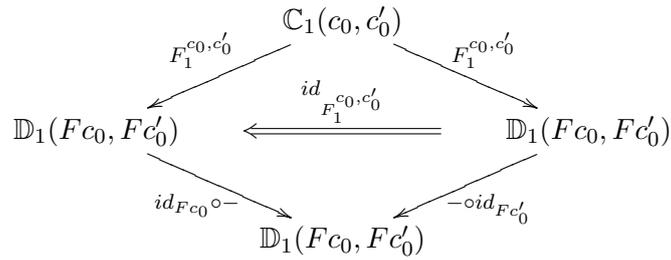
\xymatrix{
& \xxa \ar[dl]_{\xxe}\ar[dr]^{\xxf}\\
\xxb\ar@{}[rr]|(.3){}="xxb"|(.7){}="xxc"\ar[dr]_{\xxg}
&&\xxc\ar[dl]^{\xxh}\\
&\xxd \ar@{=>}"xxc";"xxb"_{\xxi}}}
$$

$id_{Fc_0}\circ -$ and $-\circ id_{Fc_0'}$ are both equal to the
identity $n-1$-functor over $\D_1(Fc_0,Fc_0')$. It is
straightforward to see that these give a 2-morphism, according to our
definition.

\begin{Proposition}
Let us fix $n$-categories $\C$ and $\D$. Morphisms between them and
2-morphisms between those form a category, with composition and
units given above.
\end{Proposition}
\begin{proof}
We will sometimes denote 1-compositions of 2-morphisms just by juxtaposition. We must prove that  composition is associative and units are neutral. To this
end, we start considering a diagram:
$$
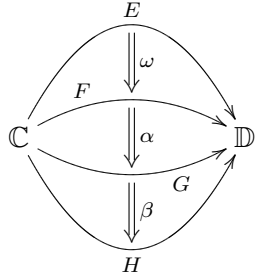
\xymatrix@C=15ex{\C \ar@/^9ex/[r]^E_{}="o1"
\ar@/^3ex/[r]^{}="o2"^(.3)F _{}="a1" \ar@/_3ex/[r]^{}="a2"_(.7)G
_{}="b1" \ar@/_9ex/[r]^{}="b2"_H &\D \ar@{=>} "o1";"o2"^{\omega}
\ar@{=>} "a1";"a2"^{\alpha} \ar@{=>} "b1";"b2"^{\beta} }
$$
We want to prove $(\omega\alpha)\beta=\omega(\alpha\beta)$.\\

For every $c_0$ in $\C_0$, by associative composition of maps
$$
[(\omega\alpha)\beta]_0(c_0)=(\omega_0c_0\circ\alpha_0c_0)\circ\beta_0c_0=
\omega_0c_0\circ(\alpha_0c_0\circ\beta_0c_0)=[\omega(\alpha\beta)]_0(c_0).
$$
Furthermore for every pair $c_0, c_0'$ in $\C_0$, associative
vertical composition of 2-morphisms of $(n-1)$-categories gives the
following diagram for both $[(\omega\alpha)\beta]_1^{c_0c_0'}$ and
$[\omega(\alpha\beta)]_1^{c_0c_0'}$
{\small$$
\xymatrix{ &&[c_0,c_0'] \ar@/_4ex/[ddll]_{E_1} \ar[ddl]^{F_1}
\ar[ddr]_{G_1} \ar@/^4ex/[ddrr]^{H_1}
\\
\\
[Ec_0,Ec_0']\ar[d]_{-\circ \omega c_0'}
&\ar@{=>}[l]_{\omega_1^{c_0,c_0'}}[Fc_0,Fc_0']\ar[dl]^{\omega c_0
\circ -} \ar[dr]_{-\circ \alpha c_0'}
&&\ar@{}[ll]|(.3){}="1"|(.7){}="2" [Gc_0,Gc_0']\ar[dl]^{\alpha c_0
\circ -} \ar[dr]_{-\circ \beta c_0'}
&\ar@{=>}[l]_{\beta_1^{c_0,c_0'}}[Hc_0,Hc_0']\ar[d]^{\beta c_0 \circ
-}
\\
  [Ec_0,Fc_0']\ar[dr]_{-\circ \alpha c_0'}
&&[Fc_0,Gc_0']\ar[dl]^{\omega c_0 \circ -} \ar[dr]_{-\circ \beta
c_0'} &&[Gc_0,Hc_0']\ar[dl]^{\alpha c_0 \circ -}
\\
&[Ec_0,Gc_0']\ar[dr]_{-\circ \beta c_0'}
&&[Fc_0,Hc_0']\ar[dl]^{\omega c_0 \circ -}
\\
&&[Ec_0,Hc_0'] \ar@{=>}"1";"2"_{\alpha_1^{c_0,c_0'}}}
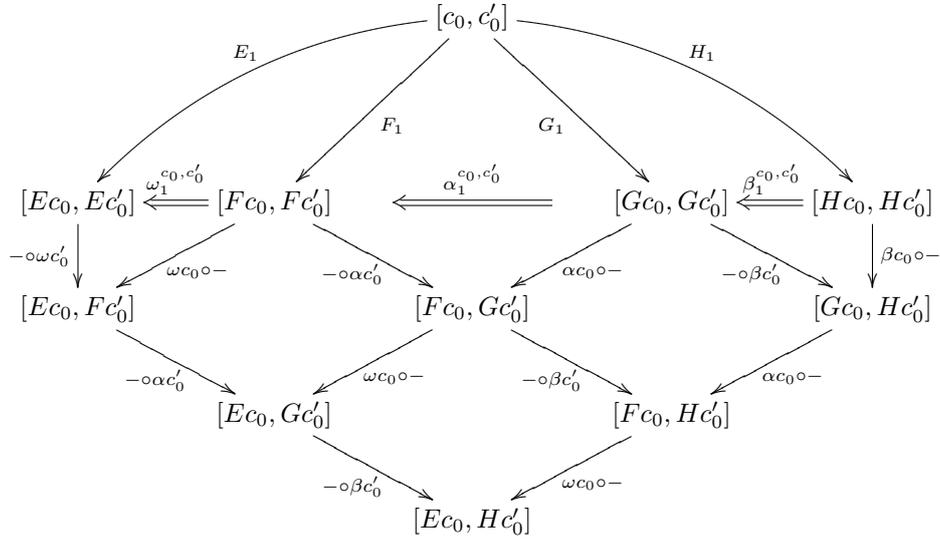
$$}
Finally, we will show that $\alpha\, id_G=\alpha$ ($id_F\, \alpha =\alpha$  is proved similarly).\\

For every $c_0$ in $\C_0$, by neutral identities of  maps
$$
[\alpha\, id_G]_0(c_0)  =  \alpha_0 c_0\, id_{Gc_0}=\alpha_0 c_0;
$$
furthermore, for every pair $c_0, c_0'$ in $\C_0$, neutral
identities for vertical composition of 2-morphisms of
$(n-1)$-categories give: {\small$$ {\def \xwa{[c_0,c_0']}
\def \xwb{[Fc_0,Fc_0']}
\def \xwc{[Gc_0,Gc_0']}
\def \xwd{[Gc_0,Gc_0']}
\def \xwe{[Fc_0,Gc_0']}
\def \xwf{[Gc_0,Gc_0']}
\def \xwg{[Fc_0,Gc_0']}
\def \xwh{F_1}
\def \xwi{G_1}
\def \xwl{G_1}
\def \xwm{-\circ\alpha c_0'}
\def \xwn{\alpha c_0\circ-}
\def \xwo{id}
\def \xwp{id}
\def \xwq{id}
\def \xwr{\alpha c_0\circ-}
\def \xws{\alpha_1^{c_0,c_0'}}
\def \xwt{id}
\raisebox{20ex}
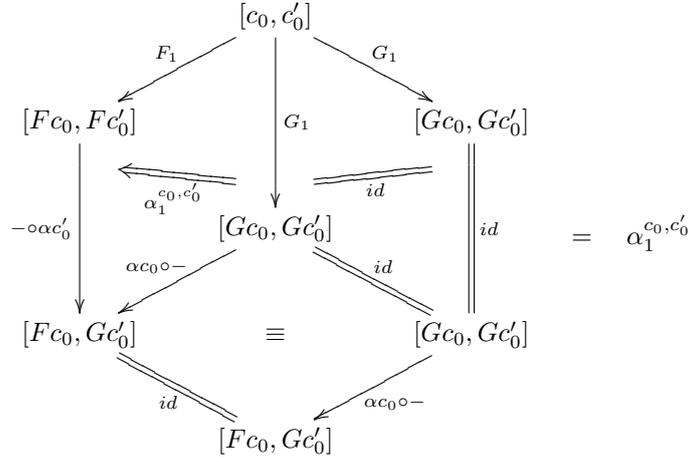
{\xymatrix{ &\xwa
\ar[dl]_{\xwh} \ar[dd]^{\xwi}|(.8){}="2" \ar[dr]^{\xwl}\\
\xwb \ar[dd]_{\xwm}|(.2){}="3"
&&\xwc\ar@{=}[dd]^{\xwp}|(.2){}="1"\\
&\xwd \ar[dl]_{\xwn} \ar@{=}[dr]^{\xwo}\\
\xwe\ar@{=}[dr]_{\xwq} &\equiv&\xwf\ar[dl]^{\xwr}\\
&\xwg \ar@{} "1";"2"|(.2){}="a1"|(.8){}="a2" \ar@{}
"2";"3"|(.2){}="b1"|(.8){}="b2" \ar@{=} "a1";"a2"^{\xwt}  \ar@{=>}
"b1";"b2"^{\xws} }}}\quad=\quad\alpha_1^{c_0,c_0'}
$$}
\end{proof}

\section{$n$\textbf{Cat}: the sesqui-categorical structure}
In the next sections we will introduce reduced left/right compositions of morphisms and 2-morphisms of n-categories, in order
to show that $n$\textbf{Cat} has a canonical sesqui-categorical structure. Notice that
$0$\textbf{Cat}=\textbf{Set} has a trivial sesqui-categorical structure (all 2-cells are identities),
while $1$\textbf{Cat}=\textbf{Cat} has a canonical 2-categorical structure,
that inherits a sesqui-categorical structure, forgetting horizontal composition
of 2-cells. Hence we may well suppose $n>1$.
\subsection{Defining reduced left-composition}

Given the situation
$$
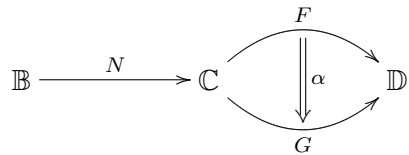
\xymatrix@C=12ex{\B\ar[r]^N&\C\ar@/^4ex/[r]^F_{}="1"\ar@/_4ex/[r]_G^{}="2"&\D\ar@{=>}"1";"2"^{\alpha}}
$$
one defines reduced horizontal composition $N\bullet^0\alpha:NF\Rightarrow NG:\B\rightarrow\D$ (or $0$-composition)
in the following way:\\

$\bullet$ for every object $b_0$ in $\B$,
$$
[N\bullet^0\alpha]_0 = \alpha_0(N(b_0)):F(N(b_0))\rightarrow G(N(b_0))
$$
$\bullet$ for every pair of objects $b_0, b_0'$ of $\B$, the diagram below describes $[N\bullet^0\alpha]_1^{b_0,b_0'}$ by means of
reduced left composition in $(n-1)$\textbf{Cat}:
$$
[N\bullet^0\alpha]_1^{b_0,b_0'}=N_1^{b_0,b_0'}\bullet^0\alpha_1^{Nb_0,Nb_0'}
$$
$$
\xymatrix@C=-6ex{
&\B_1(b_0,b_0')\ar[d]^{N_1}\ar@{.>}[ddl]_{[NF]_1}\ar@{.>}[ddr]^{[NG]_1}\\
&\C_1(Nb_0,Nb_0')\ar[dl]_{F_1}\ar[dr]^{G_1}\\
\D_1(NF(b_0),NF(b_0'))\ar[dr]_{-\circ\alpha_{Nb_0'}}
&&\D_1(NG(b_0),NG(b_0'))\ar[dl]^{\alpha_{Nb_0}\circ-}\ar@{}[ll]|(.4){}="1"|(.6){}="2"\\
&\D_1(NF(b_0),NG(b_0'))
\ar@{=>}"1";"2"_(.3){\alpha_1^{Nb_0,Nb_0'}}}
$$

To prove that these data give indeed a $2$-morphism, unit and composition axioms equations
(\ref{ntrans:units}) (\ref{ntrans:composition})  must hold.\\

To this purpose, let us chose first an object $b_0$ of $\B$ and  consider the
following chain of diagrams equalities:
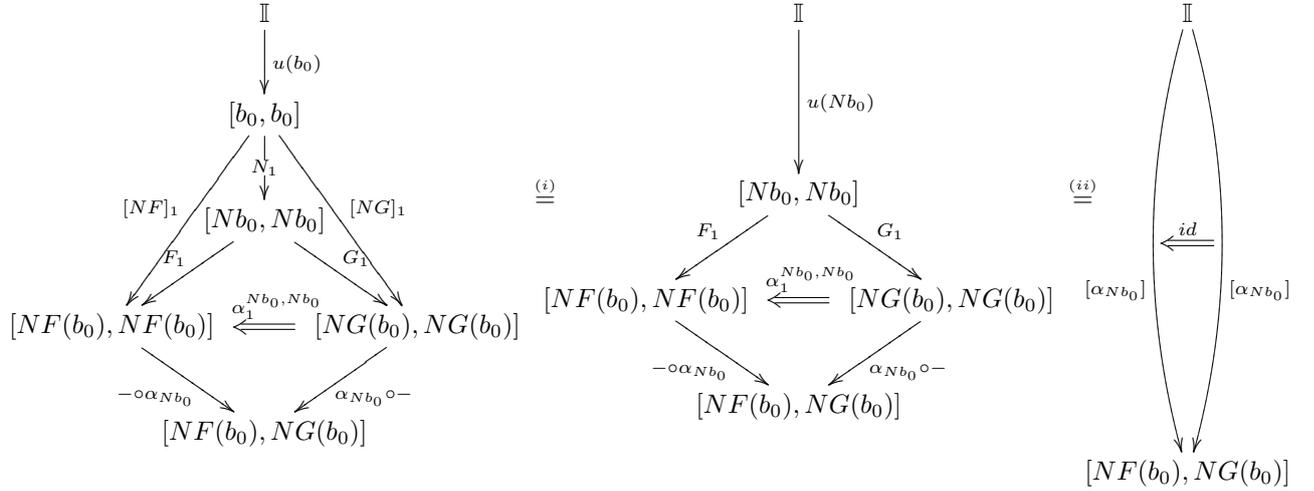
\begin{changemargin}{-10ex}{-10ex}{\small
$$
\raisebox{16ex}{\xymatrix@C=-6ex{&\I\ar[d]^{u(b_0)}\\
&[b_0,b_0]\ar[d]|{N_1}\ar[ddl]_{[NF]_1}\ar[ddr]^{[NG]_1}\\
&[Nb_0,Nb_0]\ar[dl]_{F_1}\ar[dr]^{G_1}\\
[NF(b_0),NF(b_0)]\ar[dr]_{-\circ\alpha_{Nb_0}}
&&[NG(b_0),NG(b_0)]\ar[dl]^{\alpha_{Nb_0}\circ-}\ar@{}[ll]|(.4){}="1"|(.6){}="2"\\
&[NF(b_0),NG(b_0)]
\ar@{=>}"1";"2"_(.3){\alpha_1^{Nb_0,Nb_0}}}}
\eq{i}\!\!\!\!\!\!
\raisebox{16ex}{\xymatrix@C=-6ex{&\I\ar[dd]^{u(Nb_0)}\\
{}\\
&[Nb_0,Nb_0]\ar[dl]_{F_1}\ar[dr]^{G_1}\\
[NF(b_0),NF(b_0)]\ar[dr]_{-\circ\alpha_{Nb_0}}
&&[NG(b_0),NG(b_0)]\ar[dl]^{\alpha_{Nb_0}\circ-}\ar@{}[ll]|(.4){}="1"|(.6){}="2"\\
&[NF(b_0),NG(b_0)]
\ar@{=>}"1";"2"_(.3){\alpha_1^{Nb_0,Nb_0}}}}
\eq{ii}\!\!\!\!\!\!
\raisebox{16ex}{\xymatrix@R=37ex{\I
\ar@/_3ex/[d]_(.6){[\alpha_{Nb_0}]}^{}="1"
\ar@/^3ex/[d]^(.6){[\alpha_{Nb_0}]}_{}="2"
\\
[NF(b_0),NG(b_0)] \ar@{=>}"2";"1"_{id}}}$$
}
\end{changemargin}
 $(i)$ holds by functoriality w.r.t. units (\ref{nfunct:units}), and $(ii)$ is simply functoriality w.r.t  units of $\alpha$
(\ref{ntrans:units}) for the object $Nb_0$ of $\C$.\\

Concerning composition axiom, take three objects $b_0$, $b_0'$
and $b_0''$ in $\B$, and consider the following diagram:
\begin{changemargin}{-10ex}{-10ex}{\small
$$\xymatrix@C=-10ex{
&&[b_0,b_0']\tm[b_0',b_0'']\ar[dl]_{id\tm N_1}\ar[dr]^{N_1\tm id}
\\
& [b_0,b_0']\tm[Nb_0',Nb_0'']\ar[dl]_{id\tm F_1}\ar[dd]^{id\tm G_1}
&&[Nb_0,Nb_0']\tm[b_0',b_0'']\ar[dr]^{G_1\tm id}\ar[dd]_{F_1\tm id}
\\
[b_0,b_0']\tm[NF(b_0'),NF(b_0'')]\ar[dd]_{id\tm(-\circ\alpha_{Nb_0''})}
&&&&[NG(b_0),NG(b_0')]\tm[b_0',b_0'']\ar[dd]^{(\alpha_{Nb_0})\tm id}\ar@{}[dl]|(.3){}="a1"|(.7){}="a2"
\\
&[b_0,b_0']\tm[NG(b_0'),NG(b_0'')]\ar[dl]^{\quad id\tm(\alpha_{Nb_0'}\circ-)}\ar@{}[ul]|(.3){}="b1"|(.7){}="b2"
&&[NF(b_0),NF(b_0')]\tm[b_0',b_0'']\ar[dr]_{(-\circ\alpha_{Nb_0'})\tm id\quad}
\\
[b_0,b_0']\tm[NF(b_0'),NG(b_0'')]\ar[d]_{N_1\tm id}
&&&&[NF(b_0),NG(b_0')]\tm[b_0',b_0'']\ar[d]^{id\tm N_1}
\\
[Nb_0,Nb_0']\tm[NF(b_0'),NG(b_0'')]\ar[d]_{F_1\tm id}
&&&&[NF(b_0),NG(b_0')]\tm[Nb_0',Nb_0'']\ar[d]^{id \tm G_1}
\\
[NF(b_0),NF(b_0')]\tm[NF(b_0'),NG(b_0'')]\ar[drr]_{\circ}
&&&&[NF(b_0),NG(b_0')]\tm[NG(b_0'),NG(b_0'')]\ar[dll]^{\circ}
\\
&&[NF(b_0),NG(b_0'')]
\ar@{=>}"a1";"a2"_{\alpha_1^{Nb_0,Nb_0'}\tm id\quad}
\ar@{=>}"b1";"b2"_{\quad id\tm\alpha_1^{Nb_0',Nb_0''}}}
$$
}
\end{changemargin}
By product interchange rules (see \emph{Lemma \ref{lemma:prod_interchange}}, when one of the components is an identity)
2-cells $id\times \alpha_1$ and $\alpha_1\times id$ can slide along $N_1\times id$ and $id\times N_1$ respectively,
in order to give:
\begin{changemargin}{-15ex}{-10ex}{\small
$$\xymatrix@C=-10ex{
&&[b_0,b_0']\tm[b_0',b_0'']\ar[dl]_{N_1\tm N_1}\ar[dr]^{N_1\tm N_1}
\\
& [Nb_0,Nb_0']\tm[Nb_0',Nb_0'']\ar[dl]_{id\tm F_1}\ar[dd]^{id\tm G_1}\ar@{=}[rr]
&&[Nb_0,Nb_0']\tm[Nb_0',Nb_0'']\ar[dr]^{G_1\tm id}\ar[dd]_{F_1\tm id}
\\
[Nb_0,Nb_0']\tm[NF(b_0'),NF(b_0'')]\ar[dd]_{id\tm(-\circ\alpha_{Nb_0''})}
&&&&[NG(b_0),NG(b_0')]\tm[Nb_0',Nb_0'']\ar[dd]^{(\alpha_{Nb_0})\tm id}\ar@{}[dl]|(.3){}="a1"|(.7){}="a2"
\\
&[Nb_0,Nb_0']\tm[NG(b_0'),NG(b_0'')]\ar[dl]^{\quad id\tm(\alpha_{Nb_0'}\circ-)}\ar@{}[ul]|(.3){}="b1"|(.7){}="b2"
&&[NF(b_0),NF(b_0')]\tm[Nb_0',Nb_0'']\ar[dr]_{(-\circ\alpha_{Nb_0'})\tm id\quad}
\\
[Nb_0,Nb_0']\tm[NF(b_0'),NG(b_0'')]\ar[d]_{F_1\tm id}
&&&&[NF(b_0),NG(b_0')]\tm[Nb_0',Nb_0'']\ar[d]^{id \tm G_1}
\\
[NF(b_0),NF(b_0')]\tm[NF(b_0'),NG(b_0'')]\ar[drr]_{\circ}
&&&&[NF(b_0),NG(b_0')]\tm[NG(b_0'),NG(b_0'')]\ar[dll]^{\circ}
\\
&&[NF(b_0),NG(b_0'')]
\ar@{=>}"a1";"a2"_{\alpha_1^{Nb_0,Nb_0'}\tm id\quad}
\ar@{=>}"b1";"b2"_{\quad id\tm\alpha_1^{Nb_0',Nb_0''}}}
$$
}
\end{changemargin}
now, just apply composition functoriality for $\alpha$ (\ref{ntrans:composition}) and get:
\begin{changemargin}{-10ex}{-10ex}
$$
\raisebox{16ex}{\xymatrix@C=-6ex{
&[b_0,b_0']\tm [b_0',b_0'']\ar[d]^{N_1\tm N_1}\\
&[Nb_0,Nb_0']\tm [Nb_0',Nb_0'']\ar[d]^{\circ}\\
&[Nb_0,Nb_0'']\ar[dl]_{F_1}\ar[dr]^{G_1}\\
[NF(b_0),NF(b_0'')]\ar[dr]_{-\circ\alpha_{Nb_0''}}
&&[NG(b_0),NG(b_0'')]\ar[dl]^{\alpha_{Nb_0}\circ-}\ar@{}[ll]|(.4){}="1"|(.6){}="2"\\
&[NF(b_0),NG(b_0'')]
\ar@{=>}"1";"2"_(.3){\alpha_1^{Nb_0,Nb_0''}}}}
=\raisebox{16ex}{\xymatrix@C=-6ex{
&[b_0,b_0']\tm [b_0',b_0'']\ar[d]^{\circ}\\
&[Nb_0,Nb_0'']\ar[d]^{N_1}\\
&[Nb_0,Nb_0'']\ar[dl]_{F_1}\ar[dr]^{G_1}\\
[NF(b_0),NF(b_0'')]\ar[dr]_{-\circ\alpha_{Nb_0''}}
&&[NG(b_0),NG(b_0'')]\ar[dl]^{\alpha_{Nb_0}\circ-}\ar@{}[ll]|(.4){}="1"|(.6){}="2"\\
&[NF(b_0),NG(b_0'')]
\ar@{=>}"1";"2"_(.3){\alpha_1^{Nb_0,Nb_0''}}}}
$$
\end{changemargin}
where the last equality is functoriality of $N$ w.r.t. compositions (\ref{nfunct:composition}).
And this completes the proof that left horizontal composition is well defined.

\subsection{Left-composition axioms}
Given the situation
$$
\xymatrix@C=12ex{\A\ar[r]^M&\B\ar[r]^N&\C\ar@/^4ex/[r]^F_{}="a1"\ar[r]|(.3)G^{}="a2"_{}="b1"\ar@/_4ex/[r]_H^{}="b2"&\D\ar@{=>}"a1";"a2"^{\alpha}\ar@{=>}"b1";"b2"^{\beta}}
$$
in $n$\textbf{Cat}, left-composition defined above satisfies axioms (L1) to (L4) of \emph{Proposition \ref{prop:sesqui-categories}}.\\

(L1) $$Id_{\C}\bullet^0 \alpha=\alpha$$
\begin{proof}
Let objects $c_0,c_0'$ of $\C$ be given. It is clear that
$$
[Id_{\C}\bullet^0 \alpha]_{c_0} \eq{def} \alpha_{Id_{\C}(c_0)}= \alpha_{c_0}
$$
and also that
$$
[Id_{\C}\bullet^0 \alpha]_1^{c_0,c_0'} \eq{def} [Id_{\C}]_1^{c_0,c_0'}\bullet^0 \alpha_1^{c_0,c_0'} \eq{1} Id_{\C_1(c_0,c_0')}\bullet^0 \alpha_1^{c_0,c_0'}\eq{2}\alpha_1^{c_0,c_0'}
$$
where $(1)$ comes from the definition of \emph{identity functors}, and $(2)$ is axiom (L1) for $(n-1)$\textbf{Cat}.\\
\end{proof}

(L2) $$MN\bullet^0 \alpha =M\bullet^0(N\bullet^0\alpha)$$
\begin{proof}
Let objects $a_0,a_0'$ of $\A$ be given. Then
$$
[MN\bullet^0 \alpha]_{a_0}\eq{def}\alpha_{MN(a_0)}=\alpha_{N(Ma_0)} \eq{def} [N\bullet^0\alpha]_{Ma_0} \eq{def} [M\bullet^0(N\bullet^0\alpha)]_{a_0}
$$
Furthermore,
\begin{eqnarray*}
  [MN\bullet^0\alpha]_1^{a_0,a_0'} &\eq{def}& [MN]_1^{a_0,a_0'}\bullet^0\alpha_1^{MN(a_0),MN(a_0')}\\
  &=& M_1^{a_0,a_0'}N_1^{Ma_0,Ma_0'}\bullet^0\alpha_1^{MN(a_0),MN(a_0')} \\
  &\eq{1}&  M_1^{a_0,a_0'}\bullet^0(N_1^{Ma_0,Ma_0'}\bullet^0\alpha_1^{N(Ma_0),N(Ma_0')})\\
  &\eq{def}& M_1^{a_0,a_0'}\bullet^0[N\bullet^0\alpha]_1^{Ma_0,Ma_0'} \\
  &\eq{def}& [M\bullet^0(N\bullet^0\alpha)]_1^{a_0,a_0'}
\end{eqnarray*}
where $(1)$ is axiom (L2) for $(n-1)$\textbf{Cat}.\\
\end{proof}

(L3) $$N\bullet^0 id_F=id_{NF}$$
\begin{proof}
Let objects $b_0,b_0'$ of $\B$ be given. Trivially,
$$
[N\bullet^0 id_F]_{b_0}\eq{def}[id_F]_{Nb_0}=[id_{NF}]_{b_0}
$$
and
\begin{eqnarray*}
  [N\bullet^0 id_F]_1^{b_0,b_0'} &\eq{def}& N_1^{b_0,b_0'}\bullet^0[id_F]_1^{Nb_0,Nb_0'}\eq{1}
N_1^{b_0,b_0'}\bullet^0 id_{F_1^{Nb_0,Nb_0'}}= \\
  &\eq{2}& id_{N_1^{b_0,b_0'}F_1^{Nb_0,Nb_0'}}=id_{[NF]_1^{b_0,b_0'}}\eq{def} [id_{NF}]_1^{b_0,b_0'}
\end{eqnarray*}
where $(1)$ comes from the definition of identity transformation and $(2)$ is axiom (L3) in  $(n-1)$\textbf{Cat}.\\
\end{proof}

(L4) $$N\bullet^0(\alpha\bullet^1\beta)=(N\bullet^0\alpha)\bullet^1(N\bullet^0\beta)$$
\begin{proof}
Let objects $b_0,b_0'$ of $\B$ be given. On objects:
$$
[N\bullet^0(\alpha\bullet^1\beta)]_{b_0}\eq{def}[\alpha\bullet^1\beta]_{Nb_0} =\alpha_{Nb_0}\circ \beta_{Nb_0}
\eq{def}[N\bullet^0\alpha]_{b_0}\circ [N\bullet^0\beta]_{b_0} = [(N\bullet^0\alpha)\bullet^1(N\bullet^0\beta)]_{b_0}
$$
On homs:
\begin{changemargin}{-10ex}{-10ex}
\begin{eqnarray*}
  [N\bullet^0(\alpha\bullet^1\beta)]_1^{b_0,b_0'} &\eq{def}& N_1^{b_0,b_0'}\bullet^0[\alpha\bullet^1\beta]_1^{Nb_0,Nb_0'} \\
   &\eq{1}& N_1^{b_0,b_0'}\bullet^0\Big(\big(\beta_1^{Nb_0,Nb_0'}\bullet^0(\alpha_{Nb_0}\circ-)\big)\bullet^1
   \big(\alpha_1^{Nb_0,Nb_0'}\bullet^0(-\circ\beta_1^{Nb_0,Nb_0'})\big) \Big)\\
   &\eq{2}& \big(N_1^{b_0,b_0'}\bullet^0\beta_1^{Nb_0,Nb_0'}\bullet^0(\alpha_{Nb_0}\circ-)\big)\bullet^1
   \big(N_1^{b_0,b_0'}\bullet^0\alpha_1^{Nb_0,Nb_0'}\bullet^0(-\circ\beta_1^{Nb_0,Nb_0'})\big) \\
   &\eq{def}& \big([N\bullet^0\beta]_1^{b_0,b_0'}\bullet^0([N\circ\alpha]_{b_0}\circ-)\big)\bullet^1
   \big([N\bullet^0\alpha]_1^{b_0,b_0'}\bullet^0(-\circ[N\bullet^0\beta]_{b_0'})\big) \\
   &\eq{3}& [(N\bullet^0\alpha)\bullet^1(N\bullet^0\beta)]_1^{b_0,b_0'}
\end{eqnarray*}
\end{changemargin}
where $(1)$ and $(3)$ hold by definition of vertical composites of 2-morphisms,
$(2)$ by axiom (L4) in $(n-1)$\textbf{Cat}.
\end{proof}

\subsection{Defining reduced right-composition}
Given the situation
$$
\xymatrix@C=12ex{\C\ar@/^4ex/[r]^F_{}="1"\ar@/_4ex/[r]_G^{}="2"&\D\ar[r]^L&\E\ar@{=>}"1";"2"^{\alpha}}
$$
one defines reduced horizontal composition $\alpha\bullet^0 L:FL\Rightarrow GL:\C\rightarrow\E$ (or $0$-composition)
in the following way:\\

$\bullet$ for every object $c_0$ in $\C$,
$$
[\alpha\bullet^0 L]_0 = L(\alpha_0(c_0)):L(F(c_0))\rightarrow L(G(c_0))
$$
$\bullet$ for every pair of objects $c_0, c_0'$ of $\B$, the diagram below describes $[\alpha\bullet^0 L]_1^{c_0,c_0'}$
by means of reduced right composition in $(n-1)$\textbf{Cat}:
$$
[\alpha\bullet^0 L]_1^{c_0,c_0'}=\alpha_1^{c_0,c_0'}\bullet^0L_1^{Fc_0,Gc_0'}
$$
$$
\xymatrix@C=-6ex{
&&\C_1(c_0,c_0')\ar[dl]_{F_1}\ar[dr]^{G_1}
\\
&\D_1(F(c_0),F(c_0'))\ar[dr]|{-\circ\alpha_{c_0'}}\ar@{.>}[dl]_{L_1}
&&\D_1(G(b_0),G(b_0'))\ar[dl]|{\alpha_{c_0}\circ-}\ar@{.>}[dr]^{L_1}\ar@{}[ll]|(.4){}="1"|(.6){}="2"
\\
\E_1(FL(c_0),FL(c_0'))\ar@{.>}[drr]_{-\circ L(\alpha_{c_0'})}&&\D_1(F(c_0),G(c_0'))\ar[d]^{L_1}
&&\E_1(GL(c_0),GL(c_0'))\ar@{.>}[dll]^{L(\alpha_{c_0})\circ -}
\ar@{=>}"1";"2"_(.3){\alpha_1^{c_0,c_0'}}\\
&&\E_1(FL(c_0),GL(c_0'))
}
$$

To prove that these data give indeed a $2$-morphism, unit and composition axioms equations
(\ref{ntrans:units}) (\ref{ntrans:composition}) must hold.\\

To this end, chose first an object $c_0$ of $\C$ and consider the
following chain of diagrams equalities:
{\small
$$
\raisebox{19ex}{\xymatrix@C=-6ex{
&\I\ar[d]^{u(c_0)}
\\
&[c_0,c_0]\ar[dl]_{F_1}\ar[dr]^{G_1}
\\
[Fc_0,Fc_0]\ar[dr]_{-\circ\alpha_{c_0}}
&&[Gc_0,Gc_0]\ar[dl]^{\alpha_{c_0}\circ-}\ar@{}[ll]|(.4){}="1"|(.6){}="2"
\\
&[Fc_0,Gc_0]\ar[d]^{L_1}
\\
&[FL(c_0),GL(c_0)]
\ar@{=>}"1";"2"_{\alpha_1^{c_0,c_0}}}}
\eq{i}
\raisebox{16ex}{\xymatrix{
\I\ar@/^3ex/[ddd]^{[\alpha_{c_0}]}_{}="1"
\ar@/_3ex/[ddd]_{[\alpha_{c_0}]}^{}="2"
\\
\\
\\
[Fc_0,Gc_0]\ar[d]^{L_1}
\\
[Fc_0,Gc_0]
\ar@{=>}"1";"2"_{id}}}
\eq{ii}
\raisebox{15ex}{\xymatrix{
\I\ar@/^3ex/[dddd]^{[L(\alpha_{c_0})]}_{}="1"
\ar@/_3ex/[dddd]_{[L(\alpha_{c_0})]}^{}="2"
\\
\\
\\
\\
[Fc_0,Gc_0]
\ar@{=>}"1";"2"_{id}}}
$$
}
where $(i)$ holds by unit functoriality of $\alpha$ (\ref{ntrans:units}), $(ii)$ by functoriality w.r.t. units (\ref{nfunct:units})
and by axiom (R3) of reduced right composition in the sesqui-category $(n-1)$\textbf{Cat}.\\

Concerning composition axiom, let us take three objects $c_0,c_0'$ and $c_0''$ in $\C$, and consider the following diagram:
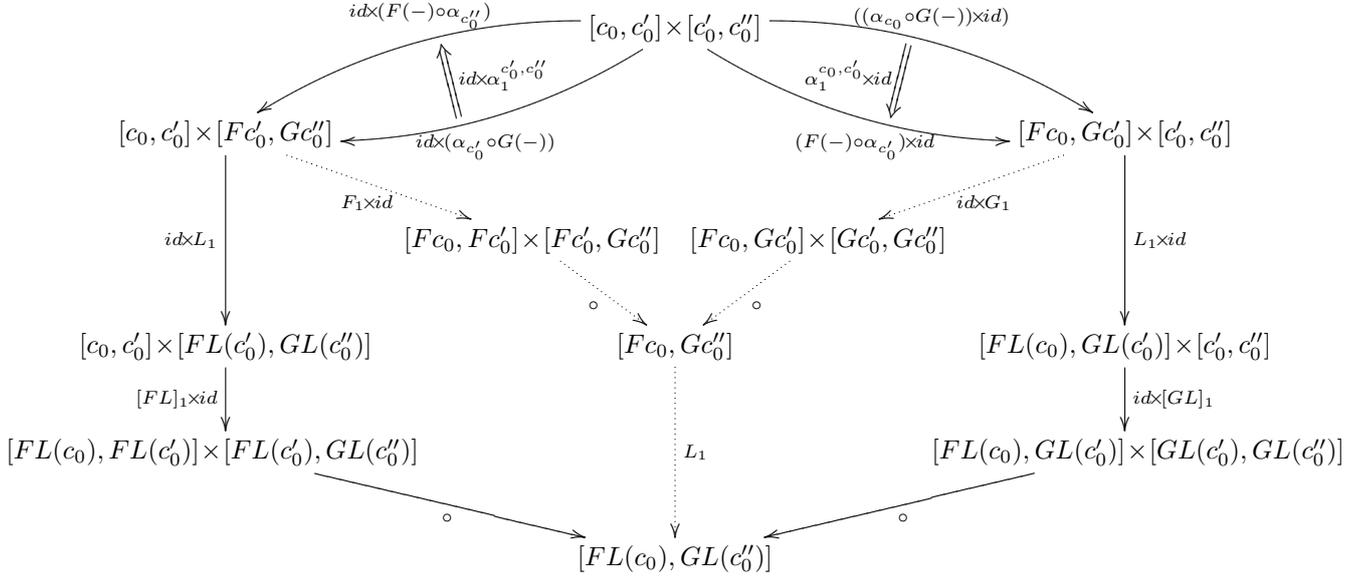
\begin{changemargin}{-15ex}{-10ex}{\small
$$
\xymatrix@C=-5ex{
&&&[c_0,c_0']\tm[c_0',c_0'']
\ar@/^4ex/[dlll]^{\ id\tm(\alpha_{c_0'}\circ G(-))}_{}="a1"
\ar@/_4ex/[dlll]_{id\tm(F(-)\circ\alpha_{c_0''})}^{}="a2"
\ar@/^4ex/[drrr]^{((\alpha_{c_0}\circ G(-))\tm id)}_{}="b1"
\ar@/_4ex/[drrr]_{(F(-)\circ\alpha_{c_0'})\tm id\ }^{}="b2"
\\
[c_0,c_0']\tm[Fc_0',Gc_0'']\ar[dd]_{id\tm L_1}\ar@{.>}[dr]_{F_1\tm id}
&&&&&&[Fc_0,Gc_0']\tm[c_0',c_0'']\ar@{.>}[dl]^{id\tm G_1}\ar[dd]^{L_1\tm id}
\\
&[Fc_0,Fc_0']\tm[Fc_0',Gc_0'']\ar@{.>}[drr]_{\circ}
&&&&[Fc_0,Gc_0']\tm[Gc_0',Gc_0'']\ar@{.>}[dll]^{\circ}
&
\\
[c_0,c_0']\tm[FL(c_0'),GL(c_0'')]\ar[d]_{[FL]_1\tm id}
&&&[Fc_0,Gc_0'']\ar@{.>}[dd]^{L_1}&&&[FL(c_0),GL(c_0')]\tm[c_0',c_0'']\ar[d]^{id\tm [GL]_1}
\\
[FL(c_0),FL(c_0')]\tm[FL(c_0'),GL(c_0'')]\quad\ar[drrr]_{\circ}
&&&&&&\quad[FL(c_0),GL(c_0')]\tm[GL(c_0'),GL(c_0'')]\ar[dlll]^{\circ}
\\
&&&[FL(c_0),GL(c_0'')]
\ar@{=>}"a1";"a2"_{id\tm\alpha_1^{c_0',c_0''}}
\ar@{=>}"b1";"b2"_{\alpha_1^{c_0,c_0'}\tm id}}
$$
}
\end{changemargin}
Here, internal dotted construction commutes with external (by product properties),  hence it can take its place
and suggests to apply composition functoriality (\ref{ntrans:composition}) for $\alpha$, in order to give
$$
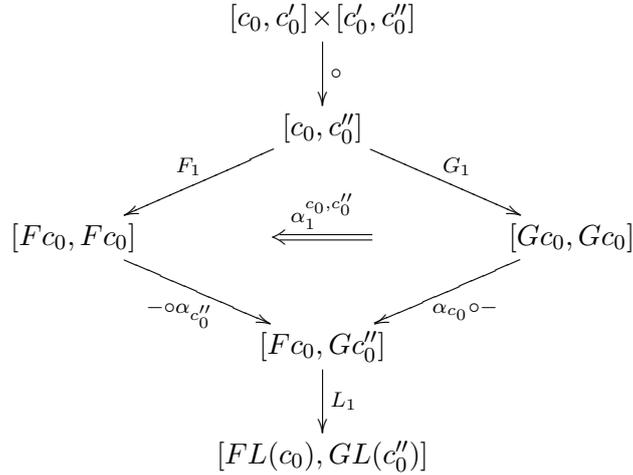
\xymatrix{
&[c_0,c_0']\tm[c_0',c_0'']\ar[d]^{\circ}
\\
&[c_0,c_0'']\ar[dl]_{F_1}\ar[dr]^{G_1}
\\
[Fc_0,Fc_0]\ar[dr]_{-\circ\alpha_{c_0''}}
&&[Gc_0,Gc_0]\ar[dl]^{\alpha_{c_0}\circ-}\ar@{}[ll]|(.4){}="1"|(.6){}="2"
\\
&[Fc_0,Gc_0'']\ar[d]^{L_1}
\\
&[FL(c_0),GL(c_0'')]
\ar@{=>}"1";"2"_{\alpha_1^{c_0,c_0''}}}
$$
and this complete the proof horizontal right composition is well defined.

\subsection{Right-composition axioms}
Given the diagram
$$
\xymatrix@C=12ex{\C\ar@/^4ex/[r]^F_{}="a1"\ar[r]|(.3)G^{}="a2"_{}="b1"\ar@/_4ex/[r]_H^{}="b2"&\D\ar[r]^L&\E\ar[r]^M&\F
\ar@{=>}"a1";"a2"^{\alpha}\ar@{=>}"b1";"b2"^{\beta}}
$$
in $n$\textbf{Cat}, right-composition defined above satisfies axioms (R1) to (R4) of \emph{Proposition \ref{prop:sesqui-categories}}.\\

(R1) $$\alpha\circ Id_{\D}=\alpha$$
\begin{proof}
Let objects $c_0,c_0'$ of $\C$ be given. It is clear that
$$
[\alpha\bullet^0 Id_{\D}]_{c_0}\eq{def} Id_{\D}(\alpha_{c_0})=\alpha_{c_0}
$$
and also that
$$
[\alpha\bullet^0 Id_{\D}]_1^{c_0,c_0'}\eq{def}\alpha_1^{c_0,c_0'}\bullet^0 [Id_{\D}]_1^{Fc_0,Gc_0'}\eq{1}
\alpha_1^{c_0,c_0'}\bullet^0 Id_{\D_1(Fc_0,Gc_0')}\eq{2}\alpha_1^{c_0,c_0'}
$$
where $(1)$ comes from the definition of \emph{identity functors}, and $(2)$ is axiom (R1) for $(n-1)$\textbf{Cat}.\\
\end{proof}

(R2) $$\alpha\bullet^0 LM=(\alpha\bullet^0 L)\bullet^0 M$$
\begin{proof}
Let objects $c_0,c_0'$ of $\C$ be given. Then
$$
[\alpha\bullet^0 LM]_{c_0}\eq{def}LM(\alpha_{c_0})=M(L(\alpha_{c_0}))\eq{def}M([\alpha\bullet^0 L]_{c_0})
\eq{def}[(\alpha\bullet^0 L)\bullet^0 M]_{c_0}
$$
Furthermore
\begin{eqnarray*}
  [\alpha\bullet^0 LM]_1^{c_0,c_0'} &\eq{def}& \alpha_1^{c_0,c_0'}\bullet^0 [LM]_1^{Fc_0,Gc_0'} \\
    &=& \alpha_1^{c_0,c_0'}\bullet^0( L_1^{Fc_0,Gc_0'} M_1^{L(Fc_0),L(Gc_0')})\\
    &\eq{1}& (\alpha_1^{c_0,c_0'}\bullet^0 L_1^{Fc_0,Gc_0'})\bullet^0 M_1^{L(Fc_0),L(Gc_0')} \\
    &\eq{def}& [\alpha\bullet^0 L]_1^{c_0,c_0'} \bullet^0 M_1^{L(Fc_0),L(Gc_0')} \\
    &\eq{def}& [(\alpha\bullet^0 L)\bullet^0 M]_1^{c_0,c_0'}
  \end{eqnarray*}
where (1) is axiom (R2) for $(n-1)$\textbf{Cat}.\\
\end{proof}

(L3) $$id_F\bullet^0 L = id_{FL}$$
\begin{proof}
Let objects $c_0,c_0'$ of $\C$ be given. Trivially,
$$
[id_F\bullet^0 L]_{c_0}\eq{def}=L([id_F]_{c_0})\eq{1}L(id_{Fc_0})\eq{2}id_{FL(c_0)}
$$
where $(1)$ holds by definition of identity transformations and $(2)$ from functoriality of $L$.
Furthermore,
\begin{eqnarray*}
  [id_F\bullet^0 L]_1^{c_0,c_0'} &\eq{def}& [id_F]_1^{c_0,c_0'}\bullet^0 L_1^{Fc_0,Fc_0'} \eq{1} id_{F_1^{c_0,c_0'}}\bullet^0 L_1^{Fc_0,Fc_0'} =\\
    &\eq{2}& id_{F_1^{c_0,c_0'}L_1^{Fc_0,Fc_0'}} = id_{[FL]_1^{c_0,c_0'}} \eq{def} [id_{FL}]_1^{c_0,c_0'}
\end{eqnarray*}
where $(1)$ comes from the definition of identity transformation and $(2)$ is axiom (R3) in  $(n-1)$\textbf{Cat}.\\
\end{proof}

(L4) $$(\alpha\bullet^1\beta)\bullet^0 L= (\alpha\bullet^0 L)\bullet^1(\beta\bullet^0 L) $$
\begin{proof}
Let objects $c_0,c_0'$ of $\C$ be given. On objects:
$$
[(\alpha\bullet^1\beta)\bullet^0 L]_{c_0} \eq{def} L([\alpha\bullet^1\beta]_{c_0}) = L(\alpha_{c_0}\circ\beta_{c_0})
=L(\alpha_{c_0})\circ L(\beta_{c_0})=[(\alpha\bullet^0 L)\bullet^1(\beta\bullet^0 L)]_{c_0}
$$
On homs:\\
$  [(\alpha\bullet^1\beta)\bullet^0 L]_1^{c_0,c0'} =$
\begin{eqnarray*}
    &\eq{def}& [\alpha\bullet^1\beta]_1^{c_0,c_0'}\bullet^0 L_1^{Fc_0,Hc_0'} \\
    &\eq{1}& \Big(\big( \beta_1^{c_0,c_0'}\bullet^0(\alpha_{c_0}\circ -)\big) \bullet^1 \big(\alpha_1^{c_0,c_0'} \bullet^0(-\circ \beta_{c_0'}) \big)\Big)\bullet^0 L_1^{Fc_0,Hc_0'} \\
    &\eq{2}& \Big(\big( \beta_1^{c_0,c_0'}\bullet^0(\alpha_{c_0}\circ -)\big)\bullet^0 L_1^{Fc_0,Hc_0'}\Big) \bullet^1 \Big(\big(\alpha_1^{c_0,c_0'} \bullet^0(-\circ \beta_{c_0'}) \big)\bullet^0 L_1^{Fc_0,Hc_0'}\Big) \\
    &\eq{3}& \Big( \beta_1^{c_0,c_0'}\bullet^0\big(L_1^{Gc_0,Hc_0'}(L(\alpha_{c_0})\circ -) \big)\Big)\bullet^1\Big(\alpha_1^{c_0,c_0'}\bullet^0\big(L_1^{Fc_0,Gc_0'}(-\circ L(\beta_{c_0'})) \big)\Big)  \\
    &\eq{4}& \Big(\big( \beta_1^{c_0,c_0'}\bullet^0 L_1^{Gc_0,Hc_0'}\big)\bullet^0(L(\alpha_{c_0})\circ -) \Big)\bullet^1\Big(\big(\alpha_1^{c_0,c_0'}\bullet^0 L_1^{Fc_0,Gc_0'}\big)\bullet^0(-\circ L(\beta_{c_0'})) \Big)\\
    &\eq{def}& \big( [\beta\bullet^0 L]_1^{c_0,c_0'}\bullet^0 ([\alpha\bullet^0 L]_{c_0}\circ -)\big)\bullet^1 \big([\alpha\bullet^0 L]_1^{c_0,c_0'}\bullet^0 (- \circ [\beta\bullet^0 L]_{c_0'})\big)\\
    &\eq{5}&[(\alpha\bullet^0 L)\bullet^1(\beta\bullet^0 L)]_1^{c_0,c_0'}
\end{eqnarray*}
where $(1)$ and $(5)$ hold by definition of vertical composites of 2-morphisms, $(2)$ by axiom (R4) in $(n-1)$\textbf{Cat},
$(3)$ by functoriality of $L$ w.r.t. $0$-composition, and $(4)$ by axiom (R2) in $(n-1)$\textbf{Cat}.\\
\end{proof}

\subsection{Whiskering axiom }
Given the situation
$$
\xymatrix@C=12ex{\B\ar[r]^N&\C\ar@/^4ex/[r]^F_{}="1"\ar@/_4ex/[r]_G^{}="2"&\D\ar[r]^L&\E\ar@{=>}"1";"2"^{\alpha}}
$$
a whiskering operation may be defined if the following equation holds:\\

(LR5)$$(N\circ\alpha)\circ L=N\circ(\alpha\circ L)$$
\begin{proof}
Let objects $b_0,b_0'$ of $\B$ be given. Then the following follows immediately from definitions
$$
[(N\bullet^0\alpha)\bullet^0 L]_{b_0}=L( [N\bullet^0\alpha]_{b_0})=L(\alpha_{Nb_0})=[\alpha\bullet^0 L]_{Nb_0}=[N\circ(\alpha\circ L)]_{b_0}
$$
Analogously, consider:
\begin{eqnarray*}
  [(N\bullet^0\alpha)\bullet^0 L]_1^{b_0,b_0'} &=& [N\bullet^0\alpha]_1^{b_0,b_0'}\bullet^0 L_1^{F(Nb_0),G(Nb_0')} \\
    &=& \Big(N_1^{b_0,b_0'}\bullet^0\alpha_1^{Nb_0,Nb_0'}\Big)\bullet^0 L_1^{F(Nb_0),G(Nb_0')} \\
    &\eq{1}& N_1^{b_0,b_0'}\bullet^0\Big(\alpha_1^{Nb_0,Nb_0'}\bullet^0 L_1^{F(Nb_0),G(Nb_0')}\Big) \\
    &=& N_1^{b_0,b_0'}\bullet^0 [\alpha\bullet^0 L]_1^{Nb_0,Nb_0'}\\
    &=& [N\bullet^0(\alpha\bullet^0 L)]_1^{b_0,b_0'}
\end{eqnarray*}
where everything comes directly from definitions, but $(1)$ that is exactly the whiskering in $(n-1)$\textbf{Cat}.
\end{proof}

\section{Products in $n$\textbf{Cat}}
In order to close the induction
on the definition of $n$\textbf{Cat}, all we need is to show that it admits finite products, according to
the 2-dimensional \emph{Universal Property \ref{UP:products}}, i.e. to show it admits binary products and terminal objects.  \\

\subsection{2-universality of categorical products}
Let two n-categories $\C$ and $\D$ be given. We know from Proposition \ref{prop:nCat_has_products} that the underlying
category $\lfloor n\mathbf{Cat}\rfloor$ admits a (standard) product of $\C$ and $\D$:
$$
\xymatrix{&\ar[dl]_{\Pi_{\C}}\C\times\D\ar[dr]^{\Pi_{\D}}\\\C&&\D}
$$
Now suppose we are given two 2-morphisms
$$
\alpha:A\Rightarrow A':\X \rightarrow \C\times\D,\quad \beta:B\Rightarrow B':\X \rightarrow \C\times\D
$$
According to \emph{Universal Property \ref{UP:products}}, what we want to prove is that there exists a unique 2-morphism
$$
\theta:T\Rightarrow T':\X\rightarrow \C\times\D
$$
such that
\begin{equation}\label{eqn:2dim_products}
\theta\bullet^0\Pi_{\C}=\alpha,\quad \theta\bullet^0\Pi_{\D}=\beta,
\end{equation}
First let us say that $T$ and $T'$ are  determined by 1-dimensional
universal property:

$T$ is such that (and univocally determined by) $\left\{\begin{array}{l}T\bullet^0\Pi_{\C}=A \\T\bullet^0\Pi_{\D}=B\end{array}\right.$,

$T'$ is such that (and univocally determined by) $\left\{\begin{array}{l}T'\bullet^0\Pi_{\C}=A' \\T'\bullet^0\Pi_{\D}=B'\end{array}\right.$.

More explicitly, for every pair of objects $x_0,x_0'$ in $\X$,
$$
T_0(x_0)=(A_0(x_0),B_0(x_0))
$$
and
$$
\xymatrix@C=16ex{
\X_1(x_0,x_0')\ar[r]^(.4){T_1^{x_0,x_0'}}\ar[dr]_(.4){\langle A_1^{x_0,x_0'},\B_1^{x_0,x_0'}\rangle }
&[\C\times\D]_1^{x_0,x_0'}((Ax_0,Bx_0),(Ax_0',Bx_0'))\ar@{=}[d]^{def}\\
&\C_1(Ax_0,ax_0')\times\D(Bx_0,Bx_0')}
$$

Similarly for $T'$: $T_0'(x_0)=(A_0'(x_0),B_0'(x_0))$
and $T_1'^{\ x_0,x_0'}=\langle A_1'{}^{x_0,x_0'},\B_1'{}^{x_0,x_0'}\rangle$

Then, $\theta=\langle\theta_0,\theta_1\rangle$ is given by:
$$
\theta_0(x_0)=(\alpha_0(x_0),\beta_0(x_0))
$$
and $\theta_1^{x_0,x_0'}$ is given by the universal property of products in
$(n-1)$\textbf{Cat}.

In fact, a suitable $\theta_1$ would fit in the following diagram:
$$
\xymatrix@C=-13ex{
&\X_1(x_0,x_0')\ar[dl]_{T_1^{x_0,x_0'}}\ar[dr]^{T_1'{}^{x_0,x_0'}}
\\
[\C\times\D]_1((Ax_0,Bx_0),(Ax_0',Bx_0'))\ar[dr]_{-\circ\theta_{x_0'}}
&&[\C\times\D]_1((A'x_0,B'x_0),(A'x_0',B'x_0'))
\ar[dl]^{\theta_{x_0}\circ-}\ar@{}[ll]|(.4){}="1"|(.6){}="2"
\\
&[\C\times\D]_1((Ax_0,Bx_0),(A'x_0',B'x_0'))\ar@{:>}"1";"2"_{?}}
$$
Spelling out the definitions, this may be written:
$$
\xymatrix@C=-13ex{
&\X_1(x_0,x_0')
\ar[dl]_{\langle A_1^{x_0,x_0'},B_1^{x_0,x_0'}\rangle\quad }
\ar[dr]^{\quad\langle A_1'{}^{x_0,x_0'},B_1'{}^{x_0,x_0'}\rangle}
\\
\C_1(Ax_0,Ax_0') \times\D_1(Bx_0,Bx_0')
\ar[dr]_{(-\circ\alpha_{x_0'})\tm(-\circ\beta_{x_0'})\quad}
&&\ \C_1(A'x_0,A'x_0') \times\D_1(B'x_0,B'x_0')
\ar[dl]^{\quad(\alpha_{x_0}\circ-)\tm(\beta_{x_0}\circ-)}\ar@{}[ll]|(.42){}="1"|(.58){}="2"
\\
&\C_1(Ax_0,A'x_0') \times\D_1(Bx_0,B'x_0')\ar@{=>}"1";"2"_{\theta_1^{x_0,x_0'}}}
$$
hence we are allowed to define
$$
\theta_1^{x_0,x_0'}=\langle\alpha_1^{x_0,x_0'},\beta_1^{x_0,x_0'}\rangle,
$$
and this choice would satisfy the universal property of products..

Concerning the first of the (\ref{eqn:2dim_products}), for every pair of objects $x_0, x_0'$ in $\X$
$$
[\theta\bullet^0 \Pi_{\C}]_0(x_0)=[\Pi_{\C}]_0(\theta_{x_0})=\pi_{\C_0}^{\C_0\tm\D_0}(\alpha_{x_0},\beta_{x_0})=\alpha_{x_0}
$$
and also
$$
[\theta \bullet^0 \Pi_{\C}]_1^{x_0,x_0'}=\theta_1^{x_0,x_0'}\bullet^0 [\Pi_{\C}]_1^{(Ax_0,Bx_0),(A'x_0',B'x_0')}=
$$
$$
=\langle\alpha_1^{x_0,x_0'},\beta_1^{x_0,x_0'}\rangle\bullet^0\Pi_{\C_1(Ax_0,A_0'x_0')}=\alpha_1^{x_0,x_0'}
$$
where everything comes directly from definitions, but the last equality which is given by the universal property defining
$\theta_1^{x_0,x_0'}$ in $(n-1)$\textbf{Cat}.
The second of the (\ref{eqn:2dim_products}) can be proved the same way.

Moreover, such a $\theta$ is unique. For, if another
$$
\eta:T\Rightarrow T':\X\rightarrow \C\times\D
$$
is such that
$$
\eta\bullet^0\Pi_{\C}=\alpha,\quad \eta\bullet^0\Pi_{\D}=\beta,
$$
then equations above determine $\eta_0$ on objects (since a map to a product, $\C_0\times\D_0$, is determined by its projections),
hence it is equal to $\theta_0$. On the other side, for objects $x_0,x_0'$ in $\X$, and $\eta_1^{x_0,x_0'}$ is given considering its
composition with projections by universality of the product in $(n-1)$\textbf{Cat}, hence it is equal to $\theta_1^{x_0,x_0'}$.\\

To conclude this section, we will see that the just defined pair $[\theta_0,\theta_1^{-,-}]$ is indeed a 2-morphism, that is, it obeys units
and composition axioms for n-transformations.\\

To this purpose, let us chose a triple of objects $x_0, x_0'$ and $x_0''$ of $\X$, and let us consider the following diagram:
\begin{changemargin}{-15ex}{-10ex}
$$
\xymatrix@C=-10ex{
&[x_0,x_0']\tm[x_0',x_0'']\ar@/_4ex/[dl]_{id\tm L^{x_0',x_0''}}^{}="id2"\ar@/^4ex/[dl]^(.3){id\tm L'{}^{x_0',x_0''}}_{}="id1"
 \ar@/_4ex/[dr]_(.3){L^{x_0',x_0''}\tm id}^{}="2id"\ar@/^4ex/[dr]^{L'{}^{x_0',x_0''}\tm id}_{}="1id"
\\
[x_0,x_0']\tm[Ax_0',A'x_0'']\tm[Bx_0',B'x_0'']\ar[d]_{\langle A_1^{x_0,x_0'},B_1^{x_0,x_0'}\rangle\tm id}
&&[Ax_0,A'x_0']\tm[Bx_0,B'x_0']\tm[x_0',x_0'']\ar[d]^{ id\tm\langle A_1'{}^{x_0',x_0''},B_1'{}^{x_0',x_0''} \rangle}
\\
[Ax_0,Ax_0']\tm[Bx_0,Bx_0']\tm[Ax_0',A'x_0'']\tm[Bx_0',B'x_0'']\ar[d]_{\tau}
&&[Ax_0,A'x_0']\tm[Bx_0,B'x_0']\tm[A'x_0',A'x_0'']\tm[B'x_0',B'x_0'']\ar[d]^{\tau}
\\
[Ax_0,Ax_0']\tm[Ax_0',A'x_0'']\tm[Bx_0,Bx_0']\tm[Bx_0',B'x_0'']\ar[dr]_{{}^{\C}\circ\times {}^{\D}\circ}
&&[Ax_0,A'x_0']\tm[A'x_0',A'x_0'']\tm[Bx_0,B'x_0']\tm[B'x_0',B'x_0'']\ar[dl]^{{}^{\C}\circ\times {}^{\D}\circ}
\\
&[Ax_0,A'x_0'']\tm[Bx_0,B'x_0'']
\ar@{=>}"id1";"id2"_{id\tm \theta_1^{x_0',x_0''}}
\ar@{=>}"1id";"2id"_{\theta_1^{x_0,x_0'}\tm id}}
$$
\end{changemargin}
where
\begin{eqnarray*}
  L^{x,y} &=& \langle A_1^{x,y},B_1^{x,y}\rangle \Big( (-\circ\alpha_{y})\times(-\circ\beta_{y})\Big) \\
  L'{}^{x,y} &=& \langle A_1'{}^{x,y},B_1'{}^{x,y}\rangle \Big( (\alpha_{x}\circ-)\times(\beta_{x}\circ-) \Big)
\end{eqnarray*}

Now, considering only the right branch of the diagram, where $\Delta$'s express diagonal morphisms,
 the following equations hold by product interchange properties:\\

$\big(\theta_1^{x_0,x_0'}\times id\big)\bullet^0 \big(id\times\langle A_1'{}^{x_0',x_0''},B_1'{}^{x_0',x_0''} \rangle \big)\ \tau\ \big({}^{\C}\!\circ\times {}^{\D}\!\circ\big)=$
\begin{eqnarray*}
   &=&  \big(\langle\alpha_1^{x_0,x_0'},\beta_1^{x_0,x_0'}\rangle\times id\big)\bullet^0 \big(id\times\langle A_1'{}^{x_0',x_0''},B_1'{}^{x_0',x_0''} \rangle \big)\ \tau\ \big({}^{\C}\!\circ\times {}^{\D}\!\circ\big)\\
   &=& (\Delta\times id)\bullet^0 \big(\alpha_1^{x_0,x_0'}\times\beta_1^{x_0,x_0'}\times id\big) \bullet^0 (id\times\Delta)(id\times A_1'{}^{x_0',x_0''}\times B_1'{}^{x_0',x_0''} )\ \tau\ \big({}^{\C}\!\circ\times {}^{\D}\!\circ\big)\\
   &=& (\Delta\times\Delta)\ \tau\bullet^0 \big(\alpha_1^{x_0,x_0'}\times id\times\beta_1^{x_0,x_0'}\times id\big)\bullet^0 (id\times A_1'{}^{x_0',x_0''}\times id \times B_1'{}^{x_0',x_0''}) \big({}^{\C}\!\circ\times {}^{\D}\!\circ\big)\\
   &=& (\Delta\times\Delta)\ \tau\bullet^0 \Big(\big(\alpha_1^{x_0,x_0'}\times A_1'{}^{x_0',x_0''}\bullet^0({}^{\C}\!\circ)\big)\times \big(\beta_1^{x_0,x_0'}\times B_1'{}^{x_0',x_0''}\bullet^0({}^{\D}\!\circ) \big)\Big)
\end{eqnarray*}
Similarly, the left branch gives:\\

$\big(id\times \theta_1^{x_0',x_0''} \big)\bullet^0 \big(\langle A_1{}^{x_0',x_0'},B_1{}^{x_0,x_0'}\times id \rangle \big)\ \tau\ \big({}^{\C}\!\circ\times {}^{\D}\!\circ\big)=$
$$
=(\Delta\times\Delta)\ \tau\bullet^0\Big(\big(A_1{}^{x_0,x_0'}\times\alpha_1^{x_0',x_0''}\bullet^0({}^{\C}\!\circ)\big)\times \big(B_1{}^{x_0,x_0'}\times\beta_1^{x_0',x_0''}\bullet^0({}^{\D}\!\circ) \big)\Big)\qquad
$$

Hence the diagram above, being the vertical composite of the two branches, may be rewritten applying rule (L4):\\

$=(\Delta\times\Delta)\ \tau\bullet^0 \Big[\Big(\big(\alpha_1^{x_0,x_0'}\times A_1'{}^{x_0',x_0''}\bullet^0({}^{\C}\!\circ)\big)\bullet^1\big(A_1{}^{x_0,x_0'}\times\alpha_1^{x_0',x_0''}\bullet^0({}^{\C}\!\circ)\big)\Big)\times$
$$
\times\Big(\big(\beta_1^{x_0,x_0'}\times B_1'{}^{x_0',x_0''}\bullet^0({}^{\D}\!\circ) \big)\bullet^1 \big(B_1{}^{x_0,x_0'}\times\beta_1^{x_0',x_0''}\bullet^0({}^{\D}\!\circ) \big)\Big)\Big]
$$
by composition axiom of $\alpha$ and $\beta$, and product properties, we get the result:
\begin{eqnarray*}
   &=& (\Delta\times\Delta)\ \tau\bullet^0 \Big((({}^{\X}\!\circ)\bullet^0\alpha_1^{x_0,x_0''})\times(({}^{\X}\!\circ)\bullet^0\beta_1^{x_0,x_0''})\Big) \\
   &=& ({}^{\X}\!\circ)\ \Delta\bullet^0(\alpha_1^{x_0,x_0''}\times\beta_1^{x_0,x_0''}) \\
   &=& ({}^{\X}\!\circ)\bullet^0\langle\alpha_1^{x_0,x_0''},\beta_1^{x_0,x_0''}\rangle\\
   &=& ({}^{\X}\!\circ)\bullet^0 \theta_1^{x_0,x_0''}
\end{eqnarray*}

\vskip3ex
Furthermore, for an object $x_0$, consider the following diagram:
$$
\xymatrix@R=10ex{
\I \ar[d]^{u(x_0)}\\
[x_0,x_0]
\ar@/^5ex/[d]^{L'{}^{x_0,x_0}}_{}="1"
\ar@/_5ex/[d]_{L^{x_0,x_0}}^{}="2"\\
[Ax_0,A'x_0]\times[Bx_0,B'x_0]
\ar@{}"1";"2"|(.2){}="a1"|(.8){}="a2"
\ar@{=>}"a1";"a2"_{\theta_1^{x_0,x_0}}}
$$
then, as for the composition axiom,
\begin{eqnarray*}
  u(x_0)\bullet^0\theta_1^{x_0,x_0} &=& u(x_0)\circ\langle\alpha_1^{x_0,x_0},\beta_1^{x_0,x_0} \rangle \\
   &=& u(x_0)\Delta\bullet^0(\alpha_1^{x_0,x_0}\times\beta_1^{x_0,x_0}) \\
   &=& \Delta u(x_0) \bullet^0(\alpha_1^{x_0,x_0}\times\beta_1^{x_0,x_0})  \\
   &=& \Delta \bullet^0 \big((u(x_0) \bullet^0\alpha_1^{x_0,x_0})\times(u(x_0) \bullet^0\beta_1^{x_0,x_0})\big) \\
   &\eq{*}& \Delta \bullet^0 \big(Id_{[\alpha_{x_0}]}\times Id_{[\beta_{x_0}]}\big) \\
   &=& Id_{\langle[\alpha_{x_0}],[\beta_{x_0}]\rangle}= Id_{[\theta_{x_0}]}
\end{eqnarray*}
where $(*)$ is given by unit axiom of $\alpha$ and $\beta$.

\section{The standard $h$-pullback in $n$\textbf{Cat}}\label{sec:h-pb}
In this section we give a construction that will be of fundamental importance for the development of
the theory. The idea is to generalize a classical homotopical construction \cite{MR0402694}
to $n$-categories, or better to $n$-groupoids, where homotopical aspects are more than a mere suggestion.\\

We start considering the following $h$-pullbacks reference diagram.

$$
\xymatrix{
\PP\ar[r]^Q\ar[d]_P
&\C\ar[d]^G
\\
\A\ar[r]_F\ar@{}[ur]|(.3){}="1"|(.7){}="2"
&\B
\ar@{=>}"1";"2"^{\varepsilon}}
$$

For n=0, classical pullback in \textbf{Set} is an instance of $h$-pullback,
 with 2-morphism $\varepsilon$ being an identity.

In fact the category of sets and maps is (seen as) the 2-trivial sesqui-category $0$\textbf{Cat},
and indeed, only condition (1) and (2) survive.
Hence, in the next sections, we will suppose integer $n>0$ been given.\\

We exhibit  a recursive construction of the standard $h$-pullback satisfying {\em Universal Property \ref{UP:h-pullbacks}}.

We will give $\PP$ in the form$(\PP_0,\PP_1^{-,-})$. \\

The set $\PP_0$ is  the following limit in \textbf{Set} (that yields indeed, also the object-components
of $F$, $G$ and $\varepsilon$):
$$
\xymatrix{
&&\PP_0\ar@{.>}[dll]_{P_0}\ar@{.>}[d]^{\varepsilon_0}\ar@{.>}[drr]^{Q_0}
\\
\A_0\ar[dr]_{F_0}
&&\B_1\ar[dl]^{s}\ar[dr]_{t}
&&\C_0\ar[dl]^{G_0}
\\
&\B_0&&\B_0
}
$$
Here $\B_1$ is the disjoint union
$$\coprod_{b_0,b_0'\in \B_0}[\B_1(b_0,b_0')]_0,$$
 and $s,t$ are \emph{source} and \emph{target} maps of 1-cells.
More explicitly,
$$
\PP_0=\{(a_0,b_1,c_0)\quad s.t.\quad a_0\in\A_0, c_0\in\C_0, b_1:Fa_0\to Gc_0\in\B_1\}
$$
$$
P_0((a_0,b_1,c_0))=a_0,\quad Q_0((a_0,b_1,c_0))=c_0,\quad \varepsilon_0((a_0,b_1,c_0))=b_1
$$
\vskip4ex
Let us fix two element of $\PP_0$:
$$
p_0=(a_0,b_1,c_0),\quad p_0'=(a_0',b_1',c_0').
$$
The hom-(n-1)category $\PP_1(p_0,p_0')$ is granted by the following h-pullback in $(n-1)$\textbf{Cat}:
$$
\xymatrix@C=10ex{
\PP_1(p_0,p_0')\ar@{.>}[rr]^{Q_1^{p_0,p_0'}}\ar@{.>}[dd]_{P_1^{p_0,p_0'}}
&&\C_1(c_0,c_0')\ar[d]^{G_1^{c_0,c_0'}}\ar@{}[ddll]|(.3){}="1"|(.7){}="2"
\\
&&\B_1(Gc_0,Gc_0')\ar[d]^{b_1\circ-}
\\
\A_1(a_0,a_0')\ar[r]_{F_1^{a_0,a_0'}}
&\B_1(Fa_0,Fa_0')\ar[r]_{-\circ b_1'}
&\B_1(Fa_0,Gc_0')
\ar@{:>}"1";"2"_{\varepsilon_1^{p_0,p_0'}}}
$$
\begin{Proposition}
The pair $(\PP_0,\PP_1^{-,-})$ yields a n-category.
\end{Proposition}
In order to prove the statement, we need to show constructions for composition and units, and
to prove that they satisfies $n$-category axioms.
\subsection{Composition}
Suppose we are given three elements of $\PP_0$
$$
p_0=(a_0,b_1,c_0),\quad p_0'=(a_0',b_1',c_0'),\quad p_0'=(a_0'',b_1'',c_0'').
$$
One defines
$$
\circ^0:\PP_1(p_0,p_0')\times\PP_1(p_0',p_0'')\rightarrow\PP_1(p_0,p_0'')
$$
by means of the universal property of pullback in $(n-1)$\textbf{Cat}.

In fact, as $\PP_1(p_0,p_0'')$ is a pullback, we can consider the four-tuple
$$
\Big(\PP_1(p_0,p_0')\times\PP_1(p_0,p_0'),\ P_1^{p_0,p_0'}\times P_1^{p_0',p_0''},\ Q_1^{p_0,p_0'}\times Q_1^{p_0',p_0''},\ \varepsilon^{p_0,p_0',p_0''}\Big)
$$
where the (n-1)-natural transformation $\varepsilon^{p_0,p_0',p_0''}$ is the composite shown below:
\begin{changemargin}{-10ex}{-10ex}{\small
\begin{equation}\label{def:epsilon123}
\xymatrix@C=-6ex{
&[p_0,p_0']\tm[a_0',a_0'']\ar[dd]_(.7){id\tm F_1}^(.8){}="s2"\ar@{.>}[dl]_{P_1\tm id}
&&\ar[ll]_{id\tm P_1}\ar[dl]^{id\tm Q_1}
[p_0,p_0']\tm[p_0',p_0''] \ar[dr]_{P_1\tm id}\ar[rr]^{Q_1\tm id}
&&[c_0,c_0']\tm[p_0',p_0'']\ar[dd]^(.7){G_1\tm id}_(.8){}="d1"\ar@{.>}[dr]^{id\tm Q_1}
\\
[a_0,_0']\tm[a_0',a_0'']\ar@{.>}[ddddd]_{\circ}
&&[p_0,p_0']\tm[c_0',c_0'']\ar[dd]^{id\tm G_1}_(.4){}="s1"
&&[a_0,a_0']\tm[p_0',p_0'']\ar[dd]_{F_1\tm id}^(.4){}="d2"
&&[c_0,c_0']\tm[c_0',c_0'']\ar@{.>}[ddddd]^{\circ}
\\
&[p_0,p_0']\tm[Fa_0',Fa_0'']\ar[dd]_{id\tm (-\circ b_1'')}
&&&&[Gc_0,Gc_0']\tm[p_0',p_0'']\ar[dd]^{(b_1\circ-)\tm id}
\\
&&[p_0,p_0']\tm[Gc_0',Gc_0'']\ar[dl]^{id\tm (b_1'\circ-)}
&&[Fa_0,Fa_0']\tm[p_0',p_0'']\ar[dr]_{(-\circ b_1')\tm id}
&&
\\
&[p_0,p_0']\tm[Fa_0',Gc_0'']\ar[d]_{P_1\tm id}
&&&&[Fa_0,Gc_0']\tm[p_0',p_0'']\ar[d]^{id\tm Q_1}
\\
&[a_0,a_0']\tm[Fa_0',Gc_0'']\ar[dr]_{F_1\tm id}
&&&&[Fa_0,Gc_0']\tm[c_0',c_0'']\ar[dl]^{id\tm G_1}
\\
[a_0,a_0'']\ar@{.>}[dr]_{F_1}
&&[Fa_0,Fa_0']\tm[Fa_0',Gc_0'']\ar[dr]_{\circ}
&&[Fa_0,Gc_0']\tm[Gc_0',Gc_0'']\ar[dl]^{\circ}
&&[c_0,c_0'']\ar@{.>}[dl]^{G_1}
\\
&[Fa_0,Fa_0'']\ar@{.>}[rr]_{-\circ b_1''}
&&[Fa_0,Gc_0'']&&
[Gc_0,Gc_0'']\ar@{.>}[ll]^{b_1\circ-}
\ar@{}"s1";"s2"|(.3){}="ss1"|(.7){}="ss2"
\ar@{}"d1";"d2"|(.3){}="dd1"|(.7){}="dd2"
\ar@{=>}"ss1";"ss2"_{id\tm \varepsilon_1^{p_0',p_0''}}
\ar@{=>}"dd1";"dd2"_{\varepsilon_1^{p_0,p_0'}\tm id}
}
\end{equation}
}\end{changemargin}
Dotted outer border is clearly equal to continuous inner border, hence, by universal property in $(n-1)$\textbf{Cat},
there exists a unique
\begin{equation}
\xymatrix{{}^{\PP}\!\circ^0: \PP_1(p_0,p_0')\times \PP_1(p_0',p_0'')\ar[r]&\PP_1(p_0,p_0'')}
\end{equation}
such that two squares below commute
\begin{equation}\label{eqn:PQ_funct}
\xymatrix{
\A_1(a_0,a_0')\times\A_1(a_0',a_0'')\ar[d]_{{}^{\A}\!\circ}
&\ar[l]_{P_1\tm P_1}\PP_1(p_0,p_0')\times\PP_1(p_0',p_0'')\ar@{.>}[d]^{{}^{\PP}\!\circ}\ar[r]^{Q_1\tm Q_1}
&\C_1(c_0,c_0')\times\C_1(c_0,c_0')\ar[d]^{{}^{\C}\!\circ}
\\
\A_1(a_0,a_0'')
&\PP_1(p_0,p_0'')\ar[l]^{P_1}\ar[r]_{Q_1}
&\C_1(c_0,c_0'')
}
\end{equation}
and
\begin{equation}\label{eqn:epsilon123_eq_comp_epsilon}
\raisebox{12ex}{
\xymatrix@R=8ex{
\PP_1(p_0,p_0')\times\PP_1(p_0',p_0'')\ar[d]^{{}^{\PP}\!\circ}\\
\PP_1(p_0,p_0'')\ar@/^4ex/[d]_{}="1"\ar@/_4ex/[d]^{}="2"
\ar@{=>}"1";"2"_{\varepsilon_1^{p_0,p_0''}}\\
\B_1(Fa_0,Gc_0'')
}}
=
\raisebox{12ex}{
\xymatrix@R=20ex{\PP_1(p_0,p_0')\times\PP_1(p_0',p_0'')
\ar@/^4ex/[d]_{}="1"\ar@/_4ex/[d]^{}="2"
\ar@{=>}"1";"2"_{\varepsilon^{p_0,p_0',p_0''}}\\
\B_1(Fa_0,Gc_0'')
}}
\end{equation}
\begin{Lemma}
Composition ${}^{\PP}\!\circ^0$ defined above is associative, {\em i.e.} the diagram below commutes in $(n-1)\mathbf{Cat}$,
for every four-tuple $(p_0, p_0',p_0'',p_0''')$ of elements of $\PP_0$:
\begin{equation}\label{lemma:pb_comp_ass}
\xymatrix@C=10ex{\PP_1(p_0,p_0')\times \PP_1(p_0',p_0'')\times \PP_1(p_0'',p_0''')
\ar[r]^(.6){id\times{}^{\PP}\circ}\ar[d]_{{}^{\PP}\circ\times id}
&\PP_1(p_0,p_0')\times \PP_1(p_0',p_0''')\ar[d]^{{}^{\PP}\circ}
\\
\PP_1(p_0,p_0'')\times \PP_1(p_0'',p_0''')\ar[r]_{{}^{\PP}\circ}
&
\PP_1(p_0,p_0''')
}
\end{equation}
\end{Lemma}
\begin{proof}

Let us consider the diagrams
$$
\xymatrix@C=12ex{
[p_0,p_0']\tm[p_0',p_0'']\tm[p_0'',p_0''']\ar@{.>}[dr]^K\ar[rr]^{Q_1\tm Q_1\tm Q_1}\ar[dd]_{P_1\tm P_1\tm P_1}
&&[c_0,c_0']\tm[c_0',c_0'']\tm[c_0'',c_0''']\ar[d]^{\Xi(\C)}
\\
&[p_0,p_0''']\ar[r]^{Q_1}\ar[d]_{P_1}
&[c_0,c_0''']\ar[d]\ar@{}[dl]|(.3){}="1"|(.7){}="2"
\\
[a_0,a_0']\tm[a_0',a_0'']\tm[a_0'',a_0''']\ar[r]_(.6){\Xi(\A)}
&[a_0,a_0''']\ar[r]
&[Fa_0,Gc_0''']
\ar@{=>}"1";"2"_{\varepsilon_1^{p_0,p_0'''}}
}
$$
where morphism $\Xi(\mathcal{X})$ can be either the composite
$({}^{\mathcal{X}}\!\circ\times id){}^{\mathcal{X}}\!\circ$ or
$(id \times{}^{\mathcal{X}}\!\circ){}^{\mathcal{X}}\!\circ$, with $\mathcal{X}$ being $\A$ or $\C$.

The idea of the proof is to use again universal property in $(n-1)$\textbf{Cat} to get unique
$$
K:\PP_1(p_0,p_0')\times \PP_1(p_0',p_0'')\times \PP_1(p_0'',p_0''')\rightarrow \PP_1(p_0,p_0''')
$$
that coincides with both composites of diagram \ref{lemma:pb_comp_ass}.
To this end, it suffices to show that the four-tuples
$$
\left(
  \begin{array}{c}
    \PP_1(p_0,p_0')\times \PP_1(p_0',p_0'')\times \PP_1(p_0'',p_0''') \\
    (P_1\times P_1\times P_1)({}^{\A}\!\circ\times id){}^{\A}\!\circ \\
    (Q_1\times Q_1\times Q_1)({}^{\C}\!\circ\times id){}^{\C}\!\circ \\
    ({}^{\PP}\!\circ\times id){}^{\PP}\!\circ\ \varepsilon_1^{p_0,p_0'''} \\
  \end{array}
\right)
$$
and
$$
\left(
  \begin{array}{c}
    \PP_1(p_0,p_0')\times \PP_1(p_0',p_0'')\times \PP_1(p_0'',p_0''') \\
    (P_1\times P_1\times P_1)(id\times {}^{\A}\!\circ){}^{\A}\!\circ \\
    (Q_1\times Q_1\times Q_1)(id\times {}^{\C}\!\circ){}^{\C}\!\circ \\
    (id\times {}^{\PP}\!\circ){}^{\PP}\!\circ\ \varepsilon_1^{p_0,p_0'''} \\
  \end{array}
\right)
$$
are equal.

First components are identical.

Equality of second components amounts to the associativity axiom for 0-composition in $\A$.

Equality of third components amounts to the associativity axiom for 0-composition in $\C$.

What remains to prove is equality of fourth components. The starting point is the 2-morphism
$({}^{\PP}\!\circ\tm id)\bullet^0{}^{\PP}\!\circ\bullet^0\varepsilon_1^{p_0,p_0'''}$, where
$$
\varepsilon_1^{p_0,p_0'''}:b_1\circ [QG]_1(-)\Rightarrow [PF]_1(-)\circ b_1'''
$$
This may be visualized as a diagram:
$$
\xymatrix@C=3ex{
&[p_0,p_0']\tm[p_0',p_0'']\tm[p_0'',p_0''']\ar[d]^{\circ\tm id}
\\
&[p_0,p_0']\tm[p_0',p_0''']\ar[d]^{\circ}
\\
&[p_0,p_0''']\ar[dl]_{P_1}\ar[dr]^{Q_1}
\\
[a_0,a_0''']\ar[d]_{F_1}^{}="2"
&&[c_0,c_0''']\ar[d]^{G_1}_{}="1"
\\
[Fa_0,Fa_0''']\ar[dr]_{-\circ b_1'''}
&&[Gc_0,Gc_0''']\ar[dl]^{b_1\circ-}
\\
&[Fc_0,Gc_0''']
\ar@{}"1";"2"|(.4){}="a1"|(.6){}="a2"
\ar@{=>}"a1";"a2"_{\varepsilon_1^{p_0,p_0'''}}}
$$
that, applying (\ref{def:epsilon123}) and (\ref{eqn:epsilon123_eq_comp_epsilon}) to $\varepsilon_1^{p_0,p_0'''}$, equals to
\begin{changemargin}{-10ex}{-10ex}{\small
$$
\xymatrix@C=-6ex{
&&[p_0,p_0']\tm[p_0',p_0'']\tm[p_0'',p_0''']\ar[d]^{\circ\tm id}
\\
[p_0,p_0'']\tm[a_0'',a_0''']\ar[dd]_(.7){id\tm F_1}^(.8){}="s2"
&&\ar[ll]_{id\tm P_1}\ar[dl]^{id\tm Q_1}
[p_0,p_0'']\tm[p_0'',p_0'''] \ar[dr]_{P_1\tm id}\ar[rr]^{Q_1\tm id}
&&[c_0,c_0'']\tm[p_0',p_0''']\ar[dd]^(.7){G_1\tm id}_(.8){}="d1"
\\
&[p_0,p_0'']\tm[c_0'',c_0''']\ar[dd]^{id\tm G_1}_(.4){}="s1"
&&[a_0,a_0'']\tm[p_0'',p_0''']\ar[dd]_{F_1\tm id}^(.4){}="d2"
\\
[p_0,p_0'']\tm[Fa_0'',Fa_0''']\ar[dd]_{id\tm (-\circ b_1''')}
&&&&[Gc_0,Gc_0'']\tm[p_0'',p_0''']\ar[dd]^{(b_1\circ-)\tm id}
\\
&[p_0,p_0'']\tm[Gc_0'',Gc_0''']\ar[dl]^{id\tm (b_1''\circ-)}
&&[Fa_0,Fa_0'']\tm[p_0'',p_0''']\ar[dr]_{(-\circ b_1'')\tm id}
&&
\\
[p_0,p_0'']\tm[Fa_0'',Gc_0''']\ar[d]_{P_1\tm id}
&&&&[Fa_0,Gc_0'']\tm[p_0'',p_0''']\ar[d]^{id\tm Q_1}
\\
[a_0,a_0'']\tm[Fa_0'',Gc_0''']\ar[dr]_{F_1\tm id}
&&&&[Fa_0,Gc_0'']\tm[c_0'',c_0''']\ar[dl]^{id\tm G_1}
\\
&[Fa_0,Fa_0'']\tm[Fa_0'',Gc_0''']\ar[dr]_{\circ}
&&[Fa_0,Gc_0'']\tm[Gc_0'',Gc_0''']\ar[dl]^{\circ}
&&
\\
&&[Fa_0,Gc_0''']&&
\ar@{}"s1";"s2"|(.3){}="ss1"|(.7){}="ss2"
\ar@{}"d1";"d2"|(.3){}="dd1"|(.7){}="dd2"
\ar@{=>}"ss1";"ss2"_{id\tm \varepsilon_1^{p_0'',p_0'''}}
\ar@{=>}"dd1";"dd2"_{\varepsilon_1^{p_0,p_0''}\tm id}
}
$$
}\end{changemargin}
Now we can apply (\ref{def:epsilon123}) again, and (\ref{eqn:epsilon123_eq_comp_epsilon}) to the right-hand side of the
diagram, to express $\varepsilon_1^{p_0,p_0''}$ in terms of $\varepsilon_1^{p_0,p_0'}$ and $\varepsilon_1^{p_0',p_0''}$.
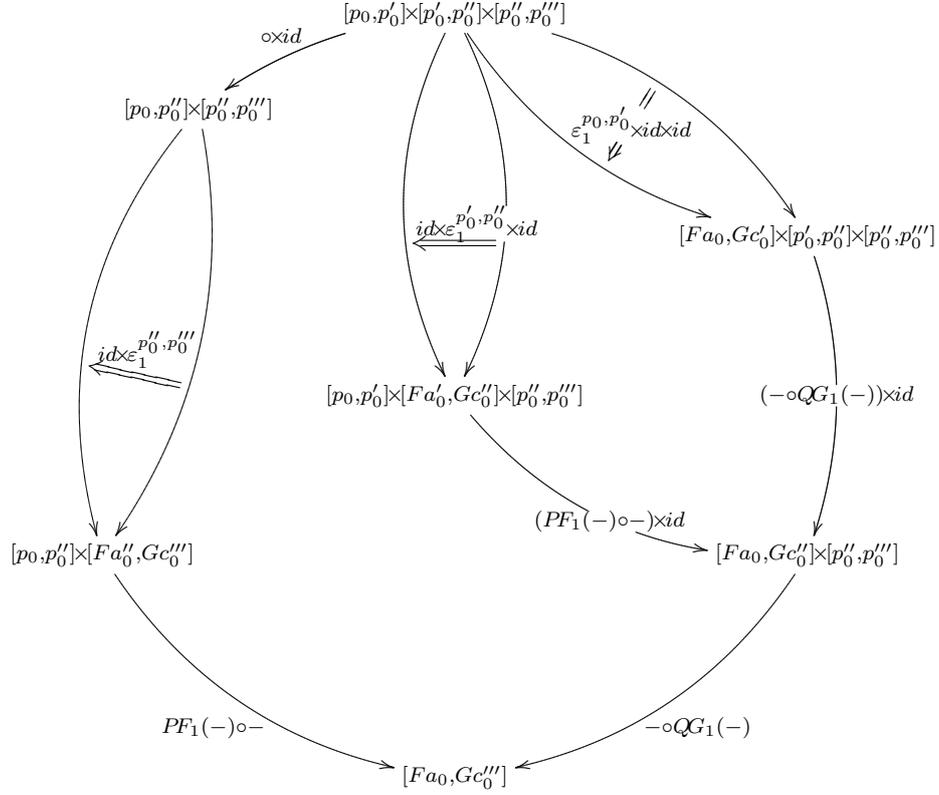
\begin{changemargin}{-1ex}{0ex}{
\def \objectstyle{\scriptstyle}
\def \labelstyle{\scriptstyle}
\def \xa{[p_0,p_0']\tm[p_0',p_0'']\tm[p_0'',p_0''']}
\def \xb{[p_0,p_0'']\tm[p_0'',p_0''']}
\def \xc{[Fa_0,Gc_0']\tm[p_0',p_0'']\tm[p_0'',p_0''']}
\def \xd{[p_0,p_0']\tm[Fa_0',Gc_0'']\tm[p_0'',p_0''']}
\def \xe{[p_0,p_0'']\tm[Fa_0'',Gc_0''']}
\def \xf{[Fa_0,Gc_0'']\tm[p_0'',p_0''']}
\def \xg{[Fa_0,Gc_0''']}
\def \xab{\circ\tm id}
\def \xact{}
\def \xacs{}
\def \xadt{}
\def \xads{}
\def \xbet{}
\def \xbes{}
\def \xcf{( -\circ Q\!G_1(-)  )\tm id}
\def \xdf{(P\!F_1(-)\circ-)\tm id}
\def \xfg{-\circ Q\!G_1(-)}
\def \xeg{P\!F_1(-)\circ-}
\def \beii{id\tm \varepsilon_1^{p_0'',p_0'''}}
\def \acii{\varepsilon_1^{p_0,p_0'}\tm id \tm id}
\def \adii{id\tm \varepsilon_1^{p_0',p_0''}\tm id \qquad}
$$
\xy 0;/r.20pc/:
(0,60)*+{\xa}="a";
(-40,45)*+{\xb}="b";
(55,25)*+{\xc}="c";
(0,0)*+{\xd}="d";
(-55,-25)*+{\xe}="e";
(55,-25)*+{\xf}="f";
(0,-60)*+{\xg}="g";
{\ar@/_7.86pt/_(.6){\xab}"a";"b"^{}="ab"};
{\ar@/_19.25pt/_(.8){\xact}"a";"c"^{}="act"};
{\ar@/^19.25pt/^(.8){\xacs}"a";"c"_{}="acs"};
{\ar@/_19.25pt/_(.8){\xadt}"a";"d"^(.6){}="adt"};
{\ar@/^19.25pt/|(.55){\adii}"a";"d"_(.6){}="ads"};
{\ar@/_23.69pt/_(.8){\xbet}"b";"e"^(.6){}="bet"};
{\ar@/^19.25pt/^(.8){\xbes}"b";"e"_(.6){}="bes"};
{\ar@/_19.25pt/|{\xdf}"d";"f"};
{\ar@/^10.91pt/|{\xcf}"c";"f"};
{\ar@/_19.25pt/_{\xeg}"e";"g"};
{\ar@/^19.25pt/^{\xfg}"f";"g"};
{\ar@{=>}|{\acii}"acs";"act"};
{\ar@{=>}_{}"ads";"adt"};
{\ar@{=>}_{\beii}"bes";"bet"};
\endxy
$$}
\end{changemargin}
Moreover, in order to shift the 2-morphism $\varepsilon_1^{p_0'',p_0'''}$ up, we apply \emph{product interchange rules}
to the left-hand side. What we get is the diagram:
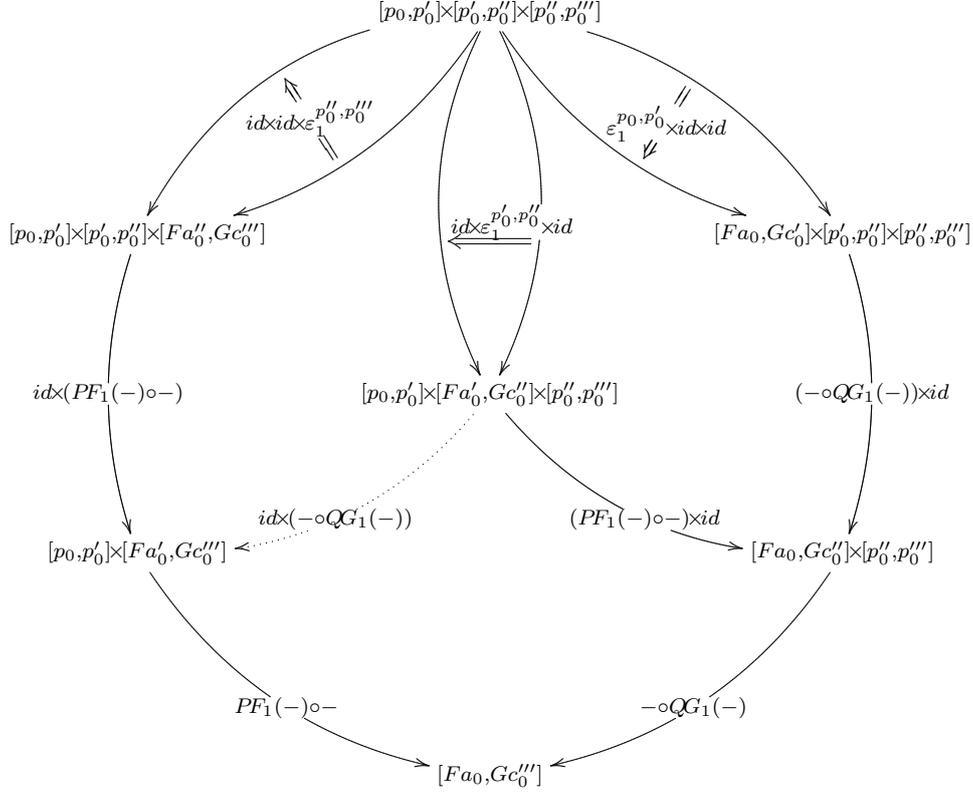
\begin{changemargin}{-1ex}{0ex}{
\def \objectstyle{\scriptstyle}
\def \labelstyle{\scriptstyle}
\def \xa{[p_0,p_0']\tm[p_0',p_0'']\tm[p_0'',p_0''']}
\def \xb{[p_0,p_0']\tm[p_0',p_0'']\tm[Fa_0'',Gc_0''']}
\def \xc{[Fa_0,Gc_0']\tm[p_0',p_0'']\tm[p_0'',p_0''']}
\def \xd{[p_0,p_0']\tm[Fa_0',Gc_0'']\tm[p_0'',p_0''']}
\def \xe{[p_0,p_0']\tm[Fa_0',Gc_0''']}
\def \xf{[Fa_0,Gc_0'']\tm[p_0'',p_0''']}
\def \xg{[Fa_0,Gc_0''']}
\def \xabt{}
\def \xabs{}
\def \xact{}
\def \xacs{}
\def \xadt{\ }
\def \xads{\ }
\def \xbe{id\tm (P\!F_1(-)\circ-)}
\def \xcf{( -\circ Q\!G_1(-)  )\tm id}
\def \xdf{(P\!F_1(-)\circ-)\tm id}
\def \xde{id\tm (-\circ Q\!G_1(-))}
\def \xfg{-\circ Q\!G_1(-)}
\def \xeg{P\!F_1(-)\circ- }
\def \abii{id\tm id\tm \varepsilon_1^{p_0'',p_0'''}}
\def \acii{\varepsilon_1^{p_0,p_0'}\tm id \tm id}
\def \adii{id\tm \varepsilon_1^{p_0',p_0''}\tm id \qquad}
$$
\xy 0;/r.20pc/:
(0,60)*+{\xa}="a";
(-55,25)*+{\xb}="b";
(55,25)*+{\xc}="c";
(0,0)*+{\xd}="d";
(-55,-25)*+{\xe}="e";
(55,-25)*+{\xf}="f";
(0,-60)*+{\xg}="g";
{\ar@/_10.91pt/|{\xbe}"b";"e"^{}="be"};
{\ar@/_19.25pt/_(.8){\xact}"a";"c"^{}="act"};
{\ar@/^19.25pt/^(.8){\xacs}"a";"c"_{}="acs"};
{\ar@/_19.25pt/|(.8){\xadt}"a";"d"^(.6){}="adt"};
{\ar@/^19.25pt/|(.55){\adii}"a";"d"_(.6){}="ads"};
{\ar@/_23.69pt/_(.8){\xabt}"a";"b"^{}="abt"};
{\ar@/^19.25pt/^(.8){\xabs}"a";"b"_{}="abs"};
{\ar@/_19.25pt/|{\xdf}"d";"f"};
{\ar@/^19.25pt/@{.>}|{\xde}"d";"e"};
{\ar@/^10.91pt/|{\xcf}"c";"f"};
{\ar@/_19.25pt/|{\xeg}"e";"g"};
{\ar@/^19.25pt/|{\xfg}"f";"g"};
{\ar@{=>}|{\acii}"acs";"act"};
{\ar@{=>}_{}"ads";"adt"};
{\ar@{=>}|{\abii}"abs";"abt"};
\endxy
$$}
\end{changemargin}
The dotted arrow fits the diagram properly, making the two regions commute. Hence the whole diagram
is perfectly symmetric, and calculations may
be carried on doing the steps in reverse order, and gain the result.
\end{proof}

\subsection{Units}
As for composition, we get unit morphisms by universal property in $(n-1)$\textbf{Cat}.

Suppose an element $p_0$ of $\PP_0$ is fixed.
$$
\xymatrix@C=10ex{
\I\ar[r]^{u(c_0)}\ar[d]_{u(a_0)}&\C_1(c_0,c_0)\ar[d]^{b_1\circ G(-)}\\
\A_1(a_0,a_0)\ar[r]_{F(-)\circ b_1}& \B_1(Fa_0,Gc_0)
}
$$
The square above commutes, hence it is an identity 2-morphism, over the same base defining $\varepsilon_1^{p_0,p_0}$,
and this implies the existence of a unique
$$
\xymatrix{{}^{\PP}\!u(p_0):\I\ar[r]&\PP_1(p_0,p_0)}.
$$
\begin{Lemma}
Units defined above are neutral w.r.t. 0-composition, {\em i.e.}, for every pair $p_0, p_0'$, in $\PP_0$, $(L)$ and $(R)$
commute.
$$
\xymatrix@C=14ex{
\ar@{}[dr]|{(R)}\PP_1(p_0,p_0'\times \I)\ar[r]^(.4){id \times u(p_0')}
&\PP_1(p_0,p_0')\times\PP_1(p_0',p_0')\ar[d]^{\circ}
\\
\ar@{}[dr]|{(L)}
\PP_1(p_0,p_0')\ar[r]|{id}\ar[d]_{\lambda}^{\cong}\ar[u]^{\rho}_{\cong}
& \PP_1(p_0,p_0')
\\
\I\times\PP_1(p_0,p_0')\ar[r]_(.4){u(p_0)\times id}
&\PP_1(p_0,p_0)\times \PP_1(p_0,p_0')\ar[u]_{\circ}
}
$$
\end{Lemma}
\begin{proof}
We show only the commutativity of $(L)$, the other being similar. Hence let us consider the composition
\begin{equation}
\xymatrix{
&[p_0,p_0']\ar[d]^{\lambda}_{\cong}\ar@/_3ex/@{.>}[ddddl]_{P_1}\ar@/^3ex/@{.>}[ddddr]^{Q_1}
\\
&\I\times[p_0,p_0']\ar[d]|{u(p_0)\times id}
\\
&[p_0,p_0]\times[p_0,p_0']\ar[d]^{\circ}
\\
&[p_0,p_0']\ar[dl]_{P_1}\ar[dr]^{Q_1}
\\
[a_0,a_0']\ar[dr]_{F(-)\circ b_1'}
&&[c_0,c_0']\ar[dl]^{b_1\circ G(-)}\ar@{}[ll]|(.4){}="1"|(.6){}="2"
\\
&[Fa_0,Gc_0']
\ar@{=>}"1";"2"_{\varepsilon_1^{p_0,p_0'}}}
\end{equation}
It is easy to see that both sides commute with dotted arrows. In fact
\begin{eqnarray*}
  \lambda({}^{\PP}u(p_0)\times id)({}^{\PP}\circ)Q_1 &\eq{i}& \lambda({}^{\PP}u(p_0)\times id)(Q_1\times Q_1) ({}^{\C}\circ)\\
   &\eq{ii}& \lambda ({}^{\C}u(c_0)\times Q_1) ({}^{\C}\circ)\\
   &\eq{iii}& Q_1
\end{eqnarray*}
where $(i)$ holds by (\ref{eqn:PQ_funct}), $(ii)$ by definition of ${}^{\PP}\!u$, $(iii)$ by the
{\em unit axiom} in $\C$. Similarly
for the left-hand side.

If one shows that composition above is equal to $\varepsilon_1^{p_0,p_0'}$, the universal property of
pullbacks implies $(L)$.
Now we can  reformulate it with the help of (\ref{def:epsilon123}) and (\ref{eqn:epsilon123_eq_comp_epsilon}):
$$
\xymatrix@C=0ex@R=10ex{
&[p_0,p_0']\ar[d]^{\lambda}_{\cong}
\\
&\I\times[p_0,p_0']\ar[d]|{u(p_0)\times id}
\\
&[p_0,p_0]\tm[p_0,p_0']\ar@/^3ex/[dl]_{}="1"\ar@/_3ex/[dl]^{}="2"\ar@{=>}"1";"2"_{id\tm\varepsilon_1^{p_0,p_0'}}
\ar[dr]^{[b_1]\tm id}
\\
[p_0,p_0]\tm[Fa_0,Gc_0']\ar[dr]_{PF_1(-)\circ-}
&&[Fa_0,Gc_0]\tm[p_0,p_0']\ar[dl]^{-\circ QG_1(-)}
\\
&[Fa_0,Gc_0']
}
$$
where we have somehow abusively replaced the identity 2-morphism on the right hand side, with its source (= target) 1-morphism.
Hence all the right-hand side, being an identity, may be cancelled. Finally, by {\em product interchange}
$$
\raisebox{18ex}{
\xymatrix{
[p_0,p_0']\ar[d]^{\lambda}_{\cong}
\\
\I\times[p_0,p_0']\ar@/^4ex/[dd]_{}="1"\ar@/_4ex/[dd]^{}="2"\ar@{=>}"1";"2"_{id\tm\varepsilon_1^{p_0,p_0'}}
\\
\\
\I\times[Fa_0,Gc_0']\ar[d]^{u(Fa_0)\circ-}
\\
[Fa_0,Gc_0']
}}
\eq{1}
\raisebox{18ex}{
\xymatrix{
[p_0,p_0']\ar@/^4ex/[dd]_{}="1"\ar@/_4ex/[dd]^{}="2"\ar@{=>}"1";"2"_{\varepsilon_1^{p_0,p_0'}}
\\
\\
[Fa_0,Gc_0']\ar[d]^{\lambda}_{\cong}
\\
\I\times[Fa_0,Gc_0']\ar[d]^{u(Fa_0)\circ-}
\\
[Fa_0,Gc_0']
}}
\eq{2}\quad
\varepsilon_1^{p_0,p_0'}
$$
where $(1)$ holds by naturality of $\lambda_{(-)}:(-)\Rightarrow\I\times(-)$, and $(2)$ by neutral identities in $\B$.
\end{proof}
\subsection{Projections and $\varepsilon$}
So far we proved that the pair $\PP=(\PP_0,\PP_1^{-,-})$ is indeed a $n$-category. In order to show it is a part of a pullback
four-tuple, we should prove that
$$
P=(P_0,P_1^{-,-}),\qquad Q=(Q_0,Q_1^{-,-}),\qquad  \varepsilon= (\varepsilon_0,\varepsilon_1^{-,-})
$$
produced in stating the definitions above, constitute respectively two n-categorie morphisms and one 2-morphism.
But this has been already proved throughout the last sections. In fact, universal definition of 0-composition above reveals that this is
just the one that makes $P$, $Q$ and $\varepsilon$ functorial. Similarly, universal definition of 0-units
reveals that these are just the ones that make $P$, $Q$ and $\varepsilon$ functorial.

\subsection{Universal property}
The final step in proving that $n$\textbf{Cat} admits $h$-pullbacks, is to show that the four-tuple
$(\PP,P,Q,\varepsilon)$ satisfies universal property of $h$-Pullbacks (\emph{UP} \ref{UP:h-pullbacks}).

To this aim, let us suppose a n-category $\X$ been given, together with morphisms and 2-morphisms
$$
M:\X\to\A,\qquad N:\X\to \C,\qquad \omega:MF\Rightarrow NG
$$
On objects, as $\PP_0$ is a limit in \textbf{Set}, it suffices to consider the cone over the same diagram defining the
latter, whose commutativity is a consequence of the very definition of $\omega$:
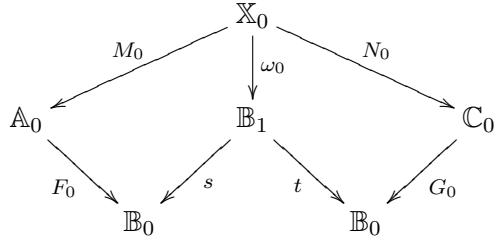
$$
\xymatrix{
&&\X_0\ar[dll]_{M_0}\ar[d]^{\omega_0}\ar[drr]^{N_0}
\\
\A_0\ar[dr]_{F_0}
&&\B_1\ar[dl]^{s}\ar[dr]_{t}
&&\C_0\ar[dl]^{G_0}
\\
&\B_0&&\B_0
}
$$
This yields a unique map $L_0:\X_0\to\PP_0$, such that:
\begin{equation}\label{pf:obj_pb}
L_0 P_0= M_0,\qquad L_0 Q_0 = N_0,\qquad L_0\varepsilon_0=\omega_0
\end{equation}

On homs, let us fix objects $x_0$ and $x_0'$ in $\PP$. By the universal property in dimension $n-1$, the four-tuple
$$
(\X_1(x_0,x_0'), M_1^{x_0,x_0'},N_1^{x_0,x_0'},\varepsilon_1^{x_0,x_0'})
$$
gives a unique morphism $L_1^{x_0,x_0'}:\X_1(x_0,x_0')\to\PP_1(Lx_0,Lx_0')$
such that
\begin{eqnarray}\label{pf:hom_pb}
\nonumber   L_1^{x_0,x_0'}\bullet^0 P_1^{Lx_0,Lx_0'}&=&M_1^{x_0,x_0'}  \\
            L_1^{x_0,x_0'}\bullet^0 Q_1^{Lx_0,Lx_0'}&=&N_1^{x_0,x_0'}  \\
\nonumber   L_1^{x_0,x_0'}\bullet^0 \varepsilon_1^{Lx_0,Lx_0'}&=&\omega_1^{x_0,x_0'}
\end{eqnarray}

\textbf{Claim}: the pair $(L_0,L_1^{-,-})$ constitutes a $n$-functor $L:\X\to\PP$.
\begin{proof} The proof is divided in two parts.\\

1. {\em Functoriality w.r.t. compositions.}

Let us fix a triple $x_0,x_0',x_0''$ of objects of $\X$. What we want to prove is the diagram below commutes:
$$
\xymatrix{
\X_1(x_0,x_0')\times\X_1(x_0',x_0'')\ar[r]^(.6){{}^\X\!\circ}\ar[d]_{L_1\times L_1}
&\X_1(x_0,x_0'')\ar[d]^{L_1}
\\
\PP_1(Lx_0,Lx_0')\times\PP_1(Lx_0',Lx_0'')\ar[r]_(.65){{}^\PP\!\circ}
&\PP_1(Lx_0,Lx_0'')
}
$$

Then, let us consider the $h$-pullback defining $\PP_1(Lx_0,Lx_0'')$. If we can show that the horizonal composition of both
composites above with $\varepsilon_1^{p_0,p_0''}$ coincide, uniqueness forces $(L_1\times L_1)\bullet^0\,{}^{\PP}\circ={}^{\X}\circ\bullet^0 L_1$.
Hence, let us follow the chain of equalities below:
\begin{changemargin}{-6ex}{6ex}
$$
\raisebox{14ex}{\xymatrix{
[x_0,x_0']\tm[x_0',x_0'']\ar[d]^{{}^\X\!\circ}
\\
[x_0,x_0'']\ar[d]^{L_1}
\\
[Lx_0,Lx_0'']
\ar@/^6ex/[dd]_{}="1"\ar@/_6ex/[dd]^{}="2"\ar@{=>}"1";"2"_{\varepsilon_1^{Lx_0,Lx_0''}}
\\
\\
[MFx_0,NGx_0'']
}}
\eq{1}
\raisebox{10ex}{\xymatrix{
[x_0,x_0']\tm[x_0',x_0'']\ar[d]^{{}^\X\!\circ}
\\
[x_0,x_0'']
\ar@/^4ex/[dd]_{}="1"\ar@/_4ex/[dd]^{}="2"\ar@{=>}"1";"2"_{\omega_1^{x_0,x_0''}}
\\
\\
[MFx_0,NGx_0'']
}}
\eq{2}
\raisebox{10ex}{\xymatrix@C=-6ex{
&[x_0,x_0']\tm[x_0',x_0'']
\ar@/^6ex/[ddl]_{}="a1"\ar@/_6ex/[ddl]^{}="a2"\ar@{=>}"a1";"a2"|{id\tm\omega_1^{x_0',x_0''}}
\ar@/^6ex/[ddr]_{}="b1"\ar@/_6ex/[ddr]^{}="b2"\ar@{=>}"b1";"b2"|{\omega_1^{x_0,x_0'}\tm id}
\\
\\
[x_0,x_0']\tm[MFx_0',NGx_0'']\ar[dr]_{MF_1(-)\circ-}
&&[MFx_0,NGx_0']\tm[x_0',x_0'']\ar[dl]^{-\circ NG_1(-)}
\\
&[MFx_0,NGx_0'']
}}
$$
\end{changemargin}
$(1)$ holds by the third equation of (\ref{pf:hom_pb}), $(2)$ by functoriality w.r.t. composition  of $\omega$,
$$
\eq{3}
\raisebox{14ex}{\xymatrix@C=-8ex{
&[x_0,x_0']\tm[x_0',x_0'']\ar[dl]_{id\tm L_1}\ar[dr]^{L_1\tm id}
\\
[x_0,x_0']\tm[Lx_0',Lx_0'']\ar@/^6ex/[dd]_{}="a1"\ar@/_6ex/[dd]^{}="a2"\ar@{=>}"a1";"a2"_{id\tm\varepsilon_1^{Lx_0',Lx_0''}}
&&[Lx_0,Lx_0']\tm[x_0',x_0'']
\ar@/^6ex/[dd]_{}="b1"\ar@/_6ex/[dd]^{}="b2"\ar@{=>}"b1";"b2"_{\varepsilon_1^{Lx_0,Lx_0'}\tm id}
\\
\\
[x_0,x_0']\tm[LPFx_0',LQGx_0'']\ar[dr]_{LPF_1(-)\circ-}
&&[LPFx_0,LQGx_0']\tm[x_0',x_0'']\ar[dl]^{-\circ LQG_1(-)}
\\
&[MFx_0,NGx_0'']
}}
$$
here, $(3)$ is a full consequence of  equations (\ref{pf:hom_pb}) and product interchange,
\begin{changemargin}{-6ex}{6ex}
$$
\eq{4}
\raisebox{10ex}{\xymatrix@C=-6ex{
&[x_0,x_0']\tm[x_0',x_0'']\ar[d]^{L_1\tm L_1}
\\
&[Lx_0,Lx_0']\tm[Lx_0',Lx_0'']
\ar@/^6ex/[ddl]_{}="a1"\ar@/_6ex/[ddl]^{}="a2"\ar@{=>}"a1";"a2"|{id\tm\varepsilon_1^{Lx_0',Lx_0''}}
\ar@/^6ex/[ddr]_{}="b1"\ar@/_6ex/[ddr]^{}="b2"\ar@{=>}"b1";"b2"|{\varepsilon_1^{Lx_0,Lx_0'}\tm id}
\\
\\
[x_0,x_0']\tm[PFx_0',QGx_0'']\ar[dr]_{PF_1(-)\circ-}
&&[PFx_0,QGx_0']\tm[x_0',x_0'']\ar[dl]^{-\circ QG_1(-)}
\\
&[PFx_0,QGx_0'']
}}
\eq{5}
\raisebox{14ex}{\xymatrix{
[x_0,x_0']\tm[x_0',x_0'']\ar[d]^{L_1 \tm L_1}
\\
[Lx_0,Lx_0']\tm[Lx_0',Lx_0'']\ar[d]^{{}^\PP\!\circ}
\\
[Lx_0,Lx_0'']
\ar@/^4ex/[dd]_{}="1"\ar@/_4ex/[dd]^{}="2"\ar@{=>}"1";"2"_{\varepsilon_1^{Lx_0,Lx_0''}}
\\
\\
[MFx_0,NGx_0'']
}}
$$
\end{changemargin}
$(4)$ is obtained sliding the $L_1$'s up, and functoriality w.r.t. composition  of $\varepsilon$ gives $(5)$.\\

2. {\em Functoriality w.r.t. units.}

Let us fix an object $x_0$ in $\X$. What we want to prove is the diagram below commutes:
$$
\xymatrix@C=16ex{\I\ar[r]^{u(x_0)}\ar[dr]_{u(Lx_0)}
&\X_1(x_0,x_0)\ar[d]^{L_1}
\\
&\PP_1(Lx_0,Lx_0)
}
$$
We proceed in a similar way. Let us consider the pullback defining $\PP_1(Lx_0,Lx_0'')$.
If we show that the horizontal composition of both
composites above with $\varepsilon_1^{Lx_0,Lx_0}$ coincide, uniqueness forces $u(x_0)\bullet^0\,L_1=u(Lx_0)$.
Hence, let us follow the chain of equalities below:
$$
\raisebox{18ex}{\xymatrix{
\I\ar[d]^{u(x_0)}
\\
[x_0,x_0]\ar[d]^{L_1}
\\
[Lx_0,Lx_0]
\ar@/^5ex/[dd]_{}="1"\ar@/_5ex/[dd]^{}="2"\ar@{=>}"1";"2"_{\varepsilon_1^{Lx_0,Lx_0}}
\\
\\
[MFx_0,NGx_0]
}}
\eq{1}
\raisebox{14ex}{\xymatrix{
\I\ar[d]^{u(x_0)}
\\
[x_0,x_0]
\ar@/^5ex/[dd]_{}="1"\ar@/_5ex/[dd]^{}="2"\ar@{=>}"1";"2"_{\omega_1^{Lx_0,Lx_0}}
\\
\\
[MFx_0,NGx_0]
}}
\eq{2}
\raisebox{14ex}{\xymatrix{
\I
\ar@/^4ex/[ddd]_{}="1"\ar@/_4ex/[ddd]^{}="2"\ar@{=>}"1";"2"_{Id_{[\omega_{x_0}]}}
\\
\\ {\phantom{[x_0,x_0']}}
\\
[MFx_0,NGx_0]
}}
\eq{3}
\raisebox{14ex}{\xymatrix{
\I\ar[d]^{u(Lx_0)}
\\
[Lx_0,Lx_0]
\ar@/^4ex/[dd]_{}="1"\ar@/_4ex/[dd]^{}="2"\ar@{=>}"1";"2"_{\varepsilon_1^{Lx_0,Lx_0}}
\\
\\
[MFx_0,NGx_0]
}}
$$
$(1)$ holds by the third of the (\ref{pf:hom_pb}), $(2)$ is just $\omega$, the  functoriality of units, and since by (\ref{pf:obj_pb}),
$\omega_{x_0}=\varepsilon_{Lx_0}$,  $(3)$ is obtained by $\varepsilon$ unit functoriality.
\end{proof}

Once we have verified that $L$ is a n-functor, equations (\ref{pf:obj_pb}) and (\ref{pf:hom_pb}) taken together are exactly conditions
1., 2. and 3. of {\em Universal Property \ref{UP:h-pullbacks}}. \\

What is still missing is uniqueness, but this is implied by the proofs. In fact let us suppose there is another
$\hat{L}\X\to\PP$ satisfying universal property. Conditions 1., 2. and 3. of (\ref{UP:h-pullbacks})
imply:
$$
\hat{L}_0 P_0= M_0,\qquad \hat{L}_0 Q_0 = N_0,\qquad \hat{L}_0\varepsilon_0=\omega_0
$$
and
\begin{eqnarray*}
 \hat{L}_1^{x_0,x_0'}\bullet^0 P_1^{\hat{L}x_0,\hat{L}x_0'}&=&M_1^{x_0,x_0'}  \\
 \hat{L}_1^{x_0,x_0'}\bullet^0 Q_1^{\hat{L}x_0,\hat{L}x_0'}&=&N_1^{x_0,x_0'}  \\
 \hat{L}_1^{x_0,x_0'}\bullet^0 \varepsilon_1^{\hat{L}x_0,\hat{L}x_0'}&=&\omega_1^{x_0,x_0'}
\end{eqnarray*}
Since $L_0$ and the $L_1$'s where determined univocally by Universal Properties of
limits in \textbf{Set} and of $h$-pullbacks in
$(n-1)$\textbf{Cat}, uniqueness of those forces
$$
\hat{L}_0=L_0\qquad\textrm{and}\qquad \hat{L}_1^{-,-}=L_1^{-,-}
$$

Hence we proved the
\begin{Theorem}
The sesqui-category $n$\textbf{\emph{Cat}} admits $h$-pullbacks.
\end{Theorem}

\subsection{Pullbacks and $h$-pullbacks}\label{def:strict_pb}
It is possible to recover the usual notion of pullback by means of a similar inductive construction: a pullback
is a  universal triple $<\Q,P,Q>$ such that $PF=QG$:
$$
\xymatrix{
\Q\ar[r]^Q\ar[d]_P
&\C\ar[d]^G
\\
\A\ar[r]_F\ar@{}[ur]|(.3){}="1"|(.7){}="2"
&\B
\ar@{}"1";"2"|{(pb)}}
$$
where $\Q_0$ is the pullback in $\mathbf{Set}$
$$
\xymatrix{
\Q_0\ar[r]^{Q_0}\ar[d]_{P_0}
&\C_0\ar[d]^{G_0}
\\
\A_0\ar[r]_{F_0}\ar@{}[ur]|(.3){}="1"|(.7){}="2"
&\B_0
\ar@{}"1";"2"|{(pb)}}
$$
and for every pair of objects $(a_0,c_0)$ and $(a_0',c_0')$ of $\Q_0$, the following pullback in $(n-1)\mathbf{Cat}$
$$
\xymatrix@C=12ex{
\Q_1\big((a_0,c_0),(a_0',c_0')\big)\ar[r]^{Q_1^{(a_0,c_0),(a_0',c_0')}}\ar[d]_{P_1^{(a_0,c_0),(a_0',c_0')}}
&\C_1(c_0,c_0')\ar[d]^{G_1^{c_0,c_0'}}
\\
\A_1(a_0,a_0')\ar[r]_{F_1^{a_0,a_0'}}\ar@{}[ur]|(.3){}="1"|(.7){}="2"
&\B_1(Fa_0=Gc_0,Fa_0'=Gc_0')
\ar@{}"1";"2"|{(pb)}}
$$
In fact this gives an $h$-pullback in the trivial (= 2-discrete) sesqui-category over the category
$\lfloor n\mathbf{Cat}\rfloor$: the triple $<\Q,P,Q>$ such that $PF=QG$ may be seen as a four-tuple
$<\Q,P,Q,id: PF\Rightarrow QG>$, and so on\dots

\chapter{$n$-Groupoids and exact sequences}\label{cha:equiv_and_ngpd}

The sesqui-category $n$\textbf{Cat} of strict and small n-categories,
defined so far, has a naturally arising notion of equivalence that
may be defined recursively. This  gives a notion of n-groupoid,
equivalent to that of Kapranov and Voevodsky in \cite{MR1130401}, as a
\emph{weakly invertible strict n-category} (see Appendix \ref{cha:appendix} for a comparison).

\section{$n$-Equivalences}
\begin{Definition}
Let n-category morphism $F:\C\to\D$   be given.

$F$ is called \emph{equivalence of n-categories} if it satisfies the following properties:\\

\framebox[1.1\width]{$n=0$}\\

$F$ is an isomorphism in \textbf{Set}.\\

\framebox[1.1\width]{$n>0$}
\begin{enumerate}
  \item\label{def:ess_surj} $F$ is essentially surjective on objects, {\em i.e.} for every object $d_0$ of $\D$,
  there exists an object $c_0$ of $\C$ and a 1-cell

  $d_1:Fc_0\to d_0$ such that for every $d_0'$  in $\C$, the morphisms
  \begin{eqnarray*}
    d_1\circ - &:& \D_1(d_0,d_0')\to \D_1(Fc_0,d_0') \\
    -\circ d_1 &:& \D_1(d_0',Fc_0)\to \D_1(d_0',d_0)
  \end{eqnarray*}
  are equivalences of (n-1)categories.
  \item for every pair $c_0,c_0'$ in $\C$,
$$
F_1^{c_0,c_0'}:\C_1(c_0,c_0')\to\D_1(Fc_0,Fc_0')
$$
is an equivalence of (n-1)categories.
\end{enumerate}
\end{Definition}
From definition above, one gets the following
\begin{Definition}\label{def:weak_inv}
A 1-cell $c_1:c_0\to c_0'$ of a n-category $\C$ is said to be \emph{weakly invertible}, or simply an \emph{equivalence}, if, for every object
$\bar{c}_0$ of $\C$, the morphisms
\begin{eqnarray*}
    c_1\circ - &:& \C_1(c_0',\bar{c}_0)\to \C_1(c_0,\bar{c}_0) \\
    -\circ c_1 &:& \C_1(\bar{c}_0,c_0)\to \C_1(\bar{c}_0,c_0')
\end{eqnarray*}
are (natural) equivalences of $(n-1)$categories.
\end{Definition}

\subsection{Inverses}

When a 1-cell is weakly invertible, then it has indeed left and right (quasi) inverses. In fact for
$c_1:c_0\to c_0'$,
$$
c_1\circ - : \C_1(c_0',c_0)\to \C_1(c_0,c_0)
$$
to be an equivalence implies that for the 1-cell $1_{c_0}:c_0\to c_0$ there exists a pair
$$
(c_1^*,\xymatrix{c_2:c_1\circ c_1^*\ar@2[r]^(.6){\sim} &1_{c_1}}),
$$
similarly for
$$
 -\circ c_1 : \C_1(c_0',c_0)\to \C_1(c_0',c_0')
$$
implies there exists a pair
$$
(c_1^{\dag},\xymatrix{c_2':c_1^{\dag}\circ c_1\ar@2[r]^(.6){\sim}& 1_{c_1'}}).
$$

\subsection{Properties}
\begin{Lemma}\label{lemma:compos_equiv}
Let n-functors $\xymatrix{\C\ar[r]^F&\D\ar[r]^G&\E}$ be given. Then
if $F$ and $G$ are equivalences, so is $F\bullet^0G$.
\end{Lemma}
\begin{proof}
The composite of isomorphisms in Set is trivially an isomorphism. Hence we may suppose $n>0$.
\begin{enumerate}
  \item for every pair $c_0,c_0'$, $[F\bullet^0 G]_1^{c_0,c_0'}$ is an equivalence. In fact
   $$
   [F\bullet^0 G]_1^{c_0,c_0'}=F_1^{c_0,c_0'}\bullet^0 G_1^{Fc_0,Fc_0'}
   $$
   two component on the right-hand side are indeed equivalences by hypothesis, and so is their composites by induction.
  \item For any object $e_0$ of $\E$ there exists a pair $(d_0,e_1:Gd_0\to e_0)$. Similarly, for any object $d_0$ in
  $\D$ there is a pair $(c_0,d_1:Fc_0\to d_0)$. Hence, for given $e_0$, those produce a pair
  $$
  (c_0,\xymatrix{G(Fc_0)\ar[r]^{Gd_1}&Gd_0\ar[r]^{e_1}&e_0})
  $$
  That left and right 0-compositions (in $\D$) with $Gd_1\circ e_1$ are equivalences, is a statement involving a composition of
  equivalences of (n-1)categories, hence given by induction. In fact, by definition
  $$
  -\circ (Gd_1\circ e_1)=(-\circ Gd_1)\bullet^0(-\circ e_1)
  $$
  and
  $$
  (Gd_1\circ e_1)\circ-=(e_1\circ-)\bullet^0(Gd_1\circ-)
  $$
    \end{enumerate}
\end{proof}

\begin{Definition}\label{def:2equiv}
A 2-morphism of n-categories $\alpha:F\Rightarrow G:\C\to\D$ is an equivalence 2-morphism, or n-natural equivalence, if\\

\framebox[1.1\width]{$n=1$}\\
$\alpha$ is a natural isomorphism.\\

\framebox[1.1\width]{$n>1$}\\
\begin{enumerate}
 \item for any object $c_0$ in $\C$, the 1-cell $\alpha_{c_0}$ is an equivalence
 \item for any pair of objects $c_0,c_0'$ in $\C$, the (n-1)transformation $\alpha_1^{c_0,c_0'}$ is an equivalence 2-morphism
\end{enumerate}
\end{Definition}

It is not precisely in the aims of this work, nevertheless it is worth mentioning  the following
\begin{Proposition}
n-categories, n-functors and n-natural equivalences form a sesqui-category, denoted $n\mathbf{Cat}_{eq}$,
\end{Proposition}
Notice that in $n\mathbf{Cat}_{eq}$, and {\em a fortiori} in $n\mathbf{Gpd}$  (t.b.d.) equivalences have
more nice properties, like they are $h$-pullback stable, have the 2of3 property and so on.\\

\begin{Definition}\label{def:h-surjective}
A morphism of $n$-categories $F:\C\to\D$ is called $h$-surjective if\\

\framebox[1.1\width]{$n=0$}\\

$F$ is a surjective map.\\[3ex]

\framebox[1.1\width]{$n>0$}
\begin{enumerate}
\item $F$ is essentially surjective on objects, i.e.

for every object $d_0$ of $\D$, there exists a pair\/ $(c_0,d_1:Lc_0\tilde{\rightarrow} d_0)$, with $d_1$ an equivalence,
\item for every pair of objects $c_0,c_0'$ of $\C$, the morphism
$$
F_1^{c_0,c_0'}: \C_1(c_0,c_0')\to\D_1(Fc_0,Fc_0')
$$
is $h$-surjective.
\end{enumerate}
\end{Definition}

\begin{Definition}
A morphism of $n$-categories $F:\C\to\D$ is called faithful if\\

\framebox[1.1\width]{$n=0$}\\

$F$ is a injective map.\\[3ex]

\framebox[1.1\width]{$n>0$}\\

for every pair $c_0,c_0'$, the (n-1)functor $F_1^{c_0,c_0'}$ is faithful.
\end{Definition}

The notion of $h$-surjective is weaker than (implied by) that of equivalence. In fact, more is true:
\begin{Proposition}\label{prop:equiv}
A   morphism of n-categories $F:\C\to\D$ is an equivalence precisely when it is faithful and $h$-surjective.
\end{Proposition}
\begin{proof}
When $n=0$ this is the characterization of bijective maps as injective plus surjective. Hence suppose $n>0$.

Let $F$ be an equivalence. Then $F$ is essentially surjective by definition. More, for every $c_0,c_0'$,
$F_1^{c_0,c_0'}$ is an equivalence in $(n-1)$\textbf{Cat}, therefore $h$-surjective by induction. Finally, the last is
also faithful by induction, and this concludes the first implication.

Conversely, let $F$ be faithful and $h$-surjective. Then it is essentially surjective by definition.
More, for every $c_0,c_0'$, $F_1^{c_0,c_0'}$'s are faithful and $h$-surjective, and inductive hypothesis implies
they are equivalences.
\end{proof}

Notice that faithfulness can be reformulated saying that the n-functor is \emph{surjective on equations}.
In fact, from the globular point of view, this is equivalent to saying that equal $n$-cells in the image
of $F$ come from equal $n$-cells of $\C$. Under this perspective, to be $h$-surjective  amounts to
being (weakly) surjective in any dimension, up to $n-1$, and last proposition says precisely that an
equivalence is (weakly) surjective on $k$-cell, with \mbox{$0\leq k \leq n$.}

\begin{Remark}
  The notions of $h$-surjective morphisms and of equivalences reduce to well known ones, when considered
  in low dimension. In fact, in dimension one those are stated explicitly in the definitions as the first
  step of the inductive process. For $n=1$ an $h$-surjective morphism
  is a functor which is  full and essentially surjective on object. Hence the notion of equivalence is
  the usual one.
\end{Remark}

Finally we state a useful Lemma, whose proof is part of the proof of {\em Lemma \ref{lemma:compos_equiv}}:

\begin{Lemma}\label{lemma:h-surj_comp}
Let $n$-functors $\xymatrix{\C\ar[r]^F&\D\ar[r]^G&\E}$ be given. Then
if $F$ and $G$ are $h$-surjective, so is $F\bullet^0G$.
\end{Lemma}

\section{$n$-Groupoids}
\begin{Definition}\label{def:n_gpd} The definition is inductive on $n$.\\

\framebox[1.1\width]{$n=0$}\\

A $0$-groupoid is a $0$-category, \emph{i.e.} a set.\\

\framebox[1.1\width]{$n>0$}\\[2ex]
A n-groupoid is a $n$-category $\C$ such that:\\

1. every 1-cell of $\C$ is an equivalence;

2. for every pair of objects $c_0,c_0'$ of $\C$ the $(n-1)$category $\C_1(c_0,c_0')$ is a $(n-1)$groupoid.\\
\end{Definition}

We denote by $n$\textbf{Gpd} the sub-sesqui-category of $n$\textbf{Cat} generated by $n$-groupoids.\\

\begin{Proposition}
For every given natural number $n$, the following is a diagram of inclusions:
$$
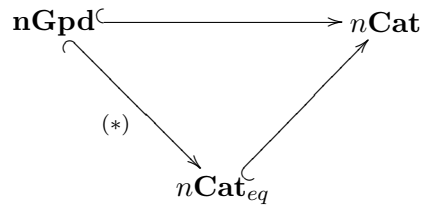
\xymatrix@R=10ex{
\mathbf{nGpd}\ar@{^(->}[rr]\ar@{_(->}[dr]_{(*)}&&n\mathbf{Cat}\\&n\mathbf{Cat}_{eq}\ar@{_(->}[ur]}
$$
\end{Proposition}
\begin{proof} The case $=1$ is well known, hence suppose $n>1$.

The only inclusion to be proved is the one marked $(*)$. To this end, it suffice to show that 2-morphisms of
n-groupoids are equivalences. But, condition 1. of {\em Lemma \ref{def:2equiv}}\/ above is automatically satisfied, as n-groupoid's
1-cells are always equivalences, condition 2. is given by induction.
\end{proof}

\begin{Proposition}
$n$\textbf{Gpd} is closed under h-pullbacks.
\end{Proposition}
\begin{proof}
For $n=0$ the result holds trivially. Hence let us suppose $n>0$. For the $h$-pullback
$$
<\PP,P,Q,\varepsilon>
$$
$$
\xymatrix{
\PP\ar[r]^Q\ar[d]_P
&\C\ar[d]^G
\\
\A\ar[r]_F\ar@{}[ur]|(.3){}="1"|(.7){}="2"
&\B
\ar@{=>}"1";"2"^{\varepsilon}}
$$
a 1-cell
$$
\xymatrix@C=2ex{
(a_0\ar[d]_{a_1}
&,
&Fa_0\ar[d]_{Fa_1}
\ar[rrr]^{b_1}
&&
&Gc_0\ar[d]^{Gc_1}
\ar@{}[dlll]|(.3){}="1"|(.7){}="2"\ar@2"1";"2"^{b_2}
&,
&c_0\ar[d]^{c_1})
\\
(a_0'
&,
&Fa_0'
\ar[rrr]_{b_1'}
&&
&Gc_0'
&,
&c_0')}
$$
induces equivalences because $a_1,c_1$ and $b_2$  do. Turning to homs, the statement  is relative to $(n-1)$ groupoids
and holds by induction.
\end{proof}

\section{The sesqui-functor $\pi_0$}
Purpose of this section is to introduce the family of sesqui-functors
$\{\pi_0^{(n)}\}_{n\in\mathbb{N^*}}$ that extends the
iso-classes functor $\mathbf{Gpd}\to\mathbf{Set}$.
\begin{Defprop}
For any integer $n>0$, there exists a classifying sesqui-functor
$$
\pi_0^{(n)}:\mathbf{nGpd} \rightarrow \mathbf{(n-1)Gpd}
$$
according to the following inductive definition.

Moreover, it commutes with finite products and it preserves equivalences.\\
\end{Defprop}

\framebox[1.1\width]{$n=1$}\\
$$
\pi_0^{(1)}:n\mathbf{Gpd} \rightarrow (n-1)\mathbf{Gpd}
$$
is the functor (= trivial sesqui-functor) $\mathbf{Gpd}\to\mathbf{Set}$ that assigns to a groupoid $\C$ the  set
$|\C|$ of isomorphism classes of objects of $\C$.

It  commutes with finite products: in fact the terminal set
$\{*\}$ is exactly the classified terminal groupoid $\mathbb{I}$, and in a product of groupoids $\C\times\D$, an isomorphism
$(c_1,d_1)$ is a pair of isomorphisms $c_1$ in $\C$ and $d_1$ in $\D$.

Finally, it sends equivalences of category in isomorphism-maps.\\[2ex]

\framebox[1.1\width]{$n>1$}

\subsection{$\pi_0$ on objects}
Let a $n$-groupoid $\C$ be given. Then, $\pi_0^{(n)}\C=([\pi_0^{(n)}\C]_0,[\pi_0^{(n)}\C]_1(-,-))$, where
\begin{itemize}
\item $[\pi_0^{(n)}\C]_0=\C_0$
\item for every pair $c_0,c_0'$ in $[\pi_0^{(n)}\C]_0$,
$$
[\pi_0^{(n)}\C]_1(c_0,c_0')=\pi_0^{(n-1)}(\C_1(c_0,c_0'))
$$
\end{itemize}
For a triple of objects $c_0,c_0',c_0''$, composition is the dotted arrow below:
$$
\xymatrix@C=16ex{
[\pi_0^{(n)}\C]_1(c_0,c_0')\times[\pi_0^{(n)}\C]_1(c_0',c_0'')
\ar@{.>}[r]^{{}^{\pi_0^{(n)}\C}_{\phantom{Q}}\!\circ}\ar@{=}[d]_{def}
&[\pi_0^{(n)}\C]_1(c_0,c_0'')\ar@{=}[dd]^{def}
\\
\pi_0^{(n-1)}(\C_1(c_0,c_0'))\times \pi_0^{(n-1)}(\C_1(c_0',c_0''))\ar@{=}[d]_{(1)}
\\
\pi_0^{(n-1)}(\C_1(c_0,c_0')\times \C_1(c_0,c_0'))\ar[r]_{\pi^{(n-1)}({}^{\C}\!\circ)}
&\pi_0^{(n-1)}(\C_1(c_0,c_0''))
}
$$
For an object $c_0$, unit morphism is the dotted arrow below:
$$
\xymatrix@C=26ex{
\Id{n}\ar@{.>}[r]^{u(c_0)}\ar@{=}[d]_{(2)}
&[\pi_0^{(n)}\C]_1(c_0,c_0)\ar@{=}[d]^{def}
\\
\pi_0^{(n-1)}(\Id{n-1})\ar[r]_{\pi_0^{(n-1)}(u(c_0))}
&\pi_0^{(n-1)}(\C_1(c_0,c_0))
}
$$
Equalities $(1)$ and $(2)$ hold because $\pi_0^{(n-1)}$ commutes with (finite) products by induction hypothesis.

\begin{Claim}
The pair $([\pi_0^{(n)}\C]_0,[\pi_0^{(n)}\C]_1^{-,-})$, with composition and units as defined above,
satisfies axioms for a \emph{(n-1)}category.

Moreover, $\pi_0^{(n)}\C$ is a \emph{(n-1)}groupoid.
\end{Claim}
\begin{proof}
We have to prove associativity and unit axioms.

Concerning associativity, let objects $c_0,c_0',c_0'',c_0'''$ be given. Let us consider the following equalities of
morphisms
$$
[\pi_0^{(n)}\C]_1(c_0,c_0')\tm[\pi_0^{(n)}\C]_1(c_0',c_0'')\tm[\pi_0^{(n)}\C]_(c_0'',c_0''')1\to[\pi_0^{(n)}\C]_1(c_0,c_0'''),
$$
\begin{eqnarray*}
  ({}^{\pi_0\C}\!\circ^{c_0,c_0',c_0''}\times id)\bullet^0  {}^{\pi_0\C}\!\circ^{c_0,c_0'',c_0'''}
    &\eq{def}&\pi_0^{(n-1)}({}^{\C}\!\circ^{c_0,c_0',c_0''}\times id)\bullet^0  \pi_0^{(n-1)}({}^{\C}\!\circ^{c_0,c_0'',c_0'''}) \\
    &\eq{1}& \pi_0^{(n-1)}(({}^{\C}\!\circ^{c_0,c_0',c_0''}\times id)\bullet^0  {}^{\C}\!\circ^{c_0,c_0'',c_0'''}) \\
    &\eq{2}& \pi_0^{(n-1)}((id\times {}^{\C}\!\circ^{c_0',c_0'',c_0'''})\bullet^0  {}^{\C}\!\circ^{c_0,c_0',c_0'''})  \\
    &\eq{3}& \pi_0^{(n-1)}(id\times {}^{\C}\!\circ^{c_0',c_0'',c_0'''})\bullet^0  \pi_0^{(n-1)}( {}^{\C}\!\circ^{c_0,c_0',c_0'''}) \\
    &\eq{def}& (id\times {}^{\pi_0\C}\!\circ^{c_0',c_0'',c_0'''})\bullet^0   {}^{\pi_0\C}\!\circ^{c_0,c_0',c_0'''}
\end{eqnarray*}
where, $(1)$ and $(3)$ hold by functoriality of $\pi_0^{(n-1)}$,
$(2)$ by associativity of 0-composition in $\C$.

Turning to left-units, for every pair $c_0,c_0'$, one has the following equalities of morphisms:
$$
[\pi_0^{(n)}\C]_1(c_0',c_0)\to [\pi_0^{(n)}\C]_1(c_0',c_0)
$$
\begin{eqnarray*}
  \lambda\bullet^0 (id\times u(c_0))\bullet^0  \circ^{c_0',c_0,c_0}
   &\eq{def}&  \pi_0^{(n-1)}(\lambda)\bullet^0 \pi_0^{(n-1)}(id\times u(c_0))\bullet^0 \pi_0^{(n-1)}(\circ^{c_0',c_0,c_0})\\
   &\eq{1}&   \pi_0^{(n-1)}(\lambda\bullet^0 (id\times u(c_0))\bullet^0 \circ^{c_0',c_0,c_0})\\
   &\eq{2}&   \pi_0^{(n-1)}(id_{\C_1(c_0',c_0)})\\
   &\eq{3}&  id_{\pi_0^{(n-1)}\C_1(c_0',c_0)}\\
   &\eq{def}&  id_{[\pi_0^{(n)}\C]_1(c_0',c_0)}
\end{eqnarray*}
where $1$ and $3$ hold by functoriality of $\pi_0^{(n-1)}$,
$(2)$ by neutrality  of 0-identities in $\C$.\\

Right units are dealt the same way.\\

Finally, in order to show that $\pi_0^{(n)}\C$ is a $(n-1)$groupoid, two facts have to be proved:

1. all 1-cells of $\pi_0^{(n)}\C$ are equivalences.

2. all homs $[\pi_0^{(n)}\C]_1(c_0,c_0')$ are $(n-2)$groupoids.\\

The first fact is an easy consequence of the very definition of $\pi_0$ on compositions. In fact, let object
$\bar{c}_0$ and 1-cell $\tilde{x}:c_0\to c_0'$ in $\pi_{(n)}\C$ be given. Then by definition,
$-\circ\tilde{x}=\pi_0^{(n-1)}(-\circ x)$, where $x=\tilde{x}$ if $n>1$, or $\{x\}_{\sim}=\tilde{x}$ if $n=1$.
The result follows, for $\pi_0^{(n-1)}$ preserves equivalences.

To prove the second statement, let us consider the $(n-2)$category $[\pi_0^{(n)}\C]_1(c_0,c_0')$. This is defined to be
$\pi_0^{(n-1)}(\C_1(c_0,c_0'))$, where $\C_1(c_0,c_0')$ is a $(n-1)$groupoid.
For $\pi_0^{(n-1)}$ the result follows  by induction.
\end{proof}

\subsection{$\pi_0$ on morphisms}

Let a n-functor $F:\C\to\D$ be given. Then, $\pi_0^{(n)}F=([\pi_0^{(n)}F]_0,[\pi_0^{(n)}F]_1^{-,-})$, where
\begin{itemize}
\item $[\pi_0^{(n)}F]_0=F_0$
\item for every pair $c_0,c_0'$ in $[\pi_0^{(n)}\C]_0$
$$
[\pi_0^{(n)}F]_1^{c_0,c_0'}=\pi_0^{(n-1)}(F_1^{c_0,c_0'})
$$
\end{itemize}
\begin{Claim}
The pair $([\pi_0^{(n)}F]_0,[\pi_0^{(n)}F]_1^{-,-})$ satisfies axioms for $(n-1)$functors.
\end{Claim}
\begin{proof}
We have to prove functoriality w.r.t. composition and units.

Concerning composition, let objects $c_0,c_0',c_0''$ be given. Let us consider the following equalities of morphisms
$$
[\pi_0^{(n)}\C]_1(c_0,c_0')\times[\pi_0^{(n)}\C]_1(c_0',c_0'')\rightarrow [\pi_0^{(n)}\D]_1(Fc_0,Fc_0'')
$$
\begin{eqnarray*}
   {}^{\pi_0^{(n)}\C}\!\circ\bullet^0  [\pi_0^{(n)}F]_1^{c_0,c_0''}
   &\eq{def}& \pi_0^{(n-1)}({}^{\C}\!\circ)\bullet^0  \pi_0^{n-1}(F_1^{c_0,c_0''})  \\
   &\eq{1}  & \pi_0^{(n-1)}({}^{\C}\!\circ\bullet^0  F_1^{c_0,c_0''})  \\
   &\eq{2}&  \pi_0^{(n-1)}((F_1^{c_0,c_0'}\times F_1^{c_0',c_0''})\bullet^0 \ {}^{\D}\!\circ)\\
   &\eq{3}&  \pi_0^{(n-1)}(F_1^{c_0,c_0'}\times F_1^{c_0',c_0''})\bullet^0 \pi_0^{(n-1)}( {}^{\D}\!\circ)\\
   &\eq{def}&  ([\pi_0^{(n)}F]_1^{c_0,c_0'}\times[\pi_0^{(n)}F]_1^{c_0',c_0''})\bullet^0 \ {}^{\pi_0^{(n)}\D}\!\circ\
\end{eqnarray*}
$(1)$ and $(3)$ are justified by the sesqui-functor $\pi_0^{(n-1)}$ preserving composition, $(2)$ is functoriality
w.r.t. composition of $F$.

Turning to identities, for every $c_0$ one has the following equalities of morphisms
$$
\Id{n}\to[\pi_0^{(n)}]_1(Fc_0,Fc_0)
$$
\begin{eqnarray*}
   u(c_0)\bullet^0 [\pi_0^{(n)}\D]_1^{c_0,c_0}
   &\eq{def}&\pi_0^{(n-1)}(u(c_0))\bullet^0  \pi_0^{(n-1)}(F_1^{c_0,c_0}) \\
   &\eq{1}&\pi_0^{(n-1)}(u(c_0)\bullet^0  F_1^{c_0,c_0})  \\
   &\eq{2}&\pi_0^{(n-1)}(u(Fc_0)  \\
   &\eq{def}&u(Fc_0)
\end{eqnarray*}
where $(1)$ holds for $\pi_0^{(n-1)}$ preserving composition, $(2)$ by functoriality w.r.t. units of $F$.
\end{proof}

\subsection{The underlying functor}
In order to be a sesqui-functor, $\pi_0$ must restrict to a functor between the underlying categories:
$$
\lfloor\pi_0^{(n)}\rfloor : \lfloor \mathbf{nGpd}\rfloor \to \lfloor\mathbf{(n-1)Gpd}\rfloor
$$
\begin{Claim}
Given the situation
$$
\xymatrix{\C\ar[r]^F&\D\ar[r]^G&\E}
$$
in \emph{$n$\textbf{Gpd}}, the assignments given above satisfy
\begin{enumerate}
  \item $\pi_0^{(n)}(F\bullet^0 G)=\pi_0^{(n)}F\bullet^0 \pi_0^{(n)}G$
  \item $\pi_0^{(n)}(id_{\C})=id_{\pi_0^{(n)}\C}$
\end{enumerate}
\end{Claim}
\begin{proof}
Let us check first number 1.
\begin{eqnarray*}
  [\pi_0^{(n)}(F\bullet^0  G)]_0
    &\eq{def}& [F\bullet^0  G]_0 \\
    &\eq{(1)}& F_0  G_0 \\
    &\eq{def}& [\pi_0^{(n)}F]_0 \bullet^0  [\pi_0^{(n)}G]_0
\end{eqnarray*}
where $(1)$ holds by  composition of n-functors.

Moreover, for objects $c_0,c_0'$,
\begin{eqnarray*}
  [\pi_0^{(n)}(F\bullet^0  G)]_1^{c_0,c_0'}
  &\eq{def}& \pi_0^{(n-1)}([F\bullet^0  G]_1^{c_0,c_0'}) \\
   &\eq{1}& \pi_0^{(n-1)}(F_1^{c_0,c_0'}\bullet^0  G_1^{Fc_0,Fc_0'}) \\
   &\eq{2}& \pi_0^{(n-1)}(F_1^{c_0,c_0'})\bullet^0 \pi_0^{(n-1)}(G_1^{Fc_0,Fc_0'}) \\
   &\eq{def}& \pi_0^{(n)}(F_1^{c_0,c_0'})\bullet^0  \pi_0^{(n)}(G_1^{Fc_0,Fc_0'})
\end{eqnarray*}
where $(1)$ holds again by composition of n-functors, and $(2)$ by functoriality of $\pi_0^{(n-1)}$.\\

Now, let us check number 2.
\begin{eqnarray*}
  [\pi_0^{(n)}(id_{\C})]_0
    &=& [id_{\C}]_0\\
    &=& id_{\C_0} \\
    &=& id_{[\pi_0^{(n)}\C]_0}\\
    &=& \Big[id_{\pi_0^{(n)}\C}\Big]_0
\end{eqnarray*}
follows straight from definitions.\\

For objects $c_0,c_0'$
\begin{eqnarray*}
  [\pi^{(n)}(id_{\C})]_1^{c_0,c_0'} &=& \pi_0^{(n-1)}([id_{\C}]_1^{c_0,c_0'}) \\
   &=& \pi_0^{(n-1)}(id_{\C_1(c_0,c_0')}) \\
   &\eq{1}& id_{\pi_0^{(n-1)}(\C_1(c_0,c_0'))} \\
   &=& id_{[\pi_0^{(n)}\C]_1(c_0,c_0')} \\
   &=& \Big[id_{\pi_0^{(n)}\C} \Big]_1^{c_0,c_0'}
\end{eqnarray*}
where every equality above comes from definitions, but $(1)$ that holds by functoriality of $\pi_0^{(n-1)}$.
\end{proof}

\subsection{$\pi_0$ on 2-morphisms}
The \emph{action} of $\pi_0^{(n)}$ on 2-morphisms is more sensible to define. Indeed, on objects and morphisms we had the
\emph{object-part}  of definitions that were mere equalities, while the homs were given by induction.
Now the situation is different, since the \emph{object-part} of a $2$-cell is a  map involving also the
1-cells of codomain $n$-groupoid. That is the reason why we must analyze carefully what happens in low dimension, in order to
start induction properly.\\

\framebox[1.1\width]{$n=2$}\\

Let us consider $\pi_0'':2\mathbf{Gpd}\to \mathbf{Gpd}$:
$$
\xymatrix@C=8ex{\C\ar@/^4ex/[r]^F_{}="1" \ar@/_4ex/[r]_G^{}="2"&\D\ar@{=>}"1";"2"^{\alpha}}
\quad\mapsto\quad
\xymatrix@C=12ex{\pi_0''\C\ar@/^4ex/[r]^{\pi_0'' F}_{}="1" \ar@/_4ex/[r]_{\pi_0'' G}^{}="2"&\pi_0''\D\ar@{=>}"1";"2"^{\pi_0''\alpha}}
$$
(in order to simplify notation, for low dimensions we use primes).

Now it is clear that $[\pi_0''\D]_1=\D_1/\sim$, where $\sim$ is the equivalence relation on the set of 1-cells of $\D$ given
by iso-2-cells. Call $p:\D_1\to\D_1/\sim$ the canonic projection onto the quotient.
Then we let
$$
[\pi_0''\alpha]_0= \alpha_0 \cdot p
$$
This is well defined, since equivalence classes in $\D_1$ respect 1-cell's sources and targets. Moreover
they are compatible with 0-composition, \emph{i.e.} $p(d_1\circ d_1')=p(d_1)\circ p(d_1')$.
In fact for every object
$c_0$ of $\C$, one has
$$
[\pi_0''\alpha]_0(c_0)= \{\alpha_0(c_0)\}_{\sim}:Fc_0=[\pi_0''F](c_0)\to [\pi_0''G](c_0)=Gc_0
$$

Hence if we choose a 1-cell $\tilde{c_1}:c_0\to c_0'$, say $\tilde{c_1}=\{c_1\}_{\sim}$, then $\pi_0''$ sends the
2-isomorphism  $\alpha_1^{c_0,c_0'}\!(c_1)$ in the equality
$$
[\pi_0''\alpha]_{c_0} \circ\{Gc_1\}_{\sim}= \{Fc_1\}_{\sim} \circ [\pi_0''\alpha]_{c_0'}
$$
This proves that $\pi_0''\alpha$ is a natural isomorphism of groupoids, i.e. a 2-morphisms in \textbf{Gpd}.\\[3ex]

\framebox[1.1\width]{$n>2$}\\

More generally, suppose we are given a 2-morphism $\alpha:F\Rightarrow G:\C\to \D$ in $n$\textbf{Gpd}.
Then we define $\pi_0^{(n)}\alpha$ as the pair $([\pi_0^{(n)}\alpha]_0,[\pi_0^{(n)}\alpha]_1^{-,-})$, where
\begin{itemize}
  \item $[\pi_0^{(n)}\alpha]_0=\alpha_0,$
  \item for every pair of objects $c_0,c_0'$ of $\C$,
$$
[\pi_0^{(n)}\alpha]_1^{c_0,c_0'}=\pi_0^{(n-1)}(\alpha_1^{c_0,c_0'}).
$$
\end{itemize}
\begin{Claim}
The pair $([\pi_0^{(n)}\alpha]_0,[\pi_0^{(n)}\alpha]_1^{-,-})$ satisfies  axioms for $(n-1)$transformations.
\end{Claim}
\begin{proof}
It is well-defined on objects, since $n>2$. Moreover, it is well-defined also on homs. In fact, given objects $c_0$ and
$c_0'$, consider the diagram:
\begin{changemargin}{-10ex}{-10ex}{\small
$$
\xymatrix@C=3ex@R=16ex{
[\pi_0^{(n)}\C]_1(c_0,c_0')\ar@{=}[dr]
\ar[rrrrr]^{[\pi_0^{(n)}G]_1^{c_0,c_0'}}_{}="x1"
\ar[ddd]_{[\pi_0^{(n)}F]_1^{c_0,c_0'}}
&&&&&[\pi_0^{(n)}\D]_1(Gc_0,Gc_0')\ar@{=}[dl]
\ar[ddd]^{[\pi_0^{(n)}\alpha]_{c_0}\circ-}
\\
&\pi_0^{(n-1)}(\C_1(c_0,c_0'))
\ar[rrr]^{}="x2"
\ar[d]_{\pi_0^{(n-1)}F_1^{c_0,c_0'}}
&&&\pi_0^{(n-1)}(\D_1(Gc_0,Gc_0'))\ar@{}[dlll]|(.3){}="1"|(.7){}="2"
\ar[d]^{\pi_0^{(n-1)}(\alpha_{c_0}\circ-)}
\\
&\pi_0^{(n-1)}(\D_1(Fc_0,Fc_0'))
\ar[rrr]_{}="y2"
&&&\pi_0^{(n-1)}(\D_1(Fc_0,Gc_0'))
\\
[\pi_0^{(n)}\D]_1(Fc_0,Fc_0')\ar@{=}[ur]
\ar[rrrrr]_{-\circ [\pi_0^{(n)}\alpha]_{c_0'}}^{}="y1"
&&&&&[\pi_0^{(n)}\D]_1(Gc_0,Gc_0')\ar@{=}[ul]
\ar@{=>}"1";"2"|{\pi_0^{(n-1)}\alpha_1^{c_0,c_0'}}
\ar@{}"x1";"x2"|(.9){\pi_0^{(n-1)}G_1^{c_0,c_0'}}
\ar@{}"y1";"y2"|(.9){\pi_0^{(n-1)}(-\circ \alpha_{c_0'})}
}
$$}
\end{changemargin}
Up and left squares are the definition of $\pi_0^{(n)}$ on n-functors (w.r.t. hom-components),
down and right commute by definition of composition in $\pi_0^{(n)}(\D)$, and of
$[\pi_0^{(n)}\alpha]_0$ above.

And that is all. In fact, coherence w.r.t. composition and units are satisfied because their diagrams-equations are
$\pi_0^{n-1}$ of corresponding diagrams-equations that hold for $\alpha$.
\end{proof}

For the same reason $\pi_0^{(n)}$
\begin{enumerate}
  \item is functorial w.r.t. vertical composition and units of 2-morphisms
  \item preserves reduced horizontal composition
\end{enumerate}
i.e. it is a sesqui-functor.
\begin{proof}
\begin{enumerate}
  \item Suppose we are given
  $$
\omega:E\Rightarrow F:\C\rightarrow \D\qquad\alpha:F\Rightarrow G:\C\rightarrow \D
  $$
in $n$\textbf{Gpd}. When $n>2$, for every object $c_0$ in $\pi_0^{(n)}\C$ one has
\begin{eqnarray*}
[\pi_0^{(n)}(\omega\alpha)]_0(c_0)
   &\eq{1}&(\omega\alpha)_0(c_0)   \\
   &\eq{2}& \omega_0c_0\circ\alpha_0c_0 \\
   &\eq{3}& [\pi_0^{(n)}(\omega)]_0(c_0)\circ[\pi_0^{(n)}(\alpha)]_0(c_0) \\
\end{eqnarray*}
where $(1)$ and $(3)$ hold by definition of $\pi_0^{(n)}$, and $(2)$ by definition of vertical composition.
When $n=2$ , one has
\begin{eqnarray*}
[\pi_0^{(n)}(\omega\alpha)]_0(c_0)
   &\eq{1}& p\big((\omega\alpha)_0(c_0)\big)   \\
   &\eq{2}& p\big(\omega_0c_0\circ\alpha_0c_0\big) \\
   &\eq{3}& p(\omega_0c_0)\circ p(\alpha_0c_0) \\
   &\eq{4}& [\pi_0^{(n)}(\omega)]_0(c_0)\circ[\pi_0^{(n)}(\alpha)]_0(c_0) \\
\end{eqnarray*}
where $(1)$ and $(4)$ hold by definition of $\pi_0^{(n)}$, $(2)$ by definition of vertical composition and
$(3)$ for $p$ is compatible with 0-composition.

Now suppose we are given objects $c_0$, $c_0'$ of $\pi_0^{(n)}\C$.
Then one has:
\begin{eqnarray*}
[\pi_0^{(n)}(\omega\alpha)]_1^{c_0,c_0'}
   &\eq{1}&\pi_0^{(n-1)}([\omega\alpha]_1^{c_0,c_0'})  \\
   &\eq{2}&\pi_0^{(n-1)}\big((\alpha_1^{c_0,c_0'}\bullet_R^0 [\omega_{c_0}\circ-])\bullet^1(\omega_1^{c_0,c_0'}\bullet_R^0 [-\circ\alpha_{c_0'}]) \big) \\
   &\eq{3}&(\pi_0^{(n-1)}\alpha_1^{c_0,c_0'}\bullet_R^0 \pi_0^{(n-1)}[\omega_{c_0}\circ-])\bullet^1(\pi_0^{(n-1)}\omega_1^{c_0,c_0'}\bullet_R^0 \pi_0^{(n-1)}[-\circ\alpha_{c_0'}]) \\
   &\eq{4}&([\pi_0^{(n)}\alpha]_1^{c_0,c_0'}\bullet_R^0 \pi_0^{(n-1)}[\omega_{c_0}\circ-])\bullet^1([\pi_0^{(n)}\omega]_1^{c_0,c_0'}\bullet_R^0 \pi_0^{(n-1)}[-\circ\alpha_{c_0'}]) \\
   &\eq{5}&([\pi_0^{(n)}\alpha]_1^{c_0,c_0'}\bullet_R^0 [[\pi_0^{(n)}\omega]_{c_0}\circ-])\bullet^1([\pi_0^{(n)}\omega]_1^{c_0,c_0'}\bullet_R^0 [-\circ[\pi_0^{(n)}\alpha]_{c_0'}])
\end{eqnarray*}
where
$(1)$ and $(4)$ hold by the definition of $\pi_0^{(n)}$ on 2-morphisms w.r.t. homs,
$(2)$ by definition of vertical composition,
$(3)$ because $\pi_0^{(n)}$ is a sesqui-functor (induction),
$(5)$ by definition of $\pi_0^{(n)}\D$ 0-composition.

Concerning units, let us consider $id_F:F\Rightarrow F$:
\begin{eqnarray*}
[\pi_0^{(n)}id_F]_0(c_0)
   &=& [id_F]_0(c_0)  \\
   &=&  id_{Fc_0}\\
   &=&  id_{[\pi_0^{(n)}F]c_0}
\end{eqnarray*}
and also
\begin{eqnarray*}
[\pi_0^{(n)}id_F]_1^{c_0,c_0'}
   &=& \pi_0^{(n-1)}([id_F]_1^{c_0,c_0'})  \\
   &=& id_{\pi_0^{(n-1)}(F_1^{c_0,c_0'})}  \\
   &=& id_{[\pi_0^{(n)}F]_1^{c_0,c_0'}}.
\end{eqnarray*}

  \item We prove the statement for reduced left-composition. Suppose 2-morphism $\alpha$ as above, and morphism
$N:\B\to\C$ be given.

When $n>2$, for any objects $b_0$ of $\B$, one has
\begin{eqnarray*}
  [\pi_0^{(n)}(N\bullet_L^0\alpha)]_0(b_0)
   &\eq{1}& [N\bullet_L^0\alpha]_0(b_0) \\
   &\eq{2}& \alpha_0(N(b_0)) \\
   &\eq{3}& [\pi_0^{(n)}\alpha]_0(N(b_0)) \\
   &\eq{4}& [\pi_0^{(n)}\alpha]_0([\pi_0^{(n)}N](b_0)) \\
   &\eq{5}& [\pi_0^{(n)}N\bullet_L^0\pi_0^{(n)}\alpha]_0(b_0)
\end{eqnarray*}
where $(1)$, $(3)$ and $(4)$ hold by definition of $\pi_0^{(n)}$, $(2)$ and $(5)$ by definition of reduced left-composition,

Moreover, let us choose two objects $b_0$ and $b_0'$ in $\B$. Then one has
\begin{eqnarray*}
[\pi_0^{(n)}(N\bullet_L^0\alpha)]_1^{b_0,b_0'}
   &\eq{1}& \pi_0^{(n-1)}([N\bullet_L^0\alpha]_1^{b_0,b_0'}) \\
   &\eq{2}& \pi_0^{(n-1)}(N_1^{b_0,b_0'}\bullet_L^0\alpha_1^{Nb_0,Nb_0'}) \\
   &\eq{3}& \pi_0^{(n-1)}(N_1^{b_0,b_0'})\bullet_L^0\pi_0^{(n-1)}(\alpha_1^{Nb_0,Nb_0'}) \\
   &\eq{4}& [\pi_0^{(n)}N]_1^{b_0,b_0'}\bullet_L^0[\pi_0^{(n)}\alpha]_1^{\pi_0^{(n)}N(b_0),\pi_0^{(n)}N(b_0')}
\end{eqnarray*}
where $(1)$ and $(4)$ hold by definition of $\pi_0^{(n)}$, $(2)$  by definition of reduced left-composition
and $(3)$ by induction hypothesis.

When $n=2$ the calculation can be carried on similarly, as we did for vertical composites above.\\

Finally, concerning reduced right-composition, the proof is similar, as the definition.
\end{enumerate}
\end{proof}
\subsection{$\pi_0$ commutes with (finite) products}
We will show that $\pi_0^{(n)}$ preserves binary products and the terminal object.
\begin{Proposition}
Let $\C$ and $\D$ be n-groupoids. Then
\begin{enumerate}
  \item $\pi_0^{(n)}(\C\times\D)=\pi_0^{(n)}\C\times\pi_0^{(n)}\D$
  \item $\pi_0^{(n)}\left(\Id{n}\right)=\Id{n-1}$
\end{enumerate}
\end{Proposition}
\begin{proof}
\begin{enumerate}
  \item Consider the following equalities:
\begin{eqnarray*}
[\pi_0^{(n)}(\C\times\D)]_0
   &\eq{1}& \pi_0^{(n-1)}([\C\times\D]_0) \\
   &\eq{2}& \pi_0^{(n-1)}(\C_0\times\D_0) \\
   &\eq{3}& \pi_0^{(n-1)}\C_0\times\pi_0^{(n-1)}\D_0 \\
   &\eq{4}& [\pi_0^{(n)}\C]_0\times[\pi_0^{(n)}\D]_0
\end{eqnarray*}
Moreover, let objects $(c_0,d_0)$ and $(c_0',d_0')$ of $\C\times\D$ be given. Then
\begin{eqnarray*}
[\pi_0^{(n)}(\C\times\D)]_1^{(c_0,d_0),(c_0',d_0')}
   &\eq{1}&  \pi_0^{(n-1)}\left([\C\times\D]_1^{(c_0,d_0),(c_0',d_0')}\right)\\
   &\eq{2}&  \pi_0^{(n-1)}\left(\C_1^{c_0,c_0'}\times\D_1^{d_0,d_0'}\right)\\
   &\eq{3}&  \pi_0^{(n-1)}(\C_1^{c_0,c_0'})\times\pi_0^{(n-1)}(\D_1^{d_0,d_0'})\\
   &\eq{4}&  [\pi_0^{(n)}\C]_1^{c_0,c_0'}\times[\pi_0^{(n)}\D]_1^{d_0,d_0'}
\end{eqnarray*}
In both cases, $(1)$ and $(4)$ follow from the definition of $\pi_0^{(n)}$, $(2)$ holds by definition of products and $(3)$ by
induction.
  \item Consider the following equalities:
\begin{eqnarray*}
\left[\pi_0^{(n)}\left(\Id{n}\right)\right]_0
   &\eq{1}& \left[\Id{n}\right]_0 \\
   &\eq{2}& \{*\}  \\
   &\eq{3}& \left[\Id{n-1}\right]_0  \\
   &\eq{4}& \left[\pi_0^{(n)}\left(\Id{n-1}\right)\right]_0
\end{eqnarray*}
where $(1)$ and $(4)$ hold by definition of $\pi_0^{(n)}$,
$(2)$ and $(3)$ hold by definition of terminal n-category.
\begin{eqnarray*}
\left[\pi_0^{(n)}\left(\Id{n}\right)\right]_1^{*,*}
   &\eq{1}& \pi_0^{(n-1)}\left(\left[\Id{n}\right]_1^{*,*}\right)  \\
   &\eq{2}& \pi_0^{(n-1)}\left(\Id{n-1}\right)  \\
   &\eq{3}& \Id{n-2}  \\
   &\eq{4}& \left[\Id{n-1}\right]_1^{*,*}
\end{eqnarray*}
where $(1)$ and $(4)$ follow from the definition of $\pi_0^{(n)}$,
 $(2)$ holds by definition of terminal n-category and $(3)$ by
induction.
\end{enumerate}
\end{proof}

\subsection{$\pi_0$ preserves equivalences}

\begin{Proposition}
Let $F:\C\to\D$ be an equivalence of n-groupoids. Then
$\pi_0^{(n)}F:\pi_0^{(n)}\C\to\pi_0^{(n)}\D$ is an equivalence
of (n-1)-groupoids.
\end{Proposition}
\begin{proof}
If $n=1$, then $F$ is an equivalence of categories, then  $\pi_0'F$ is clearly an isomorphism.
Hence we may well suppose $n>1$.
\begin{enumerate}
\item Let objects $c_0,c_0'$ of $\C$ be given. Then, by definition, $[\pi_0^{(n)}F]_1^{c_0,c_0'}=\pi_0^{(n-1)}(F_1^{c_0,c_0'})$.
This is an equivalences of (n-2)groupoids, since $F_1^{c_0,c_0'}$ is an equivalence of \makebox{(n-1)groupoids} and $\pi_0^{(n-1)}$
preserves equivalences by induction.
\item Let an object $d_0$ of $\pi_0^{(n)}\D$ be given. This is indeed an object of $\D$, hence there exists a pair
$$
(c_0,d_1:Fc_0\to d_0)
$$
with $d_1$ being an equivalence in $\D$.

Now, $c_0$ is also an object of $\pi_0^{(n)}\C$, and $d_1$ (eventually $\{d_1\}_{\sim}$ if $n=2$) is a 1-cell (hence an
equivalence) of the (n-1) groupoid $\pi_0^{(n)}\D$.
\end{enumerate}
\end{proof}

\subsection{A remark on $\pi_0$}
By the discussion in the previous section,  $\pi_0^{(n)}$ is clearly extendible to a functor
$$
\lfloor n\mathbf{Cat}\rfloor \to \lfloor (n-1)\mathbf{Cat}\rfloor
$$
on underlying categories. Difficulties arise when we try to extend it further to a sesqui-functor. In fact, even for $n=1$ our
definition fails, since a natural transformation $\alpha:F\Rightarrow G:\C\to\D$ does not imply that the
maps $\pi_0F\Rightarrow \pi_0G:\pi_0\C\to\pi_0\D$ are equal. This is the case when $\alpha$ is a natural isomorphism, {\em i.e.}
for every $c_0$ in $\C$, $\alpha_{c_0}$ is an isomorphism.

This case may be generalized in order to get a sesqui-functor on $n$-categories that remove the \emph{obstruction} by considering only
$n$-natural equivalences:
$$
\pi_0^{(n)} : n\mathbf{Cat}_{eq}\to (n-1)\mathbf{Cat}_{eq}
$$
This gives a chance to develop the theory in a more general case.

On the other hand, if one wants to keep into account \emph{all} n-transformations, one direction can be to
consider a generalization of \emph{connected-components} functor, rather then \emph{iso-classes} functor.
At the moment, we have not explored this perspective since it seems to give rise to a completely  different theory, not
consistent with low-dimensional problems we aim to generalize.

\section{The discretizer}

Purpose of this section is to introduce the family of
sesqui-functors $\{D^{(n)}\}_{n\in\mathbf{N}}$ that extends the
\emph{discrete-groupoid} functor $\mathbf{Set}\to\mathbf{Gpd}$.

\begin{Defprop}
For any integer $n>0$, there exists a $n$-discrete sesqui-functor
$$
D^{(n)}:(n-1)\mathbf{Gpd} \rightarrow n\mathbf{Gpd}
$$
according to the following recursive definition.

Moreover it commutes with finite products and it preserves equivalences
\end{Defprop}

\framebox[1.1\width]{$n=1$}\\
$$
D^{(1)}:\mathbf{Set} \rightarrow \mathbf{Gpd}
$$
is the functor (= trivial sesqui-functor) that assigns to a set $C$ the discrete
groupoid $D(C)$, {\em i.e.} with objects the elements of $C$ and only identity arrows. \\

\framebox[1.1\width]{$n>1$}\\

\subsection{$D^{(n)}$ on objects and morphisms}
Let a $(n-1)$-groupoid $\C$ be given. Then $D^{(n)}\C=([D^{(n)}\C]_0,[D^{(n)}\C]_1(-,-))$, where
\begin{itemize}
  \item $[D^{(n)}\C]_0=C_0$
  \item for every pair $c_0,c_0'$ in $[D^{(n)}\C]_0$,
  $$
  [D^{(n)}\C]_1(c_0,c_0')=D^{(n-1)}(\C_1(c_0,c_0'))
  $$
\end{itemize}
For a triple of objects $c_0,c_0',c_0''$, composition is the dotted arrow:
$$
\xymatrix@C=16ex{
[D^{(n)}\C]_1(c_0,c_0')\times[D^{(n)}\C]_1(c_0',c_0'')\ar@{.>}[r]^{{}^{D^{(n)}\C}\!\circ}\ar@{=}[d]_{def}
&[D^{(n)}\C]_1(c_0,c_0'')\ar@{=}[dd]^{def}
\\
D^{(n-1)}(\C_1(c_0,c_0'))\times D^{(n-1)}(\C_1(c_0',c_0''))\ar@{=}[d]_{(1)}
\\
D^{(n-1)}(\C_1(c_0,c_0')\times \C_1(c_0,c_0'))\ar[r]_{D^{(n-1)}({}^{\C}\!\circ)}
&D_0^{(n-1)}(\C_1(c_0,c_0''))
}
$$
For an object $c_0$, unit morphism is the dotted arrow:
$$
\xymatrix@C=26ex{
\Id{n}\ar@{.>}[r]^{u(c_0)}\ar@{=}[d]_{(2)}
&[D_0^{(n)}\C]_1(c_0,c_0)\ar@{=}[d]^{def}
\\
D_0^{(n-1)}(\Id{n-1})\ar[r]_{D_0^{(n-1)}(u(c_0))}
&D_0^{(n-1)}(\C_1(c_0,c_0))
}
$$
Equalities $(1)$ and $(2)$ hold because $D^{(n-1)}$ commutes with (finite) products by inductive hypothesis.\\

Let a $(n-1)$functor $F:\C\to\D$ be given. Then, $D^{(n)}F=([D^{(n)}F]_0,[D^{(n)}F]_1^{-,-})$, where
\begin{itemize}
\item $[D^{(n)}F]_0=F_0$
\item for every pair $c_0,c_0'$ in $[D^{(n)}\C]_0$
$$
[D^{(n)}F]_1^{c_0,c_0'}=D^{(n-1)}(F_1^{c_0,c_0'})
$$
\end{itemize}

Notice that, since definitions of $\pi_0^{(n)}$ and of $D^{(n)}$ are formally identical, proving that all above
is consistent is a matter of a syntactical substitution of the first with the second in the corresponding
proofs concerning $\pi_0^{(n)}$. Hence we have defined the following functor between
underlying categories:
$$
\lfloor D^{(n)}\rfloor : \lfloor \mathbf{(n-1)Gpd}\rfloor \to \lfloor\mathbf{nGpd}\rfloor
$$

\subsection{$D^{(n)}$ on 2-morphisms}
Unlike that of $\pi_0^{(n)}$, the definition of $D^{(n)}$ is straightforward since the beginning of induction.

Let $\alpha:F\Rightarrow G:\C\to\D$ in $n$\textbf{Gpd} be given. As usual, $D^{(n)}(\alpha)=([D^{(n)}\alpha]_0,[D^{(n)}\alpha]_1^{-,-})$
where
\begin{itemize}
\item $[D^{(n)}\alpha]_0 = \alpha_0$
\item for every pair of objects $c_0,c_0'$ in $[D\C]_0$
$$
[D^{(n)}\alpha]_1^{c_0,c_0'}=D^{(n-1)}(\alpha_1^{c_0,c_0'})
$$
\end{itemize}

\emph{Modulo} the syntactical conversion mentioned above, we can prove that $D^{(n)}$ is a sesqui-functor that commutes with products and preserves equivalences.
Hence it is well defined on $n$-groupoids.

\subsection{A remark on $D^{(n)}$}
Differently from $\pi_0^{(n)}$, the definition of $D^{(n)}$  extends with no changes to $n$-categories. In fact it
\emph{lives} more naturally in an $n$-categorical setting, and our definition is just its restriction to $n$-groupoids

\section{The adjunction $\pi_0^{(n)} \dashv D^{(n)}$}

\subsection{In low dimension}
The following adjunction has been extensively studied by category-theorist:
$$
\xymatrix@C=15ex{
\mathbf{Gpd}
\ar@/^2ex/[r]^{\pi_0}
\ar@{}[r]|{\perp}
&\mathbf{Set}
\ar@/^2ex/[l]^{D}
}
$$
Let us describe it briefly, as it will be the first step of an inductive definition for general n-groupoids.
\begin{description}
  \item[co-unit]  Given a set $S$, it is clear that $\pi_0(D(S))=S$, hence co-unit $\epsilon$ is the identity.
  \item[unit] For a groupoid $\C$, the unit
  $$
\eta=\eta_{\C}:\C\to D(\pi_0(\C))
  $$
is the {\em projection} given by $\eta_{\C}(c_0)=[c_0]_{\sim}$, and for $c_1:c_0\to c_0'$,
$\eta_{\C}(c_1)=id_{[c_0]_{\sim}}=id_{[c_0']_{\sim}}$
\end{description}

\subsection{The general setting}
Here and in the following, let an integer $n>1$ be given.

For an (n-1)groupoid $\ES$, the co-unit of the adjunction is still the identity. In fact
$\pi_0^{(n)}(D^{(n)}(\ES))=\ES$.
\begin{proof}
$$
[\pi_0^{(n)}(D^{(n)}(\ES))]_0=[D^{(n)}(\ES))]_0=\ES_0
$$
by definition, as $n>1$. On the other side, for objects $s_0,s_0'$ one has
$$
[\pi_0^{(n)}(D^{(n)}(\ES))]_1(s_0,s_0')=\pi_0^{(n-1)}([D^{(n)}(\ES)]_1(s_0,s_0'))=
\pi_0^{(n-1)}(D^{(n-1)}(\ES_1(s_0,s_0')))=\ES_1(s_0,s_0')
$$
where the last equality is precisely the induction hypothesis.

Same argument holds for morphisms.
\end{proof}

Concerning the unit of the adjunction, we state the following
\begin{Definition}
Let us fix n-groupoid $\C$. Then
$$
\eta_{\C}^{(n)}:\C\to D^{(n)}\left(\pi_0^{(n)}(\C)\right)
$$
consists of the following data:
\begin{itemize}
  \item $\left[\eta_{\C}^{(n)}\right]_0=id_{\C_0}:\C_0 \to \left[D^{(n)}\left(\pi_0^{(n)}(\C)\right)\right]_0=\C_0$
  \item for any pair $c_0,c_0'$ of objects of $\C$, $\left[\eta_{\C}^{(n)}\right]_1^{c_0,c_0'}$
  is the dotted arrow below:
$$
\xymatrix{
\C_1(c_0,c_0')\ar@{.>}[r]\ar@/_8ex/[ddr]_(.35){\eta_{\C_1(c_0,c_0')}^{(n-1)}}
&[D^{(n)}(\pi_0^{(n)}(\C))]_1(c_0,c_0')\ar@{=}[d]
\\
&D^{(n-1)}([\pi_0^{(n)}(\C)]_1(c_0,c_0'))\ar@{=}[d]
\\
&D^{(n-1)}(\pi_0^{(n-1)}(\C_1(c_0,c_0')))}
$$
\end{itemize}
\end{Definition}
\begin{Claim}\label{zz:funct}
The pair $<id_{\C_0},\eta_{\C_1(c_0,c_0')}^{(n-1)}>$ is a n-functor.
\end{Claim}
\begin{proof}
The diagrams marked $(i)$ and $(ii)$ express functoriality w.r.t. composition and units respectively, for any triple
$c_0,c_0',c_0''$ of objects of $\C$:
\begin{changemargin}{-10ex}{10ex}
$$
\xymatrix@C=15ex@R=10ex{
\C_1(c_0,c_0')\times\C_1(c_0',c_0'')
\ar[r]^{{}^{\C}\!\circ}
\ar[d]_{\eta_{\C_1(c_0,c_0')}^{(n-1)}\times\eta_{\C_1(c_0',c_0'')}^{(n-1)}}
&\C_1(c_0,c_0'')
\ar[dd]^{\eta_{\C_1(c_0,c_0'')}^{(n-1)}}
\\
D^{(n-1)}\left(\pi_0^{(n-1)}(\C_1(c_0,c_0'))\right)\times D^{(n-1)}\left(\pi_0^{(n-1)}(\C_1(c_0',c_0''))\right)
\ar@{=}[d]\ar[dr]|{{}^{D^{(n)}\pi_0^{(n)}\C}\circ}
\ar@{}[ur]|{(i)}
\\
D^{(n-1)}\left(\pi_0^{(n-1)}(\C_1(c_0,c_0')\times\C_1(c_0',c_0''))\right)
\ar[r]_{D^{(n-1)}\pi_0^{(n-1)}({}^{\C}\!\circ)}
&D^{(n-1)}\left(\pi_0^{(n-1)}(\C_1(c_0,c_0''))\right)
}
$$
\end{changemargin}
$$
\xymatrix@C=20ex@R=10ex{
\Id{n-1}
\ar[r]^{{}^{\C}u(c_0)}
\ar[dr]|{{}^{D^{(n)}\pi_0^{(n)}\C}u(c_0)}="1"
\ar@{=}[d]
&\C_1(c_0,c_0)
\ar[d]^{\eta^{(n-1)}_{\C_1(c_0,c_0)}}="2"
\\
D^{(n-1)}\left(\pi_0^{(n-1)}\left(\Id{n-1}\right)\right)
\ar[r]_{D^{(n-1)}\pi_0^{(n-1)}\left({}^{\C}u(c_0)\right)}
&D^{(n-1)}\left(\pi_0^{(n-1)}\left(\C_1(c_0,c_0) \right)\right)
\ar@{}"1";"2"|{(ii)}}
$$
Lower triangles commute  by definition. Commutativity of external diagrams
will be proved by finite induction in the following lemmas. This will imply that $(i)$ and $(ii)$
commute.
\end{proof}

First we need the following conventional
\begin{Notation}
Given the $n$-category $\C$,  for $0\leq s< m\leq n$ we write
$$
\xymatrix@C=4ex{c_m:c_s\ar@{==>}[r]&c_s'}
$$
meaning
$$
c_m:\ c_{m-1}\to c_{m-1}': \cdots :\ c_{s+1}\to c_{s+1}': c_s\to c_s'
$$
are $m$-cell, $(m-1)$-cells, \dots, $(s+1)$-cells, $s$-cells of $\C$.\\

Furthermore we inductively define
$$
\C_m(c_{m-1},c_{m-1}'):=[\C_{m-1}(c_{m-2},c_{m-2}')]_1(c_{m-1},c_{m-1}')
$$
being $\C_1(c_0,c_0')$ given by the  definition of $n$-category.
\end{Notation}

\begin{Lemma}
Given
$$
\xymatrix@C=4ex{c_{n-j-1},k_{n-j-1}:c_0\ar@{==>}[r]&c_0'},\qquad \xymatrix@C=4ex{c_{n-j-1}',k_{n-j-1}':c_0'\ar@{==>}[r]&c_0''}
$$
the following diagram commutes in $j$-\textbf{Gpd}:
\begin{changemargin}{-16ex}{-10ex}{\small
\begin{equation}\label{zzz5}
\xymatrix@C=-10ex@R=10ex{
\C_{n-j}(c_{n-j-1},k_{n-j-1})\times\C_{n-j}(c_{n-j-1}',k_{n-j-1}')
\ar[r]^{\circ^0}
\ar[d]_{\eta^{(j)}_{\C_{n-j}(\diamond,\diamond)}\times\eta^{(j)}_{\C_{n-j}(\diamond,\diamond)}}
&
\C_{n-j}(c_{n-j-1}\circ^0 c_{n-j-1}',k_{n-j-1}\circ^0 k_{n-j-1}')
\ar[dd]^{\eta^{(j)}_{\C_{n-j}(\diamond,\diamond)}}
\\
D^{(j)}\pi_0^{(j)}\left(\C_{n-j}(c_{n-j-1},k_{n-j-1})\right)\times D^{(j)}\pi_0^{(j)}\left(\C_{n-j}(c_{n-j-1}',k_{n-j-1}')\right)
\ar@{=}[dd]
\\
&
D^{(j)}\pi_0^{(j)}\left(\C_{n-j}(c_{n-j-1}\circ^0 c_{n-j-1}',k_{n-j-1}\circ^0 k_{n-j-1}') \right)
\\
D^{(j)}\pi_0^{(j)}\left(\C_{n-j}(c_{n-j-1},k_{n-j-1})\times\C_{n-j}(c_{n-j-1}',k_{n-j-1}')\right)
\ar[ur]_{D^{(j)}\pi_0^{(j)}(\circ^0)}
}
\end{equation}}
\end{changemargin}
where we write $(\diamond,\diamond)$ when substitutes are clear from the context.
\end{Lemma}
\begin{proof}
By finite induction over $j$. Diagram above for $j=1$ is the following square, of groupoids and functors
\begin{changemargin}{-16ex}{10ex}{\small
$$
\xymatrix@C=10ex@R=10ex{
\C_{n-1}(c_{n-2},k_{n-2})\times\C_{n-1}(c_{n-2}',k_{n-2}')
\ar[r]^{\circ^0}
\ar[d]_{\eta^{(1)}_{\C_{n-1}(\diamond,\diamond)}\times\eta^{(1)}_{\C_{n-1}(\diamond,\diamond)}}
&
\C_{n-1}(c_{n-2}\circ^0 c_{n-2}',k_{n-2}\circ^0 k_{n-2}')
\ar[dd]^{\eta^{(1)}_{\C_{n-1}(\diamond,\diamond)}}
\\
D^{(1)}\pi_0^{(1)}\left(\C_{n-1}(c_{n-2},k_{n-2})\right)\times D^{(1)}\pi_0^{(1)}\left(\C_{n-1}(c_{n-2}',k_{n-2}')\right)
\ar@{=}[d]
\\
D^{(1)}\pi_0^{(1)}\left(\C_{n-1}(c_{n-2},k_{n-2})\times\C_{n-1}(c_{n-2}',k_{n-2}')\right)
\ar[r]_{D^{(1)}\pi_0^{(1)}(\circ^0)}
&
D^{(1)}\pi_0^{(1)}\left(\C_{n-1}(c_{n-2}\circ^0 c_{n-2}',k_{n-2}\circ^0 k_{n-2}') \right)
}
$$}
\end{changemargin}
For this to commute, it must commute on objects and on arrows. In order to prove this, let us consider
$$
c_n:\ c_{n-1}\to k_{n-1}:\ c_{n-2}\to k_{n-2}
$$
and
$$
c_n':\ c_{n-1}'\to k_{n-1}':\ c_{n-2}'\to k_{n-2}'
$$
Moreover, let us recall that
$$
\eta_{\X}^{(1)}:=<[Id_{\X_0}]_{\sim},Id_{[dom(-)]_{\sim}}>
$$

and apply: for the pair $(c_{n-1},c_{n-1}')$, on the lower-left one has
$$
\xymatrix@C=20ex{
(c_{n-1},c_{n-1}')
\ar@{|->}[r]^{[\eta^{(1)}]_0\times[\eta^{(1)}]_0}
&
([c_{n-1}]_{\sim},[c_{n-1}']_{\sim})}
$$
$$
\xymatrix@C=20ex{
= [(c_{n-1},c_{n-1}')]_{\sim}
\ar@{|->}[r]^{[\circ^0]_{\sim}}
&
[c_{n-1}\circ^0  c_{n-1}']_{\sim}
},
$$

where on the upper-right
$$
\xymatrix@C=20ex{
(c_{n-1},c_{n-1}')
\ar@{|->}[r]^{\circ^0}
&
c_{n-1}\circ^0  c_{n-1}'
}
$$
$$
\xymatrix@C=20ex{
= c_{n-1}\circ^0  c_{n-1}'
\ar@{|->}[r]^{[\eta^{(1)}]_0}
&
[c_{n-1}\circ^0  c_{n-1}']_{\sim}
}
$$

Similarly for arrows, on the lower-left one has
$$
\xymatrix@C=20ex{
(c_{n},c_{n}')
\ar@{|->}[r]^{[\eta^{(1)}]_1\times[\eta^{(1)}]_1}
&
(Id_{[c_{n-1}]_{\sim}},Id_{[c_{n-1}']_{\sim}})}
$$
$$
\xymatrix@C=20ex{
= Id_{[(c_{n-1},c_{n-1}')]_{\sim}}
\ar@{|->}[r]^{[\circ^0]_{\sim}}
&
Id_{[c_{n-1}\circ^0  c_{n-1}']_{\sim}}
},
$$

while on the upper-right
$$
\xymatrix@C=20ex{
(c_{n},c_{n}')
\ar@{|->}[r]^{\circ^0}
&
c_{n}\circ^0  c_{n}'
}
$$
$$
\xymatrix@C=20ex{
= c_{n}\circ^0  c_{n}'
\ar@{|->}[r]^{[\eta^{(1)}]_0}
&
Id_{[c_{n-1}\circ^0  c_{n-1}']_{\sim}}
},
$$
and this completes the case $j=1$.\\

Now let us assume, as induction hypothesis, that Lemma holds for $j-1$. We will prove
it holds for $j$.\\

To this end, suppose we are given
$$
c_{n-j+1}:\ c_{n-j}\to k_{n-j}:\ c_{n-j-1}\to k_{n-j-1}
$$
and
$$
c_{n-j+1}':\ c_{n-j}'\to k_{n-j}':\ c_{n-j-1}'\to k_{n-j-1}'.
$$
\vskip4ex
For the pair $(c_{n-j},c_{n-j}')$, on the lower-left one has
$$
\xymatrix@C=20ex{
(c_{n-j},c_{n-j}')
\ar@{|->}[r]^{[\eta^{(j)}]_0\times[\eta^{(j)}]_0}
&
(c_{n-j},c_{n-j}')}
$$
$$
\xymatrix@C=20ex{
= (c_{n-j},c_{n-j}')
\ar@{|->}[r]^{\circ^0}
&
c_{n-j}\circ^0 c_{n-j}'
},
$$
\vskip4ex
where on the upper-right
$$
\xymatrix@C=20ex{
(c_{n-j},c_{n-j}')
\ar@{|->}[r]^{\circ^0}
&
c_{n-j}\circ^0  c_{n-j}'
}
$$
$$
\xymatrix@C=20ex{
= c_{n-j}\circ^0  c_{n-j}'
\ar@{|->}[r]^{[\eta^{(j)}]_0}
&
c_{n-j}\circ^0  c_{n-j}'
}
$$

For homs, let us say that
$$
\Big[\mathrm{Diagram}\ (\ref{zzz5})\ \mathrm{for}\ j\Big]_1\Big((c_{n-j},c_{n-j}'),(k_{n-j},k_{n-j}')\Big)
$$
is indeed $\Big[\mathrm{Diagram}\ (\ref{zzz5})\ \mathrm{for}\ (j-1)\Big]$. Induction completes the proof.
\end{proof}

In the same way one can prove the following
\begin{Lemma}
Given $$\xymatrix@C=4ex{c_{n-j-1}:c_0\ar@{==>}[r]&c_0},$$ the following diagram commutes in
$j$-\textbf{Gpd}:
\begin{equation}\label{zzz6}
\xymatrix@C=20ex{
\Id{j}
\ar[r]^(.4){{}^{\C}u^0(c_{0})}
\ar@{=}[d]
&\C_{n-j}(c_{n-j-1},c_{n-j-1})
\ar[d]^{\eta^j_{\C_{n-j}(c_{n-j-1},c_{n-j-1})}}
\\
D^{(j)}\pi_0^{(j)}\left(\Id{j}\right)
\ar[r]_(.4){D^{(j)}\pi_0^{(j)}\left({}^{\C}u^0(c_0)\right)}
&D^{(j)}\pi_0^{(j)}\left(\C_{n-j}(c_{n-j-1},c_{n-j-1})\right)
}
\end{equation}
\end{Lemma}

Finally we can conclude the
\begin{proof}(of {\em Claim} \ref{zz:funct}).
The first of the diagrams of {\em Claim \ref{zz:funct}}\/ amounts to diagram (\ref{zzz5}) with $j=(n-1)$, the second amounts to
diagram (\ref{zzz6}) again with $j=(n-1)$.
\end{proof}

\begin{Proposition}
$$
\eta^{(n)}:id_{n\mathbf{Gpd}}\Rightarrow  D^{(n)}\big(\pi_0^{(n)}(-)\big)
$$
is a natural transformation of sesqui-functors.
\end{Proposition}
\begin{proof}
It suffices to show that, for any $\alpha:F\Rightarrow G:\C\to\D$,
$$
\alpha \bullet^0 \eta^{(n)}_{\D} = \eta^{(n)}_{\C}\bullet^0 D^{(n)} (\pi_0^{(n)}(\alpha))
$$
$$
\xymatrix@C=20ex@R=16ex{
\C\ar[r]^{\eta_{\C}}
\ar@/_6ex/[d]_{F}="a1"
\ar@/^6ex/[d]^{G}="a2"
&
D(\pi_0(\C))
\ar@/_6ex/[d]_{D(\pi_0(F))}="b1"
\ar@/^6ex/[d]^{D(\pi_0(G))}="b2"
\\
\D\ar[r]_{\eta_{\D}}
&
D(\pi_0(\D))
\ar@{}"a1";"a2"|(.4){}="ax1"|(.6){}="ax2"
\ar@{}"b1";"b2"|(.4){}="bx1"|(.6){}="bx2"
\ar@{=>}"ax1";"ax2"^{\alpha}
\ar@{=>}"bx1";"bx2"^{D(\pi_0(\alpha))}
}
$$
By induction on $n$.\\

\framebox[1.1\width]{$n=1$}\\

The adjunction of underlying categories and functors is well known. It extends plainly to sesqui-categories:
in fact, in \textbf{Gpd}, $\alpha(c_0)$'s are isomorphisms. Hence, $D'(\pi_0'(\alpha))$ is an equality of functors.\\

\framebox[1.1\width]{$n>1$}\\

On objects, let us fix $c_0$ in $\C_0$. Then
$$
[\alpha\bullet^0 \eta^{(n)}_{\D}]_0(c_0)= \eta^{(n)}_{\D}(\alpha(c_0))=
[\eta^{(n)}_{\D_1(Fc_0,Gc_0)}]_0(\alpha (c_0))=\alpha (c_0)
$$
(or $[\alpha(c_0)]_{\sim}$ if $n=2$). On the other side,
$$
[\eta^{(n)}_{\C}\bullet^0  D^{(n)}(\pi_0^{(n)}\alpha)]_0(c_0)=[D^{(n)}(\pi_0^{(n)}\alpha)]_0([\eta^{(n)}_{\C}]_0(c_0))=
[D^{(n)}(\pi_0^{(n)}\alpha)]_0(c_0)=[\pi_0^{(n)}\alpha]_0(c_0)=\alpha(c_0)
$$
(or $[\alpha(c_0)]_{\sim}$ if $n=2$).

On homs, let us fix one more object $c_0'$. Then
$$
[\alpha\bullet^0 \eta^{(n)}_{\D}]_1^{c_0,c_0'}=\alpha_1^{c_0,c_0'}\bullet^0 [\eta^{(n)}_{\D}]_1^{Fc_0,Gc_0'}=\alpha_1^{c_0,c_0'}\bullet^0 \eta^{(n-1)}_{\D_1(Fc_0,Gc_0')}
$$
on the other side,
$$
[\eta_{\C}^{(n)}\bullet^0  D^{n}(\pi_0^{n}\alpha)]_1^{c_0,c_0'}=[\eta^{(n)_{\C}}]_1^{c_0,c_0'}\bullet^0  [D^{(n)}(\pi_0^{(n)}\alpha)]_1^{c_0,c_0'}
=\eta^{(n-1)}_{\C_1(c_0,c_0')}\bullet^0  D^{(n-1)}(\pi_0^{(n-1)}(\alpha_1^{c_0,c_0'}))
$$
Equality of the two sides
$$
\alpha_1^{c_0,c_0'}\bullet^0 \eta^{(n-1)}_{\D_1(Fc_0,Gc_0')} = \eta^{(n-1)}_{\C_1(c_0,c_0')}\bullet^0  D^{(n-1)}(\pi_0^{(n-1)}(\alpha_1^{c_0,c_0'}))
$$
is given by induction hypothesis.
\end{proof}

\begin{Theorem}
For every positive integer $n$,
$$
\pi_0^{(n)}\dashv D^{(n)}
$$
\end{Theorem}
\begin{proof}
After the discussion above, triangular identities will be proved in the following form:
$$
\xymatrix{
D^{(n)}\ar@{=>}[r]^(.4){\eta^{(n)}D^{(n)}}
\ar@{=}[dr]
& D^{(n)}\pi_0^{(n)}D^{(n)}
\ar@{=}[d]^{D^{(n)}\epsilon^{(n)}=D^{(n)}}
\\\ar@{}[ur]|(.65){(i)}&
D^{(n)}
}
\qquad
\xymatrix{
\pi_0^{(n)}\ar@{=>}[r]^(.4){\pi_0^{(n)}\eta^{(n)}}
\ar@{=}[dr]
& \pi_0^{(n)}D^{(n)}\pi_0^{(n)}
\ar@{=}[d]^{\epsilon^{(n)}\pi_0^{(n)}=\pi_0^{(n)}}
\\\ar@{}[ur]|(.65){(ii)}&
\pi_0^{(n)}
}$$
Diagram $(i)$ commutes. In fact, for a (n-1)groupoid $\ES$ one has $\eta^{(n)}_{D^{(n)}\ES}=id_{D^{(n)}\ES}$.

By induction on $n$. For $n=1$, $\ES=S$ is a set and $\eta_{D'S}'$ is given by
$$
[\eta_{D'S}']_0 : s_0 \mapsto [s_0]_{\sim}=s_0
$$
Since $D'S$ is a discrete category, it only has identity arrows, hence
$[\eta_{D'S}']_0=id$ by functoriality.

For $n>1$, by definition one has
$$
[\eta_{D^{(n)}\ES}]_0 = id_{[D^{(n)}\ES]_0}
$$
Moreover, for any pair of objects $s_0,s_0'$ of $\ES$,
\begin{eqnarray*}
[\eta_{D^{(n)}\ES}]_1^{s_0,s_0'}   &=& \eta^{(n-1)}_{[D^{(n)}\ES]_1^{s_0,s_0'}} \\
   &=&  \eta^{(n-1)}_{D^{(n-1)}(\ES_1(s_0,s_0'))}\\
   &\eq{\clubsuit}&  id_{D^{(n-1)}(\ES_1(s_0,s_0'))}\\
   &=&  id_{[D^{(n)}\ES]_1(s_0,s_0')}\\
   &=&  [id_{D^{(n)}\ES}]_1(s_0,s_0')
\end{eqnarray*}
where all equalities hold by definition, but $(\clubsuit)$ that is given by induction hypothesis.\\

Diagram $(ii)$ commutes. In fact, for a n-groupoid $\C$ one has $\pi_0^{(n)}(\eta^{(n)}_{\C})=id_{\pi_0^{(n)}\C}$.

By induction on $n$. For $n=1$, {\em i.e.} for a groupoid $\C$, $\pi_0'(\eta_{\C}')=[\pi_0'(\eta_{\C}')]_0$, and it is given by
$$
[\pi_0'(\eta_{\C}')]_0:[c_0]_{\sim}\mapsto[\pi_0']_0([\eta_{\C}']_0([c_0]_{\sim}))=\pi_0'([c_0]_{\sim})=[c_0]_{\sim}
$$

For $n>1$, by definition one has
$$
[\pi_0^{(n)}(\eta^{(n)}_{\C})]_0=id_{[\pi_0^{(n)}\C]_0}
$$

Moreover, for any pair of objects $c_0,c_0'$ of $\C$,

\begin{eqnarray*}
[\pi_0^{(n)}(\eta^{(n)}_{\C})]_1^{c_0,c_0'}   &=& \pi_0^{(n-1)}([\eta_{\C}^{(n)}]_1^{c_0,c_0'}) \\
   &=&  \pi_0^{(n-1)}(\eta^{(n)}_{\C_1(c_0,c_0')})\\
   &\eq{\clubsuit}& id_{\pi_0^{(n-1)}(\C_1(c_0,c_0'))}  \\
   &=&  id_{[\pi_0^{(n)}\C]_1(c_0,c_0')}
\end{eqnarray*}
where all equalities hold by definition, but $(\clubsuit)$ that is given by induction hypothesis.\\
\end{proof}

\section{$n$-Discrete $h$-pullbacks}
An application of the adjunction $\pi_0^{(n)} \dashv D^{(n)}$ is the following useful result, just
a special case of more general $h$-limits preservation property:
\begin{Lemma}\label{lm:discrete_h-pb}
Sesqui-functor $D^{(n)}:\ (n-1)\mathbf{Gpd}\to n\mathbf{Gpd}$ preserves $h$-pullbacks.
\end{Lemma}
\begin{proof}
We omit the superscripts being always $(n)$.
Let us consider $D$ of  the $h$-pullback $(\PP,P,Q,\phi)$ in $(n-1)\mathbf{Gpd}$
$$
\raisebox{5ex}{\xymatrix@C=12ex{
\PP\ar[r]^Q\ar[d]_P
&\C\ar[d]^G
\\
\A\ar[r]_F
\ar@{}[ur]|(.3){}="1"|(.7){}="2"\ar@2"1";"2"^{\phi}
&\B
}}\quad\xymatrix@C=8ex{\ar@{|->}[r]^{D}&{} }\quad
\raisebox{5ex}{\xymatrix@C=8ex{
D\PP\ar[r]^{DQ}\ar[d]_{DP}
&D\C\ar[d]^{DG}
\\
D\A\ar[r]_{DF}
\ar@{}[ur]|(.3){}="1"|(.7){}="2"\ar@2"1";"2"^{D\phi}
&D\B
}}$$
For the four-tuple $(\Q,M,N,\omega)$ over the base $<DF,DG>$, we can now apply
$\pi_0$
$$
\raisebox{5ex}{\xymatrix@C=8ex{
\Q\ar[r]^{N}\ar[d]_{M}
&D\C\ar[d]^{DG}
\\
D\A\ar[r]_{DF}
\ar@{}[ur]|(.3){}="1"|(.7){}="2"\ar@2"1";"2"^{\omega}
&D\B
}}
\quad\xymatrix@C=8ex{\ar@{|->}[r]^{\pi_0}&}\quad
\raisebox{5ex}{\xymatrix@C=10ex{
\pi_0\Q\ar[r]^{\pi_0N}\ar[d]_{\pi_0M}
&\C\ar[d]^{G}
\\
\A\ar[r]_{F}
\ar@{}[ur]|(.3){}="1"|(.7){}="2"\ar@2"1";"2"^{\pi_0\omega}
&\B
}}
$$
universal property for $\PP$ gives then a unique $L:\pi_0\Q\to\PP$ such that
$$
(i)\, L\bullet^0 P =\pi_0 M\quad (ii)\, L\bullet^0 Q =\pi_0 N\quad (iii)\, L\bullet^0 \phi  =\pi_0 \omega
$$
Hence the composition $\eta_{\Q}\bullet^0 DL$ witnesses the fact that $D(\PP)$ is an $h$-pullback.
In fact
\begin{eqnarray*}
   &(i)'&  \eta_{\Q}^{\phantom{\Q}}\bullet^0 DL\bullet^0 DP =    \eta_{\Q}^{\phantom{\Q}}\bullet^0 D(L\bullet^0P) \eq{i}   \eta_{\Q}^{\phantom{\Q}}\bullet^0 \pi_0M \eq{a} M\bullet^0 \eta_{D\A}^{\phantom{\Q}} \eq{b}M  \\
   &(ii)'& \eta_{\Q}^{\phantom{\Q}}\bullet^0 DL\bullet^0 DQ =    \eta_{\Q}^{\phantom{\Q}}\bullet^0 D(L\bullet^0Q) \eq{ii}  \eta_{\Q}^{\phantom{\Q}}\bullet^0 \pi_0N \eq{a} N\bullet^0 \eta_{D\C}^{\phantom{\Q}} \eq{b}N   \\
   &(iii)'&\eta_{\Q}^{\phantom{\Q}}\bullet^0 DL\bullet^0 D\phi = \eta_{\Q}^{\phantom{\Q}}\bullet^0 D(L\bullet^0\phi) \eq{iii}  \eta_{\Q}^{\phantom{\Q}}\bullet^0 \pi_0\omega \eq{a} \omega\bullet^0 \eta_{D\C}^{\phantom{\Q}} \eq{b}\omega
\end{eqnarray*}
where $(a)$ by naturality of $\eta$ and $(b)$ by triangular identities.
\end{proof}
Hence we can say that the $h$-pullback of a $n$-discrete diagram is itself $n$-discrete.

\section{Exact sequences of n-groupoids}

\subsection{Pointedness and $h$-fibers}

A pointed $n$-category  is simply a $n$-category  $\C$ with a chosen object $*_{\C}$.

A morphism of pointed $n$-categories $F:\C\to\D$ is a $n$-functor such that
$$F(*_{\C})=*_{\D}$$

A 2-morphism of pointed $n$-functors $\alpha:F\Rightarrow G:\C\to\D$ is a natural $n$-transformation
such that $\alpha_0(*_{\C})=id_{*_{\D}}:*_{\D}\to*_{\D}$.

The data described above form a sesqui-category, sub-sesqui-category of $n\mathbf{Cat}$ that we will denote
$n\mathbf{Cat}_*$. Similarly one defines the sesqui-category $n\mathbf{Gpd}_*$ of pointed n-groupoids,
sub-sesqui-category of $n\mathbf{Gpd}$.

Subscripts of the \emph{star} will be often omitted.

Notice that definitions of pointed morphism and of pointed 2-morphisms imply that $n\mathbf{Cat}_*$ and $\mathbf{nGpd}_*$
are closed under finite products and $h$-pullbacks.

\begin{Definition}
Given a morphism of n-groupoids $F:\C\to\D$, and an object $d$ of $\D$, the \emph{past} \textit{h}-fiber $\F^{(p)}$
and the \emph{future} \textit{h}-fiber $\F^{(f)}$ of $F$ over $d$ are given by the following \textit{h}-pullbacks resp.
$$
\xymatrix{
\F^{(p)}_{F,d}\ar[r]^{!}\ar[d]_E
&\I\ar[d]^{[d]}
\\
\C\ar[r]_F\ar@{}[ur]|(.3){}="1"|(.7){}="2"
&\D
\ar@{=>}"2";"1"^{\varepsilon}
}
\qquad
\xymatrix{
\F^{(f)}_{F,d}\ar[r]^E\ar[d]_{!}
&\C\ar[d]^{F}
\\
\I\ar[r]_{[d]}\ar@{}[ur]|(.3){}="1"|(.7){}="2"
&\D
\ar@{=>}"2";"1"^{\varepsilon}
}
$$
where $[d]$ is the constant functor over an object $d$, as usual.
\end{Definition}
\begin{Remark}
Distinction between future and past \textit{h}-fibers, makes more sense in $n$-categorical context. There, in fact, the $h$-fiber
may come in different tastes: the lax-version that uses lax-2-morphisms, the pseudo-version, that uses equivalence transformations.
Nevertheless, even for n-groupoids, keeping track of the direction of 2-morphisms is necessary
in order to recognize the problem of coherently choosing inverses of cells, when this problem arises.
Still, we will often omit superscripts (and subscripts), as they will be clear from diagrams.
\end{Remark}
\begin{Remark}
The notion of strict fiber, or simply fiber, is recovered by the strict pullback (see Section \ref{def:strict_pb}).
\end{Remark}

\subsection{Equivalences and $h$-surjective morphisms of n-groupoids}
The notions of \emph{\textit{h}-surjective} morphism and of \emph{equivalence} get simpler in the context of n-groupoid.
Concerning both cases, it is the notion of essential surjectivity  itself to be involved. In fact,
for a n-groupoids morphism $L:\A\to\K$, to be essentially surjective on objects amounts to the following property:
\begin{description}
  \item[Essential surjectivity (1)] for  any  object $k_0$ of $\K$, there exists a pair
  $$
(a_0,k_1:La_0\to k_0)
  $$
  with $a_0$ object of $\A$ and $k_1$ 1-cell of $\K$.
\end{description}
That is: it is no longer necessary to ask for $k_1$ to be an equivalence, since every cell in a n-groupoid
is indeed an equivalence.\\

With the notion of $h$-fiber in mind, we can further reformulate the notion of essential surjectivity.
\begin{description}
\item[Essential surjectivity (2)] for  any  object $k_0$ of $\K$, the $h$-fiber $$\F^{(p)}_{L,k_0}$$ is not empty.
\end{description}

Finally, in pointed case, fibers assume a special meaning, as referred by the following
\begin{Definition}
Let $G:\B\to\C$ be a morphism of n-groupoids. Then the fiber $\K=\F_{G,*}^{(p)}$
$$
\xymatrix{
\K\ar[r]_K \ar@/^6ex/[rr]^0_{}="1"&{\B}\ar@{}|{}="2"\ar[r]_G&\C
\ar@{=>}"1";"2"^{\kappa}}
$$
is called (past) the $h$-kernel of $G$.
\end{Definition}

We will call universal property of $h$-kernels, the universal property of pullbacks specialized for these kind of pullbacks.

\subsection{Exact sequences}

\begin{Definition}\label{def:exactness}
Let the following diagram in $\mathbf{nGpd}_*$ be given:
$$
\xymatrix{
\A\ar[r]_F \ar@/^6ex/[rr]^0_{}="1"&{\B}\ar@{}|{}="2"\ar[r]_G&\C
\ar@{=>}"1";"2"^{\varepsilon}}
$$
We call the triple $(F,\varepsilon,G)$ \emph{exact in} $\B$ if the comparison n-functor $L:\A\to\K$, given by the universal property
of the h-kernel $(\K,K,\kappa)$, is $h$-surjective.
$$
\xymatrix{
\A\ar[dd]_L="x"\ar[dr]^F\ar@/^4ex/[drr]^0_{}="a1"
\\
&\B\ar@{}="2"\ar[r]|G&\C
\\
\K\ar[ur]_K\ar@/_4ex/[urr]_0^{}="b1"
\ar@{=>}"a1";"2"^{\varepsilon}
\ar@{=>}"b1";"2"_{\kappa}
\ar@{}"x";"2"|{\equiv}
}
$$
\end{Definition}

This notion of exactness is a straightforward extensions of the notion introduced by
Vitale in \cite{Vitale2002jpaa}.

Of course it reduces to usual exactness for pointed
sets and group. Moreover it is preserved by one-point suspension and by discretization,
hence an exact sequence of groups may be considered as an exact sequence of one-point groupoids, as well
as an exact sequence of pointed discrete groupoids (with a group structure).\\

In the categorical group (pointed groupoid) situation, it has shown its usefulness
in extending  homological algebraic structures in a 1-dimensional context.\\

\subsection{$\pi_0$ preserves exactness}

In the following paragraphs we will show that, given a three-term exact sequence in $n$\textbf{Gpd}, the sesqui-functor
$\pi_0^{(n)}$ produces a three-term exact sequence in $(n-1)$\textbf{Gpd}. Preliminary Lemmas clarify
the relations between preservation of exactness and its main ingredients: $h$-surjectivity and the notion of $h$-pullback.

\begin{Lemma}\label{lemma:comparison_surjective}
Let us consider the following $h$-pullback diagram:
$$
\xymatrix{
\PP
\ar[r]^S
\ar[d]_R
&\Z
\ar[d]^H
\\
\B\ar@{}[ur]|(.3){}="1"|(.7){}="2"\ar@{=>}"1";"2"^{\varepsilon}
\ar[r]_G
&
\C
}
$$
The comparison $L:\pi_0^{(n)}\PP\to\Q$ with the $h$-pullback of $\pi_0^{(n)}(G)$ and $\pi_0^{(n)}(H)$
is $h$-surjective.
$$
\xymatrix{
\pi_0^{(n)}\PP
\ar[dr]^(.7){L}
\ar@/^4ex/[drr]^{\pi_0^{(n)}S}_(.35){}="a2"
\ar@/_4ex/[ddr]_{\pi_0^{(n)}R}^(.35){}="a1"
\\
&\Q
\ar[r]^Q
\ar[d]_P
&\pi_0^{(n)}\Z
\ar[d]^{\pi_0^{(n)}H}
\\
&\pi_0^{(n)}\B\ar@{}[ur]|(.3){}="1"|(.7){}="2"\ar@{=>}"1";"2"^{\gamma}
\ar[r]_{\pi_0^{(n)}G}
&
\pi_0^{(n)}\C
\ar@{}"a1";"a2"|(.3){}="b1"|(.7){}="b2"|(.2){\quad\pi_0^{(n)}\varepsilon}\ar@{:>}"b1";"b2"
}
$$
\end{Lemma}
\begin{proof}
By induction on $n$.\\

\framebox[1.1\width]{$n=1$}\\

The $h$-pullback $\PP$ has objects
$$
(b_0,\xymatrix{G b_0\ar[r]^{c_1}& H z_0},z_0)
$$
and arrows
$$
(b_1,=,z_1):(b_0,c_1,z_0)\to (b_0',c_1',z_0')
$$
where the ``$=$'' stays for the commutative square
$$
\xymatrix{
Gb_0\ar[r]^{c_1\phantom{'}}
\ar[d]_{Gb_1}
&Hz_0\ar[d]^{Hz_1}
\\
Gb_0'\ar[r]_{c_1'}
&Hz_0'
}
$$
Hence the set $\pi_0'(\PP)$ has elements the classes
$$
[b_0,c_1,z_0]_{\sim}.
$$
On the other side, the set $\Q$ is a usual pullback in \textbf{Set}. It has elements the
pairs $([b_0]_{\sim},[z_0]_{\sim})$ such that $\pi_0'G([b_0]_{\sim})=\pi_0'H([z_0]_{\sim})$,
{\em i.e.} $[Gb_0]_{\sim}=[Hz_0]_{\sim}$, {\em i.e.} such that there exists $c_1:Gb_0\to Hz_0$.
Then the comparison
$$
L=L_0:[b_0,c_1,z_0]_{\sim}\mapsto ([b_0]_{\sim},[z_0]_{\sim})
$$
is clearly surjective.\\

\framebox[1.1\width]{$n=2$}\\

Now the $h$-pullback $\PP$ is a 2-groupoid with objects
$$
(b_0,\xymatrix{G b_0\ar[r]^{c_1}& H z_0},z_0).
$$
Arrows are of the form
$$
(b_1,c_2,z_1):(b_0,c_1,z_0)\to (b_0',c_1',z_0')
$$
{\em i.e.}
$$
\left(
\raisebox{5ex}{\xymatrix{b_0\ar[d]_{b_1}\\b_0'}},
\raisebox{5ex}{\xymatrix{
Gb_0\ar[r]^{c_1\phantom{'}}
\ar[d]_{Gb_1}
&Hz_0\ar[d]^{Hz_1}
\\
Gb_0'\ar[r]_{c_1'}
\ar@{}[ur]|(.3){}="1"|(.7){}="2"\ar@{=>}"1";"2"^{c_2}
&Hz_0'
}},
\raisebox{5ex}{\xymatrix{z_0\ar[d]_{z_1}\\z_0'}}
\right)
$$
Finally 2-cells are of the form
$$
(b_2,=,z_2):(b_1,c_2,z_1)\Rightarrow (b_1',c_2',z_1')
$$
{\em i.e.}
$$
\left(
\raisebox{5ex}{\xymatrix{b_1\ar@2[d]_{b_2}\\b_1'}},
\raisebox{5ex}{\xymatrix{
Gb_1\circ c_1'\ar@2[r]^{c_2\phantom{'}}
\ar@2[d]_{Gb_2\circ c_1'}
&c_1\circ Hz_1\ar@2[d]^{c_1\circ Hz_2}
\\
Gb_1'\circ c_1'\ar@2[r]_{c_2'}
&c_1\circ Hz_1'
}},
\raisebox{5ex}{\xymatrix{z_1\ar@2[d]^{z_2}\\z_1'}}
\right)
$$
Therefore the groupoid $\pi_0''\PP$ has objects $(b_0,c_1,z_0)$ and arrows $[b_1,c_2,z_1]_{\sim}$.\\

On the other side, the groupoid $\Q$ has objects
$$
(b_0,\xymatrix{Gb_0\ar[r]^{[c_1]_{\sim}}&Hz_0},z_0)
$$
and arrows
$$
([b_1]_{\sim},=,[z_1]_{\sim})
$$
with $[b_1]_{\sim}:b_0\to b_0'$ in $\pi_0''\B$ and $[z_1]_{\sim}:z_0\to z_0'$ in $\pi_0''\Z$ such that
$$
\xymatrix{
Gb_0\ar[r]^{[c_1]_{\sim}}
\ar[d]_{[Gb_1]_{\sim}}
&Hz_0
\ar[d]^{[Hz_1]_{\sim}}
\\\ar@{}[ur]|{(\heartsuit)}
Gb_0'\ar[r]_{[c_1']_{\sim}}
&Hz_0'}
$$
Hence the comparison
\begin{eqnarray*}
  L &:& (b_0,c_1,z_0)\mapsto (b_0,c_1,z_0) \\
    & &  [b_1,c_2,z_1]_{\sim} \mapsto ([b_1]_{\sim},[z_1]_{\sim})
\end{eqnarray*}
is $h$-surjective. In fact it is an identity (hence strictly surjective) on objects, and full on homs.
Let us fix a  pair of objects $(b_0,c_1,z_0)$ and $(b_0',c_1',z_0')$ in the domain, and an arrow
$([b_1]_{\sim},=,[z_1]_{\sim})$ in $\Q$, where the ``$=$'' is the diagram $(\heartsuit)$ above.
Then $[c_1\circ Hz_1]_{\sim}=[Gb_1\circ c_1']_{\sim}$ if, and only if, there exists
$$
c_2:c_1\circ Hz_1 \to Gb_1\circ c_1'.
$$
In other words we get an arrow $[b_1,c_2,z_1]_{\sim}$ of $\pi_0''\PP$ that $L$ sends in
$([b_1]_{\sim},=,[z_1]_{\sim})$, {\em i.e.} $L$ is full.

\framebox[1.1\width]{$n>2$}\\

On objects, $L_0:[\pi_0^{(n)}\PP]_0\to \Q_0$ is the identity.

In fact, for big $n$   $[\pi_0^{(n)}\PP]_0=\PP_0$,
the last being the set-theoretical limit over the diagram
$$
\xymatrix{
\B_0\ar[dr]_{G_0}&&\C_1\ar[dl]^d \ar[dr]_c&&\Z_0\ar[dl]^{H_0}
\\
&\C_0&&\C_0
}
$$
Now, for $n>2$ this diagram coincides with the one defining $Q_0$:
$$
\xymatrix@C=2ex{
[\pi_0^{(n)}\B]_0\ar[dr]_{[\pi_0^{(n)}G]_0}&&[\pi_0^{(n)}\C]_1\ar[dl]^d \ar[dr]_c&&[\pi_0^{(n)}\Z]_0\ar[dl]^{[\pi_0^{(n)}H]_0}
\\
&[\pi_0^{(n)}\C]_0&&[\pi_0^{(n)}\C]_0
}
$$
Universality of limits gives $L_0=id$.

On homs, let us fix two objects $p_1=(b_0,c_1,z_0)$ and $p_1'=(b_0',c_1',z_0')$ of $[\pi_0^{(n)}\PP]_0=\PP_0$ and
compute $L_1^{p_1,p_1'}$ as usual: by means of universal property of $h$-pullbacks. In fact  the diagram
\begin{changemargin}{-10ex}{-10ex}
$$
\xymatrix{
[\pi_0^{(n)}\PP]_1(p_0,p_0')\ar[dr]^(.65){L_1^{p_0,p_0'}}
\ar@/^3.5ex/[drrr]^{[\pi_0^{(n)}S]_1}_(.3){}="1"
\ar@/_6ex/[dddr]_{[\pi_0^{(n)}R]_1}^(.2){}="2"
\\
&\Q_1(p_0,p_0')
\ar[rr]^{Q_1}  \ar[dd]_{P_1}
&& [\pi_0^{(n)}\Z]_1(z_0,z_0')
\ar[d]^{[\pi_0^{(n)}H]_1}
\ar@{}[ddll]|(.4){}="a1"|(.6){}="a2"
\\
&&&[\pi_0^{(n)}\C]_1(Hz_0,Hz_0')
\ar[d]^{c_1\circ-}
\\
&[\pi_0^{(n)}\B]_1(b_0,b_0')\ar[r]_{[\pi_0^{(n)}G]_1}
&[\pi_0^{(n)}\C]_1(Gb_0,Gb_0')\ar[r]_{-\circ c_1'}
&[\pi_0^{(n)}\C]_1(Gb_0,Hz_0')
\ar@{=>}"a1";"a2"^{\sigma}
\ar@{}"1";"2"|(.4){}="b1"|(.85){}="b2"^(.7){[\pi_0^{(n)}\varepsilon]_1}
\ar@{:>}"b1";"b2"^{}}
$$
\end{changemargin}
is the same as (and determined by)
\begin{changemargin}{-10ex}{-10ex}
$$
\xymatrix{
\pi_0^{(n-1)}(\PP_1(p_0,p_0'))\ar[dr]^(.65){L_1^{p_0,p_0'}}
\ar@/^3.5ex/[drrr]^{\pi_0^{(n-1)}S_1}_(.3){}="1"
\ar@/_6ex/[dddr]_{\pi_0^{(n-1)}R_1}^(.2){}="2"
\\
&\Q_1(p_0,p_0')
\ar[rr]^{Q_1}  \ar[dd]_{P_1}
&& [\pi_0^{(n)}\Z]_1(z_0,z_0')
\ar[d]^{\pi_0^{(n-1)}H_1}
\ar@{}[ddll]|(.4){}="a1"|(.6){}="a2"
\\
&&&\pi_0^{(n-1)}(\C_1(Hz_0,Hz_0'))
\ar[d]^{c_1\circ-}
\\
&\pi_0^{(n-1)}(\B_1(b_0,b_0'))\ar[r]_{\pi_0^{(n-1)}G_1}
&\pi_0^{(n-1)}(\C_1(Gb_0,Gb_0'))\ar[r]_{-\circ c_1'}
&\pi_0^{(n-1)}(\C_1(Gb_0,Hz_0'))
\ar@{=>}"a1";"a2"^{\sigma}
\ar@{}"1";"2"|(.4){}="b1"|(.85){}="b2"^(.7){\pi_0^{(n-1)}(\varepsilon_1)}
\ar@{:>}"b1";"b2"^{}}
$$
\end{changemargin}
This shows that $L_1^{p_0,p_0'}$ is itself a comparison between $\pi_0$ of an $h$-pullback and a $h$-pullback
of a $\pi_0$ of a diagram (of (n-1)groupoids), hence it is $h$-surjective by induction hypothesis.

In conclusion we have shown that $L=<L_0,L_1^{-,-}>$ is $h$-surjective.
\end{proof}

\begin{Lemma}\label{lemma:pi_0_surjective}
If the $n$-functor $L:\A\to\K$ is $h$-surjective, then also $\pi_0^{(n)}(L)$ is $h$-surjective.
\end{Lemma}
\begin{proof}
By induction on $n$.\\

\framebox[1.1\width]{$n=1$}\\

Let $L$ be a $h$-surjective functor between groupoids, {\em i.e.} $L$ is full and essentially surjective on objects.
Therefore for an element $[k_0]_{\sim}\in \pi_0'\K$ there exists a pair
$$
(a_0,k_1:La_0\to k_0)
$$
Hence
$$
(\pi_0'L)([a_0]_{\sim})=[La_0]_{\sim}=[k_0]_{\sim}\qquad \mathrm{in}\ \pi_0'\K.
$$

\framebox[1.1\width]{$n=2$}\\

Let $L$ be a $h$-surjective morphism between $2$-groupoids, {\em i.e.}
\begin{enumerate}
  \item for any $k_0$ there exist $(a_0,k_1:La_0\to k_0)$.
  \item for any pair $a_0,a_0'$,
  $$
  L_1^{a_0,a_0'}:\A_1(a_0,a_0')\to\K_1(La_0,La_0')
  $$
  is $h$-surjective.
\end{enumerate}
Since $[\pi_0''L]_0=L_0$, for any $k_0$ one has  $[k_1]_{\sim}:La_0\to k_0$, and this proves the first condition.

Moreover, once we fix a pair $a_0,a_0'$, by definition one has $[\pi_0''L]_1^{a_0,a_0'}=\pi_0'(L_1^{a_0,a_0'})$.
Hence it is $h$-surjective by previous case.\\

\framebox[1.1\width]{$n>2$}\\

More generally, a morphism $L$ of $n$-groupoids is $h$-surjective when conditions 1. and 2. above both hold.
Since $[\pi_0^{(n)}L]_0=L_0$ and $[\pi_0^{(n)}(\K)]_1=\pi_0^{(n-1)}(\K_1)$, condition 1. for $L$ and for $\pi_0^{(n)}L$
is the same, hence it holds.

Moreover, whence we fix a pair $a_0,a_0'$, by definition one has $[\pi_0^{(n)}L]_1^{a_0,a_0'}=\pi_0^{(n-1)}(L_1^{a_0,a_0'})$.
Hence it is $h$-surjective by induction hypothesis.
\end{proof}

Finally we are ready to state and prove the following important
\begin{Theorem}\label{thm:pi_0_preserves_exactness}
Given an exact sequence in $n\mathbf{Gpd}_*$
$$
\xymatrix@C=9ex{
\A\ar[r]_F \ar@/^6ex/[rr]^0_{}="1"&{\B}\ar@{}|{}="2"\ar[r]_G&\C
\ar@{=>}"1";"2"^{\varepsilon}}
$$
the sequence
$$
\xymatrix{
\pi_0^{(n)}\A\ar[r]_{\pi_0^{(n)}F} \ar@/^7ex/[rr]^0_{}="1"&{\pi_0^{(n)}\B}
\ar@{}|{}="2"\ar[r]_{\pi_0^{(n)}G}&\pi_0^{(n)}\C
\ar@{=>}"1";"2"^{\pi_0^{(n)}\varepsilon}}
$$
is exact in $(n-1)\mathbf{Gpd}_*$
\end{Theorem}
\begin{proof}
Let us consider the diagram
$$
\xymatrix@C=10ex{
\mathrm{Ker}\left(\pi_0^{(n)}G\right)\ar[dr]^(.8){}="y"
\ar@/^4ex/[drr]^0="x"
\\
\ar[u]^{L'}
\pi_0^{(n)}\big(\mathrm{Ker}G\big)\ar[r]
&\pi_0^{(n)}\B\ar[r]_{\pi_0^{(n)}G}\ar@{}"y";"x"|(.3){}="2"|(.7){}="1"
&\pi_0^{(n)}\C
\\
\ar[u]^{\pi_0^{(n)}L}
\pi_0^{(n)}\A\ar[ur]_{\pi_0^{(n)}F}
\ar@{=>}"1";"2"^{\kappa}
}
$$
$L$ is the comparison in $n$\textbf{Gpd}, hence $h$-surjective by hypothesis. Therefore
$\pi_0^{(n)}L$ is  $h$-surjective by {\em Lemma \ref{lemma:pi_0_surjective}}.

$L'$ is the comparison in (n-1)\textbf{Gpd}, $h$-surjective by {\em Lemma \ref{lemma:comparison_surjective}}.

Finally, their composition is again $h$-surjective by {\em Lemma \ref{lemma:h-surj_comp}}, and it is the comparison
between $\pi_0^{(n)}\A$ and the kernel of $\pi_0^{(n)}G$ by uniqueness in universal property of $h$-kernels.
\end{proof}

\section{The sesqui-functor $\pi_1$}

Purpose of this section is to introduce the family of sesqui-functors $\{\pi_1^{(n)}\}_{n\in\mathbb{N}}$
that extends the \emph{isos-of-the-point} functor $\mathbf{Gpd}_*\to\mathbf{Set}_*$ that
assigns to each pointed groupoid the pointed (hom-)set of endo-arrows
of the point.

\begin{Defprop}
For any integer $n>0$, there exists a  sesqui-functor
$$
\pi_1^{(n)}:n\mathbf{Gpd}_* \rightarrow (n-1)\mathbf{Gpd}_*
$$
contra-variant on 2-morphisms, according to the following recursive definition.
\end{Defprop}

\framebox[1.1\width]{$n=1$}\\
$\pi_1^{(1)}$is the functor (= trivial sesqui-functor) $\mathbf{Gpd}_*\to\mathbf{Set}_*$ that assigns to a pointed
groupoid $\C$ the pointed set $\C(*,*)$. It can be considered contra-variant,
since 2-morphisms in $\mathbf{Set}$ are equalities. \\

\framebox[1.1\width]{$n>1$}\\

Let a $n$-groupoid $\C$ be given. Then $\pi_1^{(n)}\C=\C_1(*,*)$.\\

Let a morphism of $n$-groupoids $F:\C\to\D$ be given. Then
$\pi_1^{(n)}F=F_1^{*,*}$.

Of course these assignments give indeed
a functor between underlying categories. In fact
$$
\pi_1^{(n)}(id_{\C})=[id_{\C}]_1^{*,*}=id_{\C_1(*,*)}=id_{\pi_1^{(n)}\C}
$$
and for every other $G:\D\to\E$, one has
$$
\pi_1^{(n)}(F\bullet^0 G)=[F\bullet^0 G]_1^{*,*}=F_1^{*,*}\bullet^0 G_1^{*,*}=\pi_1^{(n)}(F)\bullet^0 \pi_1^{(n)}(G).
$$

Let  a 2-morphism $\alpha:F\Rightarrow G:\C\to\D$ be given.
Then we define
$$
\pi_1^{(n)}(\alpha)=\alpha_1^{*,*}:\pi_1^{(n)}(G)\Rightarrow\pi_1^{(n)}(F).
$$
In fact, since compositions with identities gives identity functors,
 contra-variance is explained by the following diagram:

$$
{
\def \xxa{\C_1(*,*)}
\def \xxb{\D_1(*,*)}
\def \xxc{\D_1(*,*)}
\def \xxd{\D_1(*,*)}
\def \xxe{F_1^{*,*}}
\def \xxf{G_1^{*,*}}
\def \xxh{1_*\circ -}
\def \xxg{-\circ1_*}
\def \xxi{\alpha_1^{*,*}}
\xymatrix{
& \xxa \ar[dl]_{\xxe}\ar[dr]^{\xxf}\\
\xxb\ar@{}[rr]|(.4){}="xxb"|(.6){}="xxc"\ar[dr]_{\xxg}
&&\xxc\ar[dl]^{\xxh}\\
&\xxd
\ar@{=>}"xxc";"xxb"_{\xxi}}}
$$

In order to show that so-defined $\pi_1^{(n)}$ is indeed a
sesqui-functor, two facts have to be proved regarding $\pi_1^{(n)}$.

\begin{enumerate}
  \item it is functorial on hom-categories.
  \item it preserves reduced horizontal compositions.
\end{enumerate}

\begin{proof}
\begin{enumerate}
  \item Suppose we are given
  $$
\omega:E\Rightarrow F:\C\rightarrow \D\qquad\alpha:F\Rightarrow G:\C\rightarrow \D
  $$
in $n$\textbf{Gpd}.
Then
$$
\pi_1^{(n)}(\omega\bullet^1 \alpha)=[\omega\bullet^1 \alpha]_1^{*,*}
=\alpha_1^{*,*}\bullet^1 \omega_1^{*,*}=\pi_1^{(n)}(\alpha)\bullet^1 \pi_1^{(n)}(\omega)
$$
$$
\pi_1^{(n)}(id_F)=[id_F]_1^{*,*}=id_{F_1^{*,*}}=id_{\pi_1^{(n)}(F)}
$$
\item We prove the statement for reduced left-composition. Suppose 2-morphism $\alpha$ as above, and morphism
$N:\B\to\C$ be given. Then
$$
\pi_1^{(n)}(N\bullet^0_L\alpha)=[N\bullet^0_L\alpha]_1^{*,*}=N_1^{*,*}\bullet^0_L \alpha_1^{*,*}=\pi_1^{(n)}(N)\bullet^0_L \pi_1^{(n)}(\alpha).
$$
Finally, concerning reduced right-composition, the proof is similar, as the definition
\end{enumerate}
\end{proof}
The following proposition is a direct consequence of the definitions. Hence it needs no proof.
\begin{Proposition}
Sesqui-functor $\pi_1^{(n)}$ commutes with finite products and preserves equivalences.
\end{Proposition}
\begin{Remark}
The definition of $\pi_1^{(n)}$ given here makes it
of difficult use w.r.t. the inductive  setting developed so far, where  everything is given
as a pair, where the first component lives in \textbf{Set}, the second in $(n-1)$\textbf{Cat} (or $(n-1)$\textbf{Gpd}).
This motives a further search for a different (but equivalent) definition of
$\pi_1^{(n)}$, see {\em Corollary \ref{cor:eqv_pi_1}}.
Thereafter it will also be  shown that $\pi_1^{(n)}$ preserves exactness, as a consequence of universality
of its definition, although  this could be proved here directly.
\end{Remark}

\chapter{$3$-Morphisms of $n$-categories}\label{cha:three}

 \section{What structure for $n$\textbf{Cat}?}

So far we have shown that n-categories organizes naturally into a
sesqui-category, {\em ditto}  for $n$-groupoids. This gives a setting
to deal not only with $n$-categories and $n$-functors, but also with
their 2-morphisms, namely lax-n-transformations.

Yet the necessity of introducing 3-morphisms (lax-$n$-modifications)
takes us out of that comfortable setting, into the unknown territory of
sesqui-categorically enriched structures.

Following this suggestion, we have named the new setting {\em sesqui${}^2$-category}.
This notion  is closely related with that of \textbf{Gray}-category \cite{MR0412252,MR0371990}
(or 3D-Tas see \cite{MR1748753}) and incorporates a horizontal
dimension raising composition of 2-morphisms. In fact the set of axioms which
define the former is a subset of those defining the latter.

In order to fully justify the name chosen to denote such a structure, it
would be interesting to investigate explicitly the enrichment that generates this notion from that of
sesqui-category.

What we present here is a treatable inductive approach, comprehensive of a useful characterization
given in {\em Theorem \ref{thm:caratterizzo_sesqui2}}.

\begin{Definition}\label{def:sequi2cat}
A (small) {\em sesqui${}^2$-category} $\mathcal{C}$ consists of:\\

$\bullet$ A 3-truncated reflexive globular set $\mathcal{C}_{\bullet}$:
$$
\xymatrix{
\mathcal{C}_3
\ar@<+1.5ex>[r]^{d_2}
\ar@<-1.5ex>[r]_{c_2}
&
\ar[l]|{e_2}
\mathcal{C}_2
\ar@<+1.5ex>[r]^{d_1}
\ar@<-1.5ex>[r]_{c_1}
&
\ar[l]|{e_1}
\mathcal{C}_1
\ar@<+1.5ex>[r]^{d_0}
\ar@<-1.5ex>[r]_{c_0}
&
\ar[l]|{e_0}
\mathcal{C}_0
}
$$
with operations
$$
\bullet^m:\ \mathcal{C}_p\,\,{}_{c_m}\!\!\!\times_{d_m}\mathcal{C}_q\to\mathcal{C}_{p+q-m-1},\qquad m< \mathrm{min}(p,q)
$$
such that the following axioms hold:\\

$(i)$
For every pair $\C,\D\in\mathcal{C}_0$, the localization $\mathcal{C}(\C,\D)$ is a sesqui-category,
with\\

- object are $F,G,etc.\ \in \mathcal{C}_1(\C,\D)$\\
- for any pair of objects $F,G$, 1-cells are $\alpha,\beta,etc.\ \in\mathcal{C}_2(F,G)$\\
- for any pair of 1-cells $\alpha,\gamma:F\to G$, 2-cells are $\Lambda,\Sigma,etc.\ \in\mathcal{C}_3(\alpha,\beta)$\\

k-compositions are restrictions of $\bullet^{k+1}$-compositions:\\

- 0-composition of 1-cells of $\mathcal{C}(\C,\D)$
$$
\diamond^0:=\bullet^1:\mathcal{C}_2\,\,{}_{c_1}\!\!\!\times_{d_1}\mathcal{C}_2\to\mathcal{C}_2
$$

- left/right reduced 0-compositions of 1-cell with a 2-cell of $\mathcal{C}(\C,\D)$
$$
\diamond^0_L:=\bullet^1:\mathcal{C}_2\,\,{}_{c_1}\!\!\!\times_{d_1}\mathcal{C}_3\to\mathcal{C}_3
$$
$$
\diamond^0_R:=\bullet^1:\mathcal{C}_3\,\,{}_{c_1}\!\!\!\times_{d_1}\mathcal{C}_2\to\mathcal{C}_3
$$

- 1-compositions of 2-cells of $\mathcal{C}(\C,\D)$
$$
\diamond^1:=\bullet^2:\mathcal{C}_3\,\,{}_{c_2}\!\!\!\times_{d_2}\mathcal{C}_3\to\mathcal{C}_3
$$

$(ii)$
For every morphism $F:\C\to\D$ and objects $\B,\E$ of $\mathcal{C}$
$$
-\bullet^0 F:\ \mathcal{C}(\B,\C)\to\mathcal{C}(\B,\D)
$$
$$
F\bullet^0 -:\ \mathcal{C}(\D,\E)\to\mathcal{C}(\C,\E)
$$
are sesqui-functors.\\

\vskip4ex

$(iii)$
For every object $\C$ and objects $\B,\D$ of $\mathcal{C}$, if we denote $id_{\C}=e_0(\C)$,
$$
-\bullet^0 id_{\C}:\ \mathcal{C}(\B,\C)\to\mathcal{C}(\B,\C)
$$
$$
id_{\C}\bullet^0 -:\ \mathcal{C}(\C,\D)\to\mathcal{C}(\C,\D)
$$
are identity sesqui-functors.\\

\vskip4ex

$(iv)$  {\em (naturality axioms)}\\

For every pair of 0-composable 2-morphisms $\alpha:F\Rightarrow G:\C\to\D$ and $\beta:H\Rightarrow K:\D\to\E$\\

$(a)$
$$
\alpha\bullet^0\beta: (F\bullet^0\beta)\bullet^1(\alpha\bullet^0K)\to(\alpha\bullet_0H)\bullet^1(G\bullet^0\beta)
$$

\vskip3ex

For every 2-morphisms $\varepsilon:L\Rightarrow M:\B\to\C$ and $\beta:H\Rightarrow K:\D\to\E$,
and for every 3-morphism $\Lambda:\alpha\xymatrix@C=3ex{\ar@3[r]&}\omega:F\Rightarrow G:\C\to\D$\\

$(b)$
$$
(\alpha\bullet^0\beta)\bullet^2\Big((\Lambda\bullet^0H)\bullet^1(G\bullet^0\beta)\Big)
=\Big((F\bullet^0\beta)\bullet^1(\Lambda\bullet^0K)\Big)\bullet^2(\omega\bullet^0\beta)
$$

and\\

$(c)$
$$
\Big((L\bullet^0\Lambda)\bullet^1(\varepsilon\bullet^0G)\Big)\bullet^2(\varepsilon\bullet^0\omega)
=(\varepsilon\bullet^0\alpha)\bullet^2\Big((\varepsilon\bullet^0F)\bullet^1(M\bullet^0\Lambda)\Big)
$$

\vskip4ex

$(v)$  {\em (functoriality axioms)}\\

For every 2-morphisms $\omega:D\Rightarrow E:\B\to\C$ and $\gamma:H\Rightarrow L:\D\to\E$
and every pair of 1-composable 2-morphisms $\alpha:F\Rightarrow G:\C\to\D$ and $\beta:G\Rightarrow H:\C\to\D$\\

$(a)$
$$
(\alpha\bullet^1\beta)\bullet^0\gamma =
\Big((\alpha\bullet^0\gamma)\bullet^1(\beta\bullet^0L)\Big)
\bullet^2
\Big((\alpha\bullet^0K)\bullet^1(\beta\bullet^0\gamma)\Big)
$$

and\\

$(b)$
$$
\omega\bullet^0(\alpha\bullet^1\beta)=
\Big((\omega\bullet^0\alpha)\bullet^1(E\bullet^0\beta)\Big)
\bullet^2
\Big((D\bullet^0\alpha)\bullet^1(\omega\bullet^0\beta)\Big)
$$

\vskip4ex

$(vi)$  {\em (associativity axiom)}\\

For every 0-composable triple $x\in[\mathcal{C}(\B,\C)]_p$, $y\in[\mathcal{C}(\C,\D)]_q$
and $z\in[\mathcal{C}(\D,\E)]_r$, with \mbox{$p+q+r\leq 2$}
$$
(x\bullet^0y)\bullet^0z=x\bullet^0(y\bullet^0z)
$$

\vskip4ex

$(vii)$  {\em (identity axioms)}\\

For morphisms $E:\B\to\C$ and $H:\D\to\E$, and 2-morphism $\alpha:F\Rightarrow G:\C\to \D$,
$$
id_F\bullet^0\alpha=id_{F\bullet^0\alpha},\qquad \alpha\bullet^0 id_G= id_{\alpha\bullet^0 G}
$$
\end{Definition}

\begin{Remark}

1. Axiom $(iv)(c)$ is better understood when visualized as in the following diagram (same notation)
$$
\raisebox{0ex}{\xymatrix@C=14ex{
F\bullet^0H
\ar@2@/^10ex/[r]^{\omega\bullet^0H}_{}="a2"
\ar@2[r]|{\alpha\bullet^0H}^{}="a1"
\ar@{}"a1";"a2"|(.3){}="A1"|(.7){}="A2"
\ar@3"A1";"A2"_{\Lambda\bullet^0H}
\ar@2[d]_{F\bullet^0\beta}
&
G\bullet^0H
\ar@2[d]^{G\bullet^0\beta}
\\
\ar@{}[ur]|(.3){}="1"|(.7){}="2"
\ar@3"1";"2"^{\alpha\bullet^0\beta\ \ }
F\bullet^0K
\ar@2[r]_{\alpha\bullet^0K}
&
G\bullet^0K
}}
=
\raisebox{10ex}{\xymatrix@C=14ex{
F\bullet^0H
\ar@2[r]^{\omega\bullet^0H}
\ar@2[d]_{F\bullet^0\beta}
&
G\bullet^0H
\ar@2[d]^{G\bullet^0\beta}
\\
\ar@{}[ur]|(.3){}="1"|(.7){}="2"
\ar@3"1";"2"^{\omega\bullet^0\beta\ \ }
F\bullet^0K
\ar@2[r]|{\omega\bullet^0K}_{}="a2"
\ar@2@/_10ex/[r]_{\alpha\bullet^0K}^{}="a1"
&
G\bullet^0K
\ar@{}"a1";"a2"|(.3){}="A1"|(.7){}="A2"
\ar@3"A1";"A2"_{\Lambda\bullet^0K}
}}
$$
The same can be claimed for axiom $(iv)(b)$.\\

2. Axiom $(v)(a)$ is better understood when visualized as in the following diagram (same notation)
$$
\raisebox{5ex}{\xymatrix@C=14ex{
F\bullet^0K
\ar@2[r]^{(\alpha\bullet^1\beta)\bullet^0K}
\ar@2[d]_{F\bullet^0\gamma}
&
H\bullet^0K
\ar@2[d]^{H\bullet^0\gamma}
\\
F\bullet^0L
\ar@2[r]_{(\alpha\bullet^1\beta)\bullet^0L}
\ar@{}[ur]|(.3){}="1"|(.7){}="2"
\ar@3"1";"2"^{(\alpha\bullet^1\beta)\bullet^0\gamma\ \ }
&
H\bullet^0L}
}
=
\raisebox{5ex}{
\xymatrix@C=6ex{
F\bullet^0K
\ar@2[r]^{\alpha\bullet^0K}
\ar@2[d]_{F\bullet^0\gamma}
&
G\bullet^0K
\ar@2[r]^{\beta\bullet^0K}
\ar@2[d]|{G\bullet^0\gamma}
&
H\bullet^0K
\ar@2[d]^{H\bullet^0\gamma}
\\
F\bullet^0L
\ar@2[r]_{\alpha\bullet^0L}
\ar@{}[ur]|(.3){}="1"|(.7){}="2"
\ar@3"1";"2"^{\alpha\bullet^0\gamma\ \ }
&
G\bullet^0L
\ar@2[r]_{\beta\bullet^0L}
\ar@{}[ur]|(.3){}="1"|(.7){}="2"
\ar@3"1";"2"^{\beta\bullet^0\gamma\ \ }
&
H\bullet^0L
}}
$$

The same can be claimed for axiom $(v)(b)$.
\end{Remark}

\begin{Theorem}\label{thm:caratterizzo_sesqui2}
Let $\mathcal{C}_{\bullet}$ be a 3-truncated reflexive globular set. Then the following two statements are equivalent.
\begin{enumerate}
  \item  $\mathcal{C}$ is a (small) {\em sesqui${}^2$-category}
  \item Axioms $(i)$, $(ii)$ and $(iii)$ of\/ {\em Definition \ref{def:sequi2cat}} hold, moreover
\end{enumerate}

$(viii)$ The 2-truncation $\xymatrix{\mathcal{C}_2\ar@<+1ex>[r]^{d_1}\ar@<-1ex>[r]_{c_1}
&\ar[l]|{e_1}\mathcal{C}_1\ar@<+1ex>[r]^{d_0}\ar@<-1ex>[r]_{c_0}
&\ar[l]|{e_0}\mathcal{C}_0}$ of $\mathcal{C}_\bullet$ is a sesqui-category.\\

$(ix)$ For every 2-morphism $\alpha:F\Rightarrow G:\C\to\D$ and objects $\B,\E$ of $\mathcal{C}$
$$
-\bullet^0 \alpha:\ -\bullet^0 F\Rightarrow-\bullet^0 G:\ \mathcal{C}(\B,\C)\to\mathcal{C}(\B,\D)
$$
$$
\alpha\bullet^0 -:\ F\bullet^0 -\Rightarrow  G\bullet^0 - :\ \mathcal{C}(\D,\E)\to\mathcal{C}(\C,\E)
$$
are \emph{lax} natural transformations of sesqui-functors.\\

$(x)$ For every morphism $F:\ \C\to\D$ of $\mathcal{C}$
$$
-\bullet^0 id_F\ -\bullet^0 F\Rightarrow -\bullet^0 F
$$
$$
id_F\bullet^0 -:\ F\bullet^0 -\Rightarrow  F\bullet^0 -
$$
are identical natural transformations.\\

$(xi)$  {\em (reduced associativity axiom)}\\

For every 0-composable triple $x\in[\mathcal{C}(\B,\C)]_p$, $y\in[\mathcal{C}(\C,\D)]_q$
and $z\in[\mathcal{C}(\D,\E)]_r$, with \mbox{$p+q+r = 2$}
$$
(x\bullet^0y)\bullet^0z=x\bullet^0(y\bullet^0z)
$$
{i.e.} for 3-morphism $\Lambda$, 2-morphisms $\alpha,\beta$ and morphisms $F,G$ of $\mathcal{C}$, the following equations hold, when
composites exist:
$$
\begin{array}{cc}
(\Lambda \bullet^0 F) \bullet^0 G=\Lambda \bullet^0 (F \bullet^0 G)
&
(\alpha\bullet^0\beta)\bullet^0 F = \alpha\bullet^0(\beta\bullet^0 F)
\\
(F\bullet^0 \Lambda)\bullet^0 F = F\bullet^0(\Lambda\bullet^0 G)
&
(\alpha\bullet^0 F)\bullet^0 \beta = \alpha\bullet^0(F\bullet^0\beta)
\\
(\Lambda\bullet^0 F)\bullet^0 G=\Lambda\bullet^0(F\bullet^0 G)
&
(F\bullet^0\alpha)\bullet^0\beta = F\bullet^0(\alpha\bullet^0\beta)
\end{array}
$$
\end{Theorem}
\begin{proof}
First we prove that  $1.$ implies  $2.$.

Condition $(viii)$ is equivalent to satisfying properties $(L1)$ to $(L4)$, $(R1)$ to $(R4)$ and
$(LR5)$ of {\em Proposition \ref{prop:sesqui-categories}}. Now, $(L1)$ and $(R1)$ hold by $(iii)$, $(L2)$ and $(R2)$  by $(iv)$, $(L3)$, $(R3)$, $(L4)$ and $(R4)$
by $(ii)$, $(LR5)$ by $(vi)$.

Condition $(ix)$ holds. In fact let us recall {\em Definition \ref{def:lax_sesqui_transformation}}.
Assignment on objects (=1-cells) is given by 0-composition, naturality by $(iv)$ and functoriality by
$(v)$ (compositions) and $(vii)$ (units).

Condition $(x)$ holds too. In fact this is implied by $(ix)$ above and $(vii)$.

Finally  $(xi)$  is a subset of $(vi)$.\\

Conversely we prove that  $2.$ implies  $1.$.

Conditions $(iv)$ and $(v)$ hold by $(ix)$.

Condition $(vi)$ holds by $(xi)$ for the cases $p+q+r=2$. What is still to prove is the case $p+q+r=0$
and the case $p+q+r=1$, that are given by $(viii)$.

Finally $(vii)$ is a consequence of $(ix)$ and $(x)$.
\end{proof}
\begin{Remark}
Notice that the characterization  given by {\em Theorem \ref{thm:caratterizzo_sesqui2}} is somehow redundant.
Nevertheless its usefulness is that it makes available practical rules in order to deal with calculations
in a sesqui${}^2$-categorical environment.
\end{Remark}

\section{Lax $n$-modification}
Purpose of the rest of the chapter is to give a proof of the following
\begin{Theorem}
The sesqui-category $n\mathbf{Cat}$, endowed with 3-morphism, their compositions, whiskering and dimension raising
0-composition of 2-morphisms is a sesqui${}^2$-category.
\end{Theorem}
This is done by means of the characterization given in {\em Theorem \ref{thm:caratterizzo_sesqui2}}.

As usual the approach is genuinely inductive, starting with the well known definition of a {\em modification}
in $\mathbf{Cat}$ \cite{B94HANDBOOK}.\\

Hence suppose given an integer $n>1$. Let us consider the following situation in $n$\textbf{Cat}:
$$
\alpha,\beta:\ F\Rightarrow G:\ \C\rightarrow \D
$$
A {\em lax $n$-modification} $
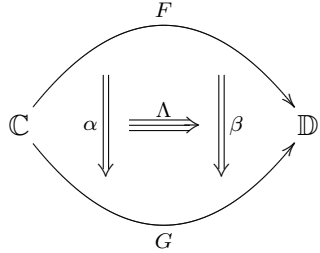
\xymatrix{\Lambda:\ \alpha\ar@3{->}[r]&\beta}$
$$
\xymatrix@C=20ex{
\C
\ar@/^8ex/[r]^{F}_(.3){}="a1"_(.7){}="b1"
\ar@/_8ex/[r]_{G}^(.3){}="a2"^(.7){}="b2"
&\D
\ar@{}"a1";"a2"|(.2){}="A1"|(.8){}="A2"
\ar@{}"b1";"b2"|(.2){}="B1"|(.8){}="B2"
\ar@{=>}"A1";"A2"_{\alpha}^{}="l1"
\ar@{=>}"B1";"B2"^{\beta}_{}="l2"
\ar@{}"l1";"l2"|(.2){}="L1"|(.8){}="L2"
\ar@3{->}"L1";"L2"^{\Lambda}
}
$$
is a pair $(\Lambda_0,\Lambda_1)$, where\\

$\bullet$ $\xymatrix{\Lambda_0:\C_0\ar[r]& {}\substack{\coprod
\\c_0\in \C_0} [\D_2(\alpha_0(c_0),\beta_0(c_0))]_0}$ is a
 map such that, for every $c_0$ in $\C_0$,
$\xymatrix{\Lambda_0(c_0):\alpha_0(c_0)\ar@{=>}[r]&\beta_0(c_0)}$.

\emph{Let us point out that subscript ``$0$'' is sometimes omitted (as in $\alpha(c_0^{\phantom{0}})$), or $c_0^{\phantom{0}}$ is itself subscripted (as in $\alpha_{c_0^{\phantom{0}}} $).\\
}
$\bullet$ (\emph{n-naturality}) for every pair of objects $c_0$, $c_0'$ of $\C$,
a 3-morphism of \mbox{$(n-1)$categories} that {\em fills} the following diagram:
$$
\xymatrix@C=10ex@R=10ex{
& \C_1(c_0,c_0')
\ar[ddl]_{F_1^{c_0,c_0'}}^(.5){}="a2"^(.8){}="b2"
\ar[ddr]^{G_1^{c_0,c_0'}}_(.5){}="a1"_(.8){}="b1"
\\
\\
\D_1(Fc_0,Fc_0')
\ar@/^6ex/[dr]^{-\circ\alpha c_0'}="x1"
\ar@/_6ex/[dr]_{-\circ\beta c_0'}="x2"
&&
\D_1(Gc_0,Gc_0')
\ar@/_6ex/[dl]_{\alpha c_0\circ-}="y1"
\ar@/^6ex/[dl]^{\beta c_0\circ-}="y2"
\\
&
\D_1(Fc_0,Gc_0')
\ar@{}"x1";"x2"|(.3){}="X1"|(.7){}="X2"
\ar@{}"y1";"y2"|(.3){}="Y1"|(.7){}="Y2"
\ar@{=>}"X1";"X2"|{-\circ\Lambda c_0'}
\ar@{=>}"Y1";"Y2"|{\Lambda c_0 \circ -}
\ar@{}"a1";"a2"|(.3){}="A1"|(.7){}="A2"
\ar@{}"b1";"b2"|(.3){}="B1"|(.7){}="B2"
\ar@{=>}"A1";"A2"_{\alpha_1^{c_0,c_0'}}^{}="z2"
\ar@{:>}"B1";"B2"^{\beta_1^{c_0,c_0'}}_{}="z1"
\ar@{}"z1";"z2"|(.25){}="Z1"|(.75){}="Z2"
\ar@3{->}"Z1";"Z2"_{\Lambda_1^{c_0,c_0'}}
}
$$
{\em i.e.}
$$
\xymatrix{
G_1^{c_0,c_0'}\bullet^0(-\circ \alpha c_0')
\ar@2[r]^{\alpha_1^{c_0,c_0'}}
\ar@2[d]_{id\bullet_0(\Lambda c_0 \circ -)}
&F_1^{c_0,c_0'}\bullet^0(-\circ\alpha c_0')
\ar@2[d]^{id\bullet_0(-\circ\Lambda c_0)}
\\
\ar@{}[ur]|(.3){}="1"|(.7){}="2"\ar@3"1";"2"^{\Lambda_1^{c_0,c_0'}}
G_1^{c_0,c_0'}\bullet^0(\beta c_0\circ -)
\ar@2[r]_{\beta_1^{c_0,c_0'}}
&F_1^{c_0,c_0'}\bullet^0(-\circ\beta c_0')
}
$$

These data must obey to {\em functoriality} axioms described by the following equations of 3-diagrams in (n-1)\textbf{Cat}:\\

$\bullet$ ({\em functoriality w.r.t. $0$-composition}) for every triple $c_0,c_0',c_0''$ of objects of $\C$\\
\begin{changemargin}{-10ex}{-10ex}
$$
\xymatrix@R=12ex@C=6ex{\ar@{=}[rr]
\C_1(c_0,c_0')\times\C_1(c_0',c_0'')
\ar@/_18ex/[dd]|(.6){}="3"|(.45){id\times \beta_1^{c_0',c_0''}  }
\ar@/_10ex/@{.>}[dd]|(.4){}="4"
\ar@/^10ex/[dd]|(.6){}="2"
\ar@/^18ex/[dd]|(.4){}="1"|(.55){\quad id\times \alpha_1^{c_0',c_0''}}
&&\C_1(c_0,c_0')\times\C_1(c_0',c_0'')
\ar@/_18ex/[dd]|(.6){}="w3"|(.4){\beta_1^{c_0,c_0'}\times id \quad }
\ar@/_10ex/@{.>}[dd]|(.4){}="w4"
\ar@/^10ex/[dd]|(.6){}="w2"
\ar@/^18ex/[dd]|(.4){}="w1"|(.6){\quad \alpha_1^{c_0,c_0}\times id}
\\
\\
\C_1(c_0,c_0')\times\D_1(Fc_0',Gc_0'')
\ar@{}"1";"2"|(.2){}="a1"|(.8){}="a2"
\ar@{}"2";"3"|(.2){}="b1"|(.8){}="b2"
\ar@{}"4";"3"|(.2){}="d1"|(.8){}="d2"
\ar@{}"1";"4"|(.2){}="c1"|(.8){}="c2"
\ar@{=>}"a1";"a2"
\ar@{=>}"b1";"b2"^{id\times F(-)\circ\Lambda c_0''}
\ar@{:>}"c1";"c2"_{id\times \Lambda c_0' \circ G(-)}
\ar@{:>}"d1";"d2"
\ar@{}"4";"2"|(.35){}="z5"|(.65){}="z6"
\ar@3"z5";"z6"_(.01){id \times\Lambda_1^{c_0',c_0''}\quad}
\ar[dr]_{F(-)\circ-}
&&\D_1(Fc_0,Gc_0')\times\C_1(c_0',c_0'')
\ar@{}"w1";"w2"|(.2){}="wa1"|(.8){}="wa2"
\ar@{}"w2";"w3"|(.2){}="wb1"|(.8){}="wb2"
\ar@{}"w4";"w3"|(.2){}="wd1"|(.8){}="wd2"
\ar@{}"w1";"w4"|(.2){}="wc1"|(.8){}="wc2"
\ar@{=>}"wa1";"wa2"
\ar@{=>}"wb1";"wb2"^{F(-)\circ \Lambda c_0'\times id}
\ar@{:>}"wc1";"wc2"_{\Lambda c_0\circ G(-)\times id}
\ar@{:>}"wd1";"wd2"
\ar@{}"w4";"w2"|(.35){}="wz5"|(.65){}="wz6"
\ar@3"wz5";"wz6"_(.4){\Lambda_1^{c_0,c_0'}\times id}
\ar[dl]^{-\circ G(-)}
\\
&\D_1(Fc_0,Gc_0'')}
$$
\end{changemargin}
\begin{equation}
=
\end{equation}
$$
\xymatrix@R=12ex@C=10ex{\C_1(c_0,c_0')\times\C_1(c_0',c_0'')\ar[d]^{\circ}\\
\C_1(c_0,c_0'')
\ar@/_18ex/[dd]|(.6){}="3"|(.45){F(-)\circ \Lambda c_0''\quad }
\ar@/_10ex/@{.>}[dd]|(.4){}="4"
\ar@/^10ex/[dd]|(.6){}="2"
\ar@/^18ex/[dd]|(.4){}="1"|(.55){\quad \Lambda c_0\circ G(-)}
\\
\\
\D_1(Fc_0,Gc_0'')
\ar@{}"2";"1"|(.2){}="a1"|(.8){}="a2"
\ar@{}"2";"3"|(.2){}="b1"|(.8){}="b2"
\ar@{}"3";"4"|(.2){}="d1"|(.8){}="d2"
\ar@{}"1";"4"|(.2){}="c1"|(.8){}="c2"
\ar@{=>}"a1";"a2"
\ar@{=>}"b1";"b2"^{\alpha_1^{c_0,c_0''}}
\ar@{:>}"c1";"c2"_{\beta_1^{c_0,c_0''}}
\ar@{:>}"d1";"d2"
\ar@{}"1";"3"|(.35){}="z5"|(.65){}="z6"
\ar@3"z5";"z6"^(.4){\Lambda_1^{c_0,c_0''}}}
$$
namely:
$$
(\Lambda_1^{c_0,c_0'}\circ G_1^{c_0',c_0''})\bullet^2 (F_1^{c_0,c_0'}\circ\Lambda_1^{c_0',c_0''}) =(-\circ-)\bullet^0\Lambda_1^{c_0,c_0''}
$$
where the 2-dimensional intersection is the 2-morphism $F(-)\circ\Lambda c_0'\circ G(-)$.\\

$\bullet$ ({\em functoriality w.r.t. units}) for every object $c_0$ of $\C$\\
\begin{changemargin}{-10ex}{-10ex}
\begin{equation}
\raisebox{20ex}{\xymatrix@R=12ex@C=10ex{\Id{n-1}\ar[d]^{u(c_0)}\\
\C_1(c_0,c_0)
\ar@/_18ex/[dd]|(.6){}="3"|(.45){F(-)\circ \Lambda c_0\quad }
\ar@/_10ex/@{.>}[dd]|(.4){}="4"
\ar@/^10ex/[dd]|(.6){}="2"
\ar@/^18ex/[dd]|(.4){}="1"|(.55){\quad \Lambda c_0\circ G(-)}
\\
\\
\D_1(Fc_0,Gc_0)
\ar@{}"2";"1"|(.2){}="a1"|(.8){}="a2"
\ar@{}"2";"3"|(.2){}="b1"|(.8){}="b2"
\ar@{}"3";"4"|(.2){}="d1"|(.8){}="d2"
\ar@{}"1";"4"|(.2){}="c1"|(.8){}="c2"
\ar@{=>}"a1";"a2"
\ar@{=>}"b1";"b2"^{\alpha_1^{c_0,c_0}}
\ar@{:>}"c1";"c2"_{\beta_1^{c_0,c_0}}
\ar@{:>}"d1";"d2"
\ar@{}"1";"3"|(.35){}="z5"|(.65){}="z6"
\ar@3"z5";"z6"^(.4){\Lambda_1^{c_0,c_0}}}
}=
\raisebox{20ex}{
\xymatrix@R=36ex@C=10ex{
\Id{n-1}
\ar@/_18ex/[d]|(.6){}="3"
\ar@/_10ex/@{.>}[d]|(.4){}="4"
\ar@/^10ex/[d]|(.6){}="2"
\ar@/^18ex/[d]|(.4){}="1"
\\
\D_1(Fc_0,Gc_0)
\ar@{}"2";"1"|(.2){}="a1"|(.8){}="a2"
\ar@{}"2";"3"|(.2){}="b1"|(.8){}="b2"
\ar@{}"3";"4"|(.2){}="d1"|(.8){}="d2"
\ar@{}"1";"4"|(.2){}="c1"|(.8){}="c2"
\ar@{=>}"a1";"a2"_(.4){[\Lambda c_0]\quad}
\ar@{=}"b1";"b2"
\ar@{:}"c1";"c2"
\ar@{:>}"d1";"d2"^(.6){\quad[\Lambda c_0]}
\ar@{}"1";"3"|(.35){}="z5"|(.65){}="z6"
\ar@3"z5";"z6"^(.4){Id}}}
\end{equation}
\end{changemargin}
namely:
$$
u(c_0)\bullet^0\Lambda_1^{c_0,c_0}=Id_{[\Lambda c_0]}
$$
\vskip4ex

We write $[\Lambda c_0]$ for the constant 2-morphism given by $\Lambda c_0$; in this case it is between constant
morphisms:
$$
[\Lambda c_0]:\ [\alpha c_0]\Rightarrow [\beta c_0]:\ \Id{n-1}\to \D_1(Fc_0,Gc_0)
$$

Notice that both functoriality axioms for 3-morphisms reduce to those for 2-morphisms, when we consider only
identity 3-morphisms (i.e. 2-morphisms {\em considered as} 3-morphisms).

In the same way functoriality axioms for 2-morphisms reduce to those for 1-morphisms, when we consider only
identity 3-morphisms (i.e. 2-morphisms {\em considered as} 3-morphisms).

\section{$n\mathbf{Cat}(\C,\D)$: the underlying category}
Here and in the following three sections we consider $n$-categories $\C$ and $\D$ be given.
We consider a sesqui-category structure over the category $n\mathbf{Cat}(\C,\D)$. As we did in defining
the sesqui-category $n\mathbf{Cat}$, we start by showing the underlying category structure. This has been already
detailed in section \ref{section:homcats}, hence it suffices to recall that:
\begin{itemize}
  \item objects of $\lfloor n\mathbf{Cat}(\C,\D)\rfloor$ are $n$-functors $\C\to\D$;
  \item arrows of $\lfloor n\mathbf{Cat}(\C,\D)\rfloor$ $n$-lax transformation between them.
\end{itemize}
Composition is $2$-morphisms $1$-composition, obvious units.

\section{$n\mathbf{Cat}(\C,\D)$: the hom-categories}
Let us fix $n$-functors $F,G:\C\to\D$.
We have to define categories $\big(n\mathbf{Cat}(\C,\D)\big)(F,G)$, or more simply
 $n\mathbf{Cat}(F,G)$.

\begin{itemize}
  \item Objects  of  $n\mathbf{Cat}(F,G)$  are 2-morphisms $\alpha:F\Rightarrow G$;
  \item Arrows $\alpha\to\beta$ are 3-morphisms of $n$-categories.
\end{itemize}

\subsection{Composition}

For 3-morphisms $\Lambda=(\Lambda_0,\Lambda_1^{-,-}):\alpha\to\beta$ and $\Sigma=(\Sigma_0,\Sigma_1^{-,-}):\beta\to\gamma$ their 2-composition
$\Lambda\bullet^2\Sigma:\alpha\to\gamma$ is given by the following data:\\

$\bullet$ ({\em on objects})
$$
[\Lambda\bullet^2\Sigma]_0 : c_0\mapsto \Lambda c_0\circ^1 \Sigma c_0
$$
{\em i.r.}
$$
\xymatrix{\alpha c_0\ar@2[r]^{\Lambda c_0}&\beta c_0\ar@2[r]^{\Sigma c_0}&\gamma c_0}
$$

$\bullet$  ({\em on homs}) For chosen objects $c_0,c_0'$ one has
$$
[\Lambda\bullet^2\Sigma]_1^{c_0,c_0'}=\Big(\big(G_1^{c_0,c_0'}\bullet^0(\Lambda c_0\circ-)\big)\bullet^1 \Sigma_1^{c_0,c_0'}\Big) \bullet^2
\Big(\Lambda_1^{c_0,c_0'}\bullet^1\big(F_1^{c_0,c_0'}\bullet^0(-\circ\Sigma c_0')\big)\Big)
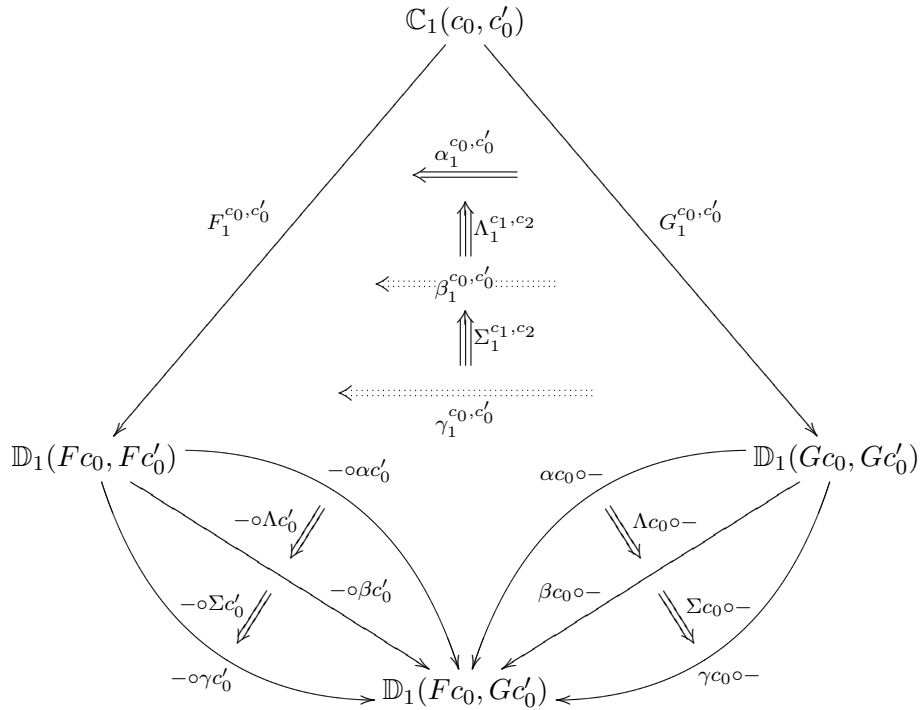
$$
$$
\xymatrix@C=15ex@R=15ex{
& \C_1(c_0,c_0')
\ar[ddl]_{F_1^{c_0,c_0'}}^(.35){}="a2"^(.6){}="b2"^(.85){}="c2"
\ar[ddr]^{G_1^{c_0,c_0'}}_(.35){}="a1"_(.6){}="b1"_(.85){}="c1"
\\
\\
\D_1(Fc_0,Fc_0')
\ar@/^8ex/[dr]^{-\circ\alpha c_0'}|{}="x1"
\ar[dr]^(.65){-\circ\beta c_0'}|{}="x2"
\ar@/_8ex/[dr]_{-\circ \gamma c_0'}|{}="x3"
&&
\D_1(Gc_0,Gc_0')
\ar@/_8ex/[dl]_{\alpha c_0\circ-}|{}="y1"
\ar[dl]_(.65){\beta c_0\circ-}|{}="y2"
\ar@/^8ex/[dl]^{\gamma c_0\circ-}|{}="y3"
\\
&
\D_1(Fc_0,Gc_0')
\ar@{}"x1";"x2"|(.2){}="X1"|(.8){}="X2"
\ar@{}"y1";"y2"|(.2){}="Y1"|(.8){}="Y2"
\ar@{}"x2";"x3"|(.2){}="X3"|(.8){}="X4"
\ar@{}"y2";"y3"|(.2){}="Y3"|(.8){}="Y4"
\ar@{=>}"X1";"X2"_{-\circ\Lambda c_0'}
\ar@{=>}"Y1";"Y2"^{\Lambda c_0 \circ -}
\ar@{=>}"X3";"X4"_{-\circ\Sigma c_0'}
\ar@{=>}"Y3";"Y4"^{\Sigma c_0 \circ -}
\ar@{}"a1";"a2"|(.3){}="A1"|(.7){}="A2"
\ar@{}"b1";"b2"|(.3){}="B1"|(.7){}="B2"
\ar@{}"c1";"c2"|(.3){}="C1"|(.7){}="C2"
\ar@{=>}"A1";"A2"_{\alpha_1^{c_0,c_0'}}^{}="z2"
\ar@{:>}"B1";"B2"|{\beta_1^{c_0,c_0'}}_{}="z1"^{}="w2"
\ar@{:>}"C1";"C2"^{\gamma_1^{c_0,c_0'}}_{}="w1"
\ar@{}"z1";"z2"|(.25){}="Z1"|(.75){}="Z2"
\ar@{}"w1";"w2"|(.25){}="W1"|(.75){}="W2"
\ar@3{->}"Z1";"Z2"_{\Lambda_1^{c_1,c_2}}
\ar@3{->}"W1";"W2"_{\Sigma_1^{c_1,c_2}}
}
$$
We can represent this also as a 2-dimensional pasting, sometimes useful in proofs:
$$
\xymatrix@C=11ex@R=10ex{
G_1^{c_0,c_0'}\bullet^0(\alpha c_0\circ -)\ar@2[r]^{id \bullet^0 (\Lambda c_0\circ -)}\ar@2[d]_{\alpha_1^{c_0,c_0'}}
&G_1^{c_0,c_0'}\bullet^0(\beta c_0 \circ -)\ar@2[r]^{id \bullet^0 (\Sigma c_0\circ -)}\ar@2[d]^{\beta_1^{c_0,c_0'}}
\ar@{}[dl]|(.3){}="a1"|(.7){}="a2"\ar@3"a1";"a2"_{\Lambda_1^{c_0,c_0'}}
&G_1^{c_0,c_0'}\bullet^0(\gamma c_0 \circ -)\ar@2[d]^{\gamma_1^{c_0,c_0'}}
\ar@{}[dl]|(.3){}="b1"|(.7){}="b2"\ar@3"b1";"b2"_{\Sigma_1^{c_0,c_0'}}
\\
F_1^{c_0,c_0'}\bullet^0(-\circ\alpha c_0')\ar@2[r]_{id \bullet^0 (-\circ\Lambda c_0')}
&F_1^{c_0,c_0'}\bullet^0(-\circ\beta c_0')\ar@2[r]_{id \bullet^0 (-\circ\Sigma c_0')}
&F_1^{c_0,c_0'}\bullet^0(-\circ\gamma c_0')
}
$$
Notice that
$$
(-\circ\Lambda c_0')\bullet^1(-\circ\Sigma c_0')=-\circ (\Lambda c_0'\circ \Sigma c_0')=-\circ [\Lambda\bullet^2\Sigma]c_0'
$$
$$
(\Lambda c_0\circ-)\bullet^1(\Sigma c_0\circ-)=(\Lambda c_0\circ \Sigma c_0)\circ - = [\Lambda\bullet^2\Sigma]c_0\circ -
$$
These data form indeed a 3-morphism. In fact let us consider the following diagram
for every triple of objects $c_0,c_0',c_0''$ of $\C$
\begin{changemargin}{-10ex}{-10ex}
$$
\xymatrix@C=20ex@R=10ex{
\alpha c_0\circ G_1^{c_0,c_0'}(-)\circ G_1^{c_0',c_0''}(-)
\ar@2[r]^{[\Lambda\bullet^2\Sigma] c_0\circ id}
\ar@2[d]_{\alpha_1^{c_0,c_0'}\circ id}
&
\gamma c_0\circ G_1^{c_0,c_0'}(-)\circ G_1^{c_0',c_0''}(-)
\ar@2[d]^{\gamma_1^{c_0,c_0'}\circ id}
\ar@{}[dl]|(.4){}="a1"|(.6){}="a2"
\ar@3"a1";"a2"_{[\Lambda\bullet^2\Sigma]_1^{c_0,c_0'}\circ id\qquad }
\\
F_1^{c_0,c_0'}(-)\circ \alpha c_0'\circ G_1^{c_0',c_0''}(-)
\ar@2[d]_{id \circ \alpha_1^{c_0',c_0''}}
\ar@2[r]|{id\circ [\Lambda\bullet^2\Sigma] c_0'\circ id}
&
F_1^{c_0,c_0'}(-)\circ \gamma c_0'\circ G_1^{c_0',c_0''}(-)
\ar@2[d]^{id \circ \gamma_1^{c_0',c_0''}}
\ar@{}[dl]|(.4){}="b1"|(.6){}="b2"
\ar@3"b1";"b2"_{id \circ [\Lambda\bullet^2\Sigma]_1^{c_0',c_0''}\qquad }
\\
F_1^{c_0,c_0'}(-)\circ F_1^{c_0',c_0''}(-)\circ\alpha c_0''
\ar@2[r]_{id\circ[\Lambda\bullet^2\Sigma] c_0''}
&
F_1^{c_0,c_0'}(-)\circ F_1^{c_0',c_0''}(-)\circ\gamma c_0''
}
$$
\end{changemargin}
by definition of $[\Lambda\bullet^2\Sigma]_1^{-,-}$ one has
\begin{changemargin}{-20ex}{-10ex}
\begin{equation}\label{xxx:vc}
\xymatrix@C=14ex@R=10ex{
\alpha c_0\circ G_1^{c_0,c_0'}(-)\circ G_1^{c_0',c_0''}(-)
\ar@2[r]^{\Lambda c_0\circ id}
\ar@2[d]_{\alpha_1^{c_0,c_0'}\circ id}
&
\beta c_0\circ G_1^{c_0,c_0'}(-)\circ G_1^{c_0',c_0''}(-)
\ar@2[r]^{\Sigma c_0\circ id}
\ar@2[d]|{\beta_1^{c_0,c_0'}\circ id}
\ar@{}[dl]|(.4){}="a1"|(.6){}="a2"
\ar@3"a1";"a2"_{\Lambda_1^{c_0,c_0'}\circ id\qquad }
&
\gamma c_0\circ G_1^{c_0,c_0'}(-)\circ G_1^{c_0',c_0''}(-)
\ar@2[d]^{\gamma_1^{c_0,c_0'}\circ id}
\ar@{}[dl]|(.4){}="b1"|(.6){}="b2"
\ar@3"b1";"b2"_{\Sigma_1^{c_0,c_0'}\circ id\qquad }
\\
F_1^{c_0,c_0'}(-)\circ \alpha c_0'\circ G_1^{c_0',c_0''}(-)
\ar@2[d]_{id \circ \alpha_1^{c_0',c_0''}}
\ar@2[r]|{id\circ \Lambda c_0'\circ id}
&
F_1^{c_0,c_0'}(-)\circ \beta c_0'\circ G_1^{c_0',c_0''}(-)
\ar@2[d]|{id \circ \beta_1^{c_0',c_0''}}
\ar@2[r]|{id\circ \Sigma c_0'\circ id}
\ar@{}[dl]|(.4){}="a3"|(.6){}="a4"
\ar@3"a3";"a4"_{id\circ \Lambda_1^{c_0',c_0''}\qquad }
&
F_1^{c_0,c_0'}(-)\circ \gamma c_0'\circ G_1^{c_0',c_0''}(-)
\ar@2[d]^{id \circ \gamma_1^{c_0',c_0''}}
\ar@{}[dl]|(.4){}="b3"|(.6){}="b4"
\ar@3"b3";"b4"_{id\circ \Sigma_1^{c_0',c_0''}\qquad }
\\
F_1^{c_0,c_0'}(-)\circ F_1^{c_0',c_0''}(-)\circ\alpha c_0''
\ar@2[r]_{id\circ\Lambda c_0''}
&
F_1^{c_0,c_0'}(-)\circ F_1^{c_0',c_0''}(-)\circ\beta c_0''
\ar@2[r]_{id\circ\Sigma c_0''}
&
F_1^{c_0,c_0'}(-)\circ F_1^{c_0',c_0''}(-)\circ\gamma c_0''
}
\end{equation}
\end{changemargin}
This diagram is unambiguous because interchange holds on separate components of product
({\em product interchange} in dimension $n-1$, with intersection the constant $[\beta c_0']$).
Hence we get
\begin{changemargin}{-15ex}{-10ex}
$$
\xymatrix@C=16ex@R=10ex{
\alpha c_0\circ G_1^{c_0,c_0''}(-\circ-)
\ar@2[r]^{\Lambda c_0\circ id}
\ar@2[d]_{\alpha_1^{c_0,c_0''}}
&
\beta c_0\circ G_1^{c_0,c_0''}(-\circ-)
\ar@2[r]^{\Sigma c_0\circ id}
\ar@2[d]|{\beta_1^{c_0,c_0''}}
\ar@{}[dl]|(.4){}="a1"|(.6){}="a2"
\ar@3"a1";"a2"_{\Lambda_1^{c_0,c_0''}\qquad }
&
\gamma c_0\circ G_1^{c_0,c_0''}(-\circ-)
\ar@2[d]^{\gamma_1^{c_0,c_0''}}
\ar@{}[dl]|(.4){}="b1"|(.6){}="b2"
\ar@3"b1";"b2"_{\Sigma_1^{c_0,c_0''}\qquad }
\\
F_1^{c_0,c_0''}(-\circ-)\circ \alpha c_0''
\ar@2[r]_{id\circ\Lambda c_0''}
&
F_1^{c_0,c_0''}(-\circ-)\circ \alpha c_0''
\ar@2[r]_{id\circ\Sigma c_0''}
&
F_1^{c_0,c_0''}(-\circ-)\circ \alpha c_0''
}
$$
\end{changemargin}
More simply for any object $c_0$ of $\C$ one has
$$
\xymatrix@C=10ex@R=10ex{
\alpha c_0\circ G_1^{c_0,c_0}(u(c_0))
\ar@2[r]^{ \Lambda c_0 \circ id}
\ar@2[d]_{\alpha_1^{c_0,c_0}}
&
\beta c_0\circ G_1^{c_0,c_0}(u(c_0))
\ar@2[r]^{ \Sigma c_0 \circ id}
\ar@2[d]|{\beta_1^{c_0,c_0}}
\ar@{}[dl]|(.4){}="a1"|(.6){}="a2"
\ar@3"a1";"a2"_{\Lambda_1^{c_0,c_0}\qquad }
&
\gamma c_0\circ G_1^{c_0,c_0}(u(c_0))
\ar@2[d]^{\gamma_1^{c_0,c_0}}
\ar@{}[dl]|(.4){}="b1"|(.6){}="b2"
\ar@3"b1";"b2"_{\Sigma_1^{c_0,c_0}\qquad }
\\
F_1^{c_0,c_0}(u(c_0))\circ \alpha c_0
\ar@2[r]_{id\circ\Lambda c_0}
&
F_1^{c_0,c_0}(u(c_0))\circ \alpha c_0
\ar@2[r]_{id\circ\Sigma c_0}
&
F_1^{c_0,c_0}(u(c_0))\circ \alpha c_0
}
$$
\begin{equation}\label{xxx:vu}
=
\end{equation}
$$
\xymatrix@C=20ex@R=10ex{
[\alpha c_0]
\ar@2[r]^{ [\Lambda c_0 ]}
\ar@{=}[d]
&
[\beta c_0]
\ar@2[r]^{ [\Sigma c_0 ]}
\ar@{=}[d]
\ar@{}[dl]|(.4){}="a1"|(.6){}="a2"
\ar@3"a1";"a2"_{id }
&
[\gamma c_0]
\ar@{=}[d]
\ar@{}[dl]|(.4){}="b1"|(.6){}="b2"
\ar@3"b1";"b2"_{id }
\\
[\alpha c_0]
\ar@2[r]_{[\Lambda c_0]}
&
[\beta c_0]
\ar@2[r]_{[\Sigma c_0]}
&
[\gamma c_0]
}
$$
where, as usual, square brackets mean {\em constant}.

\subsection{Units}

For any 2-morphism $\beta:\ F\Rightarrow G:\ \C\to\D$ its identity 3-morphisms $id_{\beta}$
is given by:\\

$\bullet$ ({\em on objects})
$$
[id_{\beta}]_0 : c_0\mapsto id_{\beta c_0}
$$

$\bullet$ ({\em on homs}) For chosen objects $c_0,c_0'$ one has
$$
[id_{\beta}]_1^{-,-} = id_{\beta_1^{-,-}}
$$
It is immediate to check that above pair is indeed a 3-morphisms.

Similarly associativity and neutral units follows from same properties for 2-cells and from diagrams (\ref{xxx:vc}) and (\ref{xxx:vu})
suitably adapted ({\em adding one more column, for what concerns associativity, trivializing one column, for what concerns units}).

\section{$n\mathbf{Cat}(\C,\D)$: the sesqui-categorical structure}
In the this section we will show that hom-categories $n\mathbf{Cat}(\C,\D)$ underly a structure of
sesqui-categories, with 2-cells provided by 3-morphisms of n-categories. To this end we define reduced left/right
1-composition of a 3-morphism with a 2-morphism, according to the following reference diagram.
$$
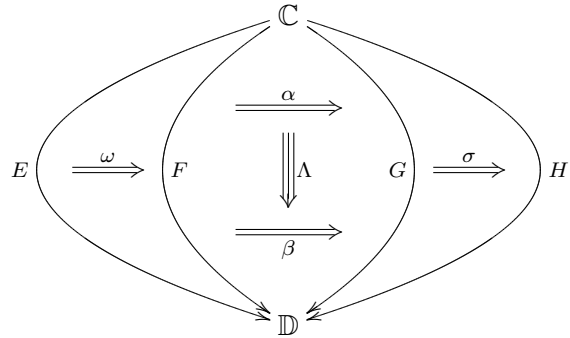
\xymatrix@R=22ex{\C
\ar@/_20ex/[d]_{E}|{}="1"
\ar@/_10ex/[d]^{F}="2"|(.3){}="u2"|(.7){}="d2"
\ar@/^10ex/[d]_{G}="3"|(.3){}="u3"|(.7){}="d3"
\ar@/^20ex/[d]^{H}|{}="4"
\\
\ar@{}"1";"2"|(.25){}="11"|(.75){}="22"
\ar@{=>}"11";"22"^{\omega}
\ar@{}"u2";"u3"|(.25){}="U2"|(.75){}="U3"
\ar@{=>}"U2";"U3"^{\alpha}_{}="l1"
\ar@{}"d2";"d3"|(.25){}="D2"|(.75){}="D3"
\ar@{=>}"D2";"D3"_{\beta}^{}="l2"
\ar@{}"3";"4"|(.25){}="33"|(.75){}="44"
\ar@{=>}"33";"44"^{\sigma}
\ar@{}"l1";"l2"|(.2){}="L1"|(.8){}="L2"
\ar@3"L1";"L2"^{\Lambda}
\D}
$$

\subsection{Reduced left-composition}

The 3-morphism
$$
\xymatrix@C=3ex{\omega\bullet^1 \Lambda:\ \omega\bullet^1\alpha\ar@3[r]
&\omega\bullet^1\beta:\ E\ar@2[r]
&G:\ \C\ar[r]
&\D}
$$
is given by the data below\\

$\bullet$ ({\em on objects}) For an object $c_0$ of $\C$
$$
[\omega\bullet^1\Lambda]_0:\ c_0\mapsto\ \raisebox{10ex}{
\xymatrix{
Ec_0\ar[d]^{\omega c_0}
\\
Fc_0
\ar@/_5ex/[d]_{\alpha c_0}^{}="1"
\ar@/^5ex/[d]^{\beta c_0}^{}="2"
\ar@{}"1";"2"|(.3){}="a1"|(.7){}="a2"
\ar@2"a1";"a2"^{\Lambda c_0}
\\
Gc_0}}
$$

$\bullet$ ({\em on homs}) For objects $c_0,c_0'$ of $\C$
\begin{eqnarray*}
  [\omega\bullet^1\Lambda]_1^{c_0,c_0'} &=& \Big(\Lambda_1^{c_0,c_0'}\bullet^0(\omega c_0\circ-)\Big)\bullet^1\Big(\omega_1^{c_0,c_0'}\bullet^0(-\circ\beta c_0')\Big) \\
   &=& (\omega c_0 \circ\Lambda_1^{c_0,c_0'})\bullet^1(\omega_1^{c_0,c_0'}\circ\beta c_0')
\end{eqnarray*}
\begin{changemargin}{-10ex}{-10ex}
$$
\xymatrix@C=6ex@R=10ex{
\omega c_0\circ\alpha c_0\circ G_1^{c_0,c_0'}(-)
\ar@{:>}[rr]^{[\omega\bullet^1\alpha]_1^{c_0,c_0'}}
\ar@2[dr]|{\omega c_0\circ \alpha_1^{c_0,c_0'}}
\ar@2[ddd]_(.45){[\omega\bullet^1\Lambda]c_0\circ G_1^{c_0,c_0'}(-)}_(.55){=\quad \omega c_0 \circ \Lambda c_0 \circ G_1^{c_0,c_0'}(-)}
&&E_1^{c_0,c_0'}(-)\circ\omega c_0'\circ\alpha c_0'
\ar@{:>}[ddd]^(.45){E_1^{c_0,c_0'}(-)\circ[\omega\bullet^1\Lambda]c_0'}^(.55){=E_1^{c_0,c_0'}(-)\circ\omega c_0'\circ\Lambda c_0'}
\\
&\omega c_0\circ F_1^{c_0,c_0'}(-)\circ\alpha c_0'
\ar@{:>}[ur]|{\omega_1^{c_0,c_0'}\circ\alpha c_0'}
\ar@2[d]^{\omega c_0\circ F_1^{c_0,c_0'}(-)\circ\Lambda c_0'}
\\
\ar@{}[ur]|(.3){}="1"|(.6){}="2"
\ar@3"1";"2"^{\omega c_0\circ\Lambda_1^{c_0,c_0'}}
&\omega c_0\circ F_1^{c_0,c_0'}(-)\circ\beta c_0'
\ar@2[dr]|{\omega_1^{c_0,c_0'}\circ\beta c_0'}
\\
\omega c_0 \circ \beta c_0\circ G_1^{c_0,c_0'}(-)
\ar@2[ur]|{\omega c_0\circ \beta_1^{c_0,c_0'}}
\ar@{:>}[rr]_{[\omega\bullet^1\beta]_1^{c_0,c_0'}}
&&E_1^{c_0,c_0'}(-)\circ\omega c_0'\circ\beta c_0'
}
$$
\end{changemargin}
The pair $<[\omega\bullet^1\Lambda]_0,[\omega\bullet^1\Lambda]_1^{-,-}>$ forms indeed a 3-morphism of n-categories.
\begin{proof}
We have to show that it satisfies composition and unit axioms. Let us begin with composition, and fix a triple
$c_0,c_0',c_0''$ of $\C$. Notice that, in order to keep diagrams in the page we denote $\circ^0$-composition
by juxtaposition, and subscripts for transformations on objects are used.
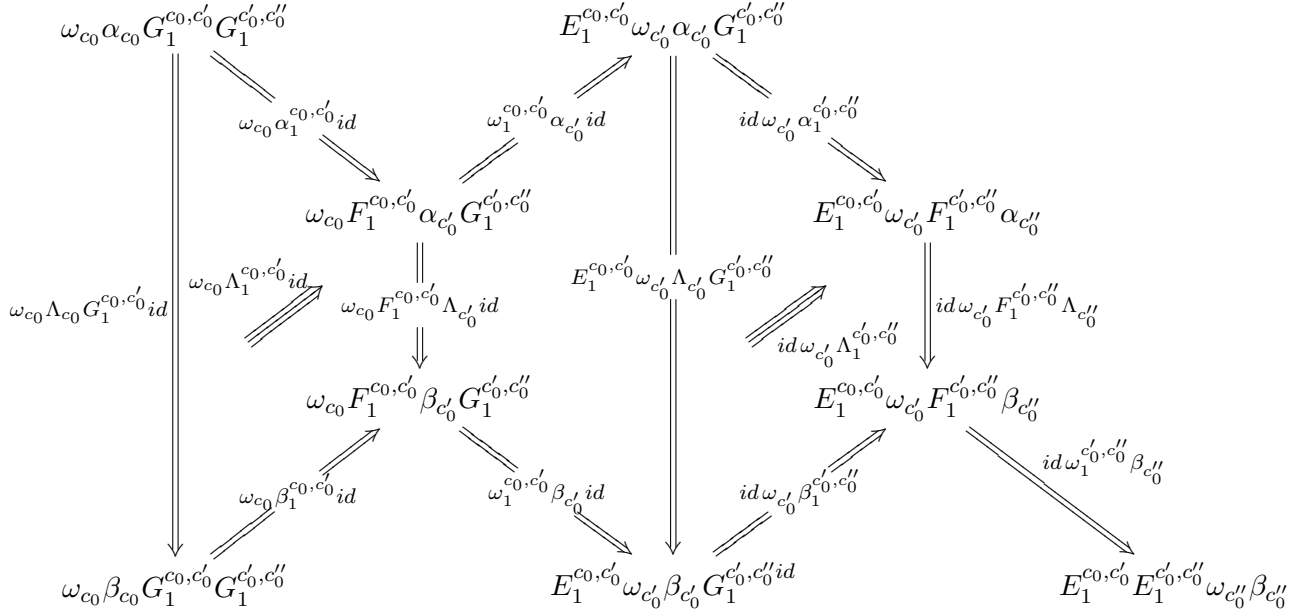
\begin{changemargin}{-10ex}{-10ex}
$$
\xymatrix@C=0ex@R=10ex{
\omega_{c_0} \alpha_{c_0}  G_1^{c_0,c_0'} G_1^{c_0',c_0''}
\ar@2[dr]|{\omega_{c_0}  \alpha_1^{c_0,c_0'} id}
\ar@2[ddd]_{\omega_{c_0}   \Lambda_{c_0}   G_1^{c_0,c_0'}  id}
&&E_1^{c_0,c_0'}  \omega_{c_0'} \alpha_{c_0'} G_1^{c_0',c_0''}
\ar@2[ddd]|(.45){E_1^{c_0,c_0'}  \omega_{c_0'} \Lambda_{c_0'}G_1^{c_0',c_0''}}
\ar@2[dr]|{id\, \omega_{c_0'}  \alpha_1^{c_0',c_0''}}
\\
&\omega_{c_0}  F_1^{c_0,c_0'}  \alpha_{c_0'} G_1^{c_0',c_0''}
\ar@2[ur]|{\omega_1^{c_0,c_0'}\alpha_{c_0'}id}
\ar@2[d]|{\omega_{c_0}  F_1^{c_0,c_0'}  \Lambda_{c_0'} id}
&&E_1^{c_0,c_0'} \omega_{c_0'}  F_1^{c_0',c_0''}  \alpha_{c_0''}
\ar@2[d]^{id\,\omega_{c_0'}  F_1^{c_0',c_0''}  \Lambda_{c_0''}}
\\
\ar@{}[ur]|(.3){}="a1"|(.6){}="a2"
\ar@3"a1";"a2"^(.7){\omega_{c_0} \Lambda_1^{c_0,c_0'} id}
&\omega_{c_0}  F_1^{c_0,c_0'}  \beta_{c_0'} G_1^{c_0',c_0''}
\ar@2[dr]|{\omega_1^{c_0,c_0'} \beta_{c_0'} id}
&
\ar@{}[ur]|(.3){}="b1"|(.6){}="b2"
\ar@3"b1";"b2"_{id\,\omega_{c_0'} \Lambda_1^{c_0',c_0''}}
&E_1^{c_0,c_0'} \omega_{c_0'}  F_1^{c_0',c_0''}  \beta_{c_0''}
\ar@2[dr]^{id\,\omega_1^{c_0',c_0''} \beta_{c_0''}}
\\
\omega_{c_0}   \beta_{c_0}  G_1^{c_0,c_0'} G_1^{c_0',c_0''}
\ar@2[ur]|{\omega_{c_0}  \beta_1^{c_0,c_0'} id}
&&E_1^{c_0,c_0'}  \omega_{c_0'} \beta_{c_0'} G_1^{c_0',c_0'' id}
\ar@2[ur]|{id\,\omega_{c_0'}  \beta_1^{c_0',c_0''}}
&&E_1^{c_0,c_0'}  E_1^{c_0',c_0''}  \omega_{c_0''} \beta_{c_0''}}
$$
\end{changemargin}
By product interchange it is clear that
\begin{changemargin}{-10ex}{-10ex}
$$
\Big(\omega_1^{c_0,c_0'}\circ\alpha_{c_0'}\circ id_{G_1^{c_0',c_0''}}\Big)\bullet^1\Big(id_{E_1^{c_0,c_0'}}\circ\omega_{c_0'}\circ\Lambda_1^{c_0',c_0''}\Big)=
\Big(\omega_{c_0}\circ F_1^{c_0,c_0'}\circ\Lambda_1^{c_0',c_0''}\Big)\bullet^1\Big(\omega_1^{c_0,c_0'}\circ id_{F_1^{c_0',c_0''}}\circ \beta_{c_0''}\Big)
$$
\end{changemargin}
and the diagram can be re-drawn
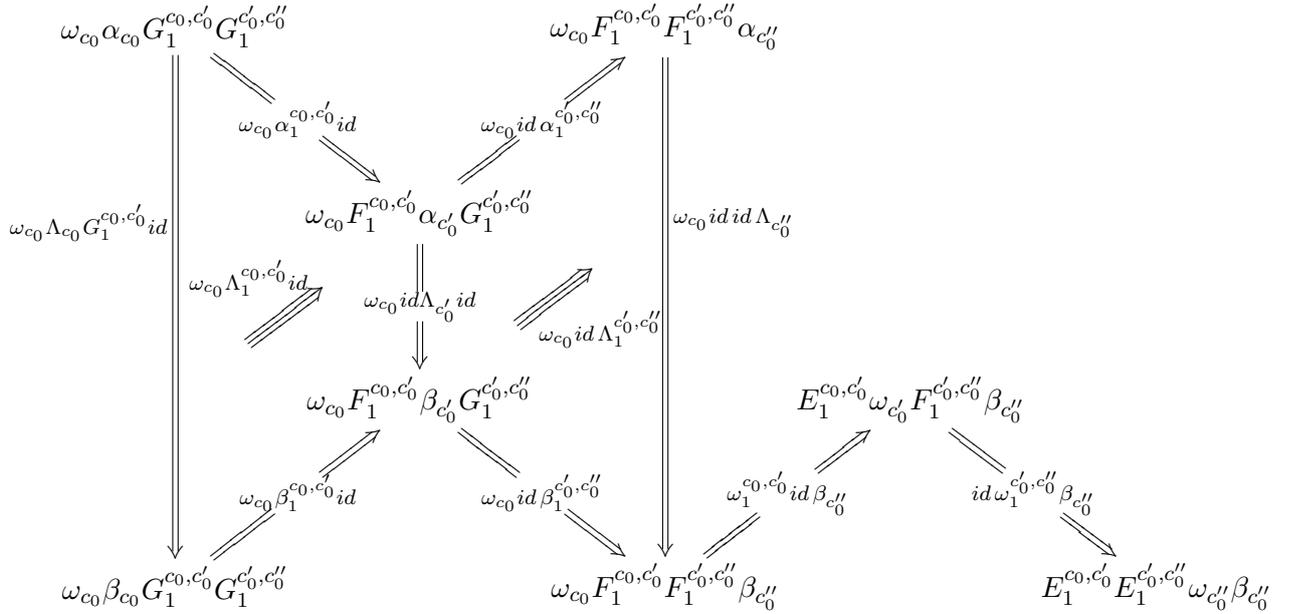
\begin{changemargin}{-10ex}{-10ex}
$$
\xymatrix@C=0ex@R=10ex{
\omega_{c_0} \alpha_{c_0}  G_1^{c_0,c_0'} G_1^{c_0',c_0''}
\ar@2[dr]|{\omega_{c_0}  \alpha_1^{c_0,c_0'} id}
\ar@2[ddd]_(.35){\omega_{c_0}   \Lambda_{c_0}   G_1^{c_0,c_0'}  id}
&&\omega_{c_0}  F_1^{c_0,c_0'}  F_1^{c_0',c_0''}  \alpha_{c_0''}
\ar@2[ddd]^(.35){\omega_{c_0}id\,id\,\Lambda_{c_0''}}\\
&\omega_{c_0}  F_1^{c_0,c_0'}  \alpha_{c_0'} G_1^{c_0',c_0''}
\ar@2[ur]|{\omega_{c_0}id\,\alpha_1^{c_0',c_0''}}
\ar@2[d]|{\omega_{c_0}  id  \Lambda_{c_0'} id}
&&
\\
\ar@{}[ur]|(.3){}="a1"|(.6){}="a2"
\ar@3"a1";"a2"^(.7){\omega_{c_0} \Lambda_1^{c_0,c_0'} id}
&\omega_{c_0}  F_1^{c_0,c_0'}  \beta_{c_0'} G_1^{c_0',c_0''}
\ar@2[dr]|{\omega_{c_0}id\,\beta_1^{c_0',c_0''}}
\ar@{}[ur]|(.4){}="b1"|(.7){}="b2"
\ar@3"b1";"b2"_(.4){\omega_{c_0}id\, \Lambda_1^{c_0',c_0''} }
&&E_1^{c_0,c_0'}    \omega_{c_0'}F_1^{c_0',c_0''} \beta_{c_0''}
\ar@2[dr]|{id\,\omega_1^{c_0',c_0''}\beta_{c_0''}}
\\
\omega_{c_0}   \beta_{c_0}  G_1^{c_0,c_0'} G_1^{c_0',c_0''}
\ar@2[ur]|{\omega_{c_0}  \beta_1^{c_0,c_0'} id}
&&\omega_{c_0}  F_1^{c_0,c_0'}   F_1^{c_0',c_0''}  \beta_{c_0''}
\ar@2[ur]|{\omega_1^{c_0,c_0'}id\,\beta_{c_0''}}
&&E_1^{c_0,c_0'}  E_1^{c_0',c_0''}  \omega_{c_0''} \beta_{c_0''}}
$$
\end{changemargin}
Now let us consider the pasting of the left-hand side of above diagram, and write it equationally:
$$
\begin{array}{c}
\Big[
\big(\Lambda_1^{c_0,c_0'}\circ G_1^{c_0',c_0''}
\big)
\bullet^0
\big(\omega_{c_0}\circ-
\big)
\Big]
\bullet^1
\Big[
\big(F_1^{c_0,c_0'}\circ \beta_1^{c_0',c_0''}
\big)
\bullet^0
\big(\omega_{c_0}\circ-
\big)
\Big]
\\ \phantom{\Big[}\bullet^2 \phantom{\Big]}\\
\Big[
\big(\alpha_1^{c_0,c_0'}\circ G_1^{c_0',c_0''}
\big)
\bullet^0
\big(\omega_{c_0}\circ-
\big)
\Big]
\bullet^1
\Big[
\big(F_1^{c_0,c_0'}\circ \Lambda_1^{c_0',c_0''}
\big)
\bullet^0
\big(\omega_{c_0}\circ-
\big)
\Big]
\end{array}
$$
By {\em whiskering interchange property} $(LRW)$ this equals to
$$
\begin{array}{c}
\Big[
\big(\Lambda_1^{c_0,c_0'}\circ G_1^{c_0',c_0''}
\big)
\bullet^1
\big(F_1^{c_0,c_0'}\circ \beta_1^{c_0',c_0''}
\big)
\Big]
\bullet^0
\big(\omega_{c_0}\circ-
\big)
\\ \phantom{\Big[}\bullet^2 \phantom{\Big]}\\
\Big[
\big(\alpha_1^{c_0,c_0'}\circ G_1^{c_0',c_0''}
\big)
\bullet^1
\big(F_1^{c_0,c_0'}\circ \Lambda_1^{c_0',c_0''}
\big)
\Big]
\bullet^0
\big(\omega_{c_0}\circ-
\big)
\end{array}
$$
and by functoriality of right 0-whiskering $(R4)''$
\begin{changemargin}{-15ex}{-15ex}
$$
\left[ \Big[
\big(\Lambda_1^{c_0,c_0'}\circ G_1^{c_0',c_0''}
\big)
\bullet^1
\big(F_1^{c_0,c_0'}\circ \beta_1^{c_0',c_0''}
\big)
\Big]^{\phantom{X}}
\!\!\!\!
\bullet^2
\Big[
\big(\alpha_1^{c_0,c_0'}\circ G_1^{c_0',c_0''}
\big)
\bullet^1
\big(F_1^{c_0,c_0'}\circ \Lambda_1^{c_0',c_0''}
\big)
\Big]
\right]
\bullet^0
\big(\omega_{c_0}\circ-
\big)
$$
\end{changemargin}
finally by functoriality w.r.t. composition of 3-morphism $\Lambda$
$$
\Lambda_1^{c_0,c_0''}\bullet^0(\omega_{c_0}\circ-)=\omega_{c_0}\circ\Lambda_1^{c_0,c_0''}
$$
Concerning composite 2-morphism on the right-hand side,
we can  apply whiskering properties and composition axiom for 2-morphisms
$$
\big(\omega_1^{c_0,c_0'}\circ F_1^{c_0,c_0''}\circ \beta_{c_0''}\big)
\bullet^1
\big(E_1^{c_0,c_0'}\circ\omega_1^{c_0',c_0''}\circ \beta_{c_0''}\big)
=\omega_1^{c_0,c_0''}\circ \beta_{c_0''}
$$
and diagram above can be re-drawn as follows
$$
\xymatrix@C=15ex@R=20ex{
\omega_{c_0}\alpha_{c_0}G_1^{c_0,c_0''}
\ar@2[r]^{\omega_{c_0}\alpha_1^{c_0,c_0''}}
\ar@2[d]_{\omega_{c_0}\Lambda_{c_0}id}
&\omega_{c_0}F_1^{c_0,c_0''}\alpha_{c_0''}
\ar@2[d]^{\omega_{c_0}id\,\Lambda_{c_0''}}
\\
\ar@{}[ur]|(.3){}="1"|(.6){}="2"
\ar@3"1";"2"^{\omega_{c_0}\Lambda_1^{c_0,c_0''}}
\omega_{c_0}\beta_{c_0}G_1^{c_0,c_0''}
\ar@2[r]_{\omega_{c_0}\beta_1^{c_0,c_0''}}
&\omega_{c_0}F_1^{c_0,c_0''}\beta_{c_0''}
\ar@2[r]_{\omega_1^{c_0,c_0''}\beta_{c_0''}}
&E_1^{c_0,c_0''}\omega_{c_0''}\beta_{c_0''}
}
$$
and this conclude the proof of composition axiom.\\

Concerning units, for an object $c_0$ of $\C$ one has
$$
\xymatrix@C=10ex@R=20ex{
\omega_{c_0}\alpha_{c_0}G_1^{c_0,c_0}(u(c_0))
\ar@2[r]^{\omega_{c_0}\alpha_1^{c_0,c_0}}
\ar@2[d]_{\omega_{c_0}\Lambda_{c_0}id}
&\omega_{c_0}F_1^{c_0,c_0}(u(c_0))\alpha_{c_0}
\ar@2[d]^{\omega_{c_0}id\,\Lambda_{c_0}}
\\
\ar@{}[ur]|(.3){}="1"|(.6){}="2"
\ar@3"1";"2"^{\omega_{c_0}\Lambda_1^{c_0,c_0}}
\omega_{c_0}\beta_{c_0}G_1^{c_0,c_0}(u(c_0))
\ar@2[r]_{\omega_{c_0}\beta_1^{c_0,c_0}}
&\omega_{c_0}F_1^{c_0,c_0}(u(c_0))\beta_{c_0}
\ar@2[r]_{\omega_1^{c_0,c_0}\beta_{c_0}}
&E_1^{c_0,c_0}(u(c_0))\omega_{c_0}\beta_{c_0}
}
$$
then by functoriality w.r.t. units of 3-morphisms of (n-1)categories
(and also by  functoriality w.r.t. units of 2-morphisms and 1-morphisms)
we get
$$
\xymatrix@C=15ex@R=20ex{
\omega_{c_0}[\alpha_{c_0}]
\ar@2[r]^{\omega_{c_0}id}
\ar@2[d]_{\omega_{c_0}[\Lambda_{c_0}]}
&\omega_{c_0}[\alpha_{c_0}]
\ar@2[d]^{\omega_{c_0}[\Lambda_{c_0}]}
\\
\ar@{}[ur]|(.3){}="1"|(.6){}="2"
\ar@3"1";"2"^{\omega_{c_0}id}
\omega_{c_0}[\beta_{c_0}]
\ar@2[r]_{\omega_{c_0}id}
&\omega_{c_0}[\beta_{c_0}]
\ar@2[r]_{\omega_{c_0}id}
&\omega_{c_0}[\beta_{c_0}]
}
$$
hence the result.
\end{proof}

\subsection{Reduced right-composition}
The 3-morphism
$$
\xymatrix@C=3ex{\Lambda\bullet^1\sigma:\ \alpha\bullet^1\sigma\ar@3[r]
&\beta\bullet^1\sigma:\ F\ar@2[r]
&H:\ \C\ar[r]
&D}
$$
is given by the data below\\

$\bullet$ ({\em on objects}) For an object $c_0$ of $\C$
$$
[\Lambda\bullet^1\sigma]_0:\ c_0\mapsto\ \raisebox{10ex}{
\xymatrix{
Fc_0
\ar@/_5ex/[d]_{\alpha c_0}^{}="1"
\ar@/^5ex/[d]^{\beta c_0}^{}="2"
\ar@{}"1";"2"|(.3){}="a1"|(.7){}="a2"
\ar@2"a1";"a2"^{\Lambda c_0}
\\
Gc_0
\ar[d]^{\sigma c_0}
\\
Hc_0}}
$$

$\bullet$ ({\em on homs}) For objects $c_0,c_0'$ of $\C$
\begin{eqnarray*}
  [\Lambda\bullet^1\sigma]_1^{c_0,c_0'} &=& \Big(\sigma_1^{c_0,c_0'}\bullet^0(\alpha c_0\circ-)\Big)
  \bullet^1\Big(\Lambda_1^{c_0,c_0'}\bullet^0(-\circ \sigma c_0')\Big) \\
   &=& (\alpha c_0\circ \sigma_1^{c_0,c_0'})\bullet^1(\Lambda_1^{c_0,c_0'}\circ \sigma c_0')
\end{eqnarray*}
\begin{changemargin}{-10ex}{-10ex}
$$
\xymatrix@C=6ex@R=10ex{
\alpha c_0\circ \sigma c_0\circ H_1^{c_0,c_0'}(-)
\ar@{:>}[rr]^{[\alpha\bullet^1\sigma]_1^{c_0,c_0'}}
\ar@2[dr]|{\alpha c_0\circ \sigma_1^{c_0,c_0'}}
\ar@{:>}[ddd]_(.45){[\Lambda\bullet^1\sigma]c_0\circ H_1^{c_0,c_0'}}_(.55){=\quad \Lambda c_0 \circ \sigma c_0\circ H_1^{c_0,c_0'}}
&&
F_1^{c_0,c_0'}(-)\circ\alpha c_0'\circ \sigma c_0'
\ar@2[ddd]^(.45){F_1^{c_0,c_0'}\circ [\Lambda\bullet^1\sigma]c_0'}^(.55){=F_1^{c_0,c_0'}\circ \Lambda c_0'\circ\sigma c_0'}
\\
&\alpha c_0\circ G_1^{c_0,c_0'}(-)\circ\sigma c_0'
\ar@2[ur]|{\alpha_1^{c_0,c_0'}\circ\sigma c_0'}
\ar@2[d]_{\Lambda c_0 \circ G_1^{c_0,c_0'}\circ \sigma c_0'}
&
\\
&\beta c_0\circ G_1^{c_0,c_0'}(-)\circ\sigma c_0'
\ar@{}[ur]|(.3){}="1"|(.6){}="2"
\ar@3"1";"2"^{\Lambda_1^{c_0,c_0'}\circ \sigma c_0'}
\ar@2[dr]|{\beta_1^{c_0,c_0'}\circ\sigma c_0'}
\\
\beta c_0\circ \sigma c_0\circ H_1^{c_0,c_0'}(-)
\ar@{:>}[ur]|{\beta c_0\circ \sigma_1^{c_0,c_0'}}
\ar@{:>}[rr]_{[\beta\bullet^1\sigma]_1^{c_0,c_0'}}
&&
F_1^{c_0,c_0'}(-)\circ\beta c_0'\circ \sigma c_0'
}
$$
\end{changemargin}
The pair $<[\Lambda\bullet^1\sigma]_0,[\Lambda\bullet^1\sigma]_1^{-,-}>$ forms indeed a 3-morphism of n-categories. The proof
is a straightforward variation of the proof for reduced right-composition above, hence it is omitted.

\subsection{Properties}
In this section we give some properties of left/right 1-composition of a 3-morphism with a 2-morphism.
They are modeled on similar properties given in the definition of a sesqui-category, and they are extremely
useful in dealing with calculations. Let us consider the diagram
$$
\xymatrix@C=5ex@R=10ex{
E'\ar@2[r]^{\omega'}
&E\ar@2[r]^{\omega}
&F\ar@2@/^8ex/[rr]^{\alpha}="1"\ar@2[rr]|{\beta}="2"\ar@2@/_8ex/[rr]_{\gamma}="3"
&&G\ar@2[r]^{\sigma}
&H\ar@2[r]^{\sigma'}
&H'
\ar@{}"1";"2"|(.3){}="l1"|(.75){}="l2"
\ar@3"l1";"l2"^{\Lambda}
\ar@{}"2";"3"|(.25){}="s1"|(.7){}="s2"
\ar@3"s1";"s2"^{\Sigma}
}
$$
as a reference for the following
\begin{Proposition}
(2-composition (i.e. vertical) composition of 3-morphisms w.r.t. (reduced) 1-composition with a 2-morphism)
$$
\begin{array}{clcl}
  (L1)' & id_F\bullet_L^1 \Lambda =\Lambda & (R1)' & \Lambda\bullet_R^1 id_G =\Lambda\\
  (L2)' & (\omega'\bullet^1\omega)\bullet_L^1 \Lambda=\omega' \bullet_L^1  (\omega\bullet_L^1\Lambda) & (R2)' & \Lambda\bullet_R^1 (\sigma\bullet^1\sigma') =(\Lambda \bullet_R^1 \sigma) \bullet_R^1 \sigma'\\
  (L3)' & \omega\bullet_L^1 id_{\alpha}= id_{\omega\alpha}& (R3)' &id_{\alpha}\bullet_R^1\sigma=id_{\alpha\sigma}\\
  (L4)' & \omega\bullet_L^1(\Lambda\bullet^2\Sigma)=(\omega\bullet_L^1\Lambda)\bullet^2(\omega\bullet_L^1\Sigma)
 & (R4)'& (\Lambda\bullet^2\Sigma)\bullet_R^1 \sigma=(\Lambda\bullet_R^1\sigma)\bullet^2(\Sigma\bullet_R^1\sigma)
\end{array}
$$
$$
(LR5)'\quad (\omega\bullet_L^1\Lambda)\bullet_R^1 \sigma= \omega\bullet_L^1(\Lambda\bullet_R^1 \sigma)
$$
\end{Proposition}
\begin{proof}
We only prove  statements $(L1)'$ to $(L4)'$ and statement $(LR5)$. Proof of properties $(R1)'$ to $(R4)'$ is just
 a straightforward variation of proof of $(L1)'$ to $(L4)'$, hence it will be omitted.

$\bullet\ (L1)'$ ({\em on objects}) Let an object $c_0$ of $\C$ be given. Then
$$
[id_F\bullet^1 \Lambda]_{c_0} = [id_F]_{c_0}\circ\Lambda_{c_0} =id_{Fc_0}\circ\Lambda_{c_0}=\Lambda_{c_0}
$$

({\em on homs}) Let objects $c_0,c_0'$ of $\C$ be given. Then

\begin{eqnarray*}
[id_F\bullet^1 \Lambda]_1^{c_0,c_0'}   &=& ([id_F]_{c_0}\circ\Lambda_1^{c_0,c_0'})\bullet^1([id_F]_1^{c_0,c_0'}\circ\beta_{c_0'})  \\
   &=&(id_{Fc_0}\circ\Lambda_1^{c_0,c_0'})\bullet^1(id_{F_1^{c_0,c_0'}}\circ\beta_{c_0'})    \\
   &\eq{\clubsuit}& \Lambda_1^{c_0,c_0'}\bullet^1 id_{F_1^{c_0,c_0'}\circ\beta_{c_0'}} \\
   &\eq{\spadesuit}&  \Lambda_1^{c_0,c_0'}
\end{eqnarray*}
where $(\clubsuit)$ holds by $(R3)$ and $(\spadesuit)$ holds by $(L1)'$ in dimension $n-1$.\\

$\bullet\ (L2)'$ ({\em on objects}) Let an object $c_0$ of $\C$ be given. Then
\begin{eqnarray*}
[(\omega'\bullet^1\omega)\bullet^1 \Lambda]_{c_0}   &=&  [\omega'\bullet^1\omega]_{c_0}\circ\Lambda_{c_0}\\
   &=&  \omega'_{c_0}\circ\omega_{c_0}\circ\Lambda_{c_0}\\
   &=&  \omega'_{c_0}\circ[\omega\bullet^1\Lambda]_{c_0}\\
   &=&  [\omega'_{c_0}\bullet^1(\omega\bullet^1\Lambda)]_{c_0}
\end{eqnarray*}

({\em on homs}) Let objects $c_0,c_0'$ of $\C$ be given. Then
$$
[(\omega'\bullet^1\omega)\bullet^1 \Lambda]_1^{c_0,c_0'}   =
$$
\begin{eqnarray*}
   &=&([\omega'\bullet^1\omega]_{c_0}\circ\Lambda_1^{c_0,c_0'})\bullet_1([\omega'\bullet^1\omega]_1^{c_0,c_0'}\circ \beta_{c_0'})\\
   &=&  (\omega'_{c_0}\circ\omega_{c_0}\circ\Lambda_1^{c_0,c_0'})\bullet_1\Big(\big((\omega'_{c_0}\circ\omega_1^{c_0,c_0'})\bullet^1 ({\omega_1'}^{c_0,c_0'}\circ\omega_{c_0'})\big)\circ \beta_{c_0'}    \Big)\\
   &=&  (\omega'_{c_0}\circ\omega_{c_0}\circ\Lambda_1^{c_0,c_0'})\bullet_1(\omega'_{c_0}\circ\omega_1^{c_0,c_0'}\circ\beta_{c_0'})\bullet^1 ({\omega_1'}^{c_0,c_0'}\circ\omega_{c_0'}\circ \beta_{c_0'})\quad(\spadesuit)\\
   &=&  \Big(\omega'_{c_0}\circ\big((\omega_{c_0}\circ\Lambda_1^{c_0,c_0'})\bullet^1(\omega_1^{c_0,c_0'}\circ\beta_{c_0'})\big)\Big)\bullet^1({\omega'_1}^{c_0,c_0'}\circ (\omega_{c_0'}\circ\beta_{c_0'}))\\
   &=&  (\omega'_{c_0}\circ[\omega\bullet^1\Lambda]_1^{c_0,c_0'})\bullet^1({\omega'_1}^{c_0,c_0'}\circ [\omega\bullet^1\beta]_{c_0'})\\
   &=&  [\omega'\bullet^1(\omega\bullet^1\Lambda)]_1^{c_0,c_0'}
\end{eqnarray*}
where expression $(\spadesuit)$ is unambiguous for the same property in dimension $n-1$.\\

$\bullet\ (L3)'$ ({\em on objects}) Let an object $c_0$ of $\C$ be given. Then quite plainly
$$
[\omega\bullet^1 id_{\alpha}]_{c_0}=\omega_{c_0}\circ [id_{\alpha}]_{c_0}=\omega_{c_0}\circ id_{\alpha_{c_0}}=
id_{\omega_{c_0}\circ\alpha_{c_0}}=id_{[\omega\bullet^1\alpha]_{c_0}}=[id_{\omega\bullet^1\alpha}]_{c_0}
$$

({\em on homs}) Let objects $c_0,c_0'$ of $\C$ be given. Then
\begin{eqnarray*}
[\omega\bullet^1 id_{\alpha}]_1^{c_0,c_0'}   &=&
(\omega_{c_0} \circ [id_{\alpha}]_1^{c_0,c_0'})\bullet^1(\omega_1^{c_0,c_0'}\circ\alpha_{c_0'})\\
   &=&[id_{\omega_{c_0} \circ \alpha}]_1^{c_0,c_0'}\bullet^1(\omega_1^{c_0,c_0'}\circ\alpha_{c_0'})  \\
   &=&(id_{\omega_{c_0} \circ \alpha_1^{c_0,c_0'}})\bullet^1(\omega_1^{c_0,c_0'}\circ\alpha_{c_0'})  \\
   &\eq{\clubsuit}&  id_{(\omega_{c_0} \circ \alpha_1^{c_0,c_0'})\bullet^1(\omega_1^{c_0,c_0'}\circ\alpha_{c_0'})}\\
   &=& id_{[\omega\bullet^1\alpha]_1^{c_0,c_0'}}  \\
   &=&[id_{\omega\bullet^1\alpha}]_1^{c_0,c_0'}
\end{eqnarray*}
where $(\clubsuit)$ holds by $(L3)'$ in dimension $n-1$.\\

$\bullet\ (L4)'$ ({\em on objects}) Let an object $c_0$ of $\C$ be given. Then
\begin{eqnarray*}
[\omega\bullet^1(\Lambda\bullet^2\Sigma)]_{c_0}
   &=& \omega_{c_0}\circ^0[\Lambda\bullet^2\Sigma]_{c_0}  \\
   &=& \omega_{c_0}\circ^0(\Lambda_{c_0}\circ^1\Sigma_{c_0})   \\
   &\eq{\heartsuit}&  (\omega_{c_0}\circ^0\Lambda_{c_0})\circ^1(\omega_{c_0}\circ^0\Sigma_{c_0})\\
   &=&[\omega\bullet^1\Lambda]_{c_0}\circ^1[\omega\bullet^1\Sigma]_{c_0}\\
   &=&[(\omega\bullet^1\Lambda)\bullet^2(\omega\bullet^1\Sigma)]_{c_0}
\end{eqnarray*}
where all equalities follow straightforward from definitions, but $(\heartsuit)$ that
is the {\em strict interchange property} of $\circ^0$ and $\circ^1$.\\

({\em on homs}) Let objects $c_0,c_0'$ of $\C$ be given. Then by definition
$$
[\omega\bullet^1(\Lambda\bullet^2\Sigma)]_1^{c_0,c_0'}=
\Big(\omega_{c_0}\circ [\Lambda\bullet^2\Sigma]_1^{c_0,c_0'}\Big)\bullet^1\Big(\omega_1^{c_0,c_0'}\circ \gamma_{c_0'}\Big)
$$
$$
=\Bigg(\omega_{c_0}\circ \Big(
   \big((\Lambda_{c_0}\circ G_1^{c_0,c_0'})\bullet^1\Sigma_1^{c_0,c_0'}\big)
   \bullet^2 \big(\Lambda_1^{c_0,c_0'}\bullet^1(F_1^{c_0,c_0'}\circ\Sigma_{c_0'} )\big)
   \Big)\Bigg)\bullet^1\Big(\omega_1^{c_0,c_0'}\circ \gamma_{c_0'}\Big)
$$
by 0-whiskering of a morphism with a 2-composition $(L4)''$ this gets
\begin{changemargin}{-15ex}{-15ex}
$$
   \Bigg(\Big(\omega_{c_0}\circ
   \big((\Lambda_{c_0}\circ G_1^{c_0,c_0'})\bullet^1\Sigma_1^{c_0,c_0'}\big)\Big)
   \bullet^2
   \Big(\omega_{c_0}\circ \big(\Lambda_1^{c_0,c_0'}\bullet^1(F_1^{c_0,c_0'}\circ\Sigma_{c_0'} )\big)
   \Big)\Bigg)\bullet^1\Big(\omega_1^{c_0,c_0'}\circ \gamma_{c_0'}\Big) \\
$$
\end{changemargin}
by 1-whiskering of a 2-morphism with a 2-composition $(R4)'$
$$
\begin{array}{c}
   \Big(\omega_{c_0}\circ
   \big((\Lambda_{c_0}\circ G_1^{c_0,c_0'})\bullet^1\Sigma_1^{c_0,c_0'}\big)\Big)\bullet^1 (\omega_1^{c_0,c_0'}\circ \gamma_{c_0'})\\
   \bullet^2 \\
   \Big(\omega_{c_0}\circ \big(\Lambda_1^{c_0,c_0'}\bullet^1(F_1^{c_0,c_0'}\circ\Sigma_{c_0'} )\big)
   \Big)\bullet^1(\omega_1^{c_0,c_0'}\circ \gamma_{c_0'})
\end{array}
$$
and by $(L4)$ and associativity of 1-composition ($LR$ whiskering property)
$$
\begin{array}{c}
   (\omega_{c_0}\circ\Lambda_{c_0}\circ G_1^{c_0,c_0'})
   \bullet^1(\omega_{c_0}\circ\Sigma_1^{c_0,c_0'})\bullet^1 (\omega_1^{c_0,c_0'}\circ \gamma_{c_0'})\\
   \bullet^2 \\
   (\omega_{c_0}\circ \Lambda_1^{c_0,c_0'})\bullet^1(\omega_{c_0}\circ F_1^{c_0,c_0'}\circ\Sigma_{c_0'} )
   \bullet^1(\omega_1^{c_0,c_0'}\circ \gamma_{c_0'})
\end{array}
$$
rearranging the terms (since interchange holds for 0-composition with a constant transformation)
$$
\begin{array}{c}
  (\omega_{c_0}\circ\Lambda_{c_0}\circ G_1^{c_0,c_0'})\bullet^1(\omega_{c_0}\circ\Sigma_1^{c_0,c_0'})\bullet^1(\omega_1^{c_0,c_0'}\circ\gamma_{c_0'})\\
  \bullet^2 \\
  {}
  (\omega_{c_0}\circ \Lambda_1^{c_0,c_0'})
  \bullet^1
  (\omega_1^{c_0,c_0'}\circ \Lambda_{c_0'})
\end{array}
$$
rearranging the terms again
$$
\begin{array}{c}
  (\omega_{c_0}\circ\Lambda_{c_0}\circ G_1^{c_0,c_0'})\bullet^1(\omega_{c_0}\circ\Sigma_1^{c_0,c_0'})\bullet^1(\omega_1^{c_0,c_0'}\circ\gamma_{c_0'})\\
  \bullet^2 \\
  {}
  (\omega_{c_0}\circ \Lambda_1^{c_0,c_0'})
  \bullet^1
  (\omega_1^{c_0,c_0'}\circ\beta_{c_0'})
  \bullet^1(E_1^{c_0,c_0'}\circ\omega_{c_0'}\circ\Sigma_{c_0'})
\end{array}
$$
just adding brackets
$$
\begin{array}{c}
  (\omega_{c_0}\circ\Lambda_{c_0}\circ G_1^{c_0,c_0'})\bullet^1\Big((\omega_{c_0}\circ\Sigma_1^{c_0,c_0'})\bullet^1(\omega_1^{c_0,c_0'}\circ\gamma_{c_0'})\Big)\\
  \bullet^2 \\
  {}\Big(
  (\omega_{c_0}\circ \Lambda_1^{c_0,c_0'})
  \bullet^1
  (\omega_1^{c_0,c_0'}\circ\beta_{c_0'})
  \Big)\bullet^1(E_1^{c_0,c_0'}\circ\omega_{c_0'}\circ\Sigma_{c_0'})
\end{array}
$$
by definition of 1-whiskering with a 3-morphism, this can be rewritten more conveniently
$$
\Big(([\omega\bullet^1\Lambda]_{c_0}\circ G_1^{c_0,c_0'})\bullet^1[\omega\bullet^1\Sigma]_1^{c_0,c_0'} \Big)
   \bullet^2
   \Big([\omega\bullet^1\Lambda]_1^{c_0,c_0'}\bullet^1(E_1^{c_0,c_0'}\circ[\omega\bullet^1\Sigma]_{c_0'})  \Big)  \\
$$
and finally, by definition of 2-composition of 3-morphisms
$$
[(\omega\bullet^1\Lambda)\bullet^2(\omega\bullet^1\Sigma)]_1^{c_0,c_0'}
$$
\vskip4ex

$\bullet\ (LR5)'$ ({\em on objects}) Let an object $c_0$ of $\C$ be given. Then
\begin{eqnarray*}
  [(\omega\bullet^1\Lambda)\bullet^1 \sigma]_{c_0} &=& [\omega\bullet^1\Lambda]_{c_0}\circ \sigma_{c_0} \\
   &=&  \omega_{c_0}\circ\Lambda_{c_0}\circ \sigma_{c_0}\\
   &=&  \omega_{c_0}\circ[\Lambda\bullet^1 \sigma]_{c_0}\\
   &=& [\omega\bullet^1(\Lambda\bullet^1 \sigma)]_{c_0}
\end{eqnarray*}
({\em on homs}) Let objects $c_0,c_0'$ of $\C$ be given. Then
$$
  [(\omega\bullet^1\Lambda)\bullet^1 \sigma]_1^{c_0,c_0'}=\qquad\qquad\qquad\qquad\qquad
$$
\vskip-4ex
\begin{eqnarray*}
   &=& ([\omega\bullet^1\alpha]_{c_0}\circ\sigma_1^{c_0,c_0'}) \bullet^1 ([\omega\bullet^1\Lambda]_1^{c_0,c_0'}\circ\sigma_{c_0'})  \\
   &=& (\omega_{c_0}\circ\alpha_{c_0}\circ\sigma_1^{c_0,c_0'}) \bullet^1 \Big(\big((\omega_{c_0}\circ\Lambda_1^{c_0,c_0'})\bullet^1(\omega_1^{c_0,c_0'}\circ\beta_{c_0'})\big)\circ\sigma_{c_0'}\Big)  \\
   &=&  (\omega_{c_0}\circ\alpha_{c_0}\circ\sigma_1^{c_0,c_0'})\bullet^1 (\omega_{c_0}\circ\Lambda_1^{c_0,c_0'}\circ\sigma_{c_0'})\bullet^1 (\omega_1^{c_0,c_0'}\circ\beta_{c_0'}\circ\sigma_{c_0'})\\
   &=& \Big(\omega_{c_0}\circ\big((\alpha_{c_0}\circ\sigma_1^{c_0,c_0'})\bullet^1(\Lambda_1^{c_0,c_0'}\circ\sigma_{c_0'})\big)\Big)\bullet^1 (\omega_1^{c_0,c_0'}\circ \beta_{c_0'}\circ\sigma_{c_0'})  \\
   &=& (\omega_{c_0}\circ [\Lambda\bullet^1\Sigma]_1^{c_0,c_0'})\bullet^1 (\omega_1^{c_0,c_0'}\circ [\beta\bullet^1\sigma]_{c_0'})  \\
   &=& [\omega\bullet^1(\Lambda\bullet^1 \sigma)]_1^{c_0,c_0'}
\end{eqnarray*}

\end{proof}

\section{0-whiskering of 3-morphisms}
In this section  we define reduced left/right
1-composition of a 3-morphism with a 1-morphism, according to the following reference diagram.
$$
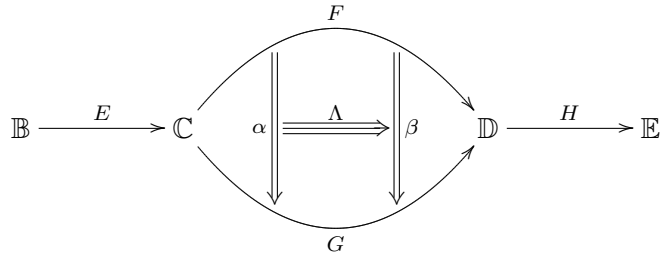
\xymatrix@C=10ex@R=14ex{
\B\ar[r]^{E}
&\C
\ar@/^8ex/[rr]_(.3){}="a1"_(.7){}="b1"^F
\ar@/_8ex/[rr]^(.3){}="a2"^(.7){}="b2"_G
&&\D\ar[r]^{H}
&\E
\ar@2"a1";"a2"_{\alpha}^{}="l1"
\ar@2"b1";"b2"^{\beta}_{}="l2"
\ar@3"l1";"l2"^{\Lambda}
}
$$

\subsection{Reduced left-composition}

The 3-morphism
$$
\xymatrix@C=3ex{
E\bullet^0\Lambda:\ E\bullet^0\alpha\ar@3[r]
&E\bullet^0\beta:\ E\bullet^0 F\ar@2[r]
&E\bullet^0 G:\ \B\ar[r]
&\D
}
$$
is given by the data below\\

$\bullet$ ({\em on objects}) For an object $b_0$ of $\B$
$$
[E\bullet^0\Lambda]_0:\ b_0\ \mapsto \ \raisebox{7ex}{
\xymatrix@R=10ex{
F(Eb_0)
\ar@/_5ex/[d]_{\alpha (Eb_0)}^{}="1"
\ar@/^5ex/[d]^{\beta (Eb_0)}_{}="2"
\\
G(Eb_0)
\ar@{}"1";"2"|(.3){}="a1"|(.7){}="a2"
\ar@2"a1";"a2"^{\Lambda(Eb_0)}
}}
$$

$\bullet$ ({\em on homs}) For objects $b_0,b_0'$ of $\B$
$$
[E\bullet^0\Lambda]_1^{b_0,b_0'}=E_1^{b_0,b_0'}(-)\bullet^0\Lambda_1^{Eb_0,Eb_0'}
$$
\begin{changemargin}{-10ex}{-10ex}
$$
\xymatrix@R=20ex@C=15ex{
\alpha(Eb_0)\circ G_1^{Eb_0,Eb_0'}\Big(E_1^{b_0,b_0'}(-)\Big)
\ar@2[r]^{\Lambda(Eb_0)\circ id}
\ar@2[d]_(.42){[E\bullet^0\alpha]_1^{b_0,b_0'}}_(.58){=\ E_1^{b_0,b_0'}\bullet^0\alpha_1^{Eb_0,Eb_0'}}
&\beta(Eb_0)\circ G_1^{Eb_0,Eb_0'}\Big(E_1^{b_0,b_0'}(-)\Big)
\ar@2[d]^(.42){[E\bullet^0\beta]_1^{b_0,b_0'}}^(.58){=\ E_1^{b_0,b_0'}\bullet^0\beta_1^{Eb_0,Eb_0'}}
\ar@{}[dl]|(.4){}="1"|(.6){}="2"
\ar@3"1";"2"_{[E\bullet^0\Lambda]_1^{b_0,b_0'}}
\\
F_1^{Eb_0,Eb_0'}\Big(E_1^{b_0,b_0'}(-)\Big)\circ\alpha(Eb_0')
\ar@2[r]_{id\circ\Lambda(Eb_0')}
&
F_1^{Eb_0,Eb_0'}\Big(E_1^{b_0,b_0'}(-)\Big)\circ\beta(Eb_0')
}
$$
\end{changemargin}

The pair $<[E\bullet^0\alpha]_0,[E\bullet^0\alpha]_1^{-,-}>$ forms indeed a 3-morphism of n-categories.
\begin{proof}
We have to show that it satisfies composition and unit axioms. Let us begin with composition, and fix a
triple $c_0,c_0',c_0''$ of $\C$. Notice that, in order to keep wide diagrams in the page we write
 $\circ^0$-composition by juxtaposition, and subscripts for transformations on objects are used.

\begin{changemargin}{-10ex}{-10ex}
$$
\xymatrix@R=20ex@C=15ex{
\alpha_{Eb_0}[EG]_1^{b_0,b_0'}[EG]_1^{b_0',b_0''}
\ar@2[r]^{\Lambda_{Eb_0} id\,id}
\ar@2[d]_{[E\bullet^0\alpha]_1^{b_0,b_0'}\,id}
&\beta_{Eb_0}  [EG]_1^{b_0,b_0'}[EG]_1^{b_0',b_0''}
\ar@2[d]^{[E\bullet^0\beta]_1^{b_0,b_0'}id}
\ar@{}[dl]|(.4){}="a1"|(.6){}="a2"
\ar@3"a1";"a2"_(.35){[E\bullet^0\Lambda]_1^{b_0,b_0'}\, id\quad}
\\
[EF]_1^{b_0,b_0'}\alpha_{Eb_0'}[EG]_1^{b_0',b_0''}
\ar@2[r]|{id\,\Lambda_{Eb_0'}id}
\ar@2[d]_{id\,[E\bullet^0\alpha]_1^{b_0',b_0''}}
&
[EF]_1^{Eb_0,Eb_0'}\beta_{Eb_0'}[EG]_1^{b_0',b_0''}
\ar@2[d]^{id\,[E\bullet^0\beta]_1^{b_0',b_0''}}
\ar@{}[dl]|(.4){}="b1"|(.6){}="b2"
\ar@3"b1";"b2"^{id\,[E\bullet^0\Lambda]_1^{b_0',b_0''}}
\\
[EF]_1^{b_0,b_0'}[EG]_1^{b_0',b_0''}\alpha_{Eb_0''}
\ar@2[r]_{id\,id\,\Lambda_{Eb_0''}}
&
[EF]_1^{Eb_0,Eb_0'}[EG]_1^{b_0',b_0''}\beta_{Eb_0''}
}
$$
\end{changemargin}
This can be written equationally
\begin{changemargin}{-10ex}{-10ex}
$$
\begin{array}{c}
\Big[
\Big(E_1^{b_0,b_0'}\bullet^0\Lambda_1^{Eb_0,Eb_0'}\Big) \circ \Big(E_1^{b_0',b_0''}\bullet^0G_1^{Eb_0',Eb_0''}\Big)
\Big]\bullet^1\Big[
\Big(E_1^{b_0,b_0'}\bullet^0F_1^{Eb_0,Eb_0'}\Big)\circ \Big(E_1^{b_0',b_0''}\bullet^0\beta_1^{Eb_0',Eb_0''}\Big)
\Big]
\\
\phantom{\Big[}\bullet^2\phantom{\Big]}
\\
\Big[
\Big(E_1^{b_0,b_0'}\bullet^0\alpha_1^{Eb_0,Eb_0'}\Big)\circ \Big(E_1^{b_0',b_0''}\bullet^0G_1^{Eb_0',Eb_0''}\Big)
\Big]\bullet^1\Big[
\Big(E_1^{b_0,b_0'}\bullet^0F_1^{Eb_0,Eb_0'}\Big)\circ \Big(E_1^{b_0',b_0''}\bullet^0\Lambda_1^{Eb_0',Eb_0''}\Big)
\Big]
\end{array}
$$
\end{changemargin}
by product interchange before $\circ^0$-composition this turns to be
\begin{changemargin}{-10ex}{-10ex}
$$
\begin{array}{c}
\Big[
\Big(E_1^{b_0,b_0'}\times E_1^{b_0',b_0''}\Big) \bullet^0 \Big(\Lambda_1^{Eb_0,Eb_0'}\circ G_1^{Eb_0',Eb_0''}\Big)
\Big]\bullet^1\Big[
\Big(E_1^{b_0,b_0'}\times E_1^{b_0',b_0''}\Big)\bullet^0 \Big(F_1^{Eb_0,Eb_0'}\circ\beta_1^{Eb_0',Eb_0''}\Big)
\Big]
\\
\phantom{\Big[}\bullet^2\phantom{\Big]}
\\
\Big[
\Big(E_1^{b_0,b_0'}\times E_1^{b_0',b_0''}\Big)\bullet^0 \Big(\alpha_1^{Eb_0,Eb_0'}\circ G_1^{Eb_0',Eb_0''}\Big)
\Big]\bullet^1\Big[
\Big(E_1^{b_0,b_0'}\times E_1^{b_0',b_0''}\Big)\bullet^0 \Big(F_1^{Eb_0,Eb_0'}\circ\Lambda_1^{Eb_0',Eb_0''}\Big)
\Big]
\end{array}
$$
\end{changemargin}
by {\em whiskering interchange property} $(LR)$ this gives
\begin{changemargin}{-10ex}{-10ex}
$$
\begin{array}{c}
\Big(E_1^{b_0,b_0'}\times E_1^{b_0',b_0''}\Big)
\bullet^0\Big[
\Big(\Lambda_1^{Eb_0,Eb_0'}\circ G_1^{Eb_0',Eb_0''}\Big)\bullet^1 \Big(F_1^{Eb_0,Eb_0'}\circ\beta_1^{Eb_0',Eb_0''}\Big)
\Big]
\\
\phantom{\Big[}\bullet^2\phantom{\Big]}
\\
\Big(E_1^{b_0,b_0'}\times E_1^{b_0',b_0''}\Big)
\bullet^0\Big[
\Big(\alpha_1^{Eb_0,Eb_0'}\circ G_1^{Eb_0',Eb_0''})\bullet^1 \Big(F_1^{Eb_0,Eb_0'}\circ\Lambda_1^{Eb_0',Eb_0''}\Big)
\Big]
\end{array}
$$
\end{changemargin}
hence by $(L4)''$
\begin{changemargin}{-15ex}{-10ex}
$$
\Big(E_1^{b_0,b_0'}\times E_1^{b_0',b_0''}\Big)\bullet^0
\left(\begin{array}{c}
\Big(\Lambda_1^{Eb_0,Eb_0'}\circ G_1^{Eb_0',Eb_0''}\Big)\bullet^1 \Big(F_1^{Eb_0,Eb_0'}\circ\beta_1^{Eb_0',Eb_0''}\Big)
\\
\bullet^2
\\
\Big(\alpha_1^{Eb_0,Eb_0'}\circ G_1^{Eb_0',Eb_0''})\bullet^1 \Big(F_1^{Eb_0,Eb_0'}\circ\Lambda_1^{Eb_0',Eb_0''}\Big)
\end{array}\right)
$$

\end{changemargin}
By composition functoriality of $\Lambda$, the second row changes to \mbox{$(-\circ-)\bullet^0\Lambda_1^{Eb_0,Eb_0''}$}
thus giving by associativity of 0-whiskering
$$
\Big(E_1^{b_0,b_0'}\circ E_1^{b_0',b_0''}\Big)\bullet^0\Lambda_1^{Eb_0,Eb_0''}=E_1^{b_0,b_0''}\bullet^0\Lambda_1^{Eb_0,Eb_0''}
$$
where the last equality follows from composition axiom for 1-morphisms.\\

Turning to units axiom, let an object $b_0$ of $\B$ be given. Then
$$
u(b_0)\bullet^0E_1^{b_0,b_0}\bullet^0\Lambda_1^{Eb_0,Eb_0}=
u(Eb_0)\bullet^0\Lambda_1^{Eb_0,Eb_0}=id_{[\Lambda(Eb_0)]}=id_{[(E\bullet^0\Lambda)(b_0)]}
$$
where the first expression is unambiguous for $\bullet^0$-associativity $(L2)''$, first equality holds
 by units axiom for 1-morphism $E$, second by units axiom for 3-morphism $\Lambda$, last is the definition.
\end{proof}

\subsection{Reduced right-composition}
The 3-morphism
$$
\xymatrix@C=3ex{
\Lambda\bullet^0H:\ \alpha\bullet^0\ar@3[r]
&\beta\bullet^0H:\ F\bullet^0 H\ar@2[r]
&G\bullet^0 H:\ \C\ar[r]
&\E
}
$$
is given by the data below\\

$\bullet$ ({\em on objects}) For an object $c_0$ of $\C$
$$
[\Lambda\bullet^0H]_0:\ c_0\ \mapsto \ \raisebox{7ex}{
\xymatrix@R=10ex{
H(Fc_0)
\ar@/_5ex/[d]_{H(\alpha c_0)}^{}="1"
\ar@/^5ex/[d]^{H(\beta c_0)}_{}="2"
\\
H(Gc_0)
\ar@{}"1";"2"|(.3){}="a1"|(.7){}="a2"
\ar@2"a1";"a2"^{H(\Lambda c_0)}
}}
$$

$\bullet$ ({\em on homs}) For objects $c_0,c_0'$ of $\C$
$$
[\Lambda\bullet^0H]_1^{c_0,c_0'}=\Lambda_1^{c_0,c_0'}\bullet^0 H_1^{Fc_0,Gc_0'}
$$
\begin{changemargin}{-15ex}{-10ex}
$$
\xymatrix@R=20ex@C=15ex{
H_1^{Fc_0,Gc_0}(\alpha_{c_0})\circ H_1^{Gc_0,Gc_0'}(G_1^{c_0,c_0'}(-))
\ar@2[r]^{H_1^{Fc_0,Gc_0}(\Lambda_{c_0})\circ id}
\ar@2[d]_{H_1^{Fc_0,Gc_0'}(\alpha_1^{c_0,c_0'})}
&
H_1^{Fc_0,Gc_0}(\beta_{c_0})\circ H_1^{Gc_0,Gc_0'}(G_1^{c_0,c_0'}(-))
\ar@2[d]^{H_1^{Fc_0,Gc_0'}(\beta_1^{c_0,c_0'})}
\ar@{}[dl]|(.4){}="1"|(.6){}="2"
\ar@3"1";"2"_{[\Lambda\bullet^0H]_1^{c_0,c_0'}}
\\
H_1^{Fc_0,Fc_0'}(F_1^{c_0,c_0'})\circ H_1^{Fc_0',Gc_0'}(\alpha_{c_0'})
\ar@2[r]_{id \circ H_1^{Fc_0',Gc_0'}(\Lambda_{c_0'})}
&
H_1^{Fc_0,Fc_0'}(F_1^{c_0,c_0'})\circ H_1^{Fc_0',Gc_0'}(\beta_{c_0'})
}
$$
\end{changemargin}
The pair $<[\Lambda\bullet^0H]_0,[\Lambda\bullet^0H]_1^{-,-}>$ forms indeed a 3-morphism of n-categories. The proof
is a straightforward variation of the proof for reduced right-composition above, hence it is omitted.

\subsection{Properties}
As we did in describing the sesqui-categorical structure for homs in $n\mathbf{Cat}$, we use again a
{\em left-and-right} approach to describe properties of the 0-whiskering of
 of a 3-morphism with a morphism. Let us consider the diagram
$$
\xymatrix@C=10ex@R=14ex{
\A\ar[r]^{E'}
&\B\ar[r]^{E}
&\C
\ar@/^8ex/[rr]_(.2){}="a1"_(.5){}="b1"_(.8){}="c1"^F
\ar@/_8ex/[rr]^(.2){}="a2"^(.5){}="b2"^(.8){}="c2"_G
&&\D\ar[r]^{H}
&\E\ar[r]^{H'}
&\F
\ar@2"a1";"a2"_{\alpha}^{}="l1"
\ar@2"b1";"b2"^(.75){\beta}_{}="l2"^{}="s1"
\ar@2"c1";"c2"^{\gamma}_{}="s2"
\ar@3"l1";"l2"^{\Lambda}
\ar@3"s1";"s2"^{\Sigma}
}
$$
as a reference for the following
\begin{Proposition}[2-composition (i.e. vertical) composition of  3-morphisms w.r.t. (reduced) 0-composition with a (1-)morphism]
$$
\begin{array}{clcl}
  (L1)'' & id_{\C}\bullet_L^0 \Lambda =\Lambda & (R1)'' & \Lambda\bullet_R^0 id_{\D} =\Lambda\\
  (L2)'' & (E'\bullet_L^0 E)\bullet_L^0 \Lambda=E' \bullet_L^0  (E\bullet_L^0\Lambda) &
  (R2)'' & \Lambda\bullet^0_R (H\bullet_R^0 H') =(\Lambda \bullet_R^0 H) \bullet_R^0 H'\\
  (L3)'' & E\bullet_L^0 id_{\alpha}= id_{E\bullet_L^0\alpha}&
  (R3)'' &id_{\alpha}\bullet_R^0H =id_{\alpha\bullet_R^0 H}\\
  (L4)'' & E\bullet_L^0(\Lambda\bullet^2\Sigma)=(E\bullet_L^0\Lambda)\bullet^2(E\bullet_L^0\Sigma)
 & (R4)''& (\Lambda\bullet^2\Sigma)\bullet_R^0 H =(\Lambda\bullet_R^0 H)\bullet^2(\Sigma\bullet_R^0 H)
\end{array}
$$
$$
(LR5)''\quad (E\bullet_L^0\Lambda)\bullet_R^0 H= E\bullet_L^0(\Lambda\bullet_R^0 H)
$$
\end{Proposition}
\begin{proof}
We prove statements $(L1)''$ to $(L4)''$ and $(LR5)$. Proofs of statements $(R1)''$ to $(R4)''$
is similar, hence it is omitted.\\

$\bullet$ $(L1)''$ ({\em on objects}) Let an object $c_0$ of $\C$ be given. Then
$$
[id_{\C}\bullet^0\Lambda]_{c_0}=\Lambda(id_{\C}(c_0))=\Lambda_{c_0}
$$

({\em on homs}) Let objects $c_0,c_0'$ of $\C$ be given. Then
$$
[id_{\C}\bullet^0\Lambda]_1^{c_0c_0'}=[id_{\C}]_1^{c_0,c_0'}\bullet^0\Lambda_1^{c_0c_0'}
=id_{\C_1(c_0,c_0')}\bullet^0\Lambda_1^{c_0c_0'}=\Lambda_1^{c_0c_0'}
$$

$\bullet$ $(L2)''$ (\em on objects) Let an object $a_0$ of $\A$ be given. Then
$$
[(E'\bullet^0 E)\bullet^0 \Lambda]_{a_0}=\Lambda_{(E'\bullet^0 E)a_0}=\Lambda_{E(E'a_0)}=
[E\bullet^0\Lambda]_{E'a_0}=[E' \bullet^0  (E\bullet^0\Lambda)]_{a_0}
$$

({\em on homs}) Let objects $a_0,a_0'$ of $\A$ be given. Then
\begin{eqnarray*}
  [(E'\bullet^0 E)\bullet^0 \Lambda]_1^{a_0,a_0'}
   &=& [E'\bullet^0 E]_1^{a_0,a_0'}\bullet^0\Lambda_1^{E(E'a_0),E(E'a_0')} \\
   &=& {E'_1}^{a_0,a_0'}\bullet^0 E_1^{E'a_0,E'a_0'}\bullet^0\Lambda_1^{E(E'a_0),E(E'a_0')} \\
   &=& {E'_1}^{a_0,a_0'}\bullet^0[E\bullet^0\Lambda]_1^{E'a_0,E'a_0'} \\
   &=& [E' \bullet^0  (E\bullet^0\Lambda)]_1^{a_0,a0'}
\end{eqnarray*}

$\bullet$ $(L3)''$ ({\em on objects}) Let an object $b_0$ of $\B$ be given. Then
$$
[E\bullet^0 id_{\alpha}]_{b_0}=[id_{\alpha}]_{Eb_0}=
id_{\alpha_{Eb_0}}=id_{[E\bullet^0\alpha]_{b_0}} = [id_{E\bullet^0\alpha}]_{b_0}
$$

({\em on homs}) Let objects $b_0,b_0'$ of $\B$ be given. Then
\begin{eqnarray*}
[E\bullet^0 id_{\alpha}]_1^{b_0,b_0'}
   &=&  E_1^{b_0,b_0'}\bullet^0 [id_{\alpha}]_1^{Eb_0,Eb_0'}\\
   &=&  E_1^{b_0,b_0'}\bullet^0 id_{\alpha_1^{Eb_0,Eb_0'}}\\
   &\eq{\clubsuit}&  id_{E_1^{b_0,b_0'}\bullet^0\alpha_1^{Eb_0,Eb_0'}} \\
   &=&  id_{[E\bullet^0\alpha]_1^{b_0,b_0'}}\\
   &=&  [id_{E\bullet^0\alpha}]_1^{b_0,b_0'}
\end{eqnarray*}
where $(\clubsuit)$ holds for the same property in dimension $n-1$.\\

$\bullet$ $(L4)''$ ({\em on objects}) Let an object $b_0$ of $\B$ be given. Then
\begin{eqnarray*}
  [E\bullet^0(\Lambda\bullet^2\Sigma)]_{b_0}
   &=&  [(\Lambda\bullet^2\Sigma]_{Eb_0}\\
   &=&  \Lambda_{Eb_0}\circ^1\Sigma_{Eb_0}\\
   &=& [E\bullet^0\Lambda]_{b_0}\circ^1[E\bullet^0\Sigma]_{b_0} \\
   &=& [(E\bullet^0\Lambda)\bullet^2(E\bullet^0\Sigma)]_{b_0}
\end{eqnarray*}

({\em on homs}) Let objects $b_0,b_0'$ of $\B$ be given. Then
\begin{eqnarray*}
[E\bullet^0(\Lambda\bullet^2\Sigma)]_1^{b_0,b_0'}
   &=&   E_1^{b_0,b_0'}\bullet^0[\Lambda\bullet^2\Sigma]_1^{Eb_0,Eb_0'}\\
   &=&   E_1^{b_0,b_0'}\bullet^0\left(
\begin{array}{c}
   (\Lambda_{Eb_0}\circ G_1^{Eb_0,Eb_0'})\bullet^1 \Sigma_1^{Eb_0,Eb_0'}\\
   \bullet^2\\
   \Lambda_1^{Eb_0,Eb_0'}\bullet^1(F_1^{Eb_0,Eb_0'}\circ\Sigma_{Eb_0'})\\
\end{array}
   \right)\\
   &\eq{\heartsuit}&   \left(
\begin{array}{c}
   E_1^{b_0,b_0'}\bullet^0\big((\Lambda_{Eb_0}\circ G_1^{Eb_0,Eb_0'})\bullet^1 \Sigma_1^{Eb_0,Eb_0'}\big)\\
   \bullet^2\\
   E_1^{b_0,b_0'}\bullet^0\big(\Lambda_1^{Eb_0,Eb_0'}\bullet^1(F_1^{Eb_0,Eb_0'}\circ\Sigma_{Eb_0'})\big)\\
\end{array}
   \right)\\
   &\eq{\spadesuit}&
\left(
\begin{array}{c}
\big(\Lambda_{Eb_0}\circ(E_1^{b_0,b_0'}\bullet^0G_1^{Eb_0,Eb_0'})\big)\bullet^1 (E_1^{b_0,b_0'}\bullet^0\Sigma_1^{Eb_0,Eb_0'})\\
\bullet^2\\
(E_1^{b_0,b_0'}\bullet^0\Lambda_1^{Eb_0,Eb_0'})\bullet^1\big((E_1^{b_0,b_0'}\bullet^0F_1^{Eb_0,Eb_0'})\circ\Sigma_{Eb_0'}\big)\\
\end{array}
\right)\\
   &=&
\left(
\begin{array}{c}
\big([E\bullet^0\Lambda]_{b_0}\circ[E\bullet^0G]_1^{b_0,b_0'}\big)\bullet^1 [E\bullet^0\Sigma]_1^{b_0,b_0'}\\
\bullet^2\\
{}[E\bullet^0\Lambda]_1^{b_0,b_0'}\bullet^1\big([E\bullet^0F]_1^{b_0,b_0'}\circ[E\bullet^0\Sigma]_{b_0'}\big)\\
\end{array}
\right)\\
   &=& [(E\bullet^0\Lambda)\bullet^2(E\bullet^0\Sigma)]_1^{b_0,b_0'}
\end{eqnarray*}
where $(\heartsuit)$ holds by the same property in dimension $n-1$, and $(\clubsuit)$ holds by
{\em whiskering interchange property}.

$\bullet$ $(LR5)''$ ({\em on objects}) Let an object $b_0$ of $\B$ be given. Then
\begin{eqnarray*}
   [(E\bullet^0\Lambda)\bullet^0 H]_{b_0}
   &=&  H([E\bullet^0\Lambda]_{b_0}=\\
   &=&  H(\Lambda_{Eb_0})\\
   &=&  [\Lambda\bullet^0 H]_{Eb_0}\\
   &=&  [E\bullet^0(\Lambda\bullet^0 H)]_{b_0}
\end{eqnarray*}

({\em on homs}) Let objects $b_0,b_0'$ of $\B$ be given. Then
\begin{eqnarray*}
   [(E\bullet^0\Lambda)\bullet^0 H]_1^{b_0,b_0'}
   &=&  [E\bullet^0\Lambda]_1^{b_0,b_0'}\bullet^0H_1^{F(Eb_0),G(Eb_0')}\\
   &=&  E_1^{b_0,b_0'}\bullet^0\Lambda_1^{Eb_0,Eb_0'} \bullet^0 H_1^{F(Eb_0),G(Eb_0')}\\
   &=&  E_1^{b_0,b_0'}\bullet^0[\Lambda\bullet^0 H]_1^{Eb_0,Eb_0'}\\
   &=&  [E\bullet^0(\Lambda\bullet^0 H)]_1^{b_0,b_0'}
\end{eqnarray*}
\end{proof}

Before switching to next section, let us give a last property that express at once functoriality
of left and right 0-composition with a morphism. To this end, let us be given also  2-morphisms
$\omega:\ M\Rightarrow F$ and $\sigma:\ G\Rightarrow N$, as represented in the diagram below
$$
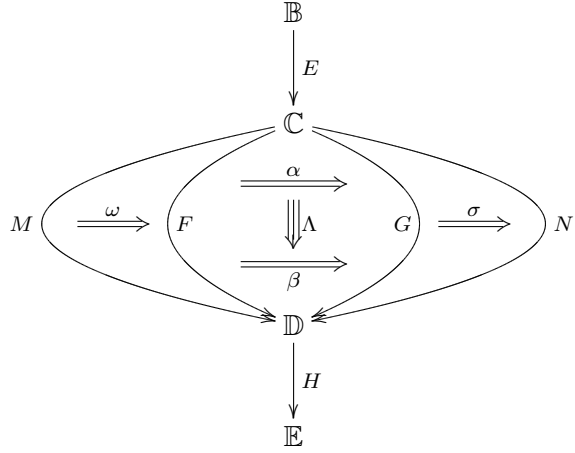
\xymatrix@R=6ex{
\B\ar[d]^{E}
\\
\C
\ar@/_20ex/[dd]_{M}|{}="1"
\ar@/_10ex/[dd]^{F}="2"|(.3){}="u2"|(.7){}="d2"
\ar@/^10ex/[dd]_{G}="3"|(.3){}="u3"|(.7){}="d3"
\ar@/^20ex/[dd]^{N}|{}="4"
\\
\\
\ar@{}"1";"2"|(.25){}="11"|(.75){}="22"
\ar@{=>}"11";"22"^{\omega}
\ar@{}"u2";"u3"|(.25){}="U2"|(.75){}="U3"
\ar@{=>}"U2";"U3"^{\alpha}_{}="l1"
\ar@{}"d2";"d3"|(.25){}="D2"|(.75){}="D3"
\ar@{=>}"D2";"D3"_{\beta}^{}="l2"
\ar@{}"3";"4"|(.25){}="33"|(.75){}="44"
\ar@{=>}"33";"44"^{\sigma}
\ar@{}"l1";"l2"|(.2){}="L1"|(.8){}="L2"
\ar@3"L1";"L2"^{\Lambda}
\D\ar[d]^{H}
\\
\E}
$$
Left/right 0-composition of a 3-morphism with a morphism satisfies also
the following property that relates 0-whiskering w.r.t. 1-whiskering:
\begin{Proposition}[Whiskering interchange property]
$$
(LRW)\qquad E\bullet^0_L\big(\omega\bullet^1_L\Lambda\bullet^1_R\sigma\big)\bullet^0_R H=
(E\bullet^0_L\omega\bullet^0_R H)\bullet^1_L(E\bullet^0_L\Lambda\bullet^0_R H)\bullet^1_R(E\bullet^0_L\sigma\bullet^0_R H)
$$
\end{Proposition}
\begin{proof}
Without loss of generality it suffices to prove the following two equalities:
\begin{eqnarray*}
   &(LRW)^1& E\bullet^0_L\big(\omega\bullet^1_L\Lambda\big)=(E\bullet^0_L\omega)\bullet^1_L(E\bullet^0_L\Lambda)\\
   &(LRW)^2& E\bullet^0_L\big(\Lambda\bullet^1_R\sigma\big)=(E\bullet^0_L\Lambda)\bullet^1_R(E\bullet^0_L\sigma)
\end{eqnarray*}

$\bullet$ $(LRW)^1$ ({\em on objects})  Let an object $b_0$ of $\B$ be given. Then
\begin{eqnarray*}
  [E\bullet^0\big(\omega\bullet^1\Lambda\big)]_{b_0}&=& [\omega\bullet^1\Lambda]_{Eb_0} \\
   &=& \omega_{Eb_0}\circ \Lambda_{Eb_0} \\
   &=&  [E\bullet^0\omega]_{b_0}\circ [E\bullet^0\Lambda]_{b_0}\\
   &=& [(E\bullet^0\omega)\bullet^1(E\bullet^0\Lambda)]_{b_0}
\end{eqnarray*}

({\em on homs}) Let objects $b_0,b_0'$ of $\B$ be given. Then

$$
   \left[E\bullet^0\big(\omega\bullet^1\Lambda\big)\right]_1^{b_0,b_0'}=\qquad\qquad\qquad
$$
\begin{eqnarray*}
   &=& E_1^{b_0,b_0'}\bullet^0[\omega\bullet^1\Lambda]_1^{Eb_0,Eb_0'} \\
   &=& E_1^{b_0,b_0'}\bullet^0\Big(\big(\omega_{Eb_0}\circ \Lambda_1^{Eb_0,Eb_0'}\big)
   \bullet^1 \big(\omega_1^{Eb_0,Eb_0'}\circ\beta_{Eb_0'} \big)\Big)\\
   &\eq{\clubsuit}&
\Big( E_1^{b_0,b_0'}\bullet^0 \big(\omega_{Eb_0}\circ \Lambda_1^{Eb_0,Eb_0'}\big)\Big)
\bullet^1
\Big( E_1^{b_0,b_0'}\bullet^0 \big(\omega_1^{Eb_0,Eb_0'}\circ\beta_{Eb_0'} \big)\Big)
    \\
   &=&
\Big( E_1^{b_0,b_0'}\bullet^0 \Lambda_1^{Eb_0,Eb_0'}\bullet^0 \big( \omega_{Eb_0}\circ-\big)\Big)
\bullet^1
\Big( E_1^{b_0,b_0'}\bullet^0 \omega_1^{Eb_0,Eb_0'}\bullet^0 \big(- \circ\beta_{Eb_0'} \big)\Big)
\\
   &=&
\Big(\omega_{Eb_0}\circ
\big(E_1^{b_0,b_0'}\bullet^0\Lambda_1^{Eb_0,Eb_0'}\big)\Big)
\bullet^1
\Big(\big(E_1^{b_0,b_0'}\bullet^0\omega_1^{Eb_0,Eb_0'}\big)\circ \beta_{Eb_0'}\Big)
\\
   &=&
\Big([E\bullet^0\omega]_{b_0}\circ [E\bullet^0\Lambda]_1^{b_0,b_0'}\Big)
\bullet^1
\Big([E\bullet^0\omega]_1^{b_0,b_0'}\circ [E\bullet^0\beta]_{b_0'}\Big)
   \\
   &=& [(E\bullet^0\omega)\bullet^1(E\bullet^0\Lambda)]_1^{b_0,b_0'}
\end{eqnarray*}
where equality $(\clubsuit)$ holds by same property in dimension $n-1$, and the following by associativity
of 0-composition.\\

$\bullet$ $(LRW)^2$ ({\em on objects}) Let an object $b_0$ of $\B$ be given. Then
\begin{eqnarray*}
   [E\bullet^0\big(\Lambda\bullet^1\sigma\big)]_{b_0}
   &=&  [\Lambda\bullet^1\sigma]_{Eb_0}\\
   &=&  \Lambda_{Eb_0}\circ\sigma_{Eb_0}\\
   &=&  [E\bullet^0\Lambda]_{b_0}\circ[E\bullet^0\sigma]_{b_0}\\
   &=& [(E\bullet^0\Lambda)\bullet^1(E\bullet^0\sigma)]_{b_0}
\end{eqnarray*}

({\em on homs}) Let objects $b_0,b_0'$ of $\B$ be given. Then
$$
[E\bullet^0\big(\Lambda\bullet^1\sigma\big)]_1^{b_0,b_0'}
$$
\begin{eqnarray*}
   &=& E_1^{b_0,b_0'}\bullet^0[\Lambda\bullet^1\sigma]_1^{Eb_0,Eb_0'}\\
   &=& E_1^{b_0,b_0'}\bullet^0
\Big(
\big(\alpha_{Eb_0}\circ\sigma_1^{Eb_0,Eb_0'}\big)
\bullet^1
\big(\Lambda_1^{Eb_0,Eb_0'}\circ\sigma_{Eb_0'}\big)
\Big)
   \\
   &=&
\Big(E_1^{b_0,b_0'}\bullet^0\sigma_1^{Eb_0,Eb_0'}\big(\alpha_{Eb_0}\circ-\big)\Big)
\bullet^1
\Big(E_1^{b_0,b_0'}\bullet^0\Lambda_1^{Eb_0,Eb_0'}\big(-\circ\sigma_{Eb_0'}\big)\Big)
    \\
   &=&
\Big(E_1^{b_0,b_0'}\bullet^0\big(\alpha_{Eb_0}\circ\sigma_1^{Eb_0,Eb_0'}\big)\Big)
\bullet^1
\Big(E_1^{b_0,b_0'}\bullet^0\big(\Lambda_1^{Eb_0,Eb_0'}\circ\sigma_{Eb_0'}\big)\Big)
    \\
   &=&
\Big(\alpha_{Eb_0}\circ\big(E_1^{b_0,b_0'}\bullet^0\sigma_1^{Eb_0,Eb_0'}\big)\Big)
\bullet^1
\Big(\big(E_1^{b_0,b_0'}\bullet^0\Lambda_1^{Eb_0,Eb_0'}\big)\circ\sigma_{Eb_0'}\Big)
    \\
   &=&
\Big([E\bullet^0\alpha]_{b_0}\circ[E\bullet^0\sigma]_1^{b_0,b_0'}\Big)
\bullet^1
\Big([E\bullet^0\Lambda]_1^{b_0,b_0'}\circ[E\bullet^0\sigma]_{b_0'}\Big)
   \\
   &=& [(E\bullet^0\Lambda)\bullet^1_(E\bullet^0\sigma)]_1^{b_0,b_0'}
\end{eqnarray*}

\end{proof}

\section{Dimension raising 0-composition of 2-morphisms}

Let two 0-intersecting 2-morphisms of n-categories be given.
$$
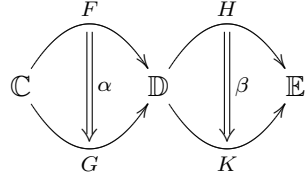
\xymatrix@C=8ex{
\C
\ar@/^5ex/[r]^{F}_{}="a1"
\ar@/_5ex/[r]_{G}^{}="a2"
\ar@{=>}"a1";"a2"^{\alpha}
&
\D
\ar@/^5ex/[r]^{H}_{}="b1"
\ar@/_5ex/[r]_{K}^{}="b2"
\ar@{=>}"b1";"b2"^{\beta}
&
\E
}
$$
It is easy to verify that in general
$$
\alpha \setminus \beta:=(F\bullet^0 \beta)\bullet^1(\alpha\bullet^0 K)\neq (\alpha\bullet^0 H)\bullet^1(G\bullet^1 \beta)=:\alpha / \beta
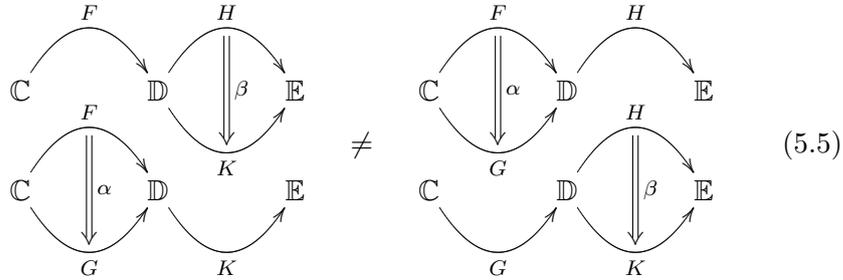
$$
\begin{equation}\label{eq:domcodstar}
\raisebox{5ex}{\xymatrix@C=8ex{
\C
\ar@/^5ex/[r]^{F}
&
\D
\ar@/^5ex/[r]^{H}_{}="b1"
\ar@/_5ex/[r]_{K}^{}="b2"
\ar@{=>}"b1";"b2"^{\beta}
&
\E
\\
\C
\ar@/^5ex/[r]^{F}_{}="a1"
\ar@/_5ex/[r]_{G}^{}="a2"
\ar@{=>}"a1";"a2"^{\alpha}
&
\D
\ar@/_5ex/[r]_{K}
&
\E
}}
\quad\neq\quad
\raisebox{5ex}{\xymatrix@C=8ex{
\C
\ar@/^5ex/[r]^{F}_{}="a1"
\ar@/_5ex/[r]_{G}^{}="a2"
\ar@{=>}"a1";"a2"^{\alpha}
&
\D
\ar@/^5ex/[r]^{H}
&
\E
\\
\C
\ar@/_5ex/[r]_{G}
&
\D
\ar@/^5ex/[r]^{H}_{}="b1"
\ar@/_5ex/[r]_{K}^{}="b2"
\ar@{=>}"b1";"b2"^{\beta}
&
\E
}}
\end{equation}
More interestingly they constitute a 3-morphism
$$
\xymatrix{
\alpha * \beta: \alpha \setminus \beta\ar@3[r]&\alpha/ \beta
}
$$
In fact, for every object $c_0$ of $\C$ one defines
$$
[\alpha * \beta]_0 : c_0\quad\mapsto\qquad \raisebox{8ex}
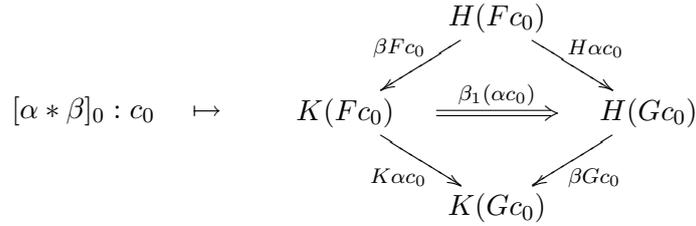
{\xymatrix@R=4ex@C=3ex{
&H(Fc_0)
\ar[dl]_{\beta Fc_0}
\ar[dr]^{H \alpha c_0}
\\
K(Fc_0)\ar@{}[rr]|(.3){}="1"|(.7){}="2"\ar@2"1";"2"^{\beta_1(\alpha c_0)}
\ar[dr]_{K\alpha c_0}
&&H(Gc_0)
\ar[dl]^{\beta Gc_0}
\\
& K(Gc_0)}}
$$
Moreover for every pair of objects $c_0,c_0'$ of $\C$ one defines
$$
[\alpha*\beta]_1^{c_0,c_0'}=\alpha_1^{c_0,c_0'}*\beta_1^{Fc_0,Gc_0'}
$$
\begin{Claim}
The pair
$$
<[\alpha*\beta]_0,[\alpha*\beta]_1^{-,-}>
$$
is indeed a 3-morphism of n-categories.
\end{Claim}
\begin{proof}
We must show that domain and codomain of $[\alpha*\beta]_1^{-,-}$ are compatible with
those of the definition of
3-morphism. Moreover the pair above must satisfy unit and composition  axioms.
In order to prove the first fact, we write the diagram that represents the 3-morphism of (n-1)categories
$$
\xymatrix{
[\alpha*\beta]_1^{c_0,c_0'}:\ \alpha_1^{c_0,c_0'}\setminus \beta_1^{Fc_0,Gc_0'}
\ar@3[r]&\alpha_1^{c_0,c_0'}/ \beta_1^{Fc_0,Gc_0'}}
$$
i.e. the composition
$$
\xymatrix@R=4ex@C=-3ex{
&[c_0,c_0']
\ar[dl]_{F_1}
\ar[dr]^{G_1}
\\
[Fc_0,Fc_0']
\ar[dr]_{-\circ \alpha c_0'}
&&
[Gc_0,Gc_0']
\ar[dl]^{\alpha c_0\circ-}
\ar@{}[ll]|(.3){}="a1"|(.7){}="a2"
\ar@2"a1";"a2"_{\alpha_1^{c_0,c_0'}}
\\
&[Fc_0,Gc_0']
\ar[dl]_{H_1}
\ar[dr]^{K_1}
\\
[H(Fc_0),H(Gc_0')]
\ar[dr]_{-\circ \beta Gc_0'}
&&
[K(Fc_0),K(Gc_0')]
\ar[dl]^{\beta Fc_0\circ-}
\ar@{}[ll]|(.3){}="b1"|(.7){}="b2"
\ar@2"b1";"b2"_{\beta_1^{Fc_0,Gc_0'}}
\\
&[H(Fc_0),K(Gc_0')]
}
$$
Its domain is computed below
$$
\xymatrix@R=10ex@C=8ex{
&[c_0,c_0']
\ar[dl]_{[FH]_1}
\ar[d]|{[GH]_1}
\ar[dr]^{\alpha c_0 \circ G_1(-)}
\\
[H(Fc_0),H(Fc_0')]
\ar[d]_{-\circ H\alpha c_0'}
&[H(Gc_0),H(Gc_0')]\ar@2[l]_{[\alpha H]_1^{c_0,c_0'}}
\ar[dl]^{H\alpha c_0\circ-}
&[Fc_0,Gc_0']
\ar[dl]_{H_1}
\ar[d]^{H_1}
\\
[H(Fc_0),H(Gc_0')]
\ar[dr]_{-\circ\beta Gc_0'}
&[H(Fc_0),H(Gc_0')]
\ar[d]_{-\circ\beta Gc_0'}
&[K(Fc_0),K(Gc_0')]\ar@2[l]_{\beta_1^{Fc_0,Gc_0'}}
\ar[dl]^{\beta Fc_0\circ-}
\\
&[H(Fc_0),K(Gc_0')]}
$$
Now, by functoriality w.r.t 0-composition, with constant left composite one has
$$
(\alpha c_0\circ-)\bullet^0 \beta_1^{Fc_0,Gc_0'}= \Big(K_1^{Gc_0,Gc_0'}\bullet^0 (\beta_1(\alpha c_0)\circ-)\Big)
\bullet^1\Big(\beta_1^{Gc_0,Gc_0'}\bullet^0(H\alpha c_0\circ -) \Big)
$$
and by definition of $*$-composition on objects,
$$
=\Big(K_1^{Gc_0,Gc_0'}\bullet^0 ([\alpha*\beta]c_0\circ-)\Big)
\bullet^1\Big(\beta_1^{Gc_0,Gc_0'}\bullet^0(H\alpha c_0\circ -) \Big)
$$
Hence we can redraw the domain
$$
\xymatrix@R=10ex@C=8ex{
&[c_0,c_0']
\ar[dl]_{[FH]_1}
\ar[d]|{[GH]_1}
\ar[dr]^{[GK]_1}
\\
[H(Fc_0),H(Fc_0')]
\ar[d]|{-\circ H\alpha c_0'}
&[H(Gc_0),H(Gc_0')]\ar@2[l]_{[\alpha H]_1^{c_0,c_0'}}
\ar[d]|{-\circ \beta Gc_0'}
\ar[dl]^{H\alpha c_0\circ-}
&[K(Gc_0),K(Gc_0')]\ar@2[l]_{[G\beta]_1^{c_0,c_0'}}
\ar[dl]|{\beta Gc_0\circ-}
\ar[d]^{K\alpha c_0\circ-}
\\
[H(Fc_0),H(Gc_0')]
\ar[dr]_{-\circ\beta Gc_0'}
&[H(Gc_0),K(Gc_0')]
\ar[d]_{-\circ\beta Gc_0'}
&[K(Fc_0),K(Gc_0')]\ar@2[l]_{[\alpha*\beta]c_0\circ -}
\ar[dl]^{\beta Fc_0\circ-}
\\
&[H(Fc_0),K(Gc_0')]}
$$
And this completes the domain-part. Concerning the codomain, the calculation is similar, but on the left
side of diagrams.\\

Turning to functoriality axioms,  we start with functoriality w.r.t. units. To this end, let us suppose an object $c_0$ of
$\C$ been given. Then
\begin{eqnarray*}
  u(c_0)\bullet^0[\alpha*\beta]_1^{c_0,c_0} &=&  u(c_0)\bullet^0(\alpha^{c_0,c_0} * \beta_1^{Fc_0,Gc_0} )\\
   &\eq{i}& (u(c_0)\bullet^0\alpha^{c_0,c_0}) * \beta_1^{Fc_0,Gc_0} \\
   &\eq{ii}& id_{[\alpha(c_0)]} * \beta_1^{Fc_0,Gc_0}   \\
   &\eq{iii}& id_{[\alpha(c_0)]\bullet^0\beta_1^{Fc_0,Gc_0}} \\
   &\eq{iv}& id_{[\beta_1^{Fc_0,Gc_0}](\alpha(c_0))} \\
   &=& id_{[\alpha*\beta]_0(c_0)}
\end{eqnarray*}
where first and last equation are definitions, $(i)$ holds by $*$-associativity, $(ii)$ by units axioms for
$\alpha$, $(iii)$ by $*$-identity property, $(iv)$ is simply the application of a 2-morphism to a costant
1-morphism.\\

To prove composition axiom, we start by fixing arbitrary objects $c_0,c_0',c_0''$ of $\C$. Then it
is easier to start from the result and back-track the chain of equalities as shown below. By $*$-associativity
one has
$$
(-\circ-)\bullet^0(\alpha_1^{c_0,c_0''}*\beta_1^{Fc_0,Gc_0''})=((-\circ-)\bullet^0\alpha_1^{c_0,c_0''})*\beta_1^{Fc_0,Gc_0''}
$$
applying composition coherence of $\alpha$
$$
\Big[
\Big(\alpha_1^{c_0,c_0'}\circ G_1^{c_0',c_0''}\Big)
\bullet^1
\Big(F_1^{c_0,c_0'}\circ \alpha_1^{c_0',c_0''}\Big)
\Big]
*\beta_1^{Fc_0,Gc_0''}
$$
by $*$-functoriality
$$
\begin{array}{c}
\Big[
\Big(\alpha_1^{c_0,c_0'}\circ G_1^{c_0',c_0''}\Big)*\beta_1^{Fc_0,Gc_0''}
\Big]
\bullet^1
\Big[
\Big(F_1^{c_0,c_0'}\circ \alpha_1^{c_0',c_0''}\Big)
\bullet^0
\Big(H_1^{Fc_0,Gc_0''}\circ\beta_{Gc_0''}\Big)
\Big]
\\ \phantom{\Big[}\!\!\!\! \bullet^2 \\
\Big[
\Big(\alpha_1^{c_0,c_0'}\circ G_1^{c_0',c_0''}\Big)
\bullet^0
\Big(\beta_{Fc_0}\circ K_1^{Fc_0,Gc_0''}\Big)
\Big]
\bullet^1
\Big[
\Big(F_1^{c_0,c_0'}\circ \alpha_1^{c_0',c_0''}\Big)*\beta_1^{Fc_0,Gc_0''}
\Big]
\end{array}
$$
that is
$$
\begin{array}{c}
\Big[
\left(\Big(\alpha_1^{c_0,c_0'}\times G_1^{c_0',c_0''}\Big)\bullet^0 (-\circ-)\right)*\beta_1^{Fc_0,Gc_0''}
\Big]
\bullet^1
\Big[
\Big(F_1^{c_0,c_0'}\circ \alpha_1^{c_0',c_0''}\Big)
\bullet^0
\Big(H_1^{Fc_0,Gc_0''}\circ\beta_{Gc_0''}\Big)
\Big]
\\ \phantom{\Big[}\!\!\!\! \bullet^2 \\
\Big[
\Big(\alpha_1^{c_0,c_0'}\circ G_1^{c_0',c_0''}\Big)
\bullet^0
\Big(\beta_{Fc_0}\circ K_1^{Fc_0,Gc_0''}\Big)
\Big]
\bullet^1
\Big[
\left(\Big(F_1^{c_0,c_0'}\times \alpha_1^{c_0',c_0''}\Big)\bullet^0 (-\circ-)\right)*\beta_1^{Fc_0,Gc_0''}
\Big]
\end{array}
$$
Then by $*$-associativity we obtain
$$
\begin{array}{c}
\Big[
\Big(\alpha_1^{c_0,c_0'}\times G_1^{c_0',c_0''}\Big)*\Big( (-\circ-)\bullet^0\beta_1^{Fc_0,Gc_0''}\Big)
\Big]
\bullet^1
\Big[
\Big(F_1^{c_0,c_0'}\circ \alpha_1^{c_0',c_0''}\Big)
\bullet^0
\Big(H_1^{Fc_0,Gc_0''}\circ\beta_{Gc_0''}\Big)
\Big]
\\ \phantom{\Big[}\!\!\!\! \bullet^2 \\
\Big[
\Big(\alpha_1^{c_0,c_0'}\circ G_1^{c_0',c_0''}\Big)
\bullet^0
\Big(\beta_{Fc_0}\circ K_1^{Fc_0,Gc_0''}\Big)
\Big]
\bullet^1
\Big[
\Big(F_1^{c_0,c_0'}\times \alpha_1^{c_0',c_0''}\Big)*\Big( (-\circ-)\bullet^0\beta_1^{Fc_0,Gc_0''}\Big)
\Big]
\end{array}
$$
applying now composition axiom of $\beta$ this turns in
\begin{changemargin}{-23ex}{-20ex}
$$
\begin{array}{c}
\Big[
\Big(\alpha_1^{c_0,c_0'}\times G_1^{c_0',c_0''}\Big)*
\Big((\beta_1^{Fc_0,Gc_0'}\circ K_1^{Gc_0',Gc_0'})\bullet^1(H_1^{Fc_0,Gc_0'}\circ\beta_1^{Gc_0',Gc_0''})\Big)
\Big]
\bullet^1
\Big[
\Big(F_1^{c_0,c_0'}\circ \alpha_1^{c_0',c_0''}\Big)
\bullet^0
\Big(H_1^{Fc_0,Gc_0''}\circ\beta_{Gc_0''}\Big)
\Big]
\\ \phantom{\Big[}\!\!\!\! \bullet^2 \\
\Big[
\Big(\alpha_1^{c_0,c_0'}\circ G_1^{c_0',c_0''}\Big)
\bullet^0
\Big(\beta_{Fc_0}\circ K_1^{Fc_0,Gc_0''}\Big)
\Big]
\bullet^1
\Big[
\Big(F_1^{c_0,c_0'}\times \alpha_1^{c_0',c_0''}\Big)*
\Big((\beta_1^{Fc_0,Fc_0'}\circ K_1^{Fc_0',Gc_0''})\bullet^1(H_1^{Fc_0,Fc_0'}\circ\beta_1^{Fc_0',Gc_0''})\Big)
\Big]
\end{array}
$$
\end{changemargin}
{\em i.e.}
\begin{changemargin}{-10ex}{-10ex}
$$
\left(
  \begin{array}{c}
\Big(\alpha_1^{c_0,c_0'}\times G_1^{c_0',c_0''}\Big)*
\Big((\beta_1^{Fc_0,Gc_0'}\circ K_1^{Gc_0',Gc_0'})\bullet^1(H_1^{Fc_0,Gc_0'}\circ\beta_1^{Gc_0',Gc_0''})\Big)
\\ \bullet^1\\
\Big(F_1^{c_0,c_0'}\circ \alpha_1^{c_0',c_0''}\Big)
\bullet^0
\Big(H_1^{Fc_0,Gc_0''}\circ\beta_{Gc_0''}\Big)
  \end{array}
\right)
$$
\begin{equation}\label{zzz:ababouz}
   \bullet^2
\end{equation}
$$
\left(
  \begin{array}{c}
\Big(\alpha_1^{c_0,c_0'}\circ G_1^{c_0',c_0''}\Big)
\bullet^0
\Big(\beta_{Fc_0}\circ K_1^{Fc_0,Gc_0''}\Big)
\\ \bullet^1\\
\Big(F_1^{c_0,c_0'}\times \alpha_1^{c_0',c_0''}\Big)*
\Big((\beta_1^{Fc_0,Fc_0'}\circ K_1^{Fc_0',Gc_0''})\bullet^1(H_1^{Fc_0,Fc_0'}\circ\beta_1^{Fc_0',Gc_0''})\Big)
  \end{array}
\right)
$$
\end{changemargin}
\vskip4ex
Let us focus our attention on the second composite (w.r.t. $\bullet^2$-composition). In fact the calculations on the first component
are precisely symmetrical. \\

Applying $*$-functoriality to this we get
\begin{changemargin}{-20ex}{-20ex}
$$
  \begin{array}{c}
\Big(\alpha_1^{c_0,c_0'}\circ G_1^{c_0',c_0''}\Big)
\bullet^0
\Big(\beta_{Fc_0}\circ K_1^{Fc_0,Gc_0''}\Big)
\\ \bullet^1\\
\left(
\begin{array}{c}
\Big(
\big(F_1^{c_0,c_0'}\times \alpha_1^{c_0',c_0''}\big)
*
\big(\beta_1^{Fc_0,Fc_0'}\circ K_1^{Fc_0',Gc_0''}\big)
\Big)
\bullet^1
\Big(
\big(F_1^{c_0,c_0'}\times F_1^{c_0',c_0''}\circ\alpha_{c_0''}\big)
\bullet^0
\big(H_1^{Fc_0,Fc_0'}\circ\beta_1^{Fc_0',Gc_0''}\big)
\Big)
\\ \bullet^2 \\
\Big(
\big(F_1^{c_0,c_0'}\times\alpha_{c_0'}\circ G_1^{c_0',c_0''}\big)
\bullet^0
\big(\beta_1^{Fc_0,Fc_0'}\circ K_1^{Fc_0',Gc_0''}\big)
\Big)
\bullet^1
\Big(
\big(F_1^{c_0,c_0'}\times \alpha_1^{c_0',c_0''}\big)
*
\big(H_1^{Fc_0,Fc_0'}\circ\beta_1^{Fc_0',Gc_0''}\big)
\Big)
\end{array}
\right)
  \end{array}
$$
\end{changemargin}

\vskip4ex
By composing on product components we notice that upper $*$-composition gives indeed an identity
\begin{changemargin}{-20ex}{-20ex}
$$
\begin{array}{c}
\Big(\alpha_1^{c_0,c_0'}\circ G_1^{c_0',c_0''}\Big)
\bullet^0
\Big(\beta_{Fc_0}\circ K_1^{Fc_0,Gc_0''}\Big)
\\ \bullet^1\\
\left(
\begin{array}{c}
\Big(
(F_1^{c_0,c_0'}\bullet^0\beta_1^{Fc_0,Fc_0'})\circ(\alpha_1^{c_0',c_0''}\bullet^0K_1^{Fc_0',Gc_0''})
\Big)
\bullet^1
\Big(
\big(F_1^{c_0,c_0'}\times F_1^{c_0',c_0''}\circ\alpha_{c_0''}\big)
\bullet^0
\big(H_1^{Fc_0,Fc_0'}\circ\beta_1^{Fc_0',Gc_0''}\big)
\Big)
\\ \bullet^2 \\
\Big(
\big(F_1^{c_0,c_0'}\times\alpha_{c_0'}\circ G_1^{c_0',c_0''}\big)
\bullet^0
\big(\beta_1^{Fc_0,Fc_0'}\circ K_1^{Fc_0',Gc_0''}\big)
\Big)
\bullet^1
\Big(
(F_1^{c_0,c_0'}\bullet^0H_1^{Fc_0,Fc_0'})\circ (\alpha_1^{c_0',c_0''} * \beta_1^{Fc_0',Gc_0''})
\Big)
\end{array}
\right)
\end{array}
$$
\end{changemargin}

\vskip4ex
hence all the middle row is an identity 2-morphism, and  the whole simplifies to the following
\begin{changemargin}{-20ex}{-20ex}
$$
  \begin{array}{c}
\Big(\alpha_1^{c_0,c_0'}\circ G_1^{c_0',c_0''}\Big)
\bullet^0
\Big(\beta_{Fc_0}\circ K_1^{Fc_0,Gc_0''}\Big)
\\ \bullet^1\\
\Big(
\big(F_1^{c_0,c_0'}\times\alpha_{c_0'}\circ G_1^{c_0',c_0''}\big)
\bullet^0
\big(\beta_1^{Fc_0,Fc_0'}\circ K_1^{Fc_0',Gc_0''}\big)
\Big)
\bullet^1
\Big(
(F_1^{c_0,c_0'}\bullet^0H_1^{Fc_0,Fc_0'})\circ (\alpha_1^{c_0',c_0''} * \beta_1^{Fc_0',Gc_0''})
\Big)
  \end{array}
$$
\end{changemargin}

\vskip4ex
by $\bullet^1$-whiskering associativity this can be rearranged
\begin{changemargin}{-20ex}{-20ex}
$$
  \begin{array}{c}
\Big(\big(\alpha_1^{c_0,c_0'}\circ G_1^{c_0',c_0''}\big)
\bullet^0
\big(\beta_{Fc_0}\circ K_1^{Fc_0,Gc_0''}\big)\Big)
\bullet^1
\Big(
\big(F_1^{c_0,c_0'}\times\alpha_{c_0'}\circ G_1^{c_0',c_0''}\big)
\bullet^0
\big(\beta_1^{Fc_0,Fc_0'}\circ K_1^{Fc_0',Gc_0''}\big)
\Big)
\\ \bullet^1\\
(F_1^{c_0,c_0'}\bullet^0H_1^{Fc_0,Fc_0'})\circ (\alpha_1^{c_0',c_0''} * \beta_1^{Fc_0',Gc_0''})
  \end{array}
$$
\end{changemargin}

\vskip4ex
By functoriality of $K$
\begin{changemargin}{-20ex}{-20ex}
$$
  \begin{array}{c}
\Big(\big(\alpha_1^{c_0,c_0'}\circ G_1^{c_0',c_0''}\big)
\bullet^0
\big(\beta_{Fc_0}\circ K_1^{Fc_0,Gc_0''}\big)\Big)
\bullet^1
\Big(
\big(F_1^{c_0,c_0'}\times G_1^{c_0',c_0''}\big)
\bullet^0
\big(\beta_1^{Fc_0,Fc_0'}\circ K_1^{Fc_0',Gc_0'}(\alpha_{c_0'})\circ K_1^{Gc_0',Gc_0''}\big)
\Big)
\\ \bullet^1\\
(F_1^{c_0,c_0'}\bullet^0H_1^{Fc_0,Fc_0'})\circ (\alpha_1^{c_0',c_0''} * \beta_1^{Fc_0',Gc_0''})
  \end{array}
$$
\end{changemargin}

\vskip4ex
by functoriality of $-\circ-$
\begin{changemargin}{-20ex}{-20ex}
$$
  \begin{array}{c}
\Big(\big(\alpha_1^{c_0,c_0'}\circ G_1^{c_0',c_0''}\big)
\bullet^0
\big(\beta_{Fc_0}\circ K_1^{Fc_0,Gc_0''}\big)\Big)
\bullet^1
\Big((F_1^{c_0,c_0'}\bullet^0\beta_1^{Fc_0,Fc_0'})\circ K_1^{Fc_0',Gc_0'}(\alpha_{c_0'})\circ (G_1^{c_0',c_0''}\bullet^0K_1^{Gc_0',Gc_0''})
\Big)
\\ \bullet^1\\
(F_1^{c_0,c_0'}\bullet^0H_1^{Fc_0,Fc_0'})\circ (\alpha_1^{c_0',c_0''} * \beta_1^{Fc_0',Gc_0''})
  \end{array}
$$
\end{changemargin}

\vskip4ex
that is

\begin{changemargin}{-20ex}{-20ex}
$$
  \begin{array}{c}
\Big(\beta_{Fc_0}\circ
\big((\alpha_1^{c_0,c_0'}\circ G_1^{c_0',c_0''})
\bullet^0
K_1^{Fc_0,Gc_0''}\big)\Big)
\bullet^1
\Big((F_1^{c_0,c_0'}\bullet^0\beta_1^{Fc_0,Fc_0'})\circ K_1^{Fc_0',Gc_0'}(\alpha_{c_0'})\circ (G_1^{c_0',c_0''}\bullet^0K_1^{Gc_0',Gc_0''})
\Big)
\\ \bullet^1\\
(F_1^{c_0,c_0'}\bullet^0H_1^{Fc_0,Fc_0'})\circ (\alpha_1^{c_0',c_0''} * \beta_1^{Fc_0',Gc_0''})
  \end{array}
$$
\end{changemargin}

\vskip4ex
again by functoriality of $K$ we can write the result as
\begin{changemargin}{-20ex}{-20ex}
$$
  \begin{array}{c}
\beta_{Fc_0}\circ
\big(
\alpha_1^{c_0,c_0'}\bullet^0 K_1^{Fc_0,Gc_0'}
\big)
\circ
\big(
G_1^{c_0',c_0''}\bullet^0K_1^{Gc_0',Gc_0''}
\big)
\\
\bullet^1
\\
(F_1^{c_0,c_0'}\bullet^0\beta_1^{Fc_0,Fc_0'})\circ K_1^{Fc_0',Gc_0'}(\alpha_{c_0'})\circ (G_1^{c_0',c_0''}\bullet^0K_1^{Gc_0',Gc_0''})
\\ \bullet^1\\
(F_1^{c_0,c_0'}\bullet^0H_1^{Fc_0,Fc_0'})\circ (\alpha_1^{c_0',c_0''} * \beta_1^{Fc_0',Gc_0''})
  \end{array}
$$
\end{changemargin}

\vskip4ex
or more simply
\begin{changemargin}{-20ex}{-20ex}
$$
  \begin{array}{c}
\left(\Big(\beta_{Fc_0}\circ
\big(
\alpha_1^{c_0,c_0'}\bullet^0 K_1^{Fc_0,Gc_0'}
\big)
\Big)
\bullet^1
\Big((F_1^{c_0,c_0'}\bullet^0\beta_1^{Fc_0,Fc_0'})\circ K_1^{Fc_0',Gc_0'}(\alpha_{c_0'})\Big)
\right)
\circ (G_1^{c_0',c_0''}\bullet^0K_1^{Gc_0',Gc_0''})
\\ \bullet^1\\
(F_1^{c_0,c_0'}\bullet^0H_1^{Fc_0,Fc_0'})\circ (\alpha_1^{c_0',c_0''} * \beta_1^{Fc_0',Gc_0''})
  \end{array}
$$
\end{changemargin}

\vskip4ex
{i.e.}
$$
\begin{array}{c}
\big[(F\bullet^0\beta)\bullet^1(\alpha\bullet^0 K)\big]_1^{c_0,c_0'}\circ \big[G\bullet^0 K\big]_1^{c_0',c_0''}
\\\phantom{\Big[} \bullet^1\\
{}\big[F \bullet^0 H\big]_1^{c_0,c_0'}\circ \big[\alpha*\beta\big]_1^{c_0',c_0''}
\end{array}
$$

\vskip4ex
Carrying on the analogous calculations on the first component, (\ref{zzz:ababouz}) equals to
$$
\begin{array}{c}
\Big(\big[
\alpha*\beta
\big]_1^{c_0,c_0'}\circ \big[
G\bullet^0K
\big]_1^{c_0',c_0''}\Big)
\bullet^1
\Big({}\big[
F\bullet^0H
\big]_1^{c_0,c_0'}\circ \big[
(F\bullet^0\beta)\bullet^1(\alpha\bullet^0K)
\big]_1^{c_0',c_0''}\Big)
\\\phantom{\Big[} \bullet^2\\
\Big(\big[(F\bullet^0\beta)\bullet^1(\alpha\bullet^0 K)\big]_1^{c_0,c_0'}\circ \big[G\bullet^0 K\big]_1^{c_0',c_0''}\Big)
\bullet^1
\Big({}\big[F \bullet^0 H\big]_1^{c_0,c_0'}\circ \big[\alpha*\beta\big]_1^{c_0',c_0''}\Big)
\end{array}
$$

\vskip4ex
that is
$$
\begin{array}{c}
\Big(\big[
\alpha*\beta
\big]_1^{c_0,c_0'}\circ \big[
G\bullet^0K
\big]_1^{c_0',c_0''}\Big)
\bullet^1
\Big({}\big[
F\bullet^0H
\big]_1^{c_0,c_0'}\circ \big[
\mathbf{dom}(\alpha*\beta)
\big]_1^{c_0',c_0''}\Big)
\\\phantom{\Big[} \bullet^2\\
\Big(\big[\mathbf{cod}(\alpha*\beta)\big]_1^{c_0,c_0'}\circ \big[G\bullet^0 K\big]_1^{c_0',c_0''}\Big)
\bullet^1
\Big({}\big[F \bullet^0 H\big]_1^{c_0,c_0'}\circ \big[\alpha*\beta\big]_1^{c_0',c_0''}\Big)
\end{array}
$$
\vskip4ex
and this conclude the proof.
\end{proof}
\begin{Remark}
We have adopted the $*$-symbol instead of the more obvious $\bullet^0$ in order to emphasize the dimension-raising
property of this composition. Nevertheless $*$-properties w.r.t. other $\bullet^0$-compositions are somehow
better understood thinking only in terms of $\bullet^0$.
\end{Remark}

\begin{Lemma}
Given the case
$$
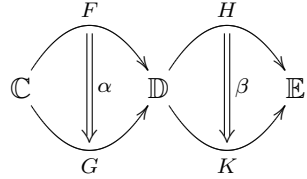
\xymatrix@C=8ex{
\C
\ar@/^5ex/[r]^{F}_{}="a1"
\ar@/_5ex/[r]_{G}^{}="a2"
\ar@{=>}"a1";"a2"^{\alpha}
&
\D
\ar@/^5ex/[r]^{H}_{}="b1"
\ar@/_5ex/[r]_{K}^{}="b2"
\ar@{=>}"b1";"b2"^{\beta}
&
\E
}
$$

If $\alpha$ is a lax natural $n$-transformation and $\beta$ is a strict natural $n$-transformation, the composition
$\alpha*\beta$ is an identity.

In this case it is possible to deal with dimension preserving 0-composition of 2-morphisms, by letting
$$
\alpha\tilde{*}\beta=\mathbf{dom}(\alpha*\beta)=\mathbf{cod}(\alpha*\beta)
$$
\end{Lemma}
\begin{proof}
Let us suppose $\alpha$ is a lax natural $n$-transformation and $\beta$ is a strict natural $n$-transformation.
Then for an object $c_0$ of $\C$, $\beta_1(\alpha_{c_0})$ is the commutative square
$\beta_{Fc_0}\circ K(\alpha_{c_0})=H(\alpha_{c_0})\circ\beta_{Gc_0}$.

Moreover once objects $c_0,c_0'$ of $\C$ are fixed, since $\beta$ is strict
$$
\beta_1^{Fc_0,Gc_0'}=id_{H_1^{Fc_0,Gc_0'}\circ\beta_{Gc_0'}}=id_{\alpha_{Fc_0}\circ K_1^{Fc_0,Gc_0'}}
$$
then
\begin{eqnarray*}
[\alpha*\beta]_1^{c_0,c_0'}
   &=& \alpha_1^{c_0,c_0'}*\beta_1^{Fc_0,Gc_0'} \\
   &=& \alpha_1^{c_0,c_0'}*id_{H_1^{Fc_0,Gc_0'}\circ\beta_{Gc_0'}} \\
   &=& id_{\alpha_1^{c_0,c_0'}\bullet^0 (H_1^{Fc_0,Gc_0'} \circ\beta_{Gc_0'})}
\end{eqnarray*}
Now the whiskering $\alpha_1^{c_0,c_0'}\bullet^0 (H_1^{Fc_0,Gc_0'} \circ\beta_{Gc_0'})$ has domain
\begin{eqnarray*}
\mathbf{dom}\big(\alpha_1^{c_0,c_0'}\bullet^0 (H_1^{Fc_0,Gc_0'} \circ\beta_{Gc_0'})\big)
   &=& (F_1^{c_0,c_0'}\circ\alpha_{c_0'})\bullet^0 (H_1^{Fc_0,Gc_0'} \circ\beta_{Gc_0'})\\
   &=& (F_1^{c_0,c_0'}\bullet^0H_1^{Fc_0,Fc_0'})\circ H(\alpha_{c_0'}) \circ\beta_{Gc_0'}\\
   &=&  [F\bullet^0H]_1^{c_0,c_0'}\circ [\alpha\tilde{*}\beta]_{c_0'}
\end{eqnarray*}
Similarly the codomain is
\begin{eqnarray*}
\mathbf{cod}\big(\alpha_1^{c_0,c_0'}\bullet^0 (H_1^{Fc_0,Gc_0'} \circ\beta_{Gc_0'})\big)
   &=& (\alpha_{c_0}\circ G_1^{c_0,c_0'})\bullet^0 (H_1^{Fc_0,Gc_0'} \circ\beta_{Gc_0'})\\
   &=& (\alpha_{c_0}\circ G_1^{c_0,c_0'})\bullet^0 (\beta_{Fc_0}\circ K_1^{Fc_0,Gc_0'})\\
   &=& \beta_{Fc_0}\circ K(\alpha_{c_0})\circ ( G_1^{c_0,c_0'} \bullet^0K_1^{Fc_0,Gc_0'})\\
   &=&[\alpha\tilde{*}\beta]_{c_0}\circ [G\bullet^0K]_1^{c_0,c_0'}
\end{eqnarray*}
hence the result is (an identity over) a 2-morphism.
\end{proof}
Importance of {\em Lemma} above is in that it allows to right-0-compose freely with constant transformations,
such as $-\circ c_2$ or $c_2\circ -$ for a 2-cell $c_2:c_1\Rightarrow c_1':c_0\to c_0'$.\\

Notice that {\em Lemma} does not hold for $\alpha$ strict and $\beta$ lax, since in this case the result is a
strict 3-morphism.

\subsection{Properties}
The following propositions conclude the description of dimension-rising composition in the sesqui${}^2$category
of strict $n$-categories.
\vskip4ex
Given the situation
$$
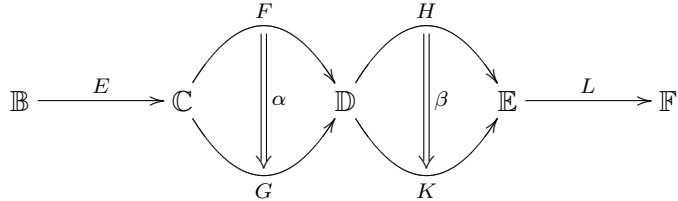
\xymatrix@C=10ex@R=10ex{
\B
\ar[r]^{E}
&
\C
\ar@/^6ex/[r]^{F}_{}="a1"
\ar@/_6ex/[r]_{G}^{}="a2"
\ar@2"a1";"a2"^{\alpha}
&
\D
\ar@/^6ex/[r]^{H}_{}="b1"
\ar@/_6ex/[r]_{K}^{}="b2"
\ar@2"b1";"b2"^{\beta}
&
\E
\ar[r]^{L}
&
\F}
$$
one has the following
\begin{Proposition}[$*$-associativity 1]
$$
\begin{array}{clcl}
  (L*A) & (E\bullet_L^0\alpha)*\beta=E\bullet^0_L(\alpha*\beta) & (R*A) & \alpha*(\beta\bullet_R^0L)=(\alpha*\beta)\bullet^0_R L
\end{array}
$$
\end{Proposition}
\begin{proof}
We prove $(L*A)$. The proof of $(R*A)$ is similar hence it is omitted. \\

({\em on objects}) Let an object $b_0$ of $\C$ be given. Then
\begin{eqnarray*}
  [(E\bullet^0\alpha)*\beta]_{b_0} &=&  \beta([E\bullet^0\alpha]_{b_0})\\
   &=& \beta(\alpha_{Eb_0}) \\
   &=&  [\alpha*\beta]_{Eb_0}\\
   &=& [E\bullet^0(\alpha*\beta)]_{b_0}
\end{eqnarray*}

({\em on homs}) Let objects $b_0,b_0'$ of $\B$ be given. Then
\begin{eqnarray*}
  [(E\bullet^0\alpha)*\beta]_1^{b_0,b_0'} &=& [E\bullet^0\alpha]_1^{b_0,b_0'}*\beta_1^{F(Eb_0),G(Eb_0')}\\
   &=&  (E_1^{b_0,b_0'}\bullet^0\alpha_1^{Eb_0,Eb_0'})*\beta_1^{F(Eb_0),G(Eb_0')}\\
   &\eq{\spadesuit}&  E_1^{b_0,b_0'}\bullet^0(\alpha_1^{Eb_0,Eb_0'}*\beta_1^{F(Eb_0),G(Eb_0')})\\
   &=&  E_1^{b_0,b_0'}\bullet^0[\alpha_1*\beta]^{Eb_0,Eb_0'}\\
   &=& [E\bullet^0(\alpha*\beta)]_1^{b_0,b_0'}
\end{eqnarray*}
where $(\spadesuit)$ holds by $*$-associativity in dimension $n-1$.
\end{proof}

\begin{Proposition}[$*$-identity]
$$
\begin{array}{clcl}
(L)&id_{E}*\alpha=id_{E\bullet_L^0\alpha}
&(R)&\alpha*id_{H}=id_{\alpha\bullet^0_R H}
\end{array}
$$
\end{Proposition}
\begin{proof}
We prove $(L)$. The proof of $(R)$ is similar hence it is omitted. \\

({\em on objects}) Let an object $b_0$ of $\B$ be given. Then the following equalities are straightforward
\begin{eqnarray*}
 [id_{E}*\alpha]_{b_0}  &=& \alpha([id_{E}]_{b_0}) \\
   &=&  \alpha(id_{Eb_0})\\
   &=&  id_{\alpha(Eb_0)}\\
   &=&  [id_{E\bullet_L^0\alpha}]_{b_0}
\end{eqnarray*}

({\em on homs}) Let objects $b_0,b_0'$ of $\B$ be given. Then
\begin{eqnarray*}
 [id_{E}*\alpha]_1^{b_0,b_0'}  &=& [id_{E}]_1^{b_0,b_0'}*\alpha_1^{Eb_0,Eb_0'}\\
   &=&  id_{E_1^{b_0,b_0'}}*\alpha_1^{Eb_0,Eb_0'} \\
   &\eq{\spadesuit}&  id_{E_1^{b_0,b_0'}\bullet^0 \alpha_1^{Eb_0,Eb_0'}} \\
   &=&[id_{E\bullet_L^0\alpha}]_1^{b_0,b_0'}
\end{eqnarray*}
where $(\spadesuit)$ is $*$-identity in dimension $n-1$.
\end{proof}

\vskip4ex

In the situation
$$
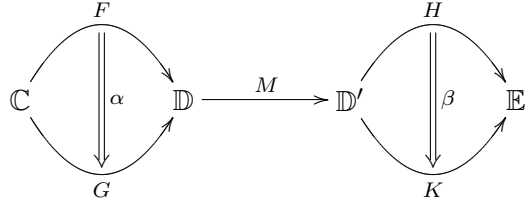
\xymatrix@C=10ex@R=10ex{
\C
\ar@/^6ex/[r]^{F}_{}="a1"
\ar@/_6ex/[r]_{G}^{}="a2"
\ar@2"a1";"a2"^{\alpha}
&
\D
\ar[r]^{M}
&
\D'
\ar@/^6ex/[r]^{H}_{}="b1"
\ar@/_6ex/[r]_{K}^{}="b2"
\ar@2"b1";"b2"^{\beta}
&
\E
}
$$
one has the following
\begin{Proposition}[$*$-associativity 2]
$$
\alpha * (M\bullet_L^0 \beta) = (\alpha\bullet_R^0 M)*\beta
$$
\end{Proposition}
\begin{proof}
({\em on objects}) Let an object $c_0$ of $\C$ be given. Then the following equalities are straightforward
\begin{eqnarray*}
  [\alpha * (M\bullet^0 \beta)]_{c_0} &=& [M\bullet^0 \beta](\alpha_{c_0}) \\
   &=& \beta\big(M(\alpha_{c_0})\big) \\
   &=& \beta\big([\alpha\bullet^0M]_{c_0} \big)\\
   &=& [(\alpha\bullet^0 M)*\beta]_{c_0}
\end{eqnarray*}

({\em on homs}) Let objects $c_0,c_0'$ of $\C$ be given. Then

\begin{eqnarray*}
\big[\alpha * (M\bullet^0 \beta)\big]_1^{c_0,c_0'}   &=& \alpha_1^{c_0,c_0'}*[M\bullet^0 \beta]_1^{Fc_0,Gc_0'} \\
   &=&  \alpha_1^{c_0,c_0'}*\Big(M_1^{Fc_0,Gc_0'}\bullet^0\beta_1^{M(Fc_0),M(Gc_0')}\Big) \\
   &\eq{\spadesuit}&  \Big(\alpha_1^{c_0,c_0'}\bullet^0M_1^{Fc_0,Gc_0'}\Big)*\beta_1^{M(Fc_0),M(Gc_0')}\\
   &=&  [\alpha\bullet^0M]_1^{c_0,c_0'}*\beta_1^{M(Fc_0),M(Gc_0')}\\
   &=&  \big[(\alpha\bullet^0 M)*\beta\big]_1^{c_0,c_0'}
\end{eqnarray*}
where $(\spadesuit)$ is $*$-associativity in dimension $n-1$.
\end{proof}
\vskip4ex
In the situation below
$$
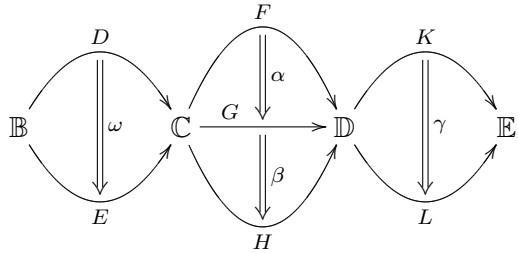
\xymatrix@C=10ex{
\B
\ar@/^6ex/[r]^{D}_{}="o1"
\ar@/_6ex/[r]_{E}^{}="o2"
\ar@2"o1";"o2"^{\omega}
&
\C
\ar@/^8ex/[r]^{F}_{}="a1"
\ar[r]^{}="a2"_{}="b1"^(.3){G}
\ar@/_8ex/[r]_{H}^{}="b2"
\ar@2"a1";"a2"^{\alpha}
\ar@2"b1";"b2"^{\beta}
&
\D
\ar@/^6ex/[r]^{K}_{}="g1"
\ar@/_6ex/[r]_{L}^{}="g2"
\ar@2"g1";"g2"^{\gamma}
&
\E
}
$$
one has the following
\begin{Proposition}[$*$-functoriality]
$$
(a)\qquad
(\alpha\bullet^1\beta)*\gamma =
\Big((\alpha*\gamma)\bullet^1(\beta\bullet^0L)\Big)
\bullet^2
\Big((\alpha\bullet^0K)\bullet^1(\beta*\gamma)\Big)
$$
$$
(b)\qquad
\omega*(\alpha\bullet^1\beta)=
\Big((\omega*\alpha)\bullet^1(E\bullet^0\beta)\Big)
\bullet^2
\Big((D\bullet^0\alpha)\bullet^1(\omega*\beta)\Big)
$$

\end{Proposition}
\begin{proof}
We prove $(a)$. The proof of $(b)$ is similar, hence it is omitted.\\

({\em on objects}) Let an object $c_0$ of $\C$ be given. Then
\begin{eqnarray*}
\left[
\begin{array}{c}
  (\alpha*\gamma)\bullet^1(\beta\bullet^0L) \\
  \bullet^2 \\
  (\alpha\bullet^0K)\bullet^1(\beta*\gamma)
\end{array}
\right]_{c_0}
&=&
\begin{array}{c}
  \left[(\alpha*\gamma)\bullet^1(\beta\bullet^0L)\right]_{c_0} \\
  \circ^1 \\
  \left[(\alpha\bullet^0K)\bullet^1(\beta*\gamma)\right]_{c_0}
\end{array}
\\
&\eq{\spadesuit}&
\left(
\begin{array}{c}
  \gamma(\alpha_{c_0})\bullet^0 L(\beta_{c_0})\\
  \circ^1 \\
  K(\alpha_{c_0})\bullet^0\gamma(\beta_{c_0})
\end{array}
\right)
  \\
   &=& \gamma(\alpha_{c_0}\circ \beta_{c_0}) \\
   &=& [(\alpha\bullet^1\beta)*\gamma]_{c_0}
\end{eqnarray*}
where all the equalities are just definitions, but $(\spadesuit)$ that is given by functoriality
w.r.t. 0-composition of $\gamma$.\\

({\em on homs}) Let objects $c_0,c_0'$ of $\C$ be given. Applying the definition of 2-composition
of 3-morphisms
$$
\left[
\begin{array}{c}
  (\alpha*\gamma)\bullet^1(\beta\bullet^0L) \\
  \bullet^2 \\
  (\alpha\bullet^0K)\bullet^1(\beta*\gamma)
\end{array}
\right]_1^{c_0,c_0'}=
$$
$$
=\begin{array}{c}
\Big[(\alpha\bullet^0 K)\bullet^1(\beta*\gamma)\Big]_1^{c_0,c_0'}\bullet^1\left([\alpha\bullet^0K]_1^{c_0,c_0'}\circ \mathbf{cod}(\beta*\gamma)\right)\\
\bullet^2\\
\left(\mathbf{dom}(\alpha*\gamma)\circ [\beta\bullet^0L]_1^{c_0,c_0'}\right)\bullet^1\Big[(\alpha*\gamma)\bullet^1(\beta\bullet^0 L)\Big]_1^{c_0,c_0'}
\end{array}
$$
by definition of whiskering of 3-morphisms and 2-morphisms
$$
\begin{array}{c}
\left([\alpha\bullet^0 K]_{c_0}\circ[\beta*\gamma]_1^{c_0,c_0'}\right)\bullet^1\left([\alpha\bullet^0K]_1^{c_0,c_0'}\circ \mathbf{cod}(\beta*\gamma)\right)\\
\bullet^2\\
\left(\mathbf{dom}(\alpha*\gamma)\circ [\beta\bullet^0L]_1^{c_0,c_0'}\right)\bullet^1\left([\alpha*\gamma]_1^{c_0,c_0'}\circ [\beta\bullet^0 L]_{c_0'}\right)
\end{array}
$$
that is
$$
\begin{array}{c}
\left([\alpha\bullet^0 K]_{c_0}\circ[\beta*\gamma]_1^{c_0,c_0'}\right)\bullet^1\left(\left([\alpha\bullet^0K]_1^{c_0,c_0'}\right)\circ K(\beta_{c_0'})\circ\gamma_{Hc_0'}\right)\\
\bullet^2\\
\left(\gamma_{Fc_0}\circ L(\alpha_{c_0})\circ [\beta\bullet^0L]_1^{c_0,c_0'}\right)\bullet^1\left([\alpha*\gamma]_1^{c_0,c_0'}\circ [\beta\bullet^0 L]_{c_0'}\right)
\end{array}
$$
by definition of $*$-composition on homs (and of 0-whiskering for 2-morphisms)
$$
\begin{array}{c}
\left(K(\alpha_{c_0})\circ\left(\beta_1^{c_0,c_0'}*\gamma_1^{Gc_0,Hc_0'}\right)\right)\bullet^1\left(\left(\alpha_1^{c_0,c_0'}\bullet^0K_1^{Fc_0,Gc_0'}\right)\circ K(\beta_{c_0'})\circ\gamma_{Hc_0'}\right)\\
\bullet^2\\
\left(\gamma_{Fc_0}\circ L(\alpha_{c_0})\circ \left(\beta_1^{c_0,c_0'}\bullet^0L_1^{Gc_0,Hc_0'}\right)\right)\bullet^1\left(\left(\alpha_1^{c_0,c_0'}*\gamma_1^{Fc_0,Gc_0'}\right)\circ L(\beta_{c_0'})\right)
\end{array}
$$
this can be rearranged
$$
\begin{array}{c}
\left(K(\alpha_{c_0})\circ\left(\beta_1^{c_0,c_0'}*\gamma_1^{Gc_0,Hc_0'}\right)\right)\bullet^1\left(\alpha_1^{c_0,c_0'}\bullet^0\left(K_1^{Fc_0,Gc_0'}\circ K(\beta_{c_0'})\circ\gamma_{Hc_0'}\right)\right)\\
\bullet^2\\
\left(\beta_1^{c_0,c_0'}\bullet^0\left(\gamma_{Fc_0}\circ L(\alpha_{c_0})\circ L_1^{Gc_0,Hc_0'}\right)\right)\bullet^1\left(\left(\alpha_1^{c_0,c_0'}*\gamma_1^{Fc_0,Gc_0'}\right)\circ L(\beta_{c_0'})\right)
\end{array}
$$
by functoriality of $L$ and $K$
$$
\begin{array}{c}
\left(K(\alpha_{c_0})\circ\left(\beta_1^{c_0,c_0'}*\gamma_1^{Gc_0,Hc_0'}\right)\right)\bullet^1\left(\left(\alpha_1^{c_0,c_0'}\circ\beta_{c_0'}\right)\bullet^0\left(K_1^{Fc_0,Hc_0'}\circ\gamma_{Hc_0'}\right)\right)\\
\bullet^2\\
\left(\left(\alpha_{c_0}\circ\beta_1^{c_0,c_0'}\right)\bullet^0\left(\gamma_{Fc_0}\circ L_1^{Fc_0,Hc_0'}\right)\right)\bullet^1\left(\left(\alpha_1^{c_0,c_0'}*\gamma_1^{Fc_0,Gc_0'}\right)\circ L(\beta_{c_0'})\right)
\end{array}
$$
by functoriality (and $*$-associativity)
$$
\begin{array}{c}
\left(\left(\alpha_{c_0}\circ\beta_1^{c_0,c_0'}\right)*\gamma_1^{Fc_0,Hc_0'}\right)\bullet^1\left(\left(\alpha_1^{c_0,c_0'}\circ\beta_{c_0'}\right)\bullet^0\left(K_1^{Fc_0,Hc_0'}\circ\gamma_{Hc_0'}\right)\right)\\
\bullet^2\\
\left(\left(\alpha_{c_0}\circ\beta_1^{c_0,c_0'}\right)\bullet^0\left(\gamma_{Fc_0}\circ L_1^{Fc_0,Hc_0'}\right)\right)\bullet^1\left(\left(\alpha_1^{c_0,c_0'}\circ\beta_{c_0'}\right)*\gamma_1^{Fc_0,Hc_0'}\right)
\end{array}
$$
and finally by $*$-functoriality
$$
\Big((\alpha_{c_0}\circ\beta_1^{c_0,c_0'})\bullet^1(\alpha_1^{c_0,c_0'}\circ\beta_{c_0'})\Big)*\gamma_1^{Fc_0,Hc_0'}
$$
that is the result:
$$
\Big[(\alpha\bullet^1\beta)*\gamma\Big]_1^{c_0,c_0'}.
$$

\end{proof}

\chapter{$h$-Pullbacks revisited and the long exact sequence}\label{cha:twopullbacks}

\section{2-dimensional $h$-pullbacks in $n$\textbf{Cat}}
We introduce here a notion of 2-dimensional h-pullback in the sesqui${}^2$-category $n$\textbf{Cat}. It
will be shown that our construction of the {\em standard} $h$-pullback of $n$-categories is an instance
of such a 2-dimensional one.\\[6ex]

In order to fix notation, let us consider the following diagram in $n$\textbf{Cat}
$$
\xymatrix{&\C\ar[d]^G\\\A\ar[r]_F&\B}
$$
A $h$-2pullback of $F$ and $G$ is a four-tuple $(\PP, P,Q,\varepsilon)$
$$
\xymatrix{
\PP\ar[r]^Q\ar[d]_P
&\C\ar[d]^G
\\
\A\ar[r]_F\ar@{}[ur]|(.3){}="1"|(.7){}="2"
&\B
\ar@{=>}"1";"2"^{\varepsilon}}
$$
that satisfies the following $2$-dimensional universal property:

\begin{UP}[$h$-2pullbacks]\label{UP:2h-pullbacks}

For any other two four-tuple

$$
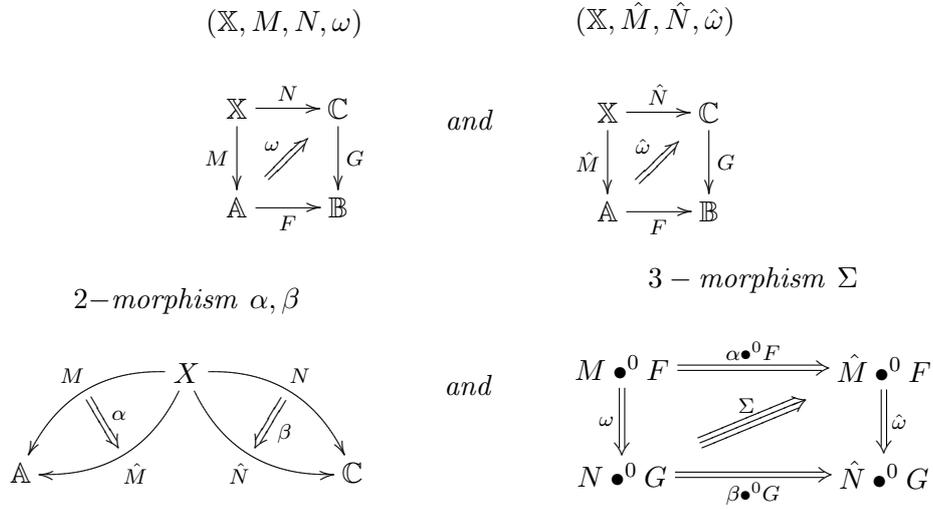
\begin{array}{c}
  (\X, M,N,\omega) \\ \\
  \xymatrix{
\X\ar[r]^N\ar[d]_M
&\C\ar[d]^G
\\
\A\ar[r]_F\ar@{}[ur]|(.3){}="1"|(.7){}="2"
&\B
\ar@{=>}"1";"2"^{\omega}}
\end{array}
\qquad\mathit{and}\qquad
\begin{array}{c}
  (\X, \hat{M},\hat{N},\hat{\omega}) \\ \\
  \xymatrix{
\X\ar[r]^{\hat{N}}\ar[d]_{\hat{M}}
&\C\ar[d]^G
\\
\A\ar[r]_F\ar@{}[ur]|(.3){}="1"|(.7){}="2"
&\B
\ar@{=>}"1";"2"^{\hat{\omega}}}
\end{array}
$$
$$
\begin{array}{c}
  2\mathit{-morphism}\ \alpha,\beta \\ \\
  \xymatrix@C=10ex{
&X
\ar@/_3ex/[dl]_{M}^{}="a1"
\ar@/^3ex/[dl]^{\hat{M}}_{}="a2"
\ar@2"a1";"a2"^{\alpha}
\ar@/^3ex/[dr]^{N}_{}="b1"
\ar@/_3ex/[dr]_{\hat{N}}^{}="b2"
\ar@2"b1";"b2"^{\beta}
\\
\A&&\C}
\end{array}
\qquad\mathit{and}\qquad
\begin{array}{c}
  3-\mathit{morphism}\ \Sigma \\ \\
 \xymatrix@C=12ex{
M\bullet^0F
\ar@2[r]^{\alpha\bullet^0F}
\ar@2[d]_{\omega}
&
\hat{M}\bullet^0F
\ar@2[d]^{\hat{\omega}}
\\
\ar@{}[ur]|(.3){}="1"|(.7){}="2"
\ar@3"1";"2"^{\Sigma}
N\bullet^0G
\ar@2[r]_{\beta\bullet^0G}
&
\hat{N}\bullet^0G
}
\end{array}
$$

there exists a unique $\lambda:\ L\Rightarrow \hat{L}:\ \X\rightarrow \PP$ such that

(UP)
\begin{enumerate}
  \item $\lambda\bullet^0P=\alpha$
  \item $\lambda\bullet^0Q=\beta$
  \item $\lambda * \varepsilon=\Sigma$
\end{enumerate}
\end{UP}
As  an immediate consequence of the definition, we state the following
\begin{Proposition}\label{prop:2up_implies_1up}
$2$-\emph{Universal Property} of $h$-2pullbacks implies $1$-dimensional one. Hence $h$-2pullbacks
are defined up to isomorphism.
\end{Proposition}
\begin{proof}
Just put $\alpha$, $\beta$ and $\Sigma$ identities.
\end{proof}

Let us notice that {\em Proposition \ref{prop:2up_implies_1up}} holds in every sesqui${}^2$-category.
More interestingly in $n$\textbf{Cat} a kind of converse to this proposition also holds.
\begin{Proposition}
Given the diagram
$$
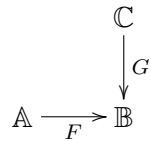
\xymatrix{&\C\ar[d]^G\\\A\ar[r]_F&\B}
$$
in $n$\textbf{Cat}, the standard $h$-pullback $<\PP,P,Q,\varepsilon>$ satisfies also {\em Universal property \ref{UP:2h-pullbacks}}.
\end{Proposition}
\begin{proof}
Firstly we remark that 1-dimensional {\em Universal Property \ref{UP:h-pullbacks}} of $h$-pullbacks
applied to the four-tuple $(\X, M,N,\omega)$ yields an $L:\X\to\PP$, while applied to $(\X, \hat{M},\hat{N},\hat{\omega})$,
 a $\hat{L}:\X\to\PP$. Those have to be \emph{domain} and \emph{co-domain} of the 2-cell  provided by
 the universal property, namely \makebox{$\lambda:\ L\Rightarrow\hat{L}$}.

We recall the constructions in order to fix notation.\\

$\bullet$ $L_0:\X_0\to\PP_0$ is the  map
$$
x_0\ \mapsto \ (Mx_0,\xymatrix{F(Mx_0)\ar[r]^{\omega(x_0)}&G(Nx_0)},Nx_0) =:p_0
$$

$\bullet$ for every pair of objects $x_0,x_0'$ of $\X$,
$$
L_1^{x_0,x_0'}:\X_1(x_0,x_0')\to \PP_1(L(x_0),L(x_0'))
$$
is given by the universal property in dimension $n-1$, and is such that
$$
L_1^{x_0,x_0'}\bullet^0 P_1^{Lx_0,Lx_0'}=M_1^{x_0,x_0'}
$$
$$
L_1^{x_0,x_0'}\bullet^0 Q_1^{Lx_0,Lx_0'}=N_1^{x_0,x_0'}
$$
$$
L_1^{x_0,x_0'}\bullet^0 \varepsilon_1^{Lx_0,Lx_0'}=\omega_1^{x_0,x_0'}
$$
The pair $L=<L_0,L_1^{-,-}>$ is a 1-morphism.

Similarly one determines $\hat{L}=<\hat{L}_0,\hat{L}_1^{-,-}>$.\\

Now we show that remaining data (namely, $\alpha$, $\beta$ and $\Sigma)$ of the hypothesis
provide a 2-morphism $\lambda:L\Rightarrow \hat{L}$ that satisfies required property.
To this end, let us consider the following assignments:\\

$\bullet$ For every object $x_0$ of $\X$,
$$
\lambda_{x_0}= (\alpha_{x_0}, \Sigma_{x_0}, \beta_{x_0} ):\ Lx_0\to\hat{L}x_0
$$
{\em i.e.}
$$
\begin{array}{c}
  \xymatrix@C=12ex{Mx_0\ar[r]^{\alpha_{x_0}}&\hat{M}x_0}
\\
    \xymatrix@C=10ex{F(Mx_0)\ar[r]^{F(\alpha_{x_0})}\ar[d]_{\omega_{x_0}}
    &F(\hat{M}x_0) \ar[d]^{\hat{\omega}_{x_0}}
    \\
    \ar@{}[ur]|(.3){}="1"|(.7){}="2"\ar@2"1";"2"^{\Sigma_{x_0}}
    G(Nx_0)\ar[r]_{G(\beta_{x_0})}&G(\hat{N})x_0}
\\
  \xymatrix@C=12ex{Nx_0\ar[r]_{\beta_{x_0}}&\hat{N}x_0}
  \end{array}
$$

$\bullet$ For every pair of objects $x_0,x_0'$ of $\X$,
$$
\xymatrix@C=-2ex{
&\X_1(x_0,x_0')
\ar[dl]_{L_1^{x_0,x_0'}}
\ar[dr]^{\hat{L}_1^{x_0,x_0'}}
\\
\PP_1(Lx_0,Lx_0')
\ar[dr]_{-\circ \lambda x_0'}
&&
\PP_1(\hat{L}x_0,\hat{L}x_0')
\ar[dl]^{\lambda x_0\circ-}
\ar@{}[ll]|(.3){}="1"|(.7){}="2"\ar@2"1";"2"_{\lambda_1^{x_0,x_0'}}
\\
&\PP_1(Lx_0,\hat{L}x_0')
}
$$
is given by the universal property for (n-1)categories.

In fact the 0-codomain of $\lambda_1^{x_0,x_0'}$, namely $\PP_1(Lx_0,\hat{L}x_0')$
is defined inductively as a $h$-2pullback in (n-1)\textbf{Cat}:
$$
\xymatrix@C=10ex@R=6ex{
\PP_1(Lx_0,\hat{L}x_0')
\ar[rr]^{Q_1^{Lx_0,\hat{L}x_0'}}
\ar[dd]_{P_1^{Lx_0,\hat{L}x_0'}}
&&\C_1(Nx_0,\hat{N}x_0')
\ar[d]^{G_1^{Nx_0,\hat{N}x_0'}}
\ar@{}[ddll]|(.4){}="1"|(.6){}="2"\ar@2"1";"2"_{\varepsilon_1^{Lx_0,\hat{L}x_0'}}
\\
&&\B_1(G(Nx_0),G(\hat{N}x_0'))
\ar[d]^{\omega x_0\circ -}
\\
\A_1(Mx_0,\hat{M}x_0')
\ar[r]_(.4){F_1^{Mx_0,\hat{M}x_0'}}
&\B_1(F(Mx_0),F(\hat{M}x_0'))
\ar[r]_{-\circ \hat{\omega}x_0'}
&\B_1(F(Mx_0),G(\hat{N}x_0'))}
$$

Over the same base are also defined
$$
\xymatrix@C=10ex@R=6ex{
\X_1(x_0,x_0')
\ar[rr]^{N_1^{x_0,x_0'}\circ\beta x_0'}
\ar[dd]_{M_1^{x_0,x_0'}\circ\alpha x_0'}
&&\C_1(Nx_0,\hat{N}x_0')
\ar[d]^{G_1^{Nx_0,\hat{N}x_0'}}
\ar@{}[ddll]|(.3){}="1"|(.7){}="2"
\ar@2"1";"2"|{\theta\ =\ (\omega_1^{x_0,x_0'}\circ G(\beta x_0'))\bullet^1([MF]_1^{x_0,x_0'}\circ \Sigma x_0')}
\\
&&\B_1(G(Nx_0),G(\hat{N}x_0'))
\ar[d]^{\omega x_0\circ -}
\\
\A_1(Mx_0,\hat{M}x_0')
\ar[r]_(.4){F_1^{Mx_0,\hat{M}x_0'}}
&\B_1(F(Mx_0),F(\hat{M}x_0'))
\ar[r]_{-\circ \hat{\omega}x_0'}
&\B_1(F(Mx_0),G(\hat{N}x_0'))}
$$
and
$$
\xymatrix@C=10ex@R=6ex{
\X_1(x_0,x_0')
\ar[rr]^{\beta x_0\circ \hat{N}_1^{x_0,x_0'}}
\ar[dd]_{\alpha x_0\circ\hat{M}_1^{x_0,x_0'}}
&&\C_1(Nx_0,\hat{N}x_0')
\ar[d]^{G_1^{Nx_0,\hat{N}x_0'}}
\ar@{}[ddll]|(.3){}="1"|(.7){}="2"
\ar@2"1";"2"|{\hat{\theta}\ = \ (\Sigma x_0\circ[\hat{N}G]_1^{x_0,x_0'})\bullet^1 (F(\alpha x_0)\circ\hat{\omega}_1^{x_0,x_0'})}
\\
&&\B_1(G(Nx_0),G(\hat{N}x_0'))
\ar[d]^{\omega x_0\circ -}
\\
\A_1(Mx_0,\hat{M}x_0')
\ar[r]_(.4){F_1^{Mx_0,\hat{M}x_0'}}
&\B_1(F(Mx_0),F(\hat{M}x_0'))
\ar[r]_{-\circ \hat{\omega}x_0'}
&\B_1(F(Mx_0),G(\hat{N}x_0'))}
$$

Moreover we can consider 2-morphisms:
$$
\alpha_1^{x_0,x_0'}:\ \alpha x_0\circ\hat{M}_1^{x_0,x_0'}\Rightarrow M_1^{x_0,x_0'}\circ\alpha x_0':\ \X_1(x_0,x_0')\to \A_1(Mx_0,\hat{M}x_0')
$$
$$
\beta_1^{x_0,x_0'}:\ \beta x_0\circ\hat{N}_1^{x_0,x_0'}\Rightarrow N_1^{x_0,x_0'}\circ\beta x_0':\ \X_1(x_0,x_0')\to \A_1(Mx_0,\hat{M}x_0')
$$
and the 3-morphism
\begin{changemargin}{-10ex}{-10ex}
$$
\xymatrix@C=14ex@R=14ex{
(\beta x_0\circ \hat{N}_1^{x_0,x_0'})\bullet^0(\omega x_0\circ G_1^{Nx_0,\hat{N}x_0'})
\ar@2[r]^{\beta_1^{x_0,x_0'}\bullet^0 id}
\ar@2[d]_{\hat{\theta}\ }
&
(N_1^{x_0,x_0'}\circ \beta x_0')\bullet^0(\omega x_0\circ G_1^{Nx_0,\hat{N}x_0'})
\ar@2[d]^{\ \theta}
\ar@{}[dl]|(.3){}="1"|(.7){}="2"
\ar@3"1";"2"|{\Sigma_1^{x_0,x_0'}}
\\
(\alpha x_0\circ \hat{M}_1^{x_0,x_0'})\bullet^0(F_1^{Mx_0,\hat{M}x_0'}\circ \hat{\omega}x_0')
\ar@2[r]_{\beta_1^{x_0,x_0'}\bullet^0 id}
&(M_1^{x_0,x_0'}\circ\alpha x_0')\bullet^0(F_1^{Mx_0,\hat{M}x_0'}\circ \hat{\omega}x_0')
}
$$
\end{changemargin}
Finally we can apply the universal property, in order to get a {\em unique} 2-morphism
$$
\lambda_1^{x_0,x_0'}:\ L_1^{x_0,x_0'}\circ \lambda x_0'\Rightarrow\lambda x_0\circ \hat{L}_1^{x_0,x_0'}
$$
such that
\begin{eqnarray}
  \lambda_1^{x_0,x_0'}\bullet^0  Q_1^{Lx_0,\hat{L}x_0'}&=&\beta_1^{x_0,x_0'}\label{eqn:2PUP1} \\
  \lambda_1^{x_0,x_0'}\bullet^0  P_1^{Lx_0,\hat{L}x_0'}&=& \alpha_1^{x_0,x_0'}\label{eqn:2PUP2} \\
  \lambda_1^{x_0,x_0'} * \varepsilon_1^{Lx_0,\hat{L}x_0'} &=& \Sigma_1^{x_0,x_0'}\label{eqn:2PUP3}
\end{eqnarray}
That the pair $\lambda=\ <\lambda_0,\lambda_1^{-,-}>$ is a 2-morphism of n-categories is proved in the
following (quite technical) {\em Lemma \ref{zzz:tech1}}.

Moreover it satisfies by construction {\em Universal Property \ref{UP:2h-pullbacks}}. In fact for any object $x_0$ of $\X$
$$
[\lambda \bullet^0 P]_{x_0} = P(\lambda_{x_0})=P\big((\alpha_{x_0}, \Sigma_{x_0}, \beta_{x_0})\big)= \alpha_{x_0}
$$
and for any pair of $x_0,x_0'$
$$
[\lambda \bullet^0P]_1^{x_0,x_0'}=\lambda_1^{x_0,x_0'}\bullet^0  P_1^{Lx_0,\hat{L}x_0'} = \alpha_1^{x_0,x_0'}
$$
thus $\lambda \bullet^0 P=\alpha$.\\

Similarly
$$
[\lambda \bullet^0 Q]_{x_0} = Q(\lambda_{x_0})=Q\big((\alpha_{x_0}, \Sigma_{x_0}, \beta_{x_0})\big)= \beta_{x_0}
$$
and
$$
[\lambda \bullet^0Q]_1^{x_0,x_0'}=\lambda_1^{x_0,x_0'}\bullet^0  Q_1^{Lx_0,\hat{L}x_0'} = \beta_1^{x_0,x_0'}
$$
thus $\lambda \bullet^0 Q=\beta$.\\

Finally
$$
[\lambda * \varepsilon]_{x_0} = \varepsilon(\lambda_{x_0})=\varepsilon\big((\alpha_{x_0}, \Sigma_{x_0}, \beta_{x_0})\big)= \Sigma_{x_0}
$$
and
$$
[\lambda *\varepsilon]_1^{x_0,x_0'}=\lambda_1^{x_0,x_0'}*  \varepsilon_1^{Lx_0,\hat{L}x_0'} = \Sigma_1^{x_0,x_0'}
$$

To conclude the proof we still need to prove uniqueness. But this will easily be achieved. Indeed the object
part of 2-morphism $\lambda$ satisfying the universal property is univocally determined by the fact that $P_0,Q_0$ and
$\varepsilon_0$ are projection, and once that is determined, uniqueness in dimension $n-1$ guaranties the
homs part.
\end{proof}

\begin{Lemma}\label{zzz:tech1}
The pair $\lambda=<\lambda_0,\lambda_1^{-,-}>$ is indeed a 2-morphism.
\end{Lemma}
\begin{proof}
We have to show that functoriality axioms for 2-morphisms are satisfied. Let us start with units axiom.\\

To this end let us fix an arbitrary object $x_0$ of $\X$. If we denote by $u(x_0)$ the identity $\Id{n-1}\to\X_1(x_0,x_0)$ then
\begin{eqnarray*}
u(x_0)\bullet^0 \lambda_1^{x_0,x_0} \bullet^0 Q_1^{Lx_0,\hat{L}x_0}   &\eq{i}&  u(x_0)\bullet^0\beta_1^{x_0,x_0} \\
   &\eq{ii}& id_{[\beta x_0]}  \\
   &\eq{iii}& id_{[Q(\lambda x_0)]} \\
   &\eq{iv}& id_{[\lambda x_0]}\bullet^0 Q_1^{Lx_0,\hat{L}x_0}
\end{eqnarray*}
where $(i)$ holds by (\ref{eqn:2PUP1}) above, $(ii)$ by unit functoriality of $\beta$, $(iii)$ by definition of
$\lambda_0$, $(iv)$ is a whiskering identity axiom. Similarly one can prove
$$
u(x_0)\bullet^0 \lambda_1^{x_0,x_0} \bullet^0 P_1^{Lx_0,\hat{L}x_0}   = id_{[\lambda x_0]}\bullet^0 P_1^{Lx_0,\hat{L}x_0}
$$
Moreover, for the formally analogous property w.r.t $*$-composition
\begin{eqnarray*}
u(x_0)\bullet^0 \lambda_1^{x_0,x_0} * \varepsilon_1^{Lx_0,\hat{L}x_0}   &=&  u(x_0)\bullet^0\Sigma_1^{x_0,x_0} \\
   &=& id_{[\Sigma x_0]}  \\
   &=& id_{[\varepsilon(\lambda x_0)]} \\
   &=& id_{[\lambda x_0]}* \varepsilon_1^{Lx_0,\hat{L}x_0}
\end{eqnarray*}
Notice that first composites are unambiguous by associativity axioms.

Calculations show that both $id_{[\lambda x_0]}$
and $u(x_0)\bullet^0 \lambda_1^{x_0,x_0}$  satisfy equations prescribed by the universal property, hence by uniqueness
they must be equal, and unit axiom is proved.\\

Turning to composition coherence, let objects $x_0,x_0',x_0''$ be given. Then\\

$$
\Big(
(\lambda_1^{x_0,x_0'}\circ \hat{L}_1^{x_0',x_0''})
\bullet^1
(L_1^{x_0,x_0'}\circ\lambda_1^{x_0',x_0''})
\Big)
\bullet^0 Q_1^{Lx_0\hat{L}x_0''}=
$$
\begin{eqnarray*}
&\eq{i}&
\Big(
(\lambda_1^{x_0,x_0'}\circ \hat{L}_1^{x_0',x_0''})
\bullet^0 Q_1^{Lx_0,\hat{L}x_0''}
\Big)
\bullet^1
\Big(
(L_1^{x_0,x_0'}\circ\lambda_1^{x_0',x_0''})
\bullet^0 Q_1^{Lx_0,\hat{L}x_0''}
\Big)
\\
&\eq{ii}&
\Big(
(\lambda_1^{x_0,x_0'}\bullet^0Q_1^{Lx_0,\hat{L}x_0'})
\circ
(\hat{L}_1^{x_0',x_0''}\bullet^0Q_1^{\hat{L}x_0',\hat{L}x_0''})
\Big)
\bullet^1
\\
&&\bullet^1\Big(
(L_1^{x_0,x_0'}\bullet^0 Q_1^{Lx_0,Lx_0'})
\circ
(\lambda_1^{x_0',x_0''}\bullet^0Q_1^{Lx_0'\hat{L}x_0''})
\Big)
\\
&\eq{iii}& (\beta_1^{x_0,x_0'}\circ \hat{N}_1^{x_0',x_0''})\bullet^1(N_1^{x_0,x_0'}\circ\beta_1^{x_0',x_0''})
\\
&\eq{iv}&(-\circ-)\bullet^0\beta_1^{x_0,x_0''}  \\
&\eq{v}&(-\circ-)\bullet^0\lambda_1^{x_0,x_0''}\bullet^0Q_1^{x_0,x_0''}
\end{eqnarray*}
where $(i)$ holds by (sesqui)functoriality of $-\bullet^0 Q_1^{Lx_0\hat{L}x_0''}$, $(ii)$ by functoriality
w.r.t. $\circ$-composition of $Q$, $(iii)$ and $(v)$ by (\ref{eqn:2PUP1}), $(iv)$ by functoriality
w.r.t. composition  of $\beta$.

Similarly one can prove
$$
\Big(
(\lambda_1^{x_0,x_0'}\circ \hat{L}_1^{x_0',x_0''})
\bullet^1
(L_1^{x_0,x_0'}\circ\lambda_1^{x_0',x_0''})
\Big)
\bullet^0 P_1^{Lx_0\hat{L}x_0''}=
(-\circ-)\bullet^0\lambda_1^{x_0,x_0''}\bullet^0P_1^{x_0,x_0''}
$$
\vskip6ex
Now let us compute
$$
\Big(
(\lambda_1^{x_0,x_0'}\circ \hat{L}_1^{x_0',x_0''})
\bullet^{1}
(L_1^{x_0,x_0'}\circ\lambda_1^{x_0',x_0''})
\Big)
*\varepsilon_1^{Lx_0,\hat{L}x_0''}
$$
by $*$-functoriality this equals to
\begin{changemargin}{-15ex}{-15ex}
\begin{equation}\label{eqn:monster}
\begin{array}{c}
\underbrace{\Big(
(\lambda_1^{x_0,x_0'}\circ \hat{L}_1^{x_0',x_0''})*\varepsilon^{Lx_0,\hat{L}x_0''}
\Big)}
\bullet^1
\Big(
(L_1^{x_0,x_0'}\circ \lambda_1^{x_0',x_0''})\bullet^0 ([PF]_1^{Lx_0,\hat{L}x_0''}\circ\hat{\omega}x_0'')
\Big)
\\ \bullet^2  \phantom{\Big[]} \\
\Big(
(\lambda_1^{x_0,x_0'}\circ \hat{L}_1^{x_0',x_0''})\bullet^0(\omega x_0\circ [QG]_1^{Lx_0,\hat{L}x_0''})
\Big)
\bullet^1
\Big(
(L_1^{x_0,x_0'}\circ \lambda_1^{x_0',x_0''})*\varepsilon_1^{Lx_0,\hat{L}x_0''}
\Big)
\end{array}
\end{equation}
\end{changemargin}
In order to simplify the expression, let us analyze first under-braced one.
This can be re-written explicitly and processed by $*$-associativity
$$
(\lambda_1^{x_0,x_0'}\circ \hat{L}_1^{x_0',x_0''})*\varepsilon^{Lx_0,\hat{L}x_0''}=
$$
\begin{eqnarray*}
   &=& \Big((\lambda_1^{x_0,x_0'}\times \hat{L}_1^{x_0',x_0''})\bullet^0(-\circ-)\Big)*\varepsilon^{Lx_0,\hat{L}x_0''} \\
   &=& (\lambda_1^{x_0,x_0'}\times \hat{L}_1^{x_0',x_0''})*\Big((-\circ-)\bullet^0\varepsilon^{Lx_0,\hat{L}x_0''}\Big) \\
   &=& (\lambda_1^{x_0,x_0'}\times \hat{L}_1^{x_0',x_0''})*\Big(
(\varepsilon_1^{Lx_0,\hat{L}x_0'}\circ [QG]_1^{\hat{L}x_0',\hat{L}x_0''})
\bullet^1
([PF]_1^{Lx_0,\hat{L}x_0'}\circ \varepsilon_1^{\hat{L}x_0',\hat{L}x_0''})
\Big)
\end{eqnarray*}
where the last is given by composition coherence of $\varepsilon$.\\

Applying again $*$-functoriality, this turns to be
\begin{changemargin}{-10ex}{-10ex}
$$\begin{array}{c}
\Big(
(\lambda_1^{x_0,x_0'}\times\hat{L}_1^{x_0',x_0''})
*
(\varepsilon_1^{Lx_0,\hat{L}x_0'}\circ[QG]_1^{\hat{L}x_0',\hat{L}x_0''})
\Big)
\bullet^1
\Big(
((L_1^{x_0,x_0'}\circ\lambda x_0')\times\hat{L}_1^{x_0',x_0''})
\bullet^0
([PF]_1^{Lx_0,\hat{L}x_0'}\circ\varepsilon_1^{\hat{L}x_0',\hat{L}x_0''})
\Big)
\\   \bullet^2    \\
\Big(
((\lambda x_0\circ \hat{L}_1^{x_0,x_0'})\times \hat{L}_1^{x_0',x_0''})
\bullet^0
(\varepsilon_1^{Lx_0,\hat{L}x_0'}\circ[QG]_1^{\hat{L}x_0',\hat{L}x_0''})
\Big)
\bullet^1
\Big(
(\lambda_1^{x_0,x_0'}\times \hat{L}_1^{x_0',x_0''})
*
([PF]_1^{Lx_0,\hat{L}x_0'}\circ\varepsilon_1^{\hat{L}x_0',\hat{L}x_0''})
\Big)
\end{array}
$$
\end{changemargin}
Now, all the second row is clearly an identity 3-morphism, being the $*$-composition on separate components
of a product, hence it can be canceled. What remains can be re-written as
$$
\Big(
(\lambda_1^{x_0,x_0'}*\varepsilon_1^{Lx_0,\hat{L}x_0'})\circ[\hat{L}QG]_1^{x_0',x_0''}
\Big)
\bullet^1
\Big(
[LPF]_1^{x_0,x_0'}\circ F(P(\lambda x_0'))\circ(\hat{L}_1^{x_0',x_0''}\bullet^0\varepsilon_1^{\hat{L}x_0',\hat{L}x_0''})
\Big)
$$
and with the help of (\ref{eqn:2PUP3}) this is simply
$$
(\Sigma_1^{x_0,x_0'}\circ[\hat{N}G]_1^{x_0',x_0''})
\bullet^1
([MF]_1^{x_0,x_0'}\circ F(P(\lambda x_0'))\circ\hat{\omega}_1^{x_0',x_0''})
$$
Substituting the first line of (\ref{eqn:monster}) becomes a triple $\bullet^1$-composition
\begin{changemargin}{-10ex}{-10ex}
\begin{equation}\label{eqn:zzz:pspsp}
\Big(\Sigma_1^{x_0,x_0'}\circ[\hat{N}G]_1^{x_0',x_0''}\Big)
\bullet^1
\Big([MF]_1^{x_0,x_0'}\circ F(P(\lambda x_0'))\circ\hat{\omega}_1^{x_0',x_0''}
\Big)
\bullet^1
\Big(
(L_1^{x_0,x_0'}\circ \lambda_1^{x_0',x_0''})\bullet^0 ([PF]_1^{Lx_0,\hat{L}x_0''}\circ\hat{\omega}x_0'')
\Big)
\end{equation}
\end{changemargin}
Now  $F(P(\lambda x_0'))=F(\alpha x_0')=[\alpha\bullet^0F]x_0'$, furthermore
$PF(L_1^{x_0,x_0'}\circ \lambda_1^{x_0',x_0''})=PF(L_1^{x_0,x_0'})\circ PF(\lambda_1^{x_0',x_0''})=
[MF]_1^{x_0,x_0'}\circ [\alpha\bullet^0F]_1^{x_0',x_0''}$.
Hence (\ref{eqn:zzz:pspsp}) is equal to
\begin{changemargin}{-10ex}{-10ex}
$$
\Big(\Sigma_1^{x_0,x_0'}\circ[\hat{N}G]_1^{x_0',x_0''}\Big)
\bullet^1
\Big([MF]_1^{x_0,x_0'}\circ [\alpha\bullet^0F]x_0'\circ\hat{\omega}_1^{x_0',x_0''}\Big)
\bullet^1
\Big(
  [MF]_1^{x_0,x_0'}\circ [\alpha\bullet^0F]_1^{x_0',x_0''}    \circ\hat{\omega}x_0''
\Big)
$$
\end{changemargin}
{\em i.e.}
$$
\Big(\Sigma_1^{x_0,x_0'}\circ[\hat{N}G]_1^{x_0',x_0''}\Big)
\bullet^1
\Big([MF]_1^{x_0,x_0'}\circ
\Big(
\big(
[\alpha\bullet^0F]_{x_0'}\circ\hat{\omega}_1^{x_0',x_0''}
\big)
\bullet^1
\big(
[\alpha\bullet^0F]_1^{x_0',x_0''}    \circ\hat{\omega}x_0''
\big)
\Big)
\Big)
$$
By definition of 1-composition of 2-morphisms this is also
$$
\Big(\Sigma_1^{x_0,x_0'}\circ[\hat{N}G]_1^{x_0',x_0''}\Big)
\bullet^1
\Big([MF]_1^{x_0,x_0'}\circ
\big[
(\alpha\bullet^0F)\bullet^1\hat{\omega}
\big]_1^{x_0',x_0''}
\Big)
$$
Symmetrical calculations can be made on the second 2-composite of (\ref{eqn:monster}), giving the composite
$$
\begin{array}{c}
\Big(\Sigma_1^{x_0,x_0'}\circ[\hat{N}G]_1^{x_0',x_0''}\Big)
\bullet^1
\Big([MF]_1^{x_0,x_0'}\circ
\big[
(\alpha\bullet^0F)\bullet^1\hat{\omega}
\big]_1^{x_0',x_0''}
\Big)
\\ \bullet^2   \\
\Big(
[\omega\bullet^1(\beta\bullet^0G)]_1^{x_0,x_0'}\circ [\hat{N}G]_1^{x_0',x_0''}
\Big)
\bullet^1
\Big(
[MF]_1^{x_0,x_0'}\circ\Sigma_1^{x_0',x_0''}
\Big)
  \end{array}
$$
By composition coherence for $\Sigma$ we get
$$
(-\circ-)\bullet^0\Sigma_1^{x_0,x_0''}
$$
and by (\ref{eqn:2PUP3}) again
$$
(-\circ-)\bullet^0\lambda_1^{x_0,x_0''}*\varepsilon_1^{Lx_0,\hat{L}x_0''}
$$

Concluding, for every choice of three objects $x_0,x_0',x_0''$ of $\X$ the following three equations hold
$$
\left\{
  \begin{array}{c}
\Big(
(\lambda_1^{x_0,x_0'}\circ \hat{L}_1^{x_0',x_0''})
\bullet^1
(L_1^{x_0,x_0'}\circ\lambda_1^{x_0',x_0''})
\Big)
\bullet^0 Q_1^{Lx_0\hat{L}x_0''}=
(-\circ-)\bullet^0\lambda_1^{x_0,x_0''}\bullet^0Q_1^{x_0,x_0''}
\\
\Big(
(\lambda_1^{x_0,x_0'}\circ \hat{L}_1^{x_0',x_0''})
\bullet^1
(L_1^{x_0,x_0'}\circ\lambda_1^{x_0',x_0''})
\Big)
\bullet^0 P_1^{Lx_0\hat{L}x_0''}=
(-\circ-)\bullet^0\lambda_1^{x_0,x_0''}\bullet^0P_1^{x_0,x_0''}
\\
\Big(
(\lambda_1^{x_0,x_0'}\circ \hat{L}_1^{x_0',x_0''})
\bullet^{1}
(L_1^{x_0,x_0'}\circ\lambda_1^{x_0',x_0''})
\Big)
*\varepsilon_1^{Lx_0,\hat{L}x_0''}
=(-\circ-)\bullet^0\lambda_1^{x_0,x_0''}*\varepsilon_1^{Lx_0,\hat{L}x_0''}
\end{array}\right.
$$

Hence  both $\Big(
(\lambda_1^{x_0,x_0'}\circ \hat{L}_1^{x_0',x_0''})
\bullet^1
(L_1^{x_0,x_0'}\circ\lambda_1^{x_0',x_0''})
\Big)
$
and $(-\circ-)\bullet^0\lambda_1^{x_0,x_0''}$  satisfy equations prescribed by universal property, hence by uniqueness
they must be equal, and composition coherence is proved.\\
\end{proof}

\section{$\Omega$ and a second definition of $\pi_1$}
We use the $h$-2pullback defined above to give an alternative description of the
sesqui-functor $\pi_1^{(n)}$. A key observation is the analogy between the hom-$(n-1)$-groupoid of
a $n$-groupoid $\C$ and the paths of a topological space.\\

Given a $n$-groupoid ($n$-category) $\C$ and two objects $c_0,c_0'$, we define
$\PP_{c_0,c_0'}(\C)$ by means of the following $h$-pullback:
\begin{equation}\label{dgm:hom-pb}
\xymatrix{
\PP_{c_0,c_0'}(\C)
\ar[r]^{!}\ar[d]_{!}
&\Id{n}\ar[d]^{[c_0']}
\\
\Id{n}\ar@{}[ur]|(.3){}="1"|(.7){}="2"
\ar[r]_{[c_0]}
&\C\ar@{=>}"1";"2"_{\varepsilon_{\C}^{c_0,c_0'}}
}
\end{equation}
This definition easily extends to  morphisms. In fact for  $F:\C\to\D$ one defines
$$
\PP_{c_0,c_0'}(F):\ \PP_{c_0,c_0'}(\C)\to\PP_{c_0,c_0'}(\D)
$$
by means of the universal property of $h$-pullbacks yielding  $\PP_{c_0,c_0'}(\D)$, for the four-tuple
$$
<\PP_{c_0,c_0'}(\C),!,!,\varepsilon_{\C}^{c_0,c_0'}\bullet^0 F>
$$
$$
\xymatrix{
\PP_{c_0,c_0'}(\C)
\ar[rrr]^{!}_(.7){}="2"
\ar[ddd]_{!}^(.7){}="1"
\ar@{}"1";"2"|(.4){}="a1"|(.4){}="a1"|(.6){}="a2"
\ar@2"a1";"a2"^{\varepsilon_{\C}^{c_0,c_0'}}
&&&\Id{n}
\ar[ddl]_{[c_0']}
\ar[ddd]^{[Fc_0']}
\\
\\
&&\C\ar[dr]^{F}
\\
\Id{n}
\ar[urr]^{[c_0]}
\ar[rrr]_{[Fc_0]}
&&&\D
}
$$
This makes $\PP_{c_0,c_0'}(-)$ ``somehow'' functorial: in fact for $H:\D\to\E$,
$$
\PP_{c_0,c_0'}(F)\bullet^0\PP_{c_0,c_0'}(H)=\PP_{c_0,c_0'}(F\bullet^0 H),\qquad \PP_{c_0,c_0'}(id_{\C})=id_{\PP_{c_0,c_0'}(\C)}
$$
\vskip3ex
Unfortunately this does not extend straightforward to 2-morphisms. In fact for a pair of parallel morphisms
$F,G:\C\to\D$, $\PP_{c_0,c_0'}(F)$ and $\PP_{c_0,c_0'}(G)$ are no longer parallel, this making it difficult
to extend $\PP_{c_0,c_0'}(-)$ to natural $n$-transformations.
\vskip3ex
Indeed in applying the same argument as for defining $\PP_{c_0,c_0'}(-)$ on morphisms, the corresponding
diagram (shown below) suggests to consider the 0-composition of 2-morphisms
$$
\varepsilon_{\C}^{c_0,c_0'}*\alpha:\xymatrix{\varepsilon_{\C}^{c_0,c_0'}\setminus\alpha\ar@3[r]&\varepsilon_{\C}^{c_0,c_0'} / \alpha}
$$
where as usual
$$
\varepsilon_{\C}^{c_0,c_0'}\setminus\alpha=([c_0]\bullet^0\alpha)\bullet^1(\varepsilon_{\C}^{c_0,c_0'}\bullet^0G)=[\alpha_{c_0}]\bullet^1(\varepsilon_{\C}^{c_0,c_0'}\bullet^0G)
$$
and
$$
\varepsilon_{\C}^{c_0,c_0'}/ \alpha=(\varepsilon_{\C}^{c_0,c_0'}\bullet^0F)\bullet^1([c_0']\bullet^0\alpha)=(\varepsilon_{\C}^{c_0,c_0'}\bullet^0F)\bullet^1[\alpha_{c_0'}]
$$
$$
\xymatrix{
\PP_{c_0,c_0'}(\C)
\ar[rrr]^{!}_(.7){}="2"
\ar[ddd]_{!}^(.7){}="1"
\ar@{}"1";"2"|(.4){}="a1"|(.6){}="a2"
\ar@2"a1";"a2"^{\varepsilon_{\C}^{c_0,c_0'}}
&&&\Id{n}
\ar[ddl]_{[c_0']}
\ar[ddd]^{[Gc_0']}
\\
\\
&&\C\ar@/^2.5ex/[dr]^(.4){G}^{}="4"\ar@/_2.5ex/[dr]_(.4){F}^{}="3"
\ar@{}"3";"4"|(.2){}="a3"|(.8){}="a4"
\ar@2"a3";"a4"^{\alpha}
\\
\Id{n}
\ar[urr]^{[c_0]}
\ar[rrr]_{[Fc_0]}
&&&\D
}
$$

Hence we can consider the four-ples
$$
\begin{array}{cc}
<\PP_{c_0,c_0'}(\C),!,!,\varepsilon_{\C}^{c_0,c_0'}\setminus\alpha> & <\PP_{c_0,c_0'}(\C),!,!,\varepsilon_{\C}^{c_0,c_0'}/ \alpha>
\\&\\
\xymatrix@C=10ex@R=10ex{
\PP_{c_0,c_0'}(\C)
\ar[r]^{!}
\ar[d]_{!}
&\Id{n}\ar[d]^{[Gc_0']}
\\
\ar@{}[ur]|(.25){}="1"|(.75){}="2"\ar@2"1";"2"|{\varepsilon_{\C}^{c_0,c_0'}\setminus\alpha}
\Id{n}
\ar[r]_{[Fc_0]}
&\D
}
&
\xymatrix@C=10ex@R=10ex{
\PP_{c_0,c_0'}(\C)
\ar[r]^{!}
\ar[d]_{!}
&\Id{n}\ar[d]^{[Gc_0']}
\\
\ar@{}[ur]|(.25){}="1"|(.75){}="2"\ar@2"1";"2"|{\varepsilon_{\C}^{c_0,c_0'}/ \alpha}
\Id{n}
\ar[r]_{[Fc_0]}
&\D}
\end{array}
$$
together with $id_{!}:\,!\Rightarrow\,!$ (taken two times) and the 3-morphism $\varepsilon_{\C}^{c_0,c_0'}*\alpha$.
Applying the universal property of $h$-2pullbacks we get a 2-morphism
$$
\PP_{c_0,c_0'}(\alpha):\ \PP_{[\alpha_{c_0}]}\circ \PP_{c_0,c_0'}(G))\Rightarrow \PP_{c_0,c_0'}(F)\circ \PP_{[\alpha_{c_0'}]}:\ \PP_{c_0,c_0'}(C)\to \PP_{Fc_0,Gc_0'}(\D)
$$
such that
\begin{eqnarray}
\nonumber  \PP_{c_0,c_0'}(\alpha)\bullet^0\,!&=&id_{!} \\
\nonumber  \PP_{c_0,c_0'}(\alpha)\bullet^0\,!&=&id_{!} \\
\PP_{c_0,c_0'}(\alpha)*\varepsilon_{\D}^{Fc_0,Gc_0'} &=& \varepsilon_{\C}^{c_0,c_0'}*\alpha\label{eqn:pi_1_on_2morph}
\end{eqnarray}
Notice that we have denoted by $\PP_{[\alpha_{c_0}]}\circ \PP_{c_0,c_0'}(G))$ and $\PP_{c_0,c_0'}(F)\circ \PP_{[\alpha_{c_0'}]}$
the morphisms obtained by applying one-dimensional the universal property to $\varepsilon_{\C}^{c_0,c_0'}\setminus\alpha$
and $\varepsilon_{\C}^{c_0,c_0'}/ \alpha$ respectively. Therefore the symbol $\circ$ involved
should be considered just a typographical suggestion. Indeed it can be shown that it is a 0-composition of morphisms, but this
would lead us far from the point. Furthermore it is inessential with respect to our purposes. For this reasons
its further developing is left to the curious reader.
\vskip5ex
Purpose of the rest of the section is to give a characterization of $\pi_1$ as a consequence of
the following
\begin{Theorem}\label{thm:pi_1}
For every n-category $\C$, and every two objects $c_0,c_0'$ in $\C$,
there exists a canonical isomorphism
$$
\mathfrak{S}_{\C}^{c_0,c_0'} :D(\C_1(c_0,c_0'))\to\PP_{c_0,c_0'}(\C)
$$
In the case of pointed n-groupoids, this gives a natural isomorphism with components
$$
\mathfrak{S}_{\C}^{*,*} :D(\pi_1(\C))\to\Omega(\C)
$$
where $\Omega(\C)=\PP_{*,*}(\C)$
\end{Theorem}
We start by making explicit $h$-pullback of (\ref{dgm:hom-pb}), but first we need
to be more precise on units.\\

\begin{Remark}
Let $\C$ be a n-category. For a fixed object $c_0$ of $\C$, let us
consider the unit (n-1)functor given by the n-category structure of $\C$:
$$
^{\C}\!u^0(c_0):\Id{n-1}\to \C_1(c_0,c_0)
$$
We can make it explicit as a pair
\begin{eqnarray*}
  [^{\C}\! u^0 (c_0) ]_0 &:& *  \mapsto id(c_0)\quad \in [\C_1(c_0,c_0')]_0 \\
  {}[^{\C}\! u^0 (c_0) ]_1 &:&\Id{n-2}\mapsto [\C_1(c_0,c_0')]_1(id(c_0),id(c_0))
\end{eqnarray*}
Now, by functoriality we get the interchange
$$[^{\C}\! u^0 (c_0) ]_1=\,{}^{\C}u^1(id(c_0))={}^{\C_1(c_0,c_0)}u^0(id(c_0))
$$
and this allows the following explicit definition:
$$
{}^{\C}u^0(c_0)=<u^{(1)}(c_0),u^{(2)}(c_0),\cdots,u^{(n)}(c_0)>
$$
where $u^{(k)}(c_0)$ is the identity $k$-cell over $c_0$.
\end{Remark}
In the rest of this section, in order to simplify notation, the n-category
$\PP_{c_0,c_0'}(\C)$ will be denoted by $\Q$.

\begin{Proposition}\label{zzz:propx}
Given the $h$-pullback of $n$-categories
\begin{equation}\label{dgm:hom-pb-Q}
\xymatrix{
\Q
\ar[r]^{!}\ar[d]_{!}
&\Id{n}\ar[d]^{[c_0']}
\\
\Id{n}\ar@{}[ur]|(.3){}="1"|(.7){}="2"
\ar[r]_{[c_0]}
&\C\ar@{=>}"1";"2"_{\varepsilon}
}
\end{equation}
the hom-$(n-k)$-category
$$
\Q_k\left(\Big(u^{(k-1)}(*),\xymatrix{c_{k-1}\ar[r]^{c_k}&c_{k-1}'},u^{(k-1)}(*)\Big),\Big(u^{(k-1)}(*),\xymatrix{c_{k-1}\ar[r]^{c'_k}&c_{k-1}'},u^{(k-1)}(*)\Big)\right)
$$
is well defined and it is given by $h$-pullback over the diagram
$$
\xymatrix{
&\Id{n-k}\ar[d]^{[c'_k]}
\\
\Id{n-k}\ar[r]_(.35){[c_k]}
&\C_k(c_{k-1},c'_{k-1})
}
$$
\end{Proposition}
\begin{proof}
By finite induction over $k$.\\

\framebox[1.1\width]{$k=1$}\\

We recall the definition of $h$-pullback:\\

$\Q_0$ is given by the limit in \textbf{Set}

$$
\raisebox{10ex}{\xymatrix{
&&\Q_0\ar@{.>}[lld]_{!}\ar@{.>}[d]^{\varepsilon_0}\ar@{.>}[drr]^{!}
\\
\{*\}=[\I]_0\ar[dr]_{[c_0]_0}
&&[\C_1]_0\ar[dl]^d\ar[dr]_c
&&[\I]_0=\{*\}\ar[dl]^{[c_0']_0}
\\
&\C_0&&\C_0}}
$$

$\Q_1((*,c_1,*),(*,c_1',*))$ is given by a $h$-pullback

$$
\raisebox{10ex}{
\xymatrix{
\Q_1(\diamond,\diamond)\ar[rr]^{!}\ar[dd]_{!}
&&\Id{n-1}\ar[d]^{[c_0']_1}\ar@{}[ddll]|(.4){}="1"|(.6){}="2"
\\
&&\C_1(c_0',c_0')\ar[d]^{c_1\circ -}
\\
\Id{n-1}\ar[r]_(.4){[c_0]_1}
&\C_1(c_0,c_0)\ar[r]_{-\circ c_1'}
&\C_1(c_0,c_0')\ar@{=>}"1";"2"^{\varepsilon_1^{\diamond,\diamond}}}}
=
\raisebox{10ex}{
\xymatrix@R=14ex{
\Q_1(\diamond,\diamond)
\ar[r]^{!}\ar[d]_{!}
&\Id{n-1}\ar[d]^{[c_1']}
\\
\Id{n-1}\ar@{}[ur]|(.4){}="1"|(.6){}="2"
\ar[r]_{[c_1]}
&\C_1(c_0,c_0')\ar@{=>}"1";"2"_{\varepsilon_1^{\diamond,\diamond}}}
}$$
with the symbol $\diamond$ when the substitute is clear from the context. \\

\framebox[1.1\width]{$k>1$}\\

Induction hypothesis gives the following definition for
$$
\Q_{k-1}\Big((*,\xymatrix{c_{k-2}\ar[r]^{c_{k-1}}&c_{k-2}'},*),(*,\xymatrix{c_{k-2}\ar[r]^{c_{k-1}'}&c_{k-2}'},*)\Big)
$$
$$
\xymatrix{
\Q_{k-1}(\diamond,\diamond)
\ar[r]^{!}\ar[d]_{!}
&\Id{n-k+1}\ar[d]^{[c_{k-1}']}
\\
\Id{n-k+1}\ar@{}[ur]|(.4){}="1"|(.6){}="2"
\ar[r]_(.35){[c_{k-1}]}
&\C_{k-1}(c_{k-2},c_{k-2}')\ar@{=>}"1";"2"_{\varepsilon_{k-1}^{\diamond,\diamond}}}
$$
More explicitly, we get the following set-theoretical limit
$$
\xymatrix@C=-4ex{
&&[\Q_{k-1}(\diamond,\diamond)]_0\ar@{.>}[lld]_{!}\ar@{.>}[d]^{[\varepsilon_{k-1}^{\diamond,\diamond}]_0}\ar@{.>}[drr]^{!}
\\
\{*\}=[\Id{n-k-1}]_0\ar[dr]_{[c_0]_0}
&&[\C_{k-1}(c_{k-2},c_{k-2}')]_1\ar[dl]^d\ar[dr]_c
&&[\Id{n-k-1}]_0=\{*\}\ar[dl]^{[c_0']_0}
\\
&[\C_{k-1}(c_{k-2},c_{k-2}')]_0&&[\C_{k-1}(c_{k-2},c_{k-2}')]_0}
$$
{\em i.e.} the set $\{(*,\xymatrix{c_{k-2}\ar[r]^{c_{k-1}}&c_{k-2}'},*)\}$.

Hence $\Q_{k}((*,\xymatrix{c_{k-1}\ar[r]^{c_{k}}&c_{k-1}'},*),(*,\xymatrix{c_{k-1}\ar[r]^{c_{k}'}&c_{k-1}'},*))$
has (inductively) well defined domain and codomain, namely $(*,\xymatrix{c_{k-1}\ar[r]^{c_{k}}&c_{k-1}'},*)$ and
$(*,\xymatrix{c_{k-1}\ar[r]^{c_{k}'}&c_{k-1}'},*)$ are legitimate objects of a $\Q_{k-1}(\diamond,\diamond)$.
By definition of $h$-pullback, we can spell it out:
\begin{changemargin}{-10ex}{10ex}
$$
\xymatrix{
\Q_k((*,c_k,*),(*,c_k',*))\ar[rr]^{!}\ar[dd]_{!}
&&\left[\Id{n-k+1}\right]_1\ar[d]^{[c_{k-1}']_1}\ar@{}[ddll]|(.4){}="1"|(.6){}="2"
\\
&&\left[\C_{k-1}(c_{k-2},c_{k-2}')\right]_1(c_{k-1'},c_{k-1}')\ar[d]^{c_k\circ -}
\\
\left[\Id{n-k+1}\right]_1\ar[r]_(.4){[c_{k-1}]_1}
&\left[\C_{k-1}(c_{k-2},c_{k-2}')\right]_1(c_{k-1},c_{k-1})\ar[r]_{-\circ c_k'}
&\left[\C_{k-1}(c_{k-2},c_{k-2}')\right]_1(c_{k-1'},c_{k-1}')
\ar@{=>}"1";"2"_{\left[\varepsilon_{k-1}^{\diamond,\diamond}\right]_1^{(*,c_k,*),(*,c_k',*)}}}
$$
\end{changemargin}
that may be rewritten
$$
\xymatrix@R=14ex@C=10ex{
\Q_k((*,c_k,*),(*,c_k',*))
\ar[r]^{!}\ar[d]_{!}
&\Id{n-k}\ar[d]^{[c_k']}
\\
\Id{n-k}\ar@{}[ur]|(.4){}="1"|(.6){}="2"
\ar[r]_{[c_k]}
&\C_k(c_{k-1},c_{k-1}')\ar@{=>}"1";"2"^{\varepsilon_k^{(*,c_k,*),(*,c_k',*)}}
}$$
\end{proof}

Proof of  \emph{Proposition \ref{zzz:propx}} gives immediately the following
\begin{Corollary}\label{zzz:cor}
The 2-morphism $\varepsilon$ is given explicitly by
$$
\varepsilon=<\varepsilon_0, [\varepsilon_1^{-,-}]_0, \dots,[\varepsilon_{n-1}^{-,-}]_0, = >
$$
where
\begin{eqnarray*}
  [\varepsilon_k^{(*,c_{k-1},*),(*,c_{k-1}',*)}]_0 : \Q_k((*,c_{k-1},*),(*,c_{k-1}',*))&\to& \C_k(c_{k-1},c_{k-1}') \\
  (*,\xymatrix{c_{k-1}\ar[r]^{c_k}& c_{k-1}'},*) &\mapsto& c_k
\end{eqnarray*}
\end{Corollary}

Next Corollary states that $h$-pullbacks along two constants are n-discrete.

\begin{Corollary}
With notation as above,
$$
D(\pi_0(\Q))=\Q
$$
\end{Corollary}
\begin{proof}
It suffices to let $k=n$ in the above. $\Q_n((*,c_n,*),(*,c_n',*))$ is given by the following
pullback in \textbf{Set}:
$$
\xymatrix{
\Q_n((*,c_n,*),(*,c_n',*))\ar[r]^(.6){!}\ar[d]_{!}
&\{*\}\ar[d]^{[c_n']}
\\
\{*\}\ar[r]_{[c_n]}
&\C_n(c_{n-1},c_{n-1}')
}
$$
Hence, $\Q_n((*,c_n,*),(*,c_n',*))=\{*\}$ if $c_n=c_n'$, the \textit{empty-set}  otherwise.
\end{proof}

\subsection{$0$-composition in $\PP_{c_0,c_0'}(\C)$}
We have to describe $0$-composition in $\Q=\PP_{c_0,c_0'}(\C)$ , being all $k$-compositions
(with $k>0$) implicit in the inductive definition of $\Q$.

Let us start with functoriality of 2-morphism $\varepsilon$ with respect to $0$-composition.
For every triple $(*,c_1,*), (*,c_1',*)$ and $(*,c_1'',*)$,  diagram
 (\ref{ntrans:composition}) may be written
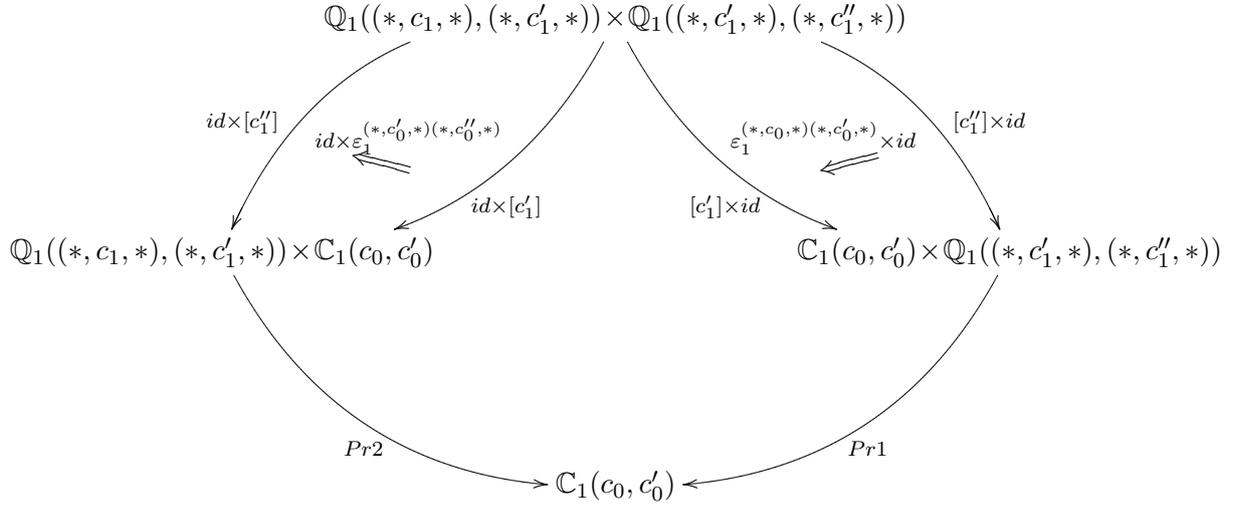
\begin{changemargin}{-10ex}{10ex}
$$
\xymatrix@C=-10ex@R=15ex{
&\Q_1((*,c_1,*), (*,c_1',*))\!\times\!\Q_1((*,c_1',*), (*,c_1'',*))
\ar@/_7ex/[dl]_(.75){id\times[c_1'']}^(.75){}="a2"   \ar@/^6ex/[dl]^(.45){id\times[c_1']}_(.45){}="a1"
\ar@/_6ex/[dr]_(.45){[c_1']\times id}^(.45){}="b2"   \ar@/^7ex/[dr]^(.75){[c_1'']\times id}_(.75){}="b1"
\\
\Q_1((*,c_1,*), (*,c_1',*))\!\times\! \C_1(c_0,c_0')\ar@/_6ex/[dr]_{Pr2}
&&\C_1(c_0,c_0')\!\times\! \Q_1((*,c_1',*), (*,c_1'',*))\ar@/^6ex/[dl]^{Pr1}
\\
&\C_1(c_0,c_0')
\ar@{}"a1";"a2"|(.35){}="aa1"|(.65){}="aa2"
\ar@{}"b1";"b2"|(.35){}="bb1"|(.65){}="bb2"
\ar@{=>}"aa1";"aa2"_{id\times\varepsilon_1^{(*,c_0',*)(*,c_0'',*)}}
\ar@{=>}"bb1";"bb2"_{\varepsilon_1^{(*,c_0,*)(*,c_0',*)}\times id}
}$$
\end{changemargin}
that reduces to
\begin{equation}\label{zzz1}
\raisebox{10ex}
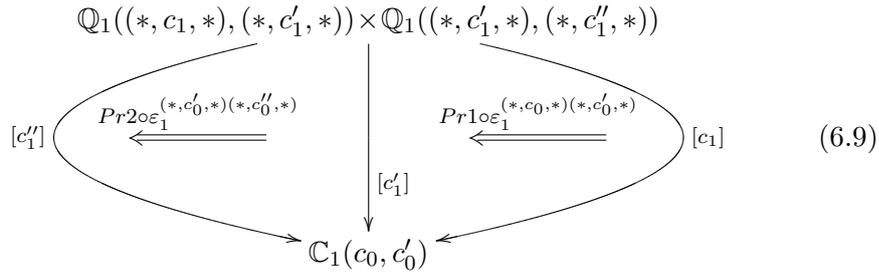
{\xymatrix@R=15ex{
\Q_1((*,c_1,*), (*,c_1',*))\!\times\!\Q_1((*,c_1',*), (*,c_1'',*))
\ar@/_25ex/[d]_{[c_1'']}="3"
\ar[d]^(.7){[c_1']}|{}="2"
\ar@/^25ex/[d]^{[c_1]}="1"
\\
\C_1(c_0,c_0')
\ar@{}"1";"2"|(.3){}="a1"|(.7){}="a2"
\ar@{}"2";"3"|(.3){}="b1"|(.7){}="b2"
\ar@{=>}"a1";"a2"_{Pr1\circ\varepsilon_1^{(*,c_0,*)(*,c_0',*)}}
\ar@{=>}"b1";"b2"_{Pr2\circ\varepsilon_1^{(*,c_0',*)(*,c_0'',*)}}
}}
\end{equation}
On the other side, the last term of equality (\ref{ntrans:composition}) may be written
\begin{equation}\label{zzz2}
\raisebox{10ex}{\xymatrix@R=15ex{
\Q_1((*,c_1,*), (*,c_1',*))\!\times\!\Q_1((*,c_1',*), (*,c_1'',*))
\ar[d]^{{}^{\Q}\circ^0}
\\
\Q_1((*,c_1,*), (*,c_1'',*))
\ar@/_10ex/[d]_{[c_1'']}="3"
\ar@/^10ex/[d]^{[c_1]}="1"
\ar@{}"1";"3"|(.3){}="b1"|(.7){}="b2"
\ar@{=>}"b1";"b2"_{\varepsilon_1^{(*,c_0,*)(*,c_0'',*)}}
\\
\C_1(c_0,c_0')}}
\end{equation}
Comparing diagrams (\ref{zzz1}) and (\ref{zzz2}),  the very definition of $0$-compositions in
$h$-pullbacks  proves the following
\begin{Proposition}\label{zzz:comp}
Let $c_0,c_0',c_0'':c_0\to c_0'$ be fixed in $\C$. Given
$$
\xymatrix@C=4ex{c_k:c_1\ar@{==>}[r]&c_1'},\quad \xymatrix@C=4ex{c_k':c_1\ar@{==>}[r]&c_1'}
$$
with $1<k\leq n$, the following equation holds:
$$
(*,c_k,*)\,{}^{\Q}\circ^0 (*,c_k,*) = (*,c_k\,{}^{\C}\circ^1  c_k',*)
$$
\end{Proposition}
\begin{Notation}
We use the notation $\xymatrix@C=4ex{c_k:c_h\ar@{==>}[r]&c_h'}$ with $h<k$,
to mean that $k$-cell $c_k$ has $h$-domain $c_h$ and $h$-codomain $c_h'$, {\em i.e.} there exist
cells $c_{k-1},c_{k-1}'\dots c_{h},c_{h}'$ such that
$$
c_k:\ \ c_{k-1}\to c_{k-1}':\ c_{k-2}\to c_{k-2}':\cdots :\ c_{h+1}\to c_{h+1}':\ c_{h}\to c_{h}'
$$
\end{Notation}
\begin{proof} (of Proposition {\ref{zzz:comp}})
By definition of vertical composition of 2-morphisms and whiskering, diagram (\ref{zzz1}) gives
$$
\left[\left[\left(Pr1\circ^0 \varepsilon_1^{(*,c_1,*),(*,c_1',*)}\right)\circ^1\left(Pr2\circ^0 \varepsilon_1^{(*,c_1',*),(*,c_1'',*)}\right)\right]_{k-1}\right]_0\big((*,c_k,*),(*,c_k',*)\big)=
$$
\begin{eqnarray*}
   &=& \left[\varepsilon_k^{(*,c_1,*),(*,c_1',*)}\right]_0 \left([Pr1]_{k-1}\big((*,c_k,*),(*,c_k',*)\big) \right){}^{\C}\!\circ^1 \cdots\\
   & & \cdots{}^{\C}\!\circ^1 \left[\varepsilon_k^{(*,c_1',*),(*,c_1'',*)}\right]_0 \left([Pr2]_{k-1}\big((*,c_k,*),(*,c_k',*)\big)  \right) \\
   &=&\left[\varepsilon_k^{(*,c_1,*),(*,c_1',*)}\right]_0\big((*,c_k,*)\big)\ {}^{\C}\!\circ^1 \left[\varepsilon_k^{(*,c_1',*),(*,c_1'',*)}\right]_0\big((*,c_k',*)\big)\\
   &=& c_k\ {}^{\C}\!\circ^1 c_k'
\end{eqnarray*}
where the last equality holds by {\em Corollary \ref{zzz:cor}}. Next, diagram (\ref{zzz2}) gives
$$
\left[\left[\left[-\ {}^{\Q}\circ^0-\right]\circ^0 \varepsilon_1^{(*,c_1,*),(*,c_1'',*)}\right]_{k-1}\right]_0\big((*,c_k,*),(*,c_k',*)\big)=
$$
\begin{eqnarray*}
   &=&  \left[\varepsilon_k^{(*,c_1,*),(*,c_1'',*)}\right]_0\big((*,c_k,*)\ {}^{\Q}\circ^0(*,c_k',*)\big)\\
   &=&  \left[\varepsilon_k^{(*,c_1,*),(*,c_1'',*)}\right]_0\big((*,\bar{c}_k,*)\big)=  \bar{c}_k
\end{eqnarray*}
Then, by comparison we get
$$
(*,c_k,*)\ {}^{\Q}\circ^0(*,c_k',*) = (*,\bar{c}_k,*)
$$
if, and only if,
$$
\bar{c}_k= c_k \ {}^{\C}\circ^1 c_k'
$$
and this completes the proof.
\end{proof}

\subsection{$0$-units in $\PP_{c_0,c_0'}(\C)$}
We have to describe $0$-units in $\Q=\PP_{c_0,c_0'}(\C)$ , being all $k$-units
(with $k>0$) implicit in the inductive definition of $\Q$.

Let us start with functoriality of 2-morphism $\varepsilon$ with respect to $0$-units, for every
$(*,c_1,*)$ in $\Q_0$ one can consider
$$
\xymatrix{{}^{\Q}u^0((*,c_1,*)):\Id{n-1}\ar[r]&\Q_1((*,c_1,*),(*,c_1,*))}
$$
Unit coherence (\ref{ntrans:units}) is then the equality
$$
\raisebox{11ex}{\xymatrix@R=10ex{
\Id{n-1}\ar[d]^{{}^{\Q}u^0((*,c_1,*))}
\\
\Q_1((*,c_1,*),(*,c_1,*))
\ar@/_10ex/[d]_{[c_1]}="2"
\ar@/^10ex/[d]^{[c_1]}="1"
\\
\C_1(c_0,c_0')
\ar@{}"1";"2"|(.3){}="a1"|(.7){}="a2"
\ar@{=>}"a1";"a2"_{\varepsilon_1^{(*,c_1,*),(*,c_1,*)}}
}}
=
\raisebox{11ex}{
\xymatrix@R=24ex{
\Id{n-1}
\ar@/_8ex/[d]_{[c_1]}="2"
\ar@/^8ex/[d]^{[c_1]}="1"
\\
\C_1(c_0,c_0')
\ar@{}"1";"2"|(.3){}="a1"|(.7){}="a2"
\ar@{=>}"a1";"a2"_{id}
}}
$$
This comparison, with the explicit description of $\varepsilon$ given in \emph{Corollary \ref{zzz:cor}},
proves the following
\begin{Proposition}\label{zzz:unts}
Let $c_0:c_0\to c_0'$ be fixed in $\C$. For $1<k\leq n$, the following equation holds
$$
\left[{}^{\Q}u^0((*,c_1,*)) \right]_k = \left(*,\left[{}^{\C}u^1(c_1)\right],*\right)
$$
\end{Proposition}

\subsection{Comparison isomorphism $\mathfrak{S}$}
What we are going to state provides an extremely powerful tool in developing the theory.
\begin{Lemma}
Let $c_0,c_0'$ be objects of an  n-category $\C$.
The assignment
$$
\mathfrak{S}_{\C}^{c_0,c_0'}=\mathfrak{S}:D(\C_1(c_0,c_0'))\to \PP_{c_0,c_0'}(\C)=\Q
$$
given explicitly by
$$
\mathfrak{S}=<\mathfrak{S}_0, \mathfrak{S}_1, \dots,\mathfrak{S}_n>
$$
with
\begin{eqnarray*}
  \mathfrak{S}_{i-1} &:=&  c_i \mapsto (*,c_i,*),\qquad i=1,2,\dots,n\\
  \mathfrak{S}_n &:=& \mathfrak{S}_{n-1}
\end{eqnarray*}
is an isomorphism of n-discrete n-categories.
\end{Lemma}
\begin{proof}
By induction on $n$.\\

\framebox[1.1\width]{$n=1$}\\
The map $\mathfrak{S}(c_1)=(*,c_1,*)$ is trivially an isomorphism between discrete categories
$D(\C_1(c_0,c_0'))$ and $\PP_{c_0,c_0'}(\C)$.\\

\framebox[1.1\width]{$n>1$}\\
Let us denote
$$
\mathfrak{S}=\ <\mathfrak{S}_0, \{\mathfrak{S}_1,\dots,\mathfrak{S}_n\}>.
$$
 In order for $\mathfrak{S}$ to be an isomorphism of
n-categories the following facts have to be checked:
\begin{enumerate}
  \item $\mathfrak{S}_0$ is an isomorphism
  \item for every pair $c_1,c_1':c_0\to c_0'$,
$$
\{\mathfrak{S}_1,\dots,\mathfrak{S}_n\}^{c_1,c_1'}
$$
  is an isomorphism of $(n-1)$-categories
  \item above data satisfy usual functoriality axioms
\end{enumerate}

1. Since $n>1$, $[D(\C_1(c_0,c_0'))]_0=[\C_1(c_0,c_0')]_0$. Yet, by {\em Proposition \ref{dgm:hom-pb}} one has
$[\PP_{c_0,c_0'}(\C)]_0=\{*\}\times[\C_1(c_0,c_0')]_0\times\{*\}$. Hence the assignment
$\mathfrak{S}_0(c_1)=(*,c_1,*)$ is clearly an isomorphism.\\

2. For any pair $c_1,c_1':c_0\to c_0'$, induction hypothesis guaranties the existence of an
isomorphism $T^{c_1,c_1'}$:
$$
\xymatrix{
\left[D\big(\C_1(c_0,c_0')\big)\right]_1(c_1,c_1')\ar@{=}[d]\ar@{.>}[r]
&\left[\PP_{c_0,c_0'}(\C)\right]_1((*,c_1,*),(*,c_1',*))\ar@{=}[d]
\\
D\left(\big[\C_1(c_0,c_0')\big]_1(c_1,c_1')\right)\ar[r]_{T^{c_1,c_1'}}
&\PP_{c_1,c_1'}(\C_1(c_0,c_0'))
}
$$
defined by
\begin{eqnarray*}
  T_{k-1}^{c_1,c_1'}(c_k)&=&(*,c_k,*),\qquad k=2,\dots n \\
  T_n^{c_1,c_1'} &=& T_{n-1}^{c_1,c_1'}
\end{eqnarray*}

Hence we let
\begin{eqnarray*}
  \mathfrak{S}_{k-1}&=&\coprod_{c_1,c_1'\in\C_1(c_0,c_0')} T_{k-1}^{c_1,c_1'},\qquad k=2,\dots n \\
  \mathfrak{S}_{n} &=& \mathfrak{S}_{n-1}
\end{eqnarray*}
so that the isomorphism $T^{c_1,c_1'}$ is exactly $\{\mathfrak{S}_1,\dots,\mathfrak{S}_n\}^{c_1,c_1'}$.\\

3. We want to prove that $<\mathfrak{S}_0,T^{c_1,c_1'}>$ is an (iso)morphism of n-categories, {\em i.e.}
it satisfies usual coherence axioms.\\

$\bullet$ Let $c_1,c_1',c_1'':c_0\to c_0'$ be given. Coherence w.r.t. composition amounts to the
commutativity of the following diagram:
\begin{changemargin}{-10ex}{10ex}
$$
\xymatrix@C=15ex{
\left[D\big(\C_1(c_0,c_0')\big)\right]_1(c_1,c_1') \times \left[D\big(\C_1(c_0,c_0')\big)\right]_1(c_1',c_1'')\ar[r]^(.6){{}^{D(\C_1(c_0,c_0'))}\circ^0}
\ar@{=}[d]
\ar@{}[dr]|{(i)}
&\left[D\big(\C_1(c_0,c_0')\big)\right]_1(c_1,c_1'')\ar@{=}[d]
\\
\ar@{}[ddr]|{(ii)}
D\big(\C_2(c_1,c_1')\times \C_2(c_1',c_1'')\big)\ar@{=}[d]\ar[r]_{D({}^{\C}\!\circ^1)}
&D\big(\C_2(c_1,c_1'')\big)\ar[dd]^{T^{c_1,c_1''}}
\\
D\big(\C_2(c_1,c_1')\big)\times D\big( \C_2(c_1',c_1'')\big)\ar[d]_{T^{c_1,c_1'}\times T^{c_1',c_1''}}
\\
\Q_1((*,c_1,*),(*,c_1',*))\times\Q_1((*,c_1',*),(*,c_1'',*))
\ar[r]_(.6){{}^{\Q}\!\circ^0}
&\Q_1((*,c_1,*),(*,c_1',*))
}$$
\end{changemargin}
Here $(i)$ commutes by definition, while $(ii)$ commutes \emph{point-wise}.
In fact, for any $k=2,\dots,n$ and for $\xymatrix@C=3ex{c_k:c_1\ar@{==>}[r]&c_1'}$ and
$\xymatrix@C=3ex{c_k':c_1'\ar@{==>}[r]&c_1''}$ {Proposition \ref{zzz:comp}} gives
$$
\xymatrix@C=15ex@R=10ex{
(c_k,c_k')\ar@{|->}[r]^{D({}^{\C}\!\circ^1)}
\ar@{|->}[d]_{T^{c_1,c_1'}\times T^{c_1',c_1''}}
&c_k{}^{\C}\!\circ^1 c_k'
\ar@{|->}[d]^{T^{c_1,c_1''}}
\\
\big((*,c_k,*),(*,c_k',*)\big)
\ar@{|->}[r]_{{}^{\C}\!\circ^0}
&(*,c_k {}^{\C}\!\circ^1 c_k',*)}
$$

$\bullet$ Let $c_1:c_0\to c:0'$ be given. Coherence w.r.t. units amounts to the
commutativity of the following diagram:
$$
\xymatrix@C=15ex{
\Id{n-1}
\ar[r]^{{}^{D(\C_1(c_0,c_0'))}u^0(c_1)}
\ar[dr]_{{}^{\C}\!u^1(c_1)}
\ar[dddr]_{{}^{\Q}\!u^0(*,c_1,*)}
&\left[D\big(C_1(c_0,c_0')\big)\right]
\ar@{=}[d]
\\
&D\big(\C_2(c_1,c_1)\big)
\ar[dd]^{T^{c_1}}
\\
\\
&\Q_1\big((*,c_1,*),(*,c_1,*)\big)
}
$$
Upper triangle commutes by definition, lower triangle commutes \emph{point-wise}.
In fact, for any $k=2,\dots,n$  {\em Proposition \ref{zzz:unts}} gives
$$
\xymatrix@C=15ex@R=10ex{
{\ *\ }
\ar@{|->}[r] ^{[{}^{\C}\!u^1(c_1)]_k}
\ar@{|->}[dr]_{[{}^{\Q}\!u^0(*,c_1,*)]_k}
&[{}^{\C}\!u^1(c_1)]_k(*)
\ar@{|->}[d]^{T^{c_1,c_1}}
\\
&(* ,[{}^{\C}\!u^1(c_1)]_k(*),*)
}
$$
Finally, in the proof we did not explicit the level of $\mathfrak{S}_n$, as n-categories considered
are $n$-discrete and sesqui-functor $D$ is a full inclusion.
\end{proof}

\subsection{Back to the Theorem}
Now that we have developed the machinery, we are able to prove the main theorem of the section.\\

\begin{proof}[Proof of\/ {\em Theorem \ref{thm:pi_1}}]
The previous Lemma guaranties precisely the existence of a canonical isomorphism of n-categories
$$
\mathfrak{S}_{\C}^{c_0,c_0'}:D(\C_1(c_0^{\phantom{I}},c_0'))\to \PP_{c_0^{\phantom{I}},c_0'}(\C)=\Q
$$
for any pair of objects $c_0,c_0'$. Further, for a n-functor $F:\C\to\D$ we get a $(c_0,c_0')$-indexed
family of commutative squares:
\begin{equation}\label{zzz3}
\xymatrix@C=15ex{
D(\C_1(c_0,c_0'))\ar[r]^{D(F_1^{c_0,c_0'})}
\ar[d]_{\mathfrak{S}_{\C}^{c_0,c_0'}}
&D(\D_1(Fc_0,Fc_0'))
\ar[d]^{\mathfrak{S}_{\D}^{Fc_0,Fc_0'}}
\\
\PP_{c_0,c_0'}(\C)\ar[r]_{\PP_{c_0,c_0'}(F)}
&\PP_{Fc_0,Fc_0'}(\D)
}
\end{equation}
We prove this by induction.

For $n=1$ it is just a diagram of discrete categories. It suffices to verify commutativity on objects.
To this end, let us choose a  $c_1:c_0\to c_0'$. Equations below complete the case:
\begin{eqnarray*}
  \mathfrak{S}_{\D}(DF(c_1)) &=& \mathfrak{S}_{\D}(Fc_1) \\
   &=& (*,Fc_1,*) \\
  \PP F(\mathfrak{S}_{\C}(c_1)) &=& \PP F (*,c_1,*) \\
   &=& (*,Fc_1,*).
\end{eqnarray*}
Hence let us consider a generic $n>1$. First we have to show that diagram (\ref{zzz3}) commutes on objects,
but this amounts exactly to what we have just shown for $n=1$.

Thus we fix $c_1,c_1':c_0\to c_0'$ and consider homs:
$$
\xymatrix@C=15ex{
[D(\C_1(c_0,c_0'))]_1(c_1,c_1')\ar[r]^{[D(F_1^{c_0,c_0'})]_1^{c_1,c_1'}}
\ar[d]_{[\mathfrak{S}_{\C}^{c_0,c_0'}]_1^{c_1,c_1'}}
&[D(\D_1(Fc_0,Fc_0'))]_1(Fc_1,Fc_1')
\ar[d]^{[\mathfrak{S}_{\D}^{Fc_0,Fc_0'}]_1^{Fc_1,Fc_1'}}
\\
[\PP_{c_0,c_0'}(\C)]_1((*,c_1,*),(*,c_1',*))\ar[r]_(.45){[\PP_{c_0,c_0'}(F)]_1^{\diamond,\diamond}}
&[\PP_{Fc_0,Fc_0'}(\D)]_1((*,Fc_1,*),(*,Fc_1',*))
}
$$
The definition of discretizer functor $D$ allows us to re-formulate the diagram as follows:
$$
\xymatrix@C=15ex{
D\big([\C_1(c_0,c_0')]_1(c_1,c_1')\big)\ar[r]^{D\big([F_1^{c_0,c_0'}]_1^{c_1,c_1'}\big)}
\ar[d]_{[\mathfrak{S}_{\C}^{c_0,c_0'}]_1^{c_1,c_1'}}
&D\big([\D_1(c_0,c_0')]_1(Fc_1,Fc_1')\big)
\ar[d]^{[\mathfrak{S}_{\D}^{Fc_0,Fc_0'}]_1^{Fc_1,Fc_1'}}
\\
[\PP_{c_0,c_0'}(\C)]_1((*,c_1,*),(*,c_1',*))\ar[r]_(.45){[\PP_{c_0,c_0'}(F)]_1^{\diamond,\diamond}}
&[\PP_{Fc_0,Fc_0'}(\D)]_1((*,Fc_1,*),(*,Fc_1',*))
}
$$
and the previous discussion turns it in
$$
\xymatrix@C=15ex{
D\big([\C_1(c_0,c_0')]_1(c_1,c_1')\big)\ar[r]^{D\big([F_1^{c_0,c_0'}]_1^{c_1,c_1'}\big)}
\ar[d]_{T_{\C}^{c_1,c_1'}}
&D\big([\D_1(c_0,c_0')]_1(Fc_1,Fc_1')\big)
\ar[d]^{T_{\D}^{Fc_1,Fc_1'}}
\\
\PP_{c_1,c_1'}\big(\C_1(c_0,c_0')\big)\ar[r]_(.45){\PP_{c_1,c_1'}(F_1^{c_0,c_0'})}
&\PP_{Fc_1,Fc_1'}\big(\D_1(Fc_0,Fc_0')\big)
}
$$
Now, as $T$'s are just $\mathfrak{S}$'s given for $n-1$, {\em i.e.}
$$
T_{\C}^{c_1,c_1'}=\mathfrak{S}_{\C_1(c_0,c_0')}^{c_1,c_1'}
\quad \mathrm{and}\quad
T_{\D}^{Fc_1,Fc_1'}=\mathfrak{S}_{\D_1(Fc_0,Fc_0')}^{Fc_1,Fc_1'}
$$
the last diagram commutes by induction hypothesis.\\

It is clear that all this restricts to n-groupoids. Moreover, in pointed case we obtain a 2-contra-variant natural isomorphism
of sesqui-functors, {\em i.e.} a strict natural transformation of sesqui-functors that reverses the direction
of 2-morphisms and in which the assignments on objects are isomorphisms:
$$
\xymatrix{
&\ar@{}|{}="2"(n-1)\mathbf{Gpd}_*\ar@/^3ex/[dr]^{D}\\
n\mathbf{Gpd}_*\ar@/^3ex/[ur]^{\pi_1}
\ar@/_10ex/[rr]_{\Omega}="1"
&&n\mathbf{Gpd}_*
\ar@{}"2";"1"|(.3){}="a1"|(.7){}="a2"
\ar@{=>}"a1";"a2"^{\mathfrak{S}}}
$$
Indeed in $n\mathbf{Gpd}_*$ ({\em i.e.} in $n\mathbf{Gpd}$ with $c_0=*=c_0'$), for 2-morphism
$\alpha: F\Rightarrow G:\C\to\D$, we can express the (strict) naturality condition
$$
\raisebox{10ex}{
\xymatrix@C=16ex@R=15ex{
D(\C_1(*,*))
\ar@/^4ex/[r]^{D(G_1^{*,*})}="a1"
\ar@/_4ex/[r]_{D(F_1^{*,*})}="a2"
\ar[d]_{\mathfrak{S}_{\C}^{*,*}}
&
D(\D_1(*,*))
\ar[d]^{\mathfrak{S}_{\D}^{*,*}}
&
\\
\PP_{*,*}(\C)
\ar@/^4ex/[r]^{\PP_{*,*}(G)}="b1"
\ar@/_4ex/[r]_{\PP_{*,*}(F)}="b2"
&
\PP_{*,*}(\D)
\ar@{}"a1";"a2"|(.3){}="A1"|(.7){}="A2"\ar@2"A1";"A2"^{D(\alpha_1^{*,*})}
\ar@{}"b1";"b2"|(.3){}="B1"|(.7){}="B2"\ar@2"B1";"B2"^{\PP_{*,*}(\alpha)}
}}
\begin{array}{c}
  D(\alpha_1^{*,*})\bullet^0\mathfrak{S}_{\D}^{*,*} \\ \\
  = \\ \\
  \mathfrak{S}_{\C}^{*,*}\bullet^0\PP_{*,*}(\alpha)
\end{array}
$$
The proof that this condition indeed holds is a corollary to the following {\em Lemma}, that therefore concludes
the proof.
\end{proof}
\begin{Lemma}\label{zzz:Lemmata}
Given the 2-morphism of n-groupoids $\alpha: F\Rightarrow G:\C\to\D$, then the following equation holds.
$$
\raisebox{10ex}{
\xymatrix@C=16ex@R=15ex{
D(\C_1(c_0,c_0'))
\ar@/^4ex/[r]^{D(\alpha_{c_0}\circ G_1^{c_0,c_0'})}="a1"
\ar@/_4ex/[r]_{D(F_1^{c_0,c_0'}\circ \alpha_{c_0'})}="a2"
\ar[d]_{\mathfrak{S}_{\C}^{c_0,c_0'}}
&
D(\D_1(Fc_0,Gc_0'))
\ar[d]^{\mathfrak{S}_{\D}^{Fc_0,Gc_0'}}
&
\\
\PP_{c_0,c_0'}(\C)
\ar@/^4ex/[r]^{\PP_{[\alpha_{c_0}]}\circ\PP_{c_0,c_0'}(G)}="b1"
\ar@/_4ex/[r]_{\PP_{c_0,c_0'}(F)\circ \PP_{[\alpha_{c_0'}]}}="b2"
&
\PP_{Fc_0,Gc_0'}(\D)
\ar@{}"a1";"a2"|(.24){}="A1"|(.76){}="A2"\ar@2"A1";"A2"|{D(\alpha_1^{c_0,c_0'})}
\ar@{}"b1";"b2"|(.24){}="B1"|(.76){}="B2"\ar@2"B1";"B2"|{\PP_{c_0,c_0'}(\alpha)}
}}
\begin{array}{c}
  D(\alpha_1^{c_0,c_0'})\bullet^0\mathfrak{S}_{\D}^{Fc_0,Gc_0'} \\ \\
  = \\ \\
  \mathfrak{S}_{\C}^{c_0,c_0'}\bullet^0\PP_{c_0,c_0'}(\alpha)
\end{array}
$$
\end{Lemma}
\begin{proof}
We will prove the {\em Lemma} by means of the universal property of $h$-2pullback defining $\PP_{Fc_0,Gc_0'}(\D)$.
To this end let us first consider the following quite trivial chain of equalities  (taken two times):
$$
D(\alpha_1^{c_0,c_0'})\bullet^0\mathfrak{S}_{\D}^{Fc_0,Gc_0'}\bullet^0\,!=
D(\alpha_1^{c_0,c_0'})\bullet^0\,!=[id_{!}]=
\mathfrak{S}_{\C}^{c_0,c_0'}\bullet^0[id_{!}]=
\mathfrak{S}_{\C}^{c_0,c_0'}\bullet^0\PP_{c_0,c_0'}(\alpha)\bullet^0\,!
$$
Less trivially we want to prove
$$
(D(\alpha_1^{c_0,c_0'})\bullet^0\mathfrak{S}_{\D}^{Fc_0,Gc_0'})*\varepsilon_{\D}^{Fc_0,Gc_0'}
= (\mathfrak{S}_{\C}^{c_0,c_0'}\bullet^0\PP_{c_0,c_0'}(\alpha))*\varepsilon_{\D}^{Fc_0,Gc_0'}
$$
By equation (\ref{eqn:pi_1_on_2morph}) this can be rewritten
$$
D(\alpha_1^{c_0,c_0'})*(\mathfrak{S}_{\D}^{Fc_0,Gc_0'}\bullet^0\varepsilon_{\D}^{Fc_0,Gc_0'})
= (\mathfrak{S}_{\C}^{c_0,c_0'}\bullet^0\varepsilon_{\C}^{c_0,c_0'})* \alpha
$$
Let us prove directly the equality on objects. To this end let us fix an arbitrary ``object''
$c_1$ of $D(\C_1^{c_0,c_0'})$. Then applying definitions we get\\

$
\left[D(\alpha_1^{c_0,c_0'})*\left(\mathfrak{S}_{\D}^{Fc_0,Gc_0'}\bullet^0\varepsilon_{\D}^{Fc_0,Gc_0'}\right)\right]_0(c_1)=
$
\begin{eqnarray*}
   &=&  \left[\mathfrak{S}_{\D}^{Fc_0,Gc_0'}\bullet^0\varepsilon_{\D}^{Fc_0,Gc_0'}\right]_1\left(\left[D(\alpha_1^{c_0,c_0'})\right]_0(c_1)\right)\\
   &=&  \left[\mathfrak{S}_{\D}^{Fc_0,Gc_0'}\bullet^0\varepsilon_{\D}^{Fc_0,Gc_0'}\right]_1\Bigg(
   \raisebox{5ex}{
   \xymatrix{
   Fc_0
   \ar@/_4ex/[d]_{\alpha_{c_0}\circ Gc_1}^{}="1"\ar@/^4ex/[d]^{Fc_1\circ\alpha_{c_0'}}_{}="2"
   \ar@2"1";"2"^{\alpha_{c_1}}\\Gc_0'}}\Bigg)\\
   &=&  \left[\varepsilon_{\D}^{Fc_0,Gc_0'}\right]_1 \left( \mathfrak{S}_{\D}^{Fc_0,Gc_0'}(\alpha_{c_1}) \right) \\
   &=&   \left[\varepsilon_{\D}^{Fc_0,Gc_0'} \right]_1 \big( (*,\alpha_{c_1},*) \big)\\
   &=&  \alpha_{c_1}  \\
   &=&  \alpha_1^{\phantom{0}}\left(\left[\varepsilon_{\C}^{c_0,c_0'}\right]_0\big((*,c_1,*)\big)\right)  \\
   &=&  \alpha_1^{\phantom{0}}\left(\left[\varepsilon_{\C}^{c_0,c_0'}\right]_0\big(\mathfrak{S}_{\C}^{c_0,c_0'}(c_1)\big)\right)  \\
   &=&  \alpha_1^{\phantom{0}}\left(\left[\mathfrak{S}_{\C}^{c_0,c_0'}\bullet^0\varepsilon_{\C}^{c_0,c_0'}\right]_0(c_1)\right)\\
   &=& \left[(\mathfrak{S}_{\C}^{c_0,c_0'}\bullet^0\varepsilon_{\C}^{c_0,c_0'})* \alpha\right]_0(c_1)
\end{eqnarray*}

On homs we will proceed by induction. Hence let us fix arbitrary ``objects'' $c_1,c_1'$ of $D(\C_1^{c_0,c_0'})$.
Then\\

$
\left[D(\alpha_1^{c_0,c_0'})*\left(\mathfrak{S}_{\D}^{Fc_0,Gc_0'}\bullet^0\varepsilon_{\D}^{Fc_0,Gc_0'}\right)\right]_1^{c_1,c_1'}=
$
\begin{eqnarray*}
   &\eq{i}&  \left[D(\alpha_1^{c_0,c_0'})\right]_1^{c_1,c_1'}*\left[\mathfrak{S}_{\D}^{Fc_0,Gc_0'}\bullet^0\varepsilon_{\D}^{Fc_0,Gc_0'}\right]_1^{\alpha_{c_0}\circ Gc_1,Fc_1'\circ\alpha_{c_0'}}\\
   &\eq{ii}&  \left[D(\alpha_1^{c_0,c_0'})\right]_1^{c_1,c_1'}*\left(\left[\mathfrak{S}_{\D}^{Fc_0,Gc_0'}\right]_1^{\alpha_{c_0}\circ Gc_1,Fc_1'\circ\alpha_{c_0'}}  \bullet^0\left[\varepsilon_{\D}^{Fc_0,Gc_0'}\right]_1^{(*,\alpha_{c_0}\circ Gc_1,*),(*,Fc_1'\circ\alpha_{c_0'},*)}\right)\\
   &\eq{iii}&  \left[D(\alpha_1^{c_0,c_0'})\right]_1^{c_1,c_1'}*\left(   \mathfrak{S}_{\D_1^{Fc_0,Gc_0'}}^{\alpha_{c_0}\circ Gc_1,Fc_1'\circ\alpha_{c_0'}}\bullet^0\varepsilon_{\D_1^{Fc_0,Gc_0'}}^{(*,\alpha_{c_0}\circ Gc_1,*),(*,Fc_1'\circ\alpha_{c_0'},*)}   \right)\\
   &\eq{iv}&  D\left(\left[\alpha_1^{c_0,c_0'}\right]_1^{c_1,c_1'}\right)*\left(   \mathfrak{S}_{\D_1^{Fc_0,Gc_0'}}^{\alpha_{c_0}\circ Gc_1,Fc_1'\circ\alpha_{c_0'}}\bullet^0\varepsilon_{\D_1^{Fc_0,Gc_0'}}^{(*,\alpha_{c_0}\circ Gc_1,*),(*,Fc_1'\circ\alpha_{c_0'},*)}   \right)\\
   &\eq{v}& \left(\mathfrak{S}_{\C_1^{c_0,c_0'}}^{c_1,c_1'} \bullet^0\varepsilon_{\C_1^{c_0,c_0'}}^{(*,c_1,*),(*,c_1',*)}\right)* \alpha_1^{c_0,c_0'} \\
   &\eq{vi}& \left(\left[\mathfrak{S}_{\C}^{c_0,c_0'}\right]_1^{c_1,c_1'} \bullet^0\left[\varepsilon_{\C}^{c_0,c_0'}\right]_1^{(*,c_1,*),(*,c_1',*)}\right)* \alpha_1^{c_0,c_0'} \\
   &\eq{vii}& \left[\mathfrak{S}_{\C}^{c_0,c_0'}\bullet^0\varepsilon_{\C}^{c_0,c_0'}\right]_1^{c_1,c_1'}* \alpha_1^{c_0,c_0'} \\
   &\eq{viii}& \left[(\mathfrak{S}_{\C}^{c_0,c_0'}\bullet^0\varepsilon_{\C}^{c_0,c_0'})* \alpha\right]_1^{c_1,c_1'}
\end{eqnarray*}
where $(i)$ and $(viii)$ hold by definition of $*$-composition, $(ii)$ and $(vii)$ by definition of 0-whiskering
for 2-morphisms, $(iii)$ and $(vi)$ by the inductive definition of $\mathfrak{S}$ and $\varepsilon$, $(iv)$ by definition of
discretizer and finally $(v)$ by induction hypothesis.\\

Uniqueness provided by the universal property completes the proof.
\end{proof}
\vskip5ex

At last the promised alternative description of $\pi_1$.
\begin{Corollary}\label{cor:eqv_pi_1}
Let $\C$ be a pointed n-groupoid. Then there exists a natural isomorphism of sesqui-functors with components
$$
\pi_0^{(n)}(\mathfrak{S}_{\C}^{*,*}):\ \pi_1^{(n)}(\C) \to \pi_0^{(n)}(\Omega(\C))
$$
\end{Corollary}
\begin{proof}
Since $\mathfrak{S}_{\C}^{*,*}$ is $n$-discrete, $\pi_0^{(n)}(\mathfrak{S}_{\C}^{*,*})$ is still an isomorphism.
\end{proof}
\begin{Remark}
From now on, as a consequence of {\em Corollary \ref{cor:eqv_pi_1}}  we will often identify
the sesqui-functors $\pi_1(-)$ and $\pi_0(\Omega(-))$.
\end{Remark}
\subsection{Final remark on $\mathfrak{S}$}
As the reader may guess, $\mathfrak{S}$ is of a richer nature than we have shown in previous sections.

As a matter of fact we have developed the theory as far as our purposes require. Nevertheless
we urge to give a hint of the big picture behind.\\

It is possible to define a sesqui-functor
$$
\tilde{\Sigma}^{(n)}=\tilde{\Sigma}:n\mathbf{Cat}\to (n+1)\mathbf{Cat}
$$
that assigns to a $n$-category $\C$ the $(n+1)$category $\tilde{\Sigma}(\C)$ with two distinguished
objects $*_0$ and $*_1$, and the following hom-$n$-categories
\begin{itemize}
  \item $[\tilde{\Sigma}(\C)]_1(*_0,*_0) =  \Id{n}$
  \item $[\tilde{\Sigma}(\C)]_1(*_0,*_1) =  \C$
  \item $[\tilde{\Sigma}( \C )]_1(*_1,*_0) =  \emptyset$
  \item $[\tilde{\Sigma}( \C )]_1(*_1,*_1) =  \Id{n}$
\end{itemize}
with trivial compositions and obvious units.\\

Then our comparison gives natural isomorphism
$$
\mathfrak{S}_{\tilde{\Sigma}(\C)}^{*_0,*_1}:D(\C)\to\PP_{*_0,*_1}(\tilde{\Sigma}(\C)).
$$
If we restrict to equivalence sesqui-categories $\mathbf{Cat}_{eq}$, this allows to represent any $n$-category as a ($\pi_0$ of a) specified $h$-pullback. Moreover
our $\mathfrak{S}_{\C}^{c_0,c_0'}$ are indeed hom-$n$-categories
\[
\left[\PP_{*_0,*_1}(\tilde{\Sigma}(\C))\right]_1\big((*,c_0,*),(*,c_0',*)\big)
\]
(this justifying the notation adopted) of this representation, and it is easy to
check that compositions and units are provided by the universal property of $h$-pullbacks.\\

All this suggests to further develop the theory in the direction of homotopy theory. In fact the construction
developed for $\PP_{-,-}(-)$ is indeed obtained by the more general product preserving \emph{path-functor}
$$
\PP:n\mathbf{Cat}\to n\mathbf{Cat}
$$
defined by the $h$-pullback
$$
\xymatrix{
\PP(\C)\
\ar[r]^{C}
\ar[d]_{D}
&\C\ar@{=}[d]
\\
\ar@{}[ur]|(.3){}="a1"|(.7){}="a2"\ar@2"a1";"a2"^{\varepsilon_{\C}}
\C\ar@{=}[r]
&\C}
$$
This is a universal 2-morphisms representor, in that it represents every 2-cell
$\alpha:F\Rightarrow G:\B\to\C$ as a functor $\tilde{\alpha}:\B\to\PP(\C)$ such that
$D(\tilde{\alpha})=F$, $C(\tilde{\alpha})=G$ and $\varepsilon_{\C}(\tilde{\alpha})=\alpha$
(compare with universal property of $h$-pullback).
Moreover the existence of a path-functor permits to obtain all $h$-pullbacks as ordinary limits. In fact it is easy to
show that the $h$-pullback over
$$
\xymatrix{\A\ar[r]^F&\B&\C\ar[l]_{G}}
$$
is indeed the usual categorical limit over
$$
\xymatrix{
\A\ar[r]^F&\B&\PP(\B)\ar[l]_D\ar[r]^C&\B&\C\ar[l]_G}
$$
All this can be made absolutely precise and algebraic by considering the notion of
\emph{cubical comonad} and related structures on it \cite{MR1478526}.

\section{Monoidal structure on $\Omega(\C)$}

Let $\C$ be a pointed $n$-groupoid. In applying  the loop sesqui-functor $\Omega$ to $\C$ one notices that
the new structure has no memory of the 0-composition in $\C$.

In fact, every cell of $\Omega(\C)$, {\em i.e.}
every cell of $\C$ with the object ``$*$'' as its 0-domain and 0-codomain is 0-composable. Hence we can recover
the forgotten structure, in order to get a \emph{many sorted} strict monoidal structure
$$
\big(\Omega(\C), \otimes, I\big)
$$
with $\otimes = \circ^0,$ and $I=(*,1_{*},*)$. Moreover cells are weakly invertible w.r.t. this structure, thus giving a
\emph{group like} structure to the $n$-groupoid $\Omega(\C)$.
\begin{Remark}
This is indeed the same as considering $\Omega(\C)$ as a $(n+1)$-groupoid with one object.
\end{Remark}
$n$-Functoriality axioms prove the following
\begin{Lemma}
The following two statements hold for $n$-groupoids, $n>0$.
\begin{enumerate}
\item Let $n$-functor $F:\C\to\D$ be given. Then $\Omega(F)$ is a strict monoidal $n$-functor.
\item Let natural $n$-transformation $\alpha : F\Rightarrow G:\C\to\D$ be given. Then $\Omega(F)$ is
a strict monoidal natural $n$-transformation.
\end{enumerate}
\end{Lemma}
Furthermore if we apply $\Omega$ once more, we get another monoidal structure on the $n$-discrete
$(n-1)$-discrete $n$-groupoid $\Omega(\Omega(\C))$. This new structure corresponds to 0-composition
of $\Omega(\C)$, {\em i.e.} 1-composition of $\C$. Functoriality of units and compositions allow us to
 apply an Heckmann-Hilton like argument, this showing that compositions coincide and are indeed commutative.\\

In terms of monoidal structures, we ca resume the previous discussion:
$$
\big(\Omega(\Omega(\C)), \otimes, I\big)
$$
is a  commutative strict monoidal structure on the $n$-groupoid $\Omega(\Omega(\C))$.

Notice that monoidal structure is automatically preserved by $\pi_0$, hence all this can be said for
sesqui-functor $\pi_1$:
\begin{Proposition}
Let $\C$ be a $n$-groupoid. Then the $(n-1)$-groupoid $\pi_1(\C)$ is naturally endowed with weakly
invertible strict monoidal structure. This structure is commutative when we consider $(\pi_1)^2(\C)$.
\end{Proposition}

\section{$\Omega$ and $\pi_1$ preserve exactness}

In the following paragraphs we will show that, given a three-term exact sequence in $n$\textbf{Gpd}, the sesqui-functor
$\Omega^{(n)}$ produces a three-term exact sequence in $n$\textbf{Gpd}. As a consequence, we obtain a similar result
for $\pi_1^{(n)}$\\

\begin{Lemma}\label{lem:omega_pres_pb}
Sesqui-functor  $\Omega^{(n)}$ preserves $h$-pullbacks.
\end{Lemma}
\begin{proof}
We consider a slightly more general setting, in order to get the proof of the statement as a consequence.
Notice that we will omit the superscripts  $(n)$.

Let us consider a $h$-pullback
$$
\xymatrix@C=12ex{
\Q\ar[r]^Q\ar[d]_P
&\C\ar[d]^G
\\
\A\ar[r]_F
\ar@{}[ur]|(.3){}="1"|(.7){}="2"\ar@2"1";"2"^{\phi}
&\B
}
$$
in $n\mathbf{Gpd}$. Moreover let us fix objects $\alpha\in\A_0$ and $\gamma\in\C_0$ such
that $F(\alpha)=\beta=G(\gamma)$. Next let us apply the universal property of $h$-2pullback to get $\PP(\phi)$ as in the diagram
$$
\xymatrix@C=12ex{
\PP_{q,q}(\Q)\ar[r]^{\PP_{}(P)}\ar[d]_{\PP_{}(Q)}
&\PP_{\alpha,\alpha}(\A)\ar[d]^{\PP_{}(F)}
\\
\PP_{\gamma,\gamma}(\C)\ar[r]_{\PP_{}(G)}
\ar@{}[ur]|(.3){}="2"|(.7){}="1"\ar@2"2";"1"_{\PP_{}(\phi)}
&\PP_{\beta,\beta}(\B)
}
$$
where $q=(\alpha,1_{\beta},\gamma)$.
Further let us consider the diagram
$$
\xymatrix@C=12ex{
\X\ar[r]^{M}\ar[d]_{N}
&\PP_{\alpha,\alpha}(\A)\ar[d]^{\PP_{}(F)}
\\
\PP_{\gamma,\gamma}(\C)\ar[r]_{\PP_{}(G)}
\ar@{}[ur]|(.3){}="2"|(.7){}="1"\ar@2"2";"1"^{\omega}
&\PP_{\beta,\beta}(\B)
}
$$
that by {\em Theorem \ref{thm:pi_1}} can be re-drawn
$$
\xymatrix@C=10ex{
\X\ar[r]^{M}\ar[d]_{N}
&D(\A_1(\alpha,\alpha))\ar[d]^{D(F_1^{\alpha,\alpha})}
\\
D(\C_1(\gamma,\gamma))\ar[r]_{D(G_1^{\gamma,\gamma})}
\ar@{}[ur]|(.3){}="2"|(.7){}="1"\ar@2"2";"1"^{\omega}
&D(\B_1(\beta,\beta))
}
$$
Applying $\pi_0$ one gets
$$
\xymatrix@C=12ex{
\pi_0(\X)\ar[r]^{\pi_0(M)}\ar[d]_{\pi_0(N)}
&\A_1(\alpha,\alpha)\ar[d]^{F_1^{\alpha,\alpha}}
\\
\C_1(\gamma,\gamma)\ar[r]_{G_1^{\gamma,\gamma}}
\ar@{}[ur]|(.3){}="2"|(.7){}="1"\ar@2"2";"1"^{\pi_0(\omega)}
&\B_1(\beta,\beta)
}
$$
Since $\Q_1(q,q)$ is defined as a $h$-pullback (see Section \ref{sec:h-pb}), the universal property yields a unique
$$
L:\pi_0(\X)\to\Q_1(q,q)
$$
 such that
$$
(i)\,   L\bullet^0 P_1^{q,q}=\pi_0(M)\quad
(ii)\,  L\bullet^0 Q_1^{q,q}=\pi_0(N)  \quad
(iii)\, L\bullet^0 \phi_1^{q,q}=\pi_0(\omega)
$$
which in turn implies that {\em there exists a unique}
$$
\xymatrix{\X\ar[r]^(.35){\eta_{\X}^{\phantom{\X}}}&D\pi_0(\X)\ar[r]^(.3){DL}&\D(\Q_1(q,q))=\PP_{q,q}(\Q)}
$$
{\em such that}
\begin{eqnarray*}
  \eta_{\X}^{\phantom{\X}}\bullet^0 DL \bullet^0 \PP(P)
  &=& \eta_{\X}^{\phantom{\X}}\bullet^0 DL \bullet^0 D P_1^{q,q} \\
  &=& \eta_{\X}^{\phantom{\X}}\bullet^0 D(L \bullet^0  P_1^{q,q}) \\
  &=& \eta_{\X}^{\phantom{\X}}\bullet^0 D(\pi_0(M)) \\
  &=& M
\end{eqnarray*}
where the last equality holds by universality of adjunctions.
Similarly one gets
$$
\eta_{\X}^{\phantom{\X}}\bullet^0 DL \bullet^0 \PP(Q)=N
$$
and
$$
\eta_{\X}^{\phantom{\X}}\bullet^0 DL \bullet^0 \PP(\phi)=\omega
$$
and this concludes the proof.
\end{proof}

\begin{Lemma}\label{lem:omega_pres_sur}
Sesqui-functor $\Omega^{(n)}$ preserves $h$-surjective morphisms.
\end{Lemma}
\begin{proof}
This is absolutely straightforward. Let $L:\K \to \A$ be an $h$-surjective morphism. Then, for
a fixed object $\kappa$ of $\K$
$$
\PP_{\kappa,\kappa}(L)=D(L_1^{\kappa,\kappa})
$$
Now $L_1^{\kappa,\kappa}$ is $h$-surjective by definition since $L$ is,  $D$ preserves trivially
$h$-surjective morphisms.
\end{proof}

\begin{Proposition}
Let the exact sequence of pointed $n$-groupoids
$$
\xymatrix{\A\ar[r]_F\ar@/^5ex/[rr]^{0}="1"&\B\ar@{}|{}="2"\ar[r]_{G}&\C
\ar@{}"1";"2"|(.3){}="a1"|(.75){}="a2"\ar@2"a1";"a2"^{\lambda}}
$$
be given. Then the sequence
$$
\xymatrix{\Omega\A\ar[r]^{\Omega F}\ar@/_5ex/[rr]_{0}="1"&\Omega\B\ar@{}|{}="2"\ar[r]^{\Omega G}&\Omega\C
\ar@{}"1";"2"|(.3){}="a2"|(.75){}="a1"\ar@2"a1";"a2"^{\Omega\lambda}}
$$
is an exact sequence of pointed $n$-groupoids.
\end{Proposition}
\begin{proof}
Let $L:\A\to\K$ be the comparison with the kernel of $G$, that is $h$-surjective by definition.
By {\em Lemma \ref{lem:omega_pres_pb}} above, $\Omega L$ is the comparison with the kernel of
$\Omega G$, and it is $h$-surjective by {\em Lemma \ref{lem:omega_pres_sur}}.
\end{proof}

Besides we get the following {\em for free}
\begin{Corollary}
Sesqui-functor $\pi_1^{(n)}$ preserves exact sequences, reversing  the direction of the 2-morphism.
\end{Corollary}

Similar results hold for the non pointed case, when we fix suitable points.

\section{Fibration sequence of a $n$-functor and the Ziqqurath of exact sequences}
Purpose of this and next sections is to establish the setting in order to get the main result.\\

Let be given a morphism of $n$-groupoids
$$
F:\ \B\to \C
$$

If we fix an object $\beta$ of $\B$, the fiber diagram over $F\beta$
\begin{equation}\label{zzz:fibra}
\xymatrix{\K\ar[r]_{K}\ar@/^5ex/[rr]^{[F\beta]}="1"&\B\ar@{}|{}="2"\ar[r]_{F}&\C
\ar@{}"1";"2"|(.3){}="a1"|(.75){}="a2"\ar@{=>}"a1";"a2"^{\varphi}}
\end{equation}
is an exact sequence of $n$-groupoids. In the following sections we will show that this produces a
canonical exact sequence
$$
\xymatrix@!C=10ex{
\PP_{\beta,\beta} (\B)\ar[r]_{\PP_{\beta,\beta} (F)}\ar@/^7ex/[rr]^{[(*,1_{\beta},\beta)]}="b1"
\ar@{}[rr]|{}="b2"&\PP_{F\!\beta,F\!\beta} (\C)\ar@/_7ex/[rr]_{[\beta]}="c1"\ar[r]_{\nabla}
\ar@{}[rr]|{}="c2"&\K\ar[r]_{K}\ar@/^7ex/[rr]^{[F\beta]}="a1"
&\B\ar@{}|{}="a2"\ar[r]_{F}&\C
\ar@{}"c1";"c2"|(.3){}="C1"|(.75){}="C2"\ar@{=}"C1";"C2"
\ar@{}"b1";"b2"|(.3){}="B1"|(.75){}="B2"\ar@{=>}"B1";"B2"^{\sigma}
\ar@{}"a1";"a2"|(.3){}="A1"|(.75){}="A2"\ar@{=>}"A1";"A2"^{\varphi}}
$$
{\em i.e.} the sequence represented above is exact in $\PP_{F\!\beta F\!\beta}(\C)$, in $\K$ and in $\B$.

\subsection{Connecting morphism $\nabla$}
Although  $\nabla$ is easily obtained by means of the universal property of $h$-pullback of $\K$, we will consider a
slightly more general situation in order to apply induction properly in the construction of the exact sequence.\\

Let us consider the (past) $h$-fiber  $\K=\F^{(p)}_{F,F\beta}$, $\beta\in \B_0$. Then for any other
$\beta'\in \B_0$, the universal property yields a
$$
\nabla=\nabla^{(p)}_{F\beta,\beta',F}: \PP_{F\!\beta,F\!\beta'} (\C)\to \F^{(p)}_{F,F\beta}=\K
$$
$$
\xymatrix{
\PP_{F\!\beta,F\!\beta'} (\C)
\ar[rr]\ar[dd]\ar@{->}[dr]_(.3){\nabla}
&\ar@{}[d]|(.4){\varepsilon_{\C}^{F\!\beta,F\!\beta'}\quad}
&\Id{n}\ar[d]^{[\beta']}
\\
\ar@{}[ur]|(.2){}="b1"|(.8){}="b2"\ar@{:>}"b1";"b2"
&\K\ar[r]^{K}\ar[d]
&\B\ar[d]^{F}
\\
\Id{n}\ar@{-}[r]_{id}
&
\Id{n}\ar[r]_{[F\beta]}
\ar@{}[ur]|(.3){}="a1"|(.7){}="a2"\ar@2"a1";"a2"_{\varphi}
&\C}
$$
By {\em Lemma \ref{lemma:h-pb_h-proj}} this is better understood as
\begin{equation}\label{diag:nabla}
\xymatrix{
\PP_{F\!\beta,F\!\beta'} (\C)
\ar[r]\ar[d]_{\nabla}
&\Id{n}\ar[d]^{[\beta']}
\\
\ar@{}[ur]|{(\dag)}
\K\ar[r]^{K}\ar[d]
&\B\ar[d]^{F}
\\
\Id{n}\ar[r]_{[F\beta]}
\ar@{}[ur]|(.3){}="a1"|(.7){}="a2"\ar@2"a1";"a2"_{\varphi}
&\C
}
\end{equation}
Where upper commuting square $(\dag)$ is a pullback, {\em i.e.} $\nabla$ is the strict fiber
of $K=K^{(p)}_{F,F\beta}$ over $\beta'$.
\begin{Proposition}
The sequence $<\nabla,id_{[\beta']},K>$ above is exact.
\end{Proposition}
\begin{proof}
What we must prove is that the comparison
$$
L:\ \PP_{F\!\beta,F\!\beta'}(\C)\to \F^{(f)}_{K,\beta'}=\mathbb{H}
$$
is $h$-surjective. We will prove the statement by induction.\\

For $n=1$ the functor
$$
D(\C_1(F\beta,F\beta'))\to\F^{(f)}_{K,\beta'}
$$
is clearly essentially surjective on objects and full. In fact it is defined on objects as
 in the general case $L_0$ below, hence it is essentially surjective on objects for the same
 proof. Moreover it is also full: chose two objects $c_1,c_1'$ in the domain, then we get objects
$$
\big((*,c_1 ,\beta'),1_{\beta'},* \big),\quad
\big((*,c_1',\beta'),1_{\beta'},* \big)
$$
An arrow between them is a 1-cell $b_1:\beta'\to\beta'$ such that
$$
(1)\ c_1'\circ Fb_1=c_1\qquad\ (2)\ b_1\circ 1_{\beta'}=1_{\beta'}
$$
Condition $(2)$ forces $b_1=1_{\beta'}$. This makes condition $(1)$ true only for $c_1=c_1'$,
{i.e.} the image hom-set is a singleton, implying that functor $L$ on homs is trivially surjective.\\

Henceforth let us suppose $n>1$. In order to fix notation we recall the $h$-fiber $\mathbb{H}$ is   a triple
$(\mathbb{H},H,\psi)$.\\

Equations
$$
L\bullet^0 H=\nabla,\qquad L\bullet^0\varphi= id_{[\beta']}
$$
give $L$ on objects:
$$
L_0:\ (*,\xymatrix{F\beta\ar[r]^{c_1}&F\beta'},*)\mapsto \big((*,\xymatrix{F\beta\ar[r]^{c_1}&F\beta'},\beta'),\xymatrix{\beta'\ar[r]^{1_{\beta'}}&\beta'},*\big)
$$
This is indeed an essentially surjective map. In fact, given
$$
h_0=\big((*,\xymatrix{F\beta\ar[r]^{c_1}&Fb_0},b_0),\xymatrix{b_0\ar[r]^{b_1}&\beta'},*\big)\in\mathbb{H}_0
$$
there exist a $\overline{p_0}$ and a $\overline{h_1}:h_0\to L(\overline{p_0})$: simply let
$$
\overline{p_0}=(*,\xymatrix{F\beta\ar[r]^{c_1}&Fb_0\ar[r]^{Fb_1}&F\beta'},*)
$$
and $\overline{h_1}$ as below
$$
\xymatrix@C=2ex{
h_0\ar[d]_{\overline{h_1}}
&=
&\big(\ (*\ar@<+1ex>@{.}[d]&,
&F\beta\ar[rr]^{c_1}\ar@{.}[d]\ar@{}[drr]|{=}
&&Fb_0\ar@{.>}[d]^{Fb_1}
&b_0\ar@{.>}[d]^{b_1}\ )&,
&b_0\ar[rr]_{b_1}\ar@{.>}[d]_{b_1}\ar@{}[drr]|{=}
&&\beta'\ar@{.}[d]&,
&{*}\ar@<-.75ex>@{.}[d]\ \big)
\\
L(\overline{p_0})
&=
&\big(\ (*&,
&F\beta\ar[rr]_{c_1 Fb_1}
&&Fb_0,
&b_0\ )&,
&b_0\ar[rr]_{b_1}
&&\beta'&,
&{*}\ \big)
}
$$
Then we fix a pair of objects $p_0,p_0'$ of $\PP_{F\!\beta,F\!\beta'}(\C)$
$$
p_0=(*,\xymatrix{F\beta\ar[r]^{c_1}&F\beta'},*),\qquad p_0'=(*,\xymatrix{F\beta\ar[r]^{c_1'}&F\beta'},*).
$$
As we have shown in proving the universal property of $h$-pullbacks, the comparison on homs
$$
L_1^{p_0,p_0'}:\ [\PP_{F\!\beta,F\!\beta'}(\C)]_1(p_0,p_0')\longrightarrow [\F^{(f)}_{K,\beta'}]_1(Fp_0,Fp_0')
$$
is given by universal property on homs ({\em i.e.} for $(n-1)$groupoids)
$$
L': \PP_{c_1,c_1'}(\C_1(F\beta,F\beta'))\longrightarrow \F^{(p)}_{K_1^{Hh_0,Hh_0'},c_1'}
$$
Now $c_1=c_1\circ F_1^{\beta',\beta'}(1_{\beta'})$, and  by definition $K_1^{Hh_0,Hh_0'}=K^{(f)}_{c_1\circ F_1^{\beta',\beta'}\ ,c_1'}$.

Hence we started with a comparison of the kind
\begin{equation}\label{xxx:rec1}
\PP_{x,Fy}(\Z)\longrightarrow \F^{(f)}_{K_{F,x}^{(p)}}
\end{equation}
and we obtain its homs part as a comparison of the kind
\begin{equation}\label{xxx:rec2}
\PP_{Fy,x}(\Z)\longrightarrow \F^{(p)}_{K_{F,x}^{(f)}}
\end{equation}
This situation  requires to check also that comparison (\ref{xxx:rec2}) is  essentially
surjective and then calculate it on homs. We obtain a comparison (\ref{xxx:rec1})
that terminates a two-level induction process and gives at once that both comparisons are
$h$-surjective.

In fact the same calculation as above shows that (\ref{xxx:rec2}) is essentially surjective, but
reverses the direction. Nevertheless this is not a serious obstruction, as all cells of an $n$-groupoid are equivalences.

In order to understand this fully, it may be interesting to make the construction explicit
for 1-cells. Let 1-cell $h_1:Lp_0\to Lp_0'$ as in the following diagram
$$
\xymatrix@C=2ex{
Lp_0\ar[d]_{h_1}
&=
&\big(\ (*\ar@<+1ex>@{.}[d]&,
&F\beta\ar[rr]^{c_1}\ar@{.}[d]
&&F\beta'\ar@{.>}[d]^{Fb_1}\ar@{}[dll]|(.3){}="a1"|(.7){}="a2"\ar@2"a1";"a2"^{c_2}
&\beta'\ar@{.>}[d]^{b_1}\ )&,
&\beta'\ar[rr]^{1_{\beta'}}\ar@{.>}[d]_{b_1}
&&\beta'\ar@{.}[d]\ar@{}[dll]|(.3){}="b1"|(.7){}="b2"\ar@2"b1";"b2"^{b_2}
&,
&{*}\ar@<-.75ex>@{.}[d]\ \big)
\\
Lp_0'
&=
&\big(\ (*&,
&F\beta\ar[rr]_{c_1'}
&&F\beta',
&\beta'\ )&,
&\beta'\ar[rr]_{1_{\beta'}}
&&\beta'&,
&{*}\ \big)
}
$$
{\em i.e.}
$$
h_1=\big(\ (1_{*},\xymatrix{c_1\circ Fb_1\ar@2^{c_2}[r]&1_{F\beta}\circ c_1'},b_1) ,\ \xymatrix{1_{\beta'}\circ 1_{\beta'}\ar@2[r]^{b_2}&b_1\circ 1_{\beta'}},\ 1_{*}\ \big)
$$
then there exist a $\overline{p_1}:p_0\to p_0'$ and a $\overline{h_2}: L(\overline{p_1})\to h_1$.
In fact such a $\overline{p_1}$ should be of the form $(*,\xymatrix{c_1\ar@2[r]^{\overline{c_2}}&c_1'},*)$.
Hence $L(\overline{p_1})$ is of the form
$$
\xymatrix@C=2ex{
Lp_0\ar[d]_{L(\overline{p_1})}
&=
&\big(\ (*\ar@<+1ex>@{.}[d]&,
&F\beta\ar[rr]^{c_1}\ar@{.}[d]
&&F\beta'\ar@{.}[d]\ar@{}[dll]|(.3){}="a1"|(.7){}="a2"\ar@2"a1";"a2"^{\overline{c_2}}
&\beta'\ar@{.}[d]\ )&,
&\beta'\ar[rr]^{1_{\beta'}}\ar@{.}[d]
&&\beta'\ar@{.}[d]\ar@{}[dll]|{=}
&,
&{*}\ar@<-.75ex>@{.}[d]\ \big)
\\
Lp_0'
&=
&\big(\ (*&,
&F\beta\ar[rr]_{c_1'}
&&F\beta',
&\beta'\ )&,
&\beta'\ar[rr]_{1_{\beta'}}
&&\beta'&,
&{*}\ \big)
}
$$
{\em i.e.}
$$
L(\overline{p_1})=\big(\ (1_{*},\xymatrix{c_1\circ 1_{F\beta'}\ar@2[r]^{\overline{c_2}}&1_{F\beta}\circ c_1'},1_{\beta'}) ,\ id_{1_{\beta'}},\ 1_{*}\ \big)
$$
Then it suffices to let
$$
\overline{c_2}=\xymatrix{c_1\ar@2[r]^(.4){c_1Fb_2}&c_1Fb_1\ar@2[r]^(.6){c_2}&c_1'}
$$
and get the wanted 2-cell:
$$
\xymatrix@C=1ex{
L(\overline{p_1})
\ar@2[d]_{}
&=
&\big(\ (1_{*}\ar@<+1ex>@{:}[d]&,
&c_1\circ 1_{F\beta'}\ar@2[rrr]^{(c_1Fb_2)c_2}\ar@{:>}[d]_{c_1Fb_2}
&&&1_{F\beta}\circ c_1'\ar@{:}[d]^{}\ar@{}[dlll]|{\equiv}
&1_{\beta'}\ar@{:>}[d]^{b_2}\ )&,
&1_{\beta'}\ar@{=}[rr]^{id_{1_{\beta'}}}\ar@{:}[d]_{}
&&1_{\beta'}\ar@{:>}[d]^{b_2}\ar@{}[dll]|{=}
&,
&1_{*}\ar@<-.75ex>@{:}[d]\ \big)
\\
h_1&=
&\big(\ (1_{*}&,
&c_1\circ Fb_1\ar@2[rrr]_{c_2}
&&&1_{F\beta}\circ c_1',
&b_1\ )&,
&1_{\beta'}\ar@2[rr]_{b_2}
&&b_1&,
&1_{*}\ \big)
}
$$
\end{proof}
\subsection{Connecting 2-morphism $\sigma$}
In order to paste diagram (\ref{diag:nabla}) with $\PP$ of the original sequence, exactness
in $\PP_{F\beta,\F\beta'}(\C)$ must be shown, that means we have to find the 2-morphism
$[0]\Rightarrow\PP(F)\bullet^0\nabla$ that realizes exactness. This is done by means of 2-dimensional
 universal property of $h$-pullbacks. \\

Let us recognize this fact in the following diagram:
\begin{equation}\label{zzz:situa}
\xymatrix@!=12ex{
\PP_{\beta,\beta'}(\B)\ar[r]^{\PP(F)}\ar[d]
&\PP_{F\beta,F\beta'}(\C)\ar[r]\ar[d]^{\nabla}
&\Id{n}\ar[d]^{[\beta']}
\\
\Id{n}\ar[r]_{[(*,1_{\beta},\beta)]}
\ar@{}[ur]|(.3){}="a1"|(.7){}="a2"\ar@{:>}"a1";"a2"_{\sigma}
&\K\ar[r]^{K}\ar[d]\ar@{}[ur]|{(pb)}
&\B\ar[d]^{F}
\\
&\Id{n}\ar[r]_{[F\beta]}
\ar@{}[ur]|(.3){}="b1"|(.7){}="b2"\ar@{=>}"b1";"b2"_{\varphi}
&\C
}
\end{equation}
\subsubsection{Construction of $\sigma$}
We can apply 2-dimensional universal property of $\K=\F^{(p)}_{F,F\beta}$ (even if $h$-pullbacks regularity is enough) to the following set of data:
2-morphisms $\raisebox{3.5ex}{\xymatrix@C=3ex@R=3ex{\PP(\B)\ar[r]^{[\beta]}\ar[d]&\B\ar[d]^F\\ \I\ar@{}[ur]|(.3){}="b1"|(.7){}="b2"\ar@{=}"b1";"b2"_{id}\ar[r]_{[F\beta]}&\C}}$
and $\raisebox{3.5ex}{\xymatrix@C=4ex@R=3ex{\PP(\B)\ar[r]^{[\beta']}\ar[d]&\B\ar[d]^F\\ \I\ar@{}[ur]|(.3){}="b1"|(.7){}="b2"\ar@2"b1";"b2"_{\varepsilon_{\B}F}\ar[r]_{[F\beta]}&\C}}$
over the base, 2-morphism $\raisebox{3.5ex}{\xymatrix@R=3.5ex{\PP(\B)\ar@/_2ex/[d]_{}^{}="1"\ar@/^2ex/[d]^{}_{}="2"\ar@{=}"1";"2"\\ \I}}$
and $\raisebox{3.5ex}{\xymatrix@R=3.5ex{\PP(\B)\ar@/_2ex/[d]_{[\beta]}^{}="1"\ar@/^2ex/[d]^{[\beta']}_{}="2"\ar@{=>}"1";"2"^{\varepsilon_{\B}}\\ \I}}$
and (identity) 3-morphism $\raisebox{3.5ex}{\xymatrix@R=3.5ex{[F\beta]\ar@{=}[r]\ar@{=}[d]&[\beta]F\ar@2[d]^{\varepsilon_{\B}F}\\[F\beta]\ar@2[r]_{\varepsilon_{\B}F}&[\beta']F\ar@{}[ul]|{\equiv}}}$,
where the second diagram is justified by the equalities
$$
\PP(F)\bullet^0\nabla\bullet^0\varphi=\PP(F)\bullet^0\varepsilon_{\C}=\varepsilon_{\B}\bullet^0F.
$$

Then there exists a unique $\sigma:[(*,1_{\beta},\beta)]\Rightarrow \PP(F)\bullet^0\nabla$ such that
$$
(i)\ \sigma\bullet^0K=\varepsilon_{\B}^{\beta,\beta'}\qquad(ii)\ \sigma*\varphi=id_{\varepsilon_{\B}^{ \beta,\beta'}\bullet^0 F}
$$
\subsubsection{Comparison morphism}
\begin{Proposition}
The triple $(\PP(F),\sigma,\nabla)$ is exact.
\end{Proposition}
\begin{proof}
In order to show that triple $(\PP(F),\sigma,\nabla)$ is exact, we must verify that comparison with (past) $h$-fibre
of $\nabla$ is $h$-surjective. Namely we construct the $h$-pullback
$$
\raisebox{5ex}{\xymatrix{
\N\ar[r]^(.34){J}\ar[d]&\PP_{F\beta,F\beta'}(\C)\ar[d]^{\nabla}\\
\I\ar@{}[ur]|(.3){}="1"|(.7){}="2"\ar@2"1";"2"^{\xi}
\ar[r]_{[(*,1_{F\beta},\beta)]}&\K}}\qquad\mathrm{where}\qquad
\begin{array}{lcl}
  \K&=&\F^{(p)}_{F,F\beta} \\ \\
  \nabla&=&\nabla^{(p)}_{F\beta,\beta',F} \\ \\
  \N&=&\F^{(p)}_{(*,1_{F\beta},\beta),\nabla}
\end{array}
$$
The universal property for $\N$ yields a unique $N:\PP_{\beta,\beta'}(\B)\to\N$ such that
\begin{equation}\label{zzz:PN1}
    (i)\ N\bullet^0J=\PP(F)\qquad\ (ii)\ N\bullet^0\xi=\sigma
\end{equation}
For $n=1$ the result is easily obtained. In fact in this case $\N$ is a discrete groupoid, and since also
$\PP_{\beta,\beta'}(\B)$ is, it suffices to check surjectivity on objects, and this is achieved in the
same manner as $h$-surjectivity for the general case.\\

Hence let us suppose $n>1$.

We will need an explicit description of $n$-functor $N=<N_0,N_1^{-,-}>$.
$$
N_0: (*,\xymatrix@C=3ex{\beta\ar[r]^{b_1}&\beta'},*)\mapsto\big(*,
(1_{*},\xymatrix@C=3ex{Fb_1\ar@2[r]^{id_{Fb_1}}&Fb_1},b_1)
,(*,Fb_1,*)
\big)
$$
Moreover for any pair of objects $(*,b_1,*)$ and $(*,b_1',*)$ in $\PP_{\beta,\beta'}(\B)$,
$$
N_1^{(*,b_1,*),(*,b_1',*)}
$$
is obtained by the universal property of $h$-pullbacks as shown by the diagram below, where \mbox{$p_0=(*,b_1,*)$} and \mbox{$p_0'=(*,b_1',*)$},
\begin{equation}\label{zzz:N1}
\xymatrix{
[\PP_{\beta,\beta'}(\B)]_1((*,b_1,*),(*,b_1',*))
\ar@/^5ex/[dr]^(.7){\qquad\ [\PP(F)]_1^{(*,b_1,*),(*,b_1',*)}}
\ar@/_12ex/[ddd]
\ar@{.>}[d]^{N_1^{(*,b_1,*),(*,b_1',*)}}
\\
\N_1(Np_0,Np_0')\ar[r]^(.3){J_1^{Np_0,Np_0'}}
\ar[dd]
&[\PP_{F\beta,F\beta'}(\C)]_1((*,Fb_1,*),(*,Fb_1',*))
\ar[d]^{[\nabla]_1^{(*,Fb_1,*),(*,Fb_1',*)}}
\\
&\K_1\big((*,Fb_1,\beta'),(*,Fb_1',\beta')\big)
\ar[d]^{(1_*,1_{Fb_1},b_1)\circ-}
\\
\ar@{}[uur]|(.3){}="2"|(.7){}="1"\ar@2"1";"2"_{\xi_1^{Np_0,Np_0'}}
\Id{n-1}
\ar[r]_(.3){[(1_*,1_{Fb_1'},b_1')]}
&\K_1\big((*,1_{F\beta},\beta),(*,Fb_1',\beta')\big)
}
\end{equation}
such that
\begin{equation}\label{zzz:PN2}
    (i)\ N_1^{\diamond,\diamond}\bullet^0J_1^{\diamond,\diamond}=[\PP(F)]_1^{\diamond,\diamond}\qquad\
    (ii)\ N_1^{\diamond,\diamond}\bullet^0\xi_1^{\diamond,\diamond}=\sigma_1^{\diamond,\diamond}
\end{equation}
\begin{Claim}
$N$ is essentially surjective on objects.
\end{Claim}

Let  an object $n_0$ of $\N$ be given
$$
n_0=\big(*, (1_{*},\xymatrix@C=3ex{Fb_1\ar@2[r]^{c_2}&c_1},b_1),(*,\xymatrix@C=3ex{F\beta\ar[r]^{c_1}&F\beta'},*)   \big)
$$
then it suffices to consider the object $p_0=(*,b_1,*)$ to get an arrow $N(p_0)\to n_0$, as suggested by the diagram below
$$
\xymatrix@C=1ex{
N(p_0)\quad=\ar@<-2ex>[d]
&
\big(*,\ar@{.}[d]
&(*,\ar@{:}[d]
&Fb_1\ar@2[rrr]^{id}\ar@{:}[d]
&&&Fb_1,\ar@{:>}[d]^{c_2}
&b_1),\ar@<-1ex>@{:}[d]
&(*,\ar@{.}[d]
&F\beta\ar@{.}[d]\ar[rrr]^{Fb_1}\ar@{}[d]
&&&F\beta',\ar@{.}[d]
&{*})\big)\ar@<-1ex>@{.}[d]
\\
n_0\quad=
&
\big(*,
&(*,
&Fb_1\ar@2[rrr]_{c_2}
\ar@{}[urrr]|{\equiv}
&&&c_1,
&b_1),
&(*,
&F\beta\ar[rrr]_{c_1}
\ar@{}[urrr]|(.3){}="2"|(.7){}="1"\ar@2"1";"2"_{c_2}
&&&F\beta',
&{*})\big)
}
$$

\begin{Claim}
For any pair $(*,b_1,*),(*,b_1',*)$, the $(n-1)$-functor $N_1^{(*,b_1,*),(*,b_1',*)}$ is $h$-surjective.
\end{Claim}

Here comes the inductive step: we want to get $N_1^{\diamond,\diamond}$ from a situation in which is
indeed a comparison itself, as in the original setting of diagram (\ref{zzz:situa}).
To this end let us consider the diagram
$$
\xymatrix@C=10ex@R=10ex{
\PP_{b_1,b_1'}\big(\B_1(\beta,\beta')\big)
\ar[r]^(.45){[\PP(F)]_1^{\diamond,\diamond}}
\ar[d]
&
\PP_{Fb_1,Fb_1'}\big(\C_1(F\beta,F\beta')\big)
\ar[r]
\ar[d]|{(1_{*},id_{Fb_1},b_1)\circ[\nabla]_1^{\diamond,\diamond}}
&\I\ar@{}[dl]|{(\dag)}
\ar[d]^{[b_1]}
\\
\I
\ar@{}[ur]|(.3){}="2"|(.7){}="1"\ar@2"1";"2"_{\sigma_1^{\diamond,\diamond}}
\ar[r]_(.3){[(1_{*},id_{Fb_1'},b_1')]}
&\K_1\big((*,1_{F\beta},\beta),(*,Fb_1',\beta')\big)
\ar[d]
\ar[r]^(.6){K_1^{\diamond,\diamond}}
&\B_1(\beta,\beta')
\ar[d]^{F_1^{\beta,\beta'}}
\\
&\I
\ar@{}[ur]|(.3){}="2"|(.7){}="1"\ar@2"1";"2"_{\varphi_1^{\diamond,\diamond}}
\ar[r]_{[Fb_1']}
&\C_1(F\beta,F\beta')
}
$$
By definition of $h$-pullback the square \framebox{$\varphi_1^{\diamond,\diamond}$} is the $h$-fiber
$\F^{(f)}_{F_1^{\beta,\beta'},Fb_1}$. Furthermore the square $(\dag)$ is a pullback by {\em Lemma \ref{lemma:h-pb_h-proj}}.
In fact this follows by the universal property, as the following equations hold:
\begin{eqnarray*}
   &(\clubsuit)& \big((1_{*},id_{Fb_1},b_1)\circ[\nabla]_1^{(*,Fb_1,*),(*,Fb_1',*)}\big)\bullet^0 K_1^{(*,1_{F\beta},\beta),(*,Fb_1',\beta')}= [b_1] \\
   &(\diamondsuit)& \big((1_{*},id_{Fb_1},b_1)\circ[\nabla]_1^{(*,Fb_1,*),(*,Fb_1',*)}\big)\bullet^0\varphi_1^{(*,1_{F\beta},\beta),(*,Fb_1',\beta')}= \varepsilon_{\C_1(F\beta,F\beta')}^{Fb_1,Fb_1'}
\end{eqnarray*}
{\em proof of }$(\clubsuit)$ {\em and of} $(\diamondsuit)$. They follow easily from the universal property in dimension $n-1$:
\begin{eqnarray*}
\big((1_{*},id_{Fb_1},b_1)\circ \nabla_1^{\diamond,\diamond}\big)\bullet^0 K_1^{\diamond,\diamond}
   &=& K\big((1_{*},id_{Fb_1},b_1)\big)\circ(\nabla_1^{\diamond,\diamond}\bullet^0K_1^{\diamond,\diamond})  \\
   &=& K\big((1_{*},id_{Fb_1},b_1)\big)\circ([\nabla\bullet^0K]_1^{\diamond,\diamond}) \\
   &\eq{*}& K\big((1_{*},id_{Fb_1},b_1)\big)\circ([\beta']_1^{\diamond,\diamond}) \\
   &=& \left[K\big((1_{*},id_{Fb_1},b_1)\big)\right] \\
   &=& [b_1]
\end{eqnarray*}
\begin{eqnarray*}
  \big((1_{*},id_{Fb_1},b_1)\circ[\nabla]_1^{\diamond,\diamond}\big)\bullet^0\varphi_1^{\diamond,\diamond}
   &\eq{i}& \Big(1_{Fb_1}\circ(\nabla_1^{\diamond,\diamond}\bullet^0 KF_1^{\diamond,\diamond}) \Big)
   \bullet^1\Big(1_{F\beta}\circ (\nabla_1^{\diamond,\diamond}\bullet^0\varphi_1^{\diamond,\diamond}) \Big) \\
   &=&  \Big(\nabla_1^{\diamond,\diamond}\bullet^0 KF_1^{\diamond,\diamond} \Big)\bullet^1\Big(\nabla_1^{\diamond,\diamond}\bullet^0\varphi_1^{\diamond,\diamond} \Big)\\
   &\eq{ii}&  \nabla_1^{\diamond,\diamond}\bullet^0\varphi_1^{\diamond,\diamond}\\
   &\eq{*}&  [\varepsilon_{\C}^{F\beta,F\beta'}]_1^{Fb_1,Fb_1'}\\
   &=&  \varepsilon_{\C_1(F\beta,F\beta')}^{Fb_1,Fb_1'}
\end{eqnarray*}
where equations marked $(*)$ hold for inductive definition of universal property of $h$-pullbacks, while $(i)$ is
composition axiom of 2-morphism, and $(ii)$ since the first component of 1-composition is a 1-morphism.

Hence
$$
(1_{*},id_{Fb_1},b_1)\circ\left[\nabla^{(p)}_{F\beta,\beta',F}\right]_1^{\diamond,\diamond}=\nabla^{(f)}_{F_1^{\beta,\beta'}b_1',b_1,F_1^{\beta,\beta'}}
$$
{\em i.e.} it is a  ``$\nabla$'' itself. \\

In other words, we wanted to prove that comparison with an $h$-fiber of
$$
\nabla^{(p)}_{F\beta,\beta',F}:\ \PP_{F\beta,F\beta'}(\C)\longrightarrow\F^{(p)}_{F,F\beta}
$$
is $h$-surjective, and we find out that is equivalent to asking:
\begin{enumerate}
\item its essential surjectivity
\item that the \emph{comparison} in dimension (n-1) w.r.t.
$$
\nabla^{(f)}_{F_1^{\beta,\beta'}b_1',b_1,F_1^{\beta,\beta'}}:\ \PP_{Fb_1,Fb_1'}\Big(\C_1(F\beta,F\beta')\Big)\longrightarrow \F^{(f)}_{F_1^{\beta,\beta'},F_1^{\beta,\beta'}b_1'}
$$
is $h$-surjective.
\end{enumerate}
In conclusion, we get the same construction as in dimension $n$, up to directions.

Since we are in a (weakly) invertible setting, we can apply induction properly and get the result.
\end{proof}

\subsection{The fibration sequence of $F$}
From now on we will consider $h$-kernels in the pointed setting  $n\mathbf{Gpd}_{*}$. Nevertheless all the constructions
 plainly apply to $h$-fibers in $n\mathbf{Gpd}$, as shown in the construction for diagram (\ref{zzz:fibra}).

Let $F:\C\to\D$  be a morphism of pointed $n$-groupoids. Specializing constructions above with
$\beta=*$ and $0=[*]$  we can exhibit the exact sequence
$$
\xymatrix@!C=6ex{
\Omega \B\ar[r]_{\Omega F}\ar@/^5ex/[rr]^{0}="b1"
&\Omega \C\ar@/_5ex/[rr]_{0}="c1"\ar@{}|{}="b2"\ar[r]_{\nabla}
&\ar@{}|{}="c2"\K\ar[r]_{K}\ar@/^5ex/[rr]^{0}="a1"
&\B\ar@{}|{}="a2"\ar[r]_{F}&\C
\ar@{}"c1";"c2"|(.3){}="C1"|(.75){}="C2"\ar@{=}"C1";"C2"
\ar@{}"b1";"b2"|(.3){}="B1"|(.75){}="B2"\ar@{=>}"B1";"B2"^{\sigma}
\ar@{}"a1";"a2"|(.3){}="A1"|(.75){}="A2"\ar@{=>}"A1";"A2"^{\varphi}}
$$
Since $\Omega$ preserves exactness, this gives another exact sequence
$$
\xymatrix@!C=6ex{
\Omega^2 \B\ar[r]^{\Omega^2 F}\ar@/_5ex/[rr]_{0}="b1"
&\Omega^2 \C\ar@/^5ex/[rr]^{0}="c1"\ar@{}|{}="b2"\ar[r]^{\Omega\nabla}
&\ar@{}|{}="c2"\Omega\K\ar[r]^{\Omega K}\ar@/_5ex/[rr]_{0}="a1"
&\Omega\B\ar@{}|{}="a2"\ar[r]^{\Omega F}&\Omega\C
\ar@{}"c1";"c2"|(.3){}="C1"|(.75){}="C2"\ar@{=}"C1";"C2"
\ar@{}"b1";"b2"|(.3){}="B1"|(.75){}="B2"\ar@{=>}"B2";"B1"^{\Omega\sigma}
\ar@{}"a1";"a2"|(.3){}="A1"|(.75){}="A2"\ar@{=>}"A2";"A1"^{\Omega\varphi}}
$$
Those can be pasted together in the seven-term exact sequence
$$
\xymatrix@!C=4ex{
\Omega^2 \B\ar[r]^{\Omega^2 F}\ar@/_5ex/[rr]_{0}="xb1"
&\Omega^2 \C\ar@/^5ex/[rr]^{0}="xc1"\ar@{}|{}="xb2"\ar[r]^{\Omega\nabla}
&\ar@{}|{}="xc2"\Omega\K\ar[r]^{\Omega K}\ar@/_5ex/[rr]_{0}="xa1"
&\Omega\B\ar@{}|{}="xa2"\ar[r]^{\Omega F}\ar@/^5ex/[rr]^{0}="b1"
&\Omega \C\ar@/_5ex/[rr]_{0}="c1"\ar@{}|{}="b2"\ar[r]_{\nabla}
&\ar@{}|{}="c2"\K\ar[r]_{K}\ar@/^5ex/[rr]^{0}="a1"
&\B\ar@{}|{}="a2"\ar[r]_{F}&\C
\ar@{}"c1";"c2"|(.3){}="C1"|(.75){}="C2"\ar@{=}"C1";"C2"
\ar@{}"b1";"b2"|(.3){}="B1"|(.75){}="B2"\ar@{=>}"B1";"B2"^{\sigma}
\ar@{}"a1";"a2"|(.3){}="A1"|(.75){}="A2"\ar@{=>}"A1";"A2"^{\varphi}
\ar@{}"xc1";"xc2"|(.3){}="xC1"|(.75){}="xC2"\ar@{=}"xC1";"xC2"
\ar@{}"xb1";"xb2"|(.3){}="xB1"|(.75){}="xB2"\ar@{=>}"xB2";"xB1"^{\Omega\sigma}
\ar@{}"xa1";"xa2"|(.3){}="xA1"|(.75){}="xA2"\ar@{=>}"xA2";"xA1"^{\Omega\varphi}}
$$
Of course the process can be iterated indefinitely in order to get a longer exact sequence, even if it
trivializes after $n$ applications.

\subsection{The Ziqqurath of a morphism of (pointed) $n$-groupoids}
A different perspective is gained by considering sesqui-functor $\pi_1$ in  place of $\Omega$.

In fact in the longer exact sequences obtained above, repeated applications of $\Omega$ give structures
which  are discrete in higher dimensional cells. Their exactness may be fruitfully investigated in
lower dimensional settings, {\em i.e.} after  repeated applications of $\pi_0$.

To this end we  state the following
\begin{Lemma}\label{lemma:pi0_commuta_pi1}
Sesqui-functor $\pi_0$ commutes with sesqui-functor $\pi_1$, {\em i.e.} for every integer $n>1$ the following
diagram is commutative
$$
\xymatrix@C=10ex{
n\mathbf{Gpd}_{*}
\ar[r]^{\pi_0^{(n)}}
\ar[d]_{\pi_1^{(n)}}
&(n-1)\mathbf{Gpd}_{*}
\ar[d]^{\pi_1^{(n-1)}}
\\
(n-1)\mathbf{Gpd}_{*}
\ar[r]_{\pi_0^{(n-1)}}
&
(n-2)\mathbf{Gpd}_{*}
}
$$
\end{Lemma}
\begin{proof}
This can be proved directly. For $n=2$ the diagram commutes trivially. Hence let us suppose $n>2$.
Let us be given a pointed $n$-groupoid $\C$, then by direct application of inductive definitions
involved one has
\begin{eqnarray*}
  [\pi_0(\pi_1(\C))]_0 &=& [\pi_0(\C_1(*,*))]_0 \\
   &=& [\C_1(*,*)]_0 \\
   &=& [\pi_0(\C_1(*,*))]_0 \\
   &=& [[\pi_0(\C)]_1(*,*)]_0 \\
   &=& [\pi_1(\pi_0(\C))]_0
\end{eqnarray*}
Moreover for any pair of ``objects'' $c_1,c_1'$ one has
\begin{eqnarray*}
[\pi_0(\pi_1(\C))]_1(c_1,c_1')   &=&  [\pi_0(\C_1(*,*))]_1(c_1,c_1')\\
   &=&  \pi_0([\C_1(*,*)]_1(c_1,c_1'))\\
   &=&  \big[\pi_0(\C_1(*,*))\big]_1(c_1,c_1')\\
   &=&  \big[[\pi_0(\C)]_1(*,*)\big]_1(c_1,c_1')\\
   &=&  \big[\pi_1(\pi_0(\C))\big]_1(c_1,c_1')
\end{eqnarray*}
Finally this extends plainly to morphisms and 2-morphisms.
\end{proof}
\begin{Remark}
In the language of loops, we can re-state {\em Lemma} above in other terms:
$$
\pi_0(\pi_0(\Omega(-)))=\pi_0(\Omega(\pi_0(-)))
$$
\end{Remark}
Let now a morphism $F:\C\to\D$ of pointed $n$-groupoids be given. Then the $h$-kernel exact sequence
$$
\xymatrix@!C=8ex{
\K\ar[r]_{K}\ar@/^7ex/[rr]^{0}="a1"
&\B\ar@{}|{}="a2"\ar[r]_{F}&\C
\ar@{}"a1";"a2"|(.3){}="A1"|(.75){}="A2"\ar@{=>}"A1";"A2"^{\varphi}}
$$
gives two exact sequences of pointed $(n-1)$-groupoids:
$$
\xymatrix@!C=5ex{
\pi_1\K\ar[r]^{\pi_1K}\ar@/_5ex/[rr]_{0}="a1"
&\pi_1\B\ar@{}|{}="a2"\ar[r]^{\pi_1F}&\pi_1\C
\ar@{}"a2";"a1"|(.3){}="A1"|(.75){}="A2"\ar@{=>}"A1";"A2"^{\pi_1\varphi}}
\qquad\xymatrix@!C=5ex{
\pi_0\K\ar[r]_{\pi_0K}\ar@/^5ex/[rr]^{0}="a1"
&\pi_0\B\ar@{}|{}="a2"\ar[r]_{\pi_0F}&\pi_0\C
\ar@{}"a1";"a2"|(.3){}="A1"|(.75){}="A2"\ar@{=>}"A1";"A2"^{\pi_0\varphi}}
$$
Those can be connected together in order to give a six term exact sequence of pointed $(n-1)$-groupoids
$$
\xymatrix@!C=5ex{
\pi_1\K\ar[r]^{\pi_1K}\ar@/_5ex/[rr]_{0}="a1"
&\pi_1\B\ar@{}|{}="a2"\ar[r]^{\pi_1F}\ar@/^5ex/[rr]^{0}="d1"
&\ar@{}|{}="d2"
\pi_1\C\ar[r]_{\Delta}
\ar@/_5ex/[rr]_{0}="c1"
\ar@{}"a2";"a1"|(.3){}="A1"|(.75){}="A2"\ar@{=>}"A1";"A2"^{\pi_1\varphi}
\ar@{}"d1";"d2"|(.3){}="D1"|(.75){}="D2"\ar@{=>}"D1";"D2"^{\delta}
&
\pi_0\K\ar[r]_{\pi_0K}\ar@/^5ex/[rr]^{0}="xa1"\ar@{}|{}="c2"
&\pi_0\B\ar@{}|{}="xa2"\ar[r]_{\pi_0F}&\pi_0\C
\ar@{}"c2";"c1"|(.3){}="C1"|(.65){}="C2"\ar@{=}"C1";"C2"
\ar@{}"xa1";"xa2"|(.3){}="xA1"|(.75){}="xA2"\ar@{=>}"xA1";"xA2"^{\pi_0\varphi}}
$$
where $\Delta=\pi_0(\nabla)$ and $\delta=\pi_0(\sigma)$. Notice that the three leftmost terms are endowed
with strict monoidal structure and weak inverses.

Applying $\pi_0$ and $\pi_1$,  we get two six-term exact sequences. Those can be pasted by {\em Lemma \ref{lemma:pi0_commuta_pi1}}
in a nine-term exact sequence of $(n-2)$-groupoids (cells to be pasted are dotted in the diagram):
$$
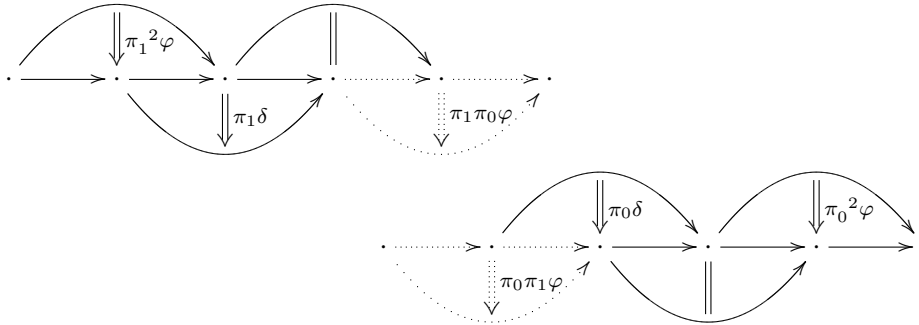
\raisebox{14ex}{\xymatrix@!C=3.5ex{
\cdot\ar[r]\ar@/^6ex/[rr]^{}="A2"
&\cdot\ar[r]\ar@{}|{}="A1"\ar@{=>}"A2";"A1"^{\pi_1^{\ 2}\varphi}
\ar@/_6ex/[rr]^{}="K2"
&\cdot\ar@{}|{}="K1"\ar@/^6ex/[rr]^{}="H2"
\ar[r]
&
\cdot\ar@{.>}[r]\ar@{..>}@/_6ex/[rr]^{}="B2"\ar@{}|{}="H1"
&\cdot\ar@{.>}[r]\ar@{}|{}="B1"\ar@{::>}"B1";"B2"^{\pi_1\pi_0\varphi}
&\cdot
\ar@{=>}"K1";"K2"^{\pi_1\delta}
\ar@{=}"H2";"H1"}}
\!\!\!\!\!\!\!\!\!\!\!\!\!\!\!\!\!\!\!\!\!\!\!\!\!\!\!\!\!\!\!\!\!\!\!\!\!\!\!
\xymatrix@!C=3.5ex{
\cdot\ar@{.>}[r]\ar@{..>}@/_6ex/[rr]_{}="A2"
&\cdot\ar@{.>}[r]\ar@{}|{}="A1"\ar@{::>}"A1";"A2"^{\pi_0\pi_1\varphi}
\ar@/^6ex/[rr]^{}="K2"
&\cdot\ar@{}|{}="K1"\ar@/_6ex/[rr]^{}="H2"
\ar[r]
&
\cdot\ar[r]\ar@/^6ex/[rr]^{}="B2"\ar@{}|{}="H1"
&\cdot\ar[r]\ar@{}|{}="B1"\ar@{=>}"B2";"B1"^{\pi_0^{\ 2}\varphi}
&\cdot
\ar@{=>}"K2";"K1"^{\pi_0\delta}
\ar@{=}"H1";"H2"}
$$
Now the three leftmost terms are endowed with a commutative strict monoidal structure and weak inverses
and the three middle terms are endowed with strict monoidal structure and weak inverses.

Iterating the process we obtain a sort of tower, a Ziqqurath, in which the lower is the level, the lower is
the dimension the longer is the length of the sequence.

$$
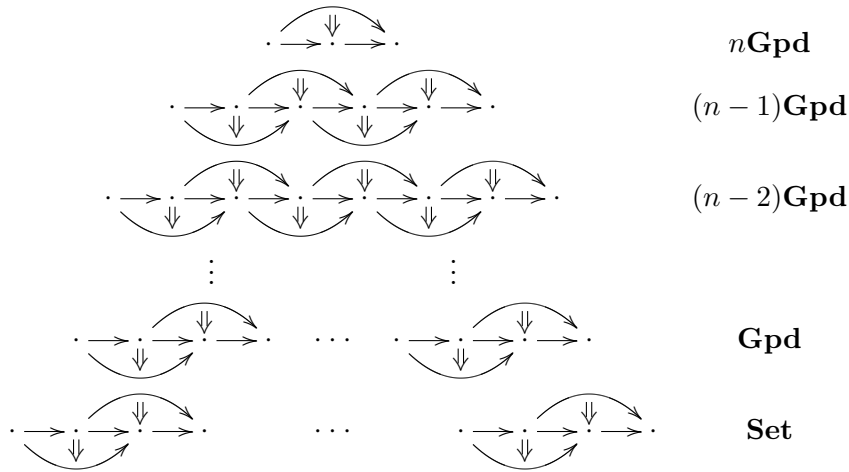
\begin{array}{cc}
\xymatrix@!C=0ex{
\cdot\ar[r]
\ar@/^3ex/[rr]_{}="1"
\ar@{}[rr]^{}="2"
&\cdot\ar[r]\ar@2"1";"2"
&\cdot
}
&n\mathbf{Gpd}
\\
\xymatrix@!C=0ex{
\cdot\ar[r]
\ar@/_3ex/[rr]^{}="a2"
\ar@{}[rr]_{}="a1"
&
\cdot\ar[r]\ar@2"a1";"a2"
\ar@/^3ex/[rr]_{}="b1"
\ar@{}[rr]^{}="b2"
&
\cdot\ar[r]\ar@2"b1";"b2"
\ar@/_3ex/[rr]^{}="xa2"
\ar@{}[rr]_{}="xa1"
&
\cdot\ar[r]\ar@2"xa1";"xa2"
\ar@/^3ex/[rr]_{}="xb1"
\ar@{}[rr]^{}="xb2"
&
\cdot\ar[r]\ar@2"xb1";"xb2"
&\cdot
}
&(n-1)\mathbf{Gpd}
\\
\xymatrix@!C=0ex{
\cdot\ar[r]
\ar@/_3ex/[rr]^{}="a2"
\ar@{}[rr]_{}="a1"
&
\cdot\ar[r]\ar@2"a1";"a2"
\ar@/^3ex/[rr]_{}="b1"
\ar@{}[rr]^{}="b2"
&
\cdot\ar[r]\ar@2"b1";"b2"
\ar@/_3ex/[rr]^{}="ya2"
\ar@{}[rr]_{}="ya1"
&
\cdot\ar[r]\ar@2"ya1";"ya2"
\ar@/^3ex/[rr]_{}="yb1"
\ar@{}[rr]^{}="yb2"
&
\cdot\ar[r]\ar@2"yb1";"yb2"
\ar@/_3ex/[rr]^{}="xa2"
\ar@{}[rr]_{}="xa1"
&
\cdot\ar[r]\ar@2"xa1";"xa2"
\ar@/^3ex/[rr]_{}="xb1"
\ar@{}[rr]^{}="xb2"
&
\cdot\ar[r]\ar@2"xb1";"xb2"
&\cdot
}
&(n-2)\mathbf{Gpd}
\\
\vdots\qquad\qquad\qquad\qquad\vdots
\\
\xymatrix@!C=0ex{
\cdot\ar[r]
\ar@/_3ex/[rr]^{}="a2"
\ar@{}[rr]_{}="a1"
&
\cdot\ar[r]\ar@2"a1";"a2"
\ar@/^3ex/[rr]_{}="b1"
\ar@{}[rr]^{}="b2"
&
\cdot\ar[r]\ar@2"b1";"b2"
&
\cdot
&
\cdot\cdot\cdot
&
\cdot\ar[r]
\ar@/_3ex/[rr]^{}="xa2"
\ar@{}[rr]_{}="xa1"
&
\cdot\ar[r]\ar@2"xa1";"xa2"
\ar@/^3ex/[rr]_{}="xb1"
\ar@{}[rr]^{}="xb2"
&
\cdot\ar[r]\ar@2"xb1";"xb2"
&\cdot
}
&\mathbf{Gpd}
\\
\xymatrix@!C=0ex{
\cdot\ar[r]
\ar@/_3ex/[rr]^{}="a2"
\ar@{}[rr]_{}="a1"
&
\cdot\ar[r]\ar@2"a1";"a2"
\ar@/^3ex/[rr]_{}="b1"
\ar@{}[rr]^{}="b2"
&
\cdot\ar[r]\ar@2"b1";"b2"
&
\cdot
&&
\cdot\cdot\cdot
&&
\cdot\ar[r]
\ar@/_3ex/[rr]^{}="xa2"
\ar@{}[rr]_{}="xa1"
&
\cdot\ar[r]\ar@2"xa1";"xa2"
\ar@/^3ex/[rr]_{}="xb1"
\ar@{}[rr]^{}="xb2"
&
\cdot\ar[r]\ar@2"xb1";"xb2"
&\cdot
}
&\mathbf{Set}
\end{array}
$$

In particular, the last row counts $3\cdot n$ terms. From left to right, there are $3\cdot n-6$ abelian groups,
3 groups and 3 pointed sets.

The row before the last counts $3\cdot(n-1)$ terms. From left to right, there are $3\cdot n-9$ strictly
commutative categorical groups, 3 categorical groups, 3 pointed groupoids. Let us observe that categorical groups
produced in this way are strict monoidal weakly invertible ones.

\appendix
\chapter{$n$-Groupoids, comparing definitions}\label{cha:appendix}
\section{$n\mathbf{Cat}$: the globular approach}
In this section we compare the classical globular definition of
$n$-category with the inductively
enriched one presented up to here.\\

The following definition is freely adapted from the one  presented
in \cite{BHch81}. It is essentially the same presented also in \cite{MR1130401}, and is indeed equivalent
to that of \cite{MR920944}.

\begin{Definition}\label{def:glob_n_cat}
A $n$-category is a reflexive $n$-truncated globular set
$$
(\mathcal{C}_{\bullet}=\{\mathcal{C}_i\}_{i=0,\cdots,n},
\{s_i,t_i:\mathcal{C}_{i+1}\to\mathcal{C}_i\}_{i}, \{
e_{i+1}:\mathcal{C}_i\to\mathcal{C}_{i+1}\}_i)
$$

We will often (ab)use the notations
\begin{eqnarray*}
s_k:\mathcal{C}_i\to\mathcal{C}_k   &\mathrm{meaning}&  s_{i-1}\cdot s_{i-2}\cdot\cdots\cdot s_{k},\\
t_k:\mathcal{C}_i\to\mathcal{C}_k   &\mathrm{meaning}&  t_{i-1}\cdot t_{i-2}\cdot\cdots\cdot t_{k},\\
e_i:\mathcal{C}_k\to\mathcal{C}_i   &\mathrm{meaning}& e_{k+1}\cdot
e_{k+2}\cdot\cdots\cdot e_{i}.
\end{eqnarray*}

$\mathcal{C}_{\bullet}$ is endowed with operations ($m<i$)
$$
\star^{m}:\mathcal{C}_i\, {}_{t_m}\!\!\times_{s_m}
\mathcal{C}_i\to\mathcal{C}_i
$$
such that:
\begin{enumerate}
\item for all $c,c'\in\mathcal{C}_{\bullet}$ and $m\leq k$,

$$
s_k(c\star^m c')= \left\{
  \begin{array}{ll}
    s_k (c) & \mathrm{if}\ m=k  \\
    s_k(c)\star^m s_k(c') &\mathrm{if}\  m<k \\
  \end{array}
\right.
$$

$$
t_k(c\star^m c')= \left\{
  \begin{array}{ll}
    t_k (c') & \mathrm{if}\ m=k  \\
    t_k(c)\star^m t_k(c') &\mathrm{if}\  m<k \\
  \end{array}
\right.
$$

\item for all $c\in\mathcal{C}_i$,  all $m$,
$$
e_{i}(s_m(c))\star^m c=c=c\star^m e_{i}(t_m(c))
$$

\item for all $c,c'\in\mathcal{C}$, all $m$,
$$
e(c\star^m c')=e(c)\star^m e(c')
$$

\item for all $c,c',c''\in\mathcal{C}$, all $m$,
$$
c\star^m(c'\star^m c'')=(c\star^mc')\star^m c''
$$

\item for all $c,c',d,d'\in\mathcal{C}$, all $p<q$,
$$
(c\star^q c')\star^p(d\star^q d')=(c\star^p d)\star^q(c'\star^p d')
$$
\end{enumerate}
In order to avoid confusion, this will be called {\em globular n-category}.\\

\end{Definition}

Given a $n$-category $\C$, this defines a globular $n$-category
$\mathcal{C}$.

In fact it suffices to let
$$
\mathcal{C}_0=\C_0
$$
and for every $0<i\leq n$,
$$
\mathcal{C}_i=\bigcup_{c_{i-1},c_{i-1}'\in\mathcal{C}_{i-1}}[\C_i(c_{i-1},c_{i-1}')]_0
$$
Sources targets and identities are obtained composing the following
ones:
$$
s_{i}:\mathcal{C}_{i+1}\to\mathcal{C}_i,\quad
t_{i}:\mathcal{C}_{i+1}\to\mathcal{C}_i,\quad
e_{i+1}:\mathcal{C}_{i}\to\mathcal{C}_{i+1}
$$
where
$$
s_{i-1}(c_i:\xymatrix{c_{i-1}\ar[r]|(.5){i}&c_{i-1}'})=c_{i-1}
$$
$$
t_{i-1}(c_i:\xymatrix{c_{i-1}\ar[r]|(.5){i}&c_{i-1}'})=c_{i-1}'
$$
$$
e_{i+1}(c_i)=\left[{}^{\C_i(c_{i-1},c_{i-1}')}\!u(c_i)\right](*):\xymatrix@C=10ex{c_{i}\ar[r]|(.45){i+1}&c_{i}}
$$
A simple calculation shows that this forms a reflexive $n$-truncated globular set.\\

Let now suppose we are given a pair $(c_i,c_i')\in\mathcal{C}_i\,
{}_{t_m}\!\!\times_{s_m} \mathcal{C}_i$.

This means that $c_i$ is a cell of $\C_{m+1}(c_m,c_m')$ and that
$c_i'$ is a cell of $\C_{m+1}(c_m',c_m'')$, for certain
$c_m,c_m',c_m''$. Then we can easily define the composition
$$
c_i\star^m c_i'= c_i\circ^m c_i'
$$
where $\circ^m$ is the $m$-composition  morphism
\begin{equation}\label{www:comp}
(-)\ {}^{\C}\!\circ^m_{c_m,c_m',c_m''}\ (-)
:\C_{m+1}(c_m,c_m')\times\C_{m+1}(c_m',c_m'')\to \C_{m+1}(c_m,c_m'')
\end{equation}
This can be seen as $0$-composition. In fact, by the inductive
definition of a $n$-category, there exist $c_{m-1},c_{m-1}'$ such
that $\xymatrix{c_m,c_m',c_m'':c_{m-1}\ar[r]|(.55){m}&c_{m-1}'}$.
Then ${}^{\C}\!\circ^m$ is indeed
$$
(-)\ {}^{\C_m(c_{m-1},c_{m-1}')}\!\circ^0_{c_m,c_m',c_m''}\ (-)
:\C_{m+1}(c_m,c_m')\times\C_{m+1}(c_m',c_m'')\to \C_{m+1}(c_m,c_m'')
$$
where, as usual, the various $\C_{m+1}(x,y)$ are indeed the short form for $[\C_m(c_{m-1},c_{m-1}')]_1(x,y)$.\\

These data satisfy axioms for a globular $n$-category.
\begin{proof}
The proof is divided into five parts, according to the five axioms.
\begin{enumerate}
  \item
The statement will be proved by (finite) induction over $k$, for a
fixed $m$.

The base of the induction is given by definition:
$$
s_{n-1}(c_n\star^{n-1}
c_n')=s_{n-1}(\xymatrix{c_{n-1}\ar[r]^{c_n}&c_{n-1}'\ar[r]^{c_n'}&c_{n-1}''})=c_{n-1}.
$$

Now suppose  $k\geq m$. Then for every $i>k$ one has
\begin{eqnarray*}
s_k(c_i\star^m c_i')   &\eq{i}&  s_k(s_{k+1}(c_i\star^m c_i'))\\
   &\eq{ii}&  s_k(s_{k+1}(c_i)\star^m s_{k+1}(c_i')) \\
   &\eq{iii}&  s_k(c_{k+1}\star^m c_{k+1}')
\end{eqnarray*}
where $(i)$ holds by definition, $(ii)$ by induction, $(iii)$ is
just a typographical substitution.

Indeed what we mean with the expression ``$c_{k+1}\star^m
c_{k+1}'$'' is the image under the morphism (\ref{www:comp}) of the
pair $(c_{k+1}, c_{k+1}')$.  This is inductively defined on homs,
hence for strictly $k>m$ one can make it explicit:
$$
(c_{k+1},c_{k+1}'):\xymatrix@C=10ex{(c_k,\bar{c}_k)\ar[r]|{k+1}&(c_k',\bar{c}_k')}
$$
This is a $(k-m)$-cell of $\C_m(c_m,c_m')\times \C_m(c_m',c_m'')$,
and its image under $[{}^{\C}\!\circ^m_{c_m,c_m',c_m''}]$ is indeed
its image under
$$
[{}^{\C}\!\circ^m_{c_m,c_m',c_m''}]_{k-m}^{(c_k,\bar{c}_k),(c_k',\bar{c}_k')}.
$$
Then functoriality over the two components of a product implies
$$
s_k(c_{k+1}\star^m c_{k+1}')=s_k(c_{k+1})\star^m s_k(c_{k+1}')=
s_k(c_i)\star^m s_k(c_i').
$$
Differently for $k=m$,
$$
s_{k}(c_{k+1}\star^{k}
c_{k+1}')=s_{k}(\xymatrix{c_{k}\ar[r]^{c_{k+1}}&c_{k}'\ar[r]^{c_{k+1}'}&c_{k}''})=c_{k}.
$$
\vskip3ex The analogous statement relative to targets is dealt
similarly.

  \item First we observe that for every $c_k:\xymatrix{c_{k-1}\ar[r]|k&c_{k-1}'}$, $k< n$, functoriality w.r.t. units
  forces the morphism ${}^{\C_k(c_{k-1},c_{k-1}')}u$ to satisfy the equation expressed by the following diagram
$$
\xymatrix@C=30ex@R=10ex{
\I\ar[r]^{id}\ar[dr]|{[{}^{\C_{k+1}(c_k,c_k)}\!u(1_{c_k})}
&\I\ar[d]|{{}^{\C_k(c_{k-1},c_{k-1}'}\!u(c_k)]_1^{*,*}}
\\
& \C_{k+1}(c_k,c_k) }
$$
that is
$$
[{}^{\C_k(c_{k,1},c_{k-1}')}\!u(c_k)]_1^{*,*}={}^{\C_{k+1}(c_k,c_k)}
\! u(1_{c_k})
$$
This fact inductively extends and gives relations between units. In
particular for $k<i$ this implies
\begin{eqnarray*}
e_i(c_k)&=&e_i(e_{i-1}(\cdots e_{k+1}(e_k(c_k))\cdots))\\
 &=&e_i(e_{i-1}(\cdots e_{k+1}([{}^{\C_k(c_{k-1},c_{k-1}')}\!u(c_k)])\cdots))\\
 &=&e_i(e_{i-1}(\cdots e_{k+1}(1_{c_k})\cdots))\\
 &=&e_i(e_{i-1}(\cdots [{}^{\C_{k+1}(c_k,c_k)} \! u(1_{c_k})]\cdots))\\
 &=&e_i(e_{i-1}(\cdots [{}^{\C_k(c_{k-1},c_{k-1}')}\!u(c_k)]_1\cdots))\\
 &&\qquad\vdots\\
 &=&[{}^{\C_k(c_{k-1},c_{k-1}')}\!u(c_k)]_{i-k}
\end{eqnarray*}
Hence we can calculate
$$
e_i(s_k(c_i)) \star^m c_i =  [{}^{\C_k(c_{k-1},c_{k-1}')}\!u(
s_k(c_i) )]_{i-k} \circ^m  c_i =c_i
$$
as a direct consequence of neutral $m$-units. \vskip3ex The
analogous statement relative to targets is dealt similarly.

  \item For $m<i<n$ we want to prove
$$
e_{i+1}(c_i\star^m \bar{c}_i)=e_{i+1}(c_i) \star^m
e_{i+1}(\bar{c}_i)
$$
or more simply
$$
1_{(c_{i}\circ^m \bar{c}_{i})}=1_{c_{i}}\circ^m 1_{\bar{c}_{i}}.
$$
This is just functoriality w.r.t. units of the morphism
$(-)\circ^m (-)$. In fact if $c_i:c_{i-1}\to c_{i-1}'$ and
$\bar{c}_i:\bar{c}_{i-1}\to \bar{c}_{i-1}'$,
 $m$-composition sends the identity
$$
(1_{c_{i}}, 1_{\bar{c}_{i}}):
\xymatrix@C=10ex{(c_{i},\bar{c}_{i})\ar[r]|{i+1}&(c_{i},\bar{c}_{i})}
$$
of $\C_{i}(c_{i-1},c_{i-1}')\times
\C_{i}(\bar{c}_{i-1},\bar{c}_{i-1}')$ to the identity
$$
1_{(c_{i}\circ^m \bar{c}_{i})}:\xymatrix@C=10ex{(c_{i}\circ^m
\bar{c}_{i})\ar[r]|{i+1}&(c_{i}\circ^m \bar{c}_{i})}
$$
of
$\C_{i}(c_{i-1}\circ^m\bar{c}_{i-1},c_{i-1}'\circ^m\bar{c}_{i-1}')$.

\item Let $m<i\leq n$. We have to prove the equality
$$
c_i\star^m (c_i'\star^m c_i'')=(c_i\star^m c_i')\star^m c_i''
$$
{\em i.e.}
$$
c_i\circ^m (c_i'\circ^m c_i'')=(c_i\circ^m c_i')\circ^m c_i''
$$
This holds by associativity of $m$-composition.

  \item
Let us suppose $p<q<i\leq n$. We have to prove the equality
$$
(c_i\star^q \bar{c}_i)\star^p(d_i\star^q \bar{d}_i)=(c_i\star^p
d_i)\star^q(\bar{c}_i\star^p \bar{d}_i)
$$
{\em i.e.}
$$
(c_i\circ^q \bar{c}_i)\circ^p(d_i\circ^q \bar{d}_i)=(c_i\circ^p
d_i)\circ^q(\bar{c}_i\circ^p \bar{d}_i)
$$
To this end let us fix notation:
\begin{eqnarray*}
 c_i        &:& \cdots \xymatrix@C=8ex{c_q\ar[r]|{q+1}&c_q'}:\cdots \xymatrix@C=8ex{c_p\ar[r]|{p+1}&c_p'} \\
 \bar{c}_i  &:& \cdots \xymatrix@C=8ex{c_q'\ar[r]|{q+1}&c_q''}:\cdots \xymatrix@C=8ex{c_p\ar[r]|{p+1}&c_p'} \\
 d_i        &:& \cdots \xymatrix@C=8ex{d_q\ar[r]|{q+1}&d_q'}:\cdots \xymatrix@C=8ex{c_p'\ar[r]|{p+1}&c_p''}\\
 \bar{d}_i  &:& \cdots \xymatrix@C=8ex{d_q'\ar[r]|{q+1}&d_q''}:\cdots \xymatrix@C=8ex{c_p'\ar[r]|{p+1}&c_p''}
\end{eqnarray*}
The $p$-composition
$$
\circ^p:\
\C_{p+1}(c_p,c_p')\times\C_{p+1}(c_p',c_p'')\to\C_{p+1}(c_p,c_p'')
$$
is functorial w.r.t. all $q$-compositions, with $p<q$. Indeed
$q$-compositions in the product
$\C_{p+1}(c_p,c_p')\times\C_{p+1}(c_p',c_p'')$ are products of
$q$-compositions in the components. Hence the pair $(c_i\circ^q
\bar{c}_i,d_i\circ^q \bar{d}_i)$ is really a composition
$(c_i,d_i)\circ^q(\bar{c}_i,\bar{d}_i)$ and its image under
$(-)\circ^p(-)$, namely $(c_i\circ^q \bar{c}_i)\circ^p(d_i\circ^q
\bar{d}_i)$, must be equal to the q-composition of the images of the
$p$-composites $c_i\circ^p d_i$ and $\bar{c}_i\circ^p \bar{d}_i$,
namely  $(c_i\circ^p d_i)\circ^q(\bar{c}_i\circ^p \bar{d}_i)$.

\end{enumerate}
\end{proof}

\emph{Vice-versa} a globular $n$-category $\mathcal{C}$ univocally
defines a $n$-category $\C$.
\begin{proof}[\emph{(}Idea of a\emph{)} proof.]
The degenerate case $n=0$ gives immediately a ($0$-truncated globular) set.

For $n=1$ {\em Definition \ref{def:glob_n_cat}} is precisely the definition of a category  as a $1$-truncated globular set.

So let us suppose $n>1$. We define a $n$-category $\C$ in the following way.

$\C_0$ is the set $\mathcal{C}_0$ of $\mathcal{C}$.

For every pair of elements $c_0,c_0'$ of $\C_0$, we can consider the $(n-1)$-truncated globular set
$\mathcal{C}(c_0,c_0')$, where
$$
[\mathcal{C}(c_0,c_0')]_0=\{c_1\in\mathcal{C}_1\ \mathrm{s.t.}\ s(c_1)=c_0, t(c_1)=c_0' \}
$$
and inductively
$$
[\mathcal{C}(c_0,c_0')]_i=\{c_{i+1}\in\mathcal{C}_{i+1}\ \mathrm{s.t.}\ s(c_{i+1})\in[\mathcal{C}(c_0,c_0')]_{i-1}, t(c_{i+1})\in[\mathcal{C}(c_0,c_0')]_{i-1}\}
$$
$k$-Sources, $k$-targets, $k$-starting identities and $k$-compositions maps, with $k=1,\dots n$, restrict properly to $\mathcal{C}(c_0,c_0')$, hence
it is a globular $(n-1)$-category. By induction hypothesis hence a $(n-1)$-category $\C_1(c_0,c_0')$.

Moreover $0$-composition defines $(n-1)$-functors
$$
\circ^0: \mathcal{C}(c_0,c_0')\times \mathcal{C}(c_0',c_0'') \to\mathcal{C}(c_0,c_0'')
$$
for every triple $c_0,c_0',c_0''$, and $0$-starting identities define $(n-1)$-functors
$$
u^0(c_0):\I\to\mathcal{C}(c_0,c_0)
$$
These data form indeed a $n$-category.
\end{proof}

\section{The groupoid condition}
Our notion of $n$-groupoid corresponds, modulo the conversions recalled above, to the notion of
\emph{ $n$-groupoid} of Kapranov and Voevodsky in \cite{MR1130401}. According to their definition a
$n$-groupoid is a globular strict-$n$-category which satisfies a so-called \emph{groupoid-condition}. This
basically says that every equation of the kind $cx=c'$ or $yc=c'$ is (weakly) solvable, when the equation
makes sense. For sake of completeness this condition is recalled below. \\
\begin{Definition}[Kapranov and Voevodsky, {\em Definition 1.1} \cite{MR1130401}]\label{def:KV}
A $n$-category $\mathcal{C}$ is called a $n$-groupoid if for all  $i<k\leq n$ the following conditions hold
\begin{description}
  \item[(GR${}_{i,k}',\ i<k-1$)]
For each $a\in \mathcal{C}_{i+1}$, $b\in\mathcal{C}_k$ $u,v\in\mathcal{C}_{k-1}$ such that
$$
s_i(a)=t_i(u)=t_i(v),\quad u\star^i a=s_{k-1}(b),\quad  v\star^i a=t_{k-1}(b)
$$
there exist $x\in\mathcal{C}_k$, $\phi\in\mathcal{C}_{k+1}$ such that
$$
s_k(\phi)= x\star^i a,\quad  t_k(\phi) =b,\quad s_{k-1}(x)=u, \quad s_{k-1}(x)=v.
$$
  \item[(GR${}_{k-1,k}'$)]
For each $a\in \mathcal{C}_{k}$, $b\in\mathcal{C}_k$  such that
$$t_{k-1}(a)=t_{k-1}(b)$$
there exist $x\in\mathcal{C}_k$, $\phi\in\mathcal{C}_{k+1}$ such that
$$
s_k(\phi)= x\star^{k-1} a,\quad  t_k(\phi) =b.
$$

  \item[(GR${}_{i,k}'',\ i<k-1$)]
For each $a\in \mathcal{C}_{i+1}$, $b\in\mathcal{C}_k$ $u,v\in\mathcal{C}_{k-1}$ such that
$$
s_i(a)=t_i(u)=t_i(v),\quad a\star^i u=s_{k-1}(b),\quad  a\star^i v=t_{k-1}(b)
$$
there exist $x\in\mathcal{C}_k$, $\phi\in\mathcal{C}_{k+1}$ such that
$$
s_k(\phi)= a\star^i x, \quad  t_k(\phi) =b,\quad s_{k-1}(x)=u, \quad s_{k-1}(x)=v.
$$

  \item[(GR${}_{k-1,k}''$)]
For each $a\in \mathcal{C}_{k}$, $b\in\mathcal{C}_k$  such that
$$s_{k-1}(a)=s_{k-1}(b)$$
there exist $x\in\mathcal{C}_k$, $\phi\in\mathcal{C}_{k+1}$ such that
$$
s_k(\phi)= a\star^{k-1} x,\quad  t_k(\phi) =b.
$$

\end{description}
\end{Definition}
The notion recalled above is a generalization of the definition by Street in \cite{MR920944} which included only axioms for
inverses (GR${}_{k-1,k}'$) and (GR${}_{k-1,k}''$). Moreover it is a wider generalization of the definition by Brown
and Higgins in \cite{BHch81} in which such axioms where taken in a strict form ($\phi=id$).\\

Notice that Kapranov and Voevodsky motivate the new axioms for they would  ensure not only the existence
of (weak) inverses. In fact they claim (but not prove) that all four axioms together imply the existence of a coherent system
of such.\\

That our definition is equivalent to \ref{def:KV} above is a corollary to Simpson accurate analysis developed
in \cite{simpson-1998}, where he uses the Tamsamani's approach  for the treatment of a groupoid condition
for weak $n$-categories \cite{Tam96}.

This is resumed in the following inductive theorem-definition

\begin{Theorem}[Simpson, {\em Theorem 2.1} \cite{simpson-1998}] Fix $n<\infty$.\\

\emph{\textbf{I. Groupoids}}
Suppose $\mathcal{C}$ is a $[$globular $]$ strict $n$-category. The following three conditions are equivalent
(and in this case we say that $\mathcal{C}$ is a strict $n$-groupoid).

(1) $\mathcal{C}$ is a $n$ groupoid in the sense of Kapranov and Voevodsky {\em (Definition \ref{def:KV})};

(2) for all $x,y \in \mathcal{C}_0$, $\mathcal{C}(x,y)$ is a strict $(n-1)$-groupoid, and for any 1-cell
$f:x\to y$ in $\mathcal{C}$, the two $[$families of $]$ morphisms of $[$left and right $]$ compositions with $f$ are equivalences
of strict $(n-1)$-groupoids;

(3) for all $x,y\in\mathcal{C}_0$, $\mathcal{C}(x,y)$ is a strict $(n-1)$-groupoid, and
$\tau_{\leq 1}\mathcal{C}$ is a groupoid.\\

\emph{\textbf{II. Truncation}} If $\mathcal{C}$ is a strict $n$-groupoid, then define $\tau_{\leq k}\mathcal{C}$
to be the strict $k$-category whose $i$-cells are those of $\mathcal{C}$ for $i<k$, and whose $k$-cells
are the equivalence classes of $k$-cells of $\mathcal{C}$ under the equivalence relation that two are equivalent if
there is a $(k+1)$-cell joining them. The fact that this is an equivalence relation is a statement about
$(n-k)$-groupoids. The set $\tau_{\leq 0}\mathcal{C}$ will also be denoted $\pi_0\mathcal{C}$. The truncation
is again a $k$-groupoid and for $n$-groupoids $\mathcal{C}$ the truncation coincide with the operation defined
in \cite{MR1130401}.\\

\emph{\textbf{III. Equivalence}}
A morphism  $F:\mathcal{C}\to\mathcal{D}$ of strict $n$-groupoids is said to be an equivalence if the following
equivalent conditions are satisfied:

(a) (this is the definition in \cite{MR1130401}) $F$ induces an isomorphism $\pi_0\mathcal{C}\to\pi_0\mathcal{D}$, and
for every object $c\in\mathcal{C}$ $F$ induces isomorphisms $\pi_i(\mathcal{C},c)\to\pi_i(\mathcal{D},F(c))$, where these
homotopy groups are defined in \cite{MR1130401};

(b) $F$ induces a surjection $\pi_0\mathcal{C}\to\pi_0\mathcal{D}$ and for every pair of objects
$x,y\in\mathcal{C}$ $F$ induces an equivalence of $(n-1)$-groupoids
$\mathcal{C}(x,y)\to\mathcal{D}(F(x),F(y))$;

(c) if $u,v$ are i-cells in $\mathcal{C}$ sharing source and target, and if
$r:\xymatrix{F(u)\ar[r]|(.45){i+1}&F(v)}$ is a $(i+1)$-cell in $\mathcal{D}$, there exists an $(i+1)$-cell
$t:\xymatrix{u\ar[r]|(.45){i+1}&v}$ of $\mathcal{C}$ and a $(i+2)$-cell $\xymatrix{F(t)\ar[r]|{i+2}&r}$ in $\mathcal{D}$
(this includes the limiting cases $i=-1$ where $u$ and $v$ are not specified, and $i=n-1,n$ where ``$(n+1)$-cell''
means equality and ``$(n+2)$-cells'' are not specified).
\end{Theorem}

Finally we can compare these conditions with our groupoid condition. \\

We have already discussed about the correspondence between
globular and enriched versions of $n$-category.

Hence we can focus on characterizing conditions.\\

Formally our {\em Definition \ref{def:n_gpd}} amounts precisely to condition \textbf{I.}(2) above. What is to be
checked is then the notion of equivalence. This is done by {\em Proposition~\ref{prop:equiv}} that is the
inductive version of \textbf{III.}(c).\\

The notion of truncation is not directly involved, nevertheless it can be recovered by successive applications
of the sesqui-functor $\pi_0$.

\section{System of adjoint inverses}

\vskip3ex
In order to be more precise about the choices of inverses, only for this section,
we give a sharper definition of equivalence that take into account directions of cells.
\begin{Definition}
Let n-category morphism $F:\C\to\D$   be given.\\

$\bullet$ $F$ is called \emph{equivalence of n-categories} if it satisfies the following properties:\\

\framebox[1.1\width]{$n=0$}\\

$F$ is an isomorphism.\\

\framebox[1.1\width]{$n>0$}
\begin{enumerate}
  \item  for every object $d_0$ of $\D$,
  there exists an object $c_0$ of $\C$ and a 1-cell
  $d_1:d_0\to Fc_0$ such that for every $d_0'$  in $\C$, the morphism
$$
    d_1\circ - : \D_1(d_0,d_0')\to \D_1(Fc_0,d_0')
$$
is an equivalence of $(n-1)$-categories, and the morphism
$$
    -\circ d_1 : \D_1(d_0',Fc_0)\to \D_1(d_0',d_0)
$$
is a co-equivalence of $(n-1)$-categories.
  \item for every pair $c_0,c_0'$ in $\C$,
$$
F_1^{c_0,c_0'}:\C_1(c_0,c_0')\to\D_1(Fc_0,Fc_0')
$$
is an equivalence of $(n-1)$-categories.
\end{enumerate}

$\bullet$ $F$ is called \emph{co-equivalence of n-categories} if it satisfies the following properties:\\

\framebox[1.1\width]{$n=0$}\\

$F$ is an isomorphism.\\

\framebox[1.1\width]{$n>0$}
\begin{enumerate}
  \item  for every object $d_0$ of $\D$,
  there exists an object $c_0$ of $\C$ and a 1-cell
  $d_1:Fc_0\to d_0$ such that for every $d_0'$  in $\C$, the morphism
$$
    d_1\circ - : \D_1(d_0,d_0')\to \D_1(Fc_0,d_0')
$$
is a co-equivalence of $(n-1)$-categories, and the morphism
$$
    -\circ d_1 : \D_1(d_0',Fc_0)\to \D_1(d_0',d_0)
$$
is an equivalence of $(n-1)$-categories.
  \item for every pair $c_0,c_0'$ in $\C$,
$$
F_1^{c_0,c_0'}:\C_1(c_0,c_0')\to\D_1(Fc_0,Fc_0')
$$
is a co-equivalence of $(n-1)$-categories.
\end{enumerate}
\end{Definition}

This suggests to improve {\em Definition \ref{def:weak_inv}}:
\begin{Definition}
A 1-cell $c_1:c_0\to c_0'$ of a n-category $\C$ is said to be \emph{weakly invertible}, or simply an \emph{equivalence}, if, for every object
$\bar{c}_0$ of $\C$, the morphism
$$
    c_1\circ - : \C_1(c_0',\bar{c}_0)\to \C_1(c_0,\bar{c}_0) \\
$$
is an equivalence of $(n-1)$-categories, and the morphism
$$
    -\circ c_1 : \C_1(\bar{c}_0,c_0)\to \C_1(\bar{c}_0,c_0')
$$
is a co-equivalence of $(n-1)$-categories.
\end{Definition}
The dual definition for a \emph{weakly co-invertible} 1-cell.\\

When a 1-cell is weakly invertible, then it has indeed left and right weak-inverses. In fact for
$c_1:c_0\to c_0'$,
$$
c_1\circ - : \C_1(c_0',c_0)\to \C_1(c_0,c_0)
$$
being an equivalence implies that for the 1-cell $1_{c_0}:c_0\to c_0$ there exists a pair
$$
(c_1^*,\xymatrix{i_2:1_{c_1}\ar@2[r]^(.5){\sim} &c_1\circ c_1^*}),
$$
similarly for
$$
 -\circ c_1 : \C_1(c_0',c_0)\to \C_1(c_0',c_0')
$$
being a  co-equivalence implies there exists a pair
$$
(c_1^{\dag},\xymatrix{e_2:c_1^{\dag}\circ c_1\ar@2[r]^(.6){\sim}& 1_{c_1'}}).
$$

Left and right inverses are indeed equivalent: following a classical group-theoretical
argument,  the 1-composition
$$
\xymatrix{
c_1^{\dag}=c_1^{\dag}\circ 1_{c_1}
\ar@2[r]^{c_1^{\dag}\circ i_2}
&c_1^{\dag}\circ 1_{c_1}\circ c_1^*
\ar@2[r]^{e_2\circ c_1^*}
&1_{c_1'}\circ c_1^*=c_1^*}
$$
witnesses the equivalence.\\

Let us suppose we have chosen a system of inverse in a $n$-groupoid $\C$, {\em i.e.} for every  $k$-cell
$c_k$ we can exhibit a weak inverse  $c_k^*$ and equivalences
$$
i_{c_k}:\xymatrix{1\ar@2[r]&c_k\circ^{k-1}c_k^*,\qquad e_{c_k}}:\xymatrix{c_k^*\circ^{k-1}c_k\ar@2[r]&1}
$$
Then it is always possible to get a system of adjoint inverses, according to the following
\begin{Defprop} By induction over $n$.\\

In $0\mathbf{Gpd}$ and in $1\mathbf{Gpd}$ inverses are unique.\\

So let us suppose $n>1$.
Let us suppose further we have chosen systems of adjoint inverses in the hom-$(n-1)$-groupoids. Then the following
procedure gives a system of adjoint inverses in $\C$.

For every four-tuple $(c_1,c_1^*,i_{c_1},e_{c_1})$ one defines a new four-tuple
$(c_1,c_1^*,i_{c_1}',e_{c_1}')$ where $e_{c_1}'=e_{c_1}$ and $i_{c_1}'$ is given by the 1-composition (0-composition is juxtaposition)
$$
\xymatrix@C=10ex{
1_{c_0}\ar@2[r]^{i_{c_1}}
&c_1 c_1^*\ar@2[r]^{c_1 e_{c_1}^* c_1^*}
&c_1 c_1^* c_1 c_1^*\ar@2[r]^{i_{c_1}^*\, c_1 c_1^*}
&c_1 c_1^*}
$$
\end{Defprop}

Then for every 1-cell $c_1:c_0\to c_0'$ one has the  triangular equivalences:
$$
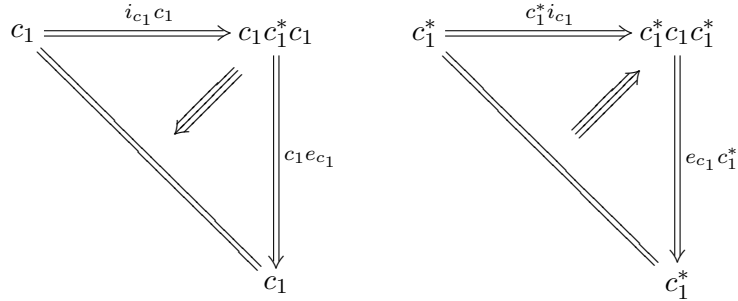
\xymatrix@!=15ex{
c_1\ar@2[r]^{i_{c_1}c_1}\ar@{=}[dr]^{}="x"
&
c_1c_1^*c_1\ar@2[d]^{c_1 e_{c_1}}\ar@{}|{}="y"\ar@{}"y";"x"|(.3){}="a1"|(.8){}="a2"
\\
&c_1
\ar@3"a1";"a2"}
\qquad
\xymatrix@!=15ex{
c_1^*\ar@2[r]^{c_1^* i_{c_1}}\ar@{=}[dr]^{}="x"
&
c_1^*c_1c_1^*\ar@2[d]^{e_{c_1}c_1^*}\ar@{}|{}="y"\ar@{}"y";"x"|(.3){}="a1"|(.8){}="a2"
\\
&c_1^*
\ar@3"a2";"a1"}
$$

\begin{proof}
It suffices to paste 3-cells in the following two diagrams, where hexagons commute quite trivially, rectangles
are identities. Concerning triangles, induction hypothesis
provide adjoint inverses for higher dimensional cells used there.

$$
\xymatrix@!C=14ex@R=5ex{
c_1\ar@2[r]^{i\, c_1}\ar@{-}[d]
&c_1c_1^*c_1\ar@2[r]^{c_1 e^* c_1^*c_1}\ar@{-}[d]
&c_1 c_1^* c_1 c_1^* c_1\ar@2[r]^{i^* c_1 c_1^* c_1}\ar@{-}[d]
&c_1 c_1^* c_1\ar@2[r]^{c_1 e}
&c_1\ar@{-}[d]
\\
c_1\ar@2[r]^{i\, c_1}\ar@{-}[d]
&c_1c_1^*c_1\ar@2[r]^{c_1 e^* c_1^*c_1}\ar@{-}[d]
&c_1 c_1^* c_1 c_1^* c_1\ar@2[r]^{c_1 c_1^* c_1 e}
&c_1 c_1^* c_1 \ar@2[r]^{i^* c_1 }\ar@{-}[d]
&c_1 \ar@{-}[d]
\\
c_1\ar@2[r]^{i\, c_1}\ar@{-}[ddrr]
&c_1c_1^*c_1\ar@2[r]^{c_1 e}\ar@{-}[dr]
&c_1 \ar@2[r]^{c_1 e^*}
&c_1 c_1^* c_1 \ar@2[r]^{i^* c_1 }\ar@{-}[dl]
&c_1 \ar@{-}[ddll]
\\
&&c_1 c_1^* c_1\ar@{}[u]|(.3){}="a1"|(.7){}="a2"&&
\\
&&c_1\ar@{}[u]|(.3){}="b1"|(.7){}="b2"&&
\ar@3"a2";"a1"^{c_1\, \eta_e^*}
\ar@3"b2";"b1"^{ \eta_i^*\, c_1}}
$$

{\small
$$
\xymatrix@!C=8ex@R=4ex{
&c_1^*\ar@2[r]^{c_1^* i}\ar@{-}[dl]
&\bullet\ar@2[r]^{c_1^* c_1 e^* c_1^*}\ar@{-}[dl]
&\bullet\ar@2[r]^{c_1^* c_1  c_1^* i^*}\ar@{-}[dl]
&c_1^* c_1  c_1^*\ar@2[r]^{e c_1^*}\ar@{-}[dl]\ar@{-}[dr]
&c_1^*\ar@{-}[dr]
\\
\bullet\ar@2[r]^{c_1^* i}\ar@{-}[d]
&\bullet\ar@2[r]^{c_1^* c_1 e^* c_1^*}\ar@{-}[d]\ar@{-}[d]
&\bullet\ar@2[r]^{c_1^* i^* c_1  c_1^* }\ar@{-}[d]
&\bullet\ar@2[r]^{c_1^* i^*}
&\bullet\ar@2[r]^{c_1^* i}\ar@{-}[d]\ar@{}[u]|(.3){}="c1"|(.7){}="c2"
&\bullet\ar@2[r]^{e c_1^*}\ar@{-}[d]
&\bullet\ar@{-}[d]
\\
\bullet\ar@2[r]^{c_1^* i}\ar@{-}[dr]
&\bullet\ar@2[r]^{c_1^* c_1 e^* c_1^*}\ar@{-}[dr]
&\bullet\ar@2[r]^{c_1^*  c_1  c_1^*  i^*}\ar@{-}[dr]
&\bullet\ar@2[r]^{c_1^* i^*}\ar@{-}[dr]
&\bullet\ar@2[r]^{c_1^* i}
&\bullet\ar@2[r]^{e c_1^*}\ar@{-}[dl]
&\bullet\ar@{-}[dl]
\\
&\bullet\ar@2[r]^(.4){c_1^* i}\ar@{-}[d]
&\bullet\ar@2[r]^(.4){c_1^* c_1 e^* c_1^*}\ar@{-}[d]
&\bullet\ar@2[r]^(.4){c_1^*  c_1  c_1^*  i^*}\ar@{-}[d]
&\bullet\ar@2[r]^{e c_1^*}\ar@{}[u]|(.3){}="d1"|(.7){}="d2"
&\bullet\ar@{-}[d]
\\
&\bullet\ar@2[r]^{c_1^* i}\ar@{-}[d]
&\bullet\ar@2[r]^{c_1^* c_1 e^* c_1^*}\ar@{-}[d]
&\bullet\ar@2[r]^{e c_1^* c_1 c_1^*}
&\bullet\ar@2[r]^{c_1^* i_1^*}\ar@{-}[d]
&\bullet\ar@{-}[d]
\\
&\bullet\ar@2[r]^{c_1^* i}\ar@{-}[ddrr]
&\bullet\ar@2[r]^{e c_1^*}\ar@{-}[dr]
&\bullet\ar@2[r]^{e^* c_1^*}
&\bullet\ar@2[r]^{c_1^* i_1^*}\ar@{-}[dl]
&\bullet\ar@{-}[ddll]
\\
&
&
&c_1^* c_1  c_1^*\ar@{}[u]|(.3){}="a1"|(.7){}="a2"
&
\\
&&&c_1^*\ar@{}[u]|(.3){}="b1"|(.7){}="b2"
\ar@3"c1";"c2"_{c_1^*\varepsilon_i}
\ar@3"d1";"d2"_{c_1^*\varepsilon_i^*}
\ar@3"a1";"a2"_{\eta_e\,c_1^*}
\ar@3"b1";"b2"_{c_1^*\,\eta_i}}
$$
}
\end{proof}

\bibliographystyle{alpha}
\addcontentsline{toc}{chapter}{Bibliography}
\bibliography{Tesi}

\backmatter
\end{document}